\documentclass[a4paper,twoside,11pt,reqno]{amsart}

\usepackage[utf8]{inputenc}

\usepackage{amsmath, amssymb, amsthm, amsfonts, amsopn, mathrsfs}

\usepackage{tikz-cd}
\usepackage{graphicx}
\usepackage[all]{xy}

\usepackage{geometry}
\usepackage{float}
\usepackage{algorithm}
\usepackage{algpseudocode}

\usepackage[bottom]{footmisc}

\usepackage{dsfont}
\usepackage{longtable}
\usepackage{booktabs}

\usepackage{enumitem}

\usepackage[pdfpagelabels]{hyperref}
\hypersetup{
  colorlinks=true,
  breaklinks=true,
  bookmarksopen=true,
  bookmarksnumbered=true,
  pdfpagemode=UseOutlines,
  plainpages=false
}

\newcommand{\dA}{{\mathcal A}}

\newcommand{\dC}{{\mathcal C}}

\newcommand{\dE}{{\mathcal E}}
\newcommand{\dF}{{\mathcal F}}

\newcommand{\dL}{{\mathcal L}}
\newcommand{\dM}{{\mathcal M}}
\newcommand{\dN}{{\mathcal N}}
\newcommand{\dO}{{\mathcal O}}

\newcommand{\dS}{{\mathcal S}}

\newcommand{\dV}{{\mathcal V}}

\newcommand{\fA}{{\mathfrak A}}
\newcommand{\fB}{{\mathfrak B}}
\newcommand{\fC}{{\mathfrak C}}

\newcommand{\fE}{{\mathfrak E}}

\newcommand{\fI}{{\mathfrak I}}

\newcommand{\fM}{{\mathfrak M}}
\newcommand{\fN}{{\mathfrak N}}
\newcommand{\fO}{{\mathfrak O}}
\newcommand{\fP}{{\mathfrak P}}

\newcommand{\fV}{{\mathfrak V}}

\newcommand{\fa}{{\mathfrak a}}

\newcommand{\fc}{{\mathfrak c}}
\newcommand{\fd}{{\mathfrak d}}

\newcommand{\fl}{{\mathfrak l}}

\newcommand{\fn}{{\mathfrak n}}
\newcommand{\fo}{{\mathfrak o}}
\newcommand{\fp}{{\mathfrak p}}
\newcommand{\fq}{{\mathfrak q}}
\newcommand{\fr}{{\mathfrak r}}
\newcommand{\fs}{{\mathfrak s}}
\newcommand{\ft}{{\mathfrak t}}

\newcommand{\bfC}{{\mathbf C}}

\newcommand{\bfF}{{\mathbf F}}

\newcommand{\bfH}{{\mathbf H}}

\newcommand{\bfP}{{\mathbf P}}
\newcommand{\bfQ}{{\mathbf Q}}
\newcommand{\bfR}{{\mathbf R}}

\newcommand{\bfZ}{{\mathbf Z}}

\newcommand{\lag}{\left\langle}
\newcommand{\rag}{\right\rangle}
\newcommand{\ncl}[1]{\langle\!\langle #1 \rangle\!\rangle}

\newcommand{\mif}{\mid}

\DeclareMathOperator{\Aut}{Aut}

\DeclareMathOperator{\im}{im}

\DeclareMathOperator{\SO}{SO}

\DeclareMathOperator{\PGammaL}{P{\Gamma}L}

\DeclareMathOperator{\GL}{GL}
\DeclareMathOperator{\PGL}{PGL}
\DeclareMathOperator{\PSL}{PSL}
\DeclareMathOperator{\SL}{SL}

\DeclareMathOperator{\Char}{char}

\DeclareMathOperator{\Comm}{Comm}

\DeclareMathOperator{\Cl}{Cl}

\DeclareMathOperator{\Ram}{Ram}

\DeclareMathOperator{\Nrd}{\rm Nrd}
\DeclareMathOperator{\trd}{\rm trd}

\DeclareMathOperator{\Gal}{Gal}

\newcommand{\an}{{\rm an}}

\newcommand{\ab}{{\rm ab}}

\DeclareMathOperator{\Sym}{\rm Sym}

\makeatletter

\newcommand{\Rmnum}[1]{\expandafter\@slowromancap\romannumeral #1@}
\makeatother

\newtheorem{thm}{Theorem}

\newtheorem{prop}[thm]{Proposition}
\newtheorem{lem}[thm]{Lemma}

\newtheorem{cor}[thm]{Corollary}

\newtheorem{thmABC}{Theorem}

\theoremstyle{definition}

\theoremstyle{remark}
\newtheorem{rmk}[thm]{Remark}
\newtheorem{nota}[thm]{Notation}

\newenvironment{pro*}[1][Proof]{{\it{#1:}} }{}
\newenvironment{pro**}[1][]{{\it{#1}} }{\hfill $\square$}

\numberwithin{equation}{section}

\title{Non-archimedean uniformization of fake quadrics}
\author{Jonah Mendel and Alina Vdovina}
\date{\today}

\begin{document}

\begin{abstract}
    We construct the first $p$-adically uniformized fake quadrics in characteristic zero, analogous to the fake quadrics in positive characteristic of Stix--Vdovina and Rungtanapirom. We find these surfaces using an algorithm, implemented in Magma, that computes fundamental domains and presentations for quaternionic semilinear $S$-arithmetic groups. With the same methods, we construct other interesting $S$-arithmetic lattices acting on products of two or more trees. We classify the torsion-free arithmetic lattices of Euler characteristic one acting on a product of two trees in characteristic zero. We also give infinite families of torsion-free lattices that act simply transitively on the vertices of a product of two or three trees.
\end{abstract}

\maketitle

\bibliographystyle{amsalpha}

\section{Introduction}
\label{sec:introduction}

Amongst groups acting on products of trees, arithmetic groups are very rare. However, arithmetic groups acting on products of trees are a significant source of inspiration for this subject, and have beautiful applications across different areas of mathematics, including geometric group theory, computer science, number theory, and algebraic geometry.

Recall that a \emph{fake quadric} over a field $k$ is a minimal, regular, smooth, projective surface $X/k$ of general type whose Chern numbers $c_1^2 = 8$ and $c_2 = 4$ agree with those of $\bfP^1_k \times \bfP^1_k$. When $k = \bfC$, these are exactly the minimal surfaces of general type with the same rational cohomology as the quadric. We say that $X$ is \emph{reducible} if it is finitely covered by a product of curves, and \emph{irreducible} otherwise. By Hirzebruch proportionality \cite{hirzebruch-2013-gesammelte-abhandlungen-collected-papers-i-19511962}, if $\Gamma$ is a torsion-free group acting properly and cocompactly on $\bfH \times \bfH$ by holomorphic isometries, and $\Gamma$ has vanishing first Betti number and holomorphic Euler characteristic $1$, then $M_\Gamma = \Gamma\backslash(\bfH \times \bfH)$ is a complex fake quadric. 

The first example of a reducible complex fake quadric $M_\Gamma$ is due to Beauville \cite[Exercise X.13(4)]{beauville-1996-complex-algebraic-surfaces}, and they were completely classified by Bauer--Catanese--Grunewald--Pignatelli \cite{bauer-catanese-and-grunewald-2008-the-classification-of-surfaces-with-p-g-q-0-isogenous-to-a-product-of-curves, bauer-catanese-grunewald-and-pignatelli-2012-quotients-of-products-of-curves-new-surfaces-with-p-g-0-and-their-fundamental-groups}. The reducible fake quadrics were generalized to Beauville surfaces, which is an extensive area of research; see, for example, \cite{fairbairn-2019-beauville-p-groups-a-survey, barker-boston-peyerimhoff-and-vdovina-2015-an-infinite-family-of-2-groups-with-mixed-beauville-structures, stix-and-vdovina-2016-series-of-p-groups-with-beauville-structure, skipper-and-thillaisundaram-2026-ggs-groups-acting-on-trees-of-growing-degrees}.

If $M_\Gamma$ is irreducible, then it follows from Margulis' arithmeticity theorem \cite{margulis-1991-discrete-subgroups-of-semisimple-lie-groups} that $\Gamma$ is commensurable with the unit group of a quaternion order, so we call $M_\Gamma$ \emph{quaternionic}. Quaternionic irreducible fake quadrics exist by work of Shavel, Takeuchi, and D\v{z}ambi\'c \cite{shavel-1978-a-class-of-algebraic-surfaces-of-general-type-constructed-from-quaternion-algebras,takeuchi-1990-algebraic-surfaces-derived-from-unit-groups-of-quaternion-algebras,dzambic-2014-fake-quadrics-from-irreducible-lattices-acting-on-the-product-of-upper-half-planes}. In \cite{linowitz-stover-and-voight-2015-commensurability-classes-of-fake-quadrics}, Linowitz--Stover--Voight (LSV) classify all maximal arithmetic groups that can contain the fundamental group of an irreducible fake quadric. This is analogous to the work of Prasad and Yeung for fake projective planes \cite{prasad-and-yeung-2005-fake-projective-planes,prasad-and-yeung-2005-fake-projective-planes-addendum}. However, while every complex fake projective plane is a ball quotient by Yau \cite{yau-1977-calabi-s-conjecture-and-some-new-results-in-algebraic-geometry}, Szab\'o recently constructed a complex fake quadric that is not aspherical, and hence not uniformized by $\bfH \times \bfH$ \cite{szabo-2026-a-fake-quadric-that-is-not-a-bidisk-quotient}. 

The first non-archimedean fake quadrics were constructed by Stix and the second author in \cite{stix-and-vdovina-2017-simply-transitive-quaternionic-lattices-of-rank-2-over-f-q-t-and-a-non-classical-fake-quadric}. Specifically, they give explicit torsion-free irreducible lattices 
$$
    \Gamma \leq \PGL_2\big(\bfF_q((t))\big) \times \PGL_2\big(\bfF_q((t))\big)
$$
for all odd prime powers $q$ such that $\Gamma$ acts simply transitively on the associated product of Bruhat--Tits trees. Adapting Mumford's theory of non-archimedean uniformization \cite{mumford-1979-an-algebraic-surface-with-k-ample-k2-9-pg-q-0}, they construct an irreducible surface over $\bfF_q((t))$ of general type that is uniformized by $\Gamma$, and in the case where $q = 3$, the surface is an irreducible (non-classical) fake quadric. Rungtanapirom also constructs a torsion-free lattice
$$
    \Gamma \leq \PGL_2\big(\bfF_2((t))\big) \times \PGL_2\big(\bfF_2((t))\big)
$$
which uniformizes a fake quadric over $\bfF_2((t))$ \cite{rungtanapirom-2018-quaternionic-arithmetic-lattices-of-rank-2-and-a-fake-quadric-in-characteristic-2}.

In this paper, we study lattices acting on products of trees in characteristic 0. Specifically, let $B/K$ be a totally definite quaternion algebra over a totally real number field $K$ with discriminant $\fd = \fd(B)$, and let $\fO$ be a hereditary order of $B$. Let $\fs$ be a square-free ideal of $\fo_K$ coprime to $\fd$ (and the level of $\fO$). In particular, $B$ splits at each prime dividing $\fs$. Then the normalizer $\Phi_{\fO,\fs}$ of the $\fo_K[1/\fs]$-order $\fO[1/\fs]$ in $B^\times$ (modulo scalars) is a lattice in
$$
    G_\fs^* = \prod_{\fp \mid \fs} \PGL_2(K_\fp),
$$
so $\Phi_{\fO,\fs}$ acts properly and cocompactly on the product of Bruhat--Tits trees $T_\fs$. By Borel and Prasad \cite{borel-and-prasad-1989-finiteness-theorems-for-discrete-subgroups-of-bounded-covolume-in-semi-simple-groups}, every lattice commensurable with $\Phi_{\fO,\fs}$ is contained in $\Phi_{\fO',\fs}$ for some hereditary order $\fO'$ of $B$ (Theorem \ref{thm:phi-maximal}). When the number of prime factors $\omega(\fs)$ of $\fs$ is $\geq 2$, Margulis' arithmeticity theorem \cite{margulis-1991-discrete-subgroups-of-semisimple-lie-groups} shows that every irreducible lattice in $G_\fs^*$ arises in this way, possibly from a quaternion algebra over another number field with the same completions. Hence every maximal irreducible lattice in $G_\fs^*$ is of the form $\Phi_{\fO,\fs}$. 

However, the lattices $\Phi_{\fO,\fs}$ are contained in finite extensions coming from Galois automorphisms that also act on $T_\fs$. If $K_0$ is the fixed field of $\Aut_\bfQ(K)$, then 
$$
    \Psi_{\fO,\fs} = \{g \in \Aut_{K_0}(B) : g(\fO[1/\fs]) = \fO[1/\fs]\}
$$
is a lattice in the \emph{semilinear group}
$$
    G_\fs = \left(\prod_{\fp\mid\fs} \PGammaL_2(K_\fp)\right) \rtimes \Sym(\fs).
$$
Here, $\PGammaL_2(K_\fp)$ is the semilinear group over $K_\fp$ and $\Sym(\fs)$ permutes the isomorphic factors. See Sections \ref{sec:quaternion-algebras} and \ref{sec:arithmetic-groups} for more details. Note that $G_\fs^*$ has finite index in $G_\fs$, and $G_\fs$ acts on the product of trees $T_\fs$. Semilinear arithmetic groups such as these are known to experts, but we could not find a systematic account of them in the literature. In the setting of fake quadrics, they first appear in the work of LSV \cite{linowitz-stover-and-voight-2015-commensurability-classes-of-fake-quadrics}. We give the background on semilinear arithmetic groups that we need in Section \ref{sec:arithmetic-groups}, following the treatment of LSV closely.

As in \cite{stix-and-vdovina-2017-simply-transitive-quaternionic-lattices-of-rank-2-over-f-q-t-and-a-non-classical-fake-quadric}, our lattices uniformize certain $p$-adic varieties. We now describe this process. Suppose that every prime dividing $\fs$ lies over the same rational prime $p$. We call such an ideal \emph{eligible} over $p$. Fix a finite extension $L/\bfQ_p$, a completion $\bfC_p$ of the algebraic closure of $L$, and embeddings $i_\fp : K_\fp \to L$ such that $i_\fp(K_\fp) = i_\fq(K_\fq)$ if $K_\fp \cong K_\fq$. Then 
$$
    H_\fs = \left(\prod_{\fp\mid\fs} \PGL_2(K_\fp)\right) \rtimes \Sym(\fs)
$$
acts on the product of Drinfeld upper half planes, 
$$
    \Omega_\fs = \prod_{\fp \mid \fs} \Omega_\fp \otimes_{K_\fp,i_\fp} L, \qquad \Omega_\fp = (\bfP^1_{K_\fp})^\an \smallsetminus \bfP^1(K_\fp),
$$
via $L$-analytic maps \cite[Theorem 1]{berkovich-1995-the-automorphism-group-of-the-drinfeld-half-plane}.  In \cite{varshavsky-1998-p-adic-uniformization-of-unitary-shimura-varieties,varshavsky-1998-p-adic-uniformization-of-unitary-shimura-varieties-ii}, Varshavsky generalizes the work of \v{C}erednik \cite{cerednik-1976-uniformization-of-algebraic-curves-by-discrete-arithmetic-subgroups-of-pgl2-kw-with-compact-quotients} and Mumford \cite{mumford-1979-an-algebraic-surface-with-k-ample-k2-9-pg-q-0} on $p$-adic uniformization of Shimura varieties. We will need the following special case of his work.

\begin{prop}\cite{varshavsky-1998-p-adic-uniformization-of-unitary-shimura-varieties-ii}
    \label{prop:uniformization}
    Let $K$, $\fs$, and $L$ be as above, and assume that $\fs$ is eligible with $r = \omega(\fs) \geq 2$. Let $\Gamma \leq H_\fs$ be an irreducible lattice that acts freely on $\Omega_\fs(\bfC_p)$. Then there is a unique (up to isomorphism) smooth projective variety $Y_\Gamma/L$ with ample canonical bundle, and an isomorphism of $L$-analytic spaces,
    $$
        (Y_\Gamma)^\an \cong \Gamma \backslash \Omega_\fs.
    $$
    The structure sheaf $\dO_{Y_\Gamma}$ has Euler characteristic $\chi(\dO_{Y_\Gamma}) = \chi(\Gamma)$ and $h^1(\dO_{Y_\Gamma}) = 0$, and for every partition $I$ of $r$, the Chern numbers of $Y_\Gamma$ are
    $$
        c_I(Y_\Gamma) = \chi(\Gamma) c_I((\bfP^1)^r).
    $$
    Furthermore, for every embedding $L \hookrightarrow \bfC$, the resulting complex variety $Y_{\Gamma,\bfC}$ satisfies $(Y_{\Gamma,\bfC})^\an \cong \Lambda_\Gamma\backslash \bfH^r$ for an irreducible, cocompact, torsion-free lattice $\Lambda_\Gamma$.

    In particular, if $r = 2$ and $\chi(\Gamma) = 1$, then $Y_\Gamma$ is a fake quadric over $L$, and $Y_{\Gamma,\bfC}$ is a quaternionic irreducible fake quadric, and thus appears in the classification of LSV.
\end{prop}

We can also deduce the commensurability class of $\Lambda_\Gamma$. See Section \ref{sec:uniformization} for details. Observe that if $\Gamma$ is torsion-free, then it acts freely on $\Omega_\fs$. However, the converse is not true (see Sections \ref{sec:uniformization} and \ref{sec:examples}). We now state our main results.

\begin{thmABC}
    \label{thm:A}
    Let $K$, $B$, and $\fs$ be as above with $\fs$ eligible and $r = \omega(\fs) = 2$. Let $\Gamma \le H_\fs$ be a lattice in the commensurability class of $B$ with $\chi(\Gamma) = 1$ that acts freely on $\Omega_\fs(\bfC_p)$. Then $\Gamma$ is in one of two commensurability classes described below.
    \begin{center}
        \setlength{\tabcolsep}{16pt}
        \begin{tabular}{ccc}
            \toprule
            $K$ & $\fd(B)$ & $p$ \\
            \midrule
            $\bfQ(\sqrt{7})$ & $\fp_2\fp_7$ & $3$ \\
            $\bfQ(\sqrt{6})$ & $\fp_2\fp_3$ & $5$ \\
            \bottomrule
        \end{tabular}
    \end{center}
    In each case, there are exactly seven such lattices up to conjugacy in $G_\fs$, and each lattice contains torsion. 
    
    For each such $\Gamma$ and each embedding $L \hookrightarrow \bfC$, the surface $Y_{\Gamma,\bfC}$ is an irreducible complex fake quadric uniformized by $\bfH \times \bfH$. Furthermore, the commensurability class of $\Lambda_\Gamma$ is independent of the choice of embedding.
\end{thmABC}

To our knowledge, these are the first $p$-adically uniformized fake quadrics in characteristic 0. The fact that all the lattices contain torsion contrasts with \cite{stix-and-vdovina-2017-simply-transitive-quaternionic-lattices-of-rank-2-over-f-q-t-and-a-non-classical-fake-quadric} and \cite{rungtanapirom-2018-quaternionic-arithmetic-lattices-of-rank-2-and-a-fake-quadric-in-characteristic-2}, where the lattices are torsion-free, but parallels the work of Allcock and Kato on fake projective planes, discussed below. The $5$-adic lattice in Theorem \ref{thm:A} is also pathological in the following sense. Suppose that $\fs = \fp\fq$ with $N\fp = N\fq = q$ odd. If $\Gamma \leq G_\fs$ is torsion-free, then $(q-1)^2/4 \mid \chi(\Gamma)$ (see Lemma \ref{lem:euler-vertex-count}).  Thus, 5-adic fake quadrics must be uniformized by lattices with torsion.

Our constructions have precedent in the literature. In \cite{allcock-and-kato-2013-the-densest-lattices-in-pgl3-q2,allcock-and-kato-2017-a-fake-projective-plane-via-2-adic-uniformization-with-torsion}, Allcock and Kato construct a lattice in $\PGL_3(\bfQ_2)$ with torsion that uniformizes a fake projective plane. Their examples are commensurable with the torsion-free lattice constructed by Mumford in \cite{mumford-1979-an-algebraic-surface-with-k-ample-k2-9-pg-q-0}. It was observed by Mumford that if $\Gamma \leq \PGL_3(\bfQ_p)$ is torsion-free and uniformizes a fake projective plane, then $p = 2$. Allcock and Kato ask whether there can be a lattice with torsion in $\PGL_3(\bfQ_p)$ with $p > 2$ that uniformizes a fake projective plane. Our $5$-adic construction answers the fake quadric analogue of their question.

The proof of Theorem \ref{thm:A} rests on the classification of Linowitz--Stover--Voight \cite{linowitz-stover-and-voight-2015-commensurability-classes-of-fake-quadrics}. It follows from Proposition \ref{prop:uniformization} that for each maximal lattice in the LSV census, there are finitely many maximal lattices $\Psi$ which can contain a lattice uniformizing a fake quadric. We then compute a fundamental domain in the associated product of Bruhat--Tits trees for each $\Psi$, and look for $\Gamma \leq \Psi$ acting freely on $\Omega_\fs$ with bounded Euler characteristic.

We give a general algorithm to compute fundamental domains of $\Psi_{\fO,\fs}$ for hereditary orders $\fO$. Our algorithm builds on several earlier algorithms. First, Franc and Masdeu compute fundamental domains in the Bruhat--Tits tree for $p$-unit groups of Eichler orders in definite rational quaternion algebras \cite{franc-and-masdeu-2014-computing-fundamental-domains-for-the-bruhat-tits-tree-for-gl2-qp-p-adic-automorphic-forms-and-the-canonical-embedding-of-shimura-curves}. Masdeu and Torrents then extend their method to totally real fields \cite{masdeu-and-torrents-2025-fundamental-domains-for-quaternionic-s-arithmetic-groups-over-totally-real-fields}, and Arenas-Carmona extends their method to global function fields \cite{arenas-carmona-2014-computing-quaternion-quotient-graphs-via-representations-of-orders}. On products of trees, \cite{chinburg-friedlander-howe-kosters-singh-stover-zhang-and-ziegler-2015-presentations-for-quaternionic-s-unit-groups} gives presentations of $S$-unit groups of orders in class number one, definite quaternion algebras over $\bfQ$. The first author and Yu  \cite{mendel-and-yu-2026-s-arithmetic-groups-acting-simply-transitively-on-products-of-bruhat-tits-trees} recently give a different algorithm for products of trees which handles class number one quaternion algebras over totally real fields based on the ideas of Rungtanapirom--Stix--Vdovina \cite{stix-and-vdovina-2017-simply-transitive-quaternionic-lattices-of-rank-2-over-f-q-t-and-a-non-classical-fake-quadric,rungtanapirom-stix-and-vdovina-2019-infinite-series-of-quaternionic-1-vertex-cube-complexes-the-doubling-construction-and-explicit-cubical-ramanujan-complexes}. Vast gives a similar algorithm  in \cite{vast-2026-new-families-of-vertex-transitive-quaternionic-lattices-in-products-of-trees}.

Sch\"onnenbeck \cite{schonnenbeck-2017-computing-s-unit-groups-of-orders} develops a general framework for presenting $S$-unit groups of orders, and computes some explicit examples. Our algorithm follows the same general idea as \cite{schonnenbeck-2017-computing-s-unit-groups-of-orders}, but it is specialized to the quaternionic semilinear groups $\Psi$. Our algorithm is implemented in Magma \cite{bosma-cannon-and-playoust-1997-the-magma-algebra-system-i-the-user-language} and relies heavily on its robust library of functions for working with quaternion algebras, much of which rests on the work of Voight \cite{voight-2005-quadratic-forms-and-quaternion-algebras-algorithms-and-arithmetic} and Kirschmer--Voight \cite{kirschmer-and-voight-2010-algorithmic-enumeration-of-ideal-classes-for-quaternion-orders}.

With our algorithm, we also construct other interesting lattices of small Euler characteristic. In particular, we obtain an Euler characteristic 2 surface over $\bfQ_2$.

\begin{thmABC}
    \label{thm:B}
    For each row of Table \ref{tab:thm-b}, there is a lattice $\Gamma \le H_\fs$ in the commensurability class determined by $K$ and $B$ with the listed properties that acts freely on $\Omega_\fs(\bfC_p)$. In particular, there exists $Y/\bfQ_2$ as in Proposition \ref{prop:uniformization} with $\chi(\dO_Y) = 2$.
\end{thmABC}

\begin{table}[ht]
    \centering
    \renewcommand{\arraystretch}{1.3}
    \setlength{\tabcolsep}{12pt}
    \begin{tabular}{lllccc}
        \toprule
        $K$ & $\fd(B)$ & $p$ & $\chi(\Gamma)$ & torsion-free & linear \\
        \midrule
        $\bfQ(\sqrt{17})$ & $\fo_K$ & $2$ & $2$ & no & yes \\
        $\bfQ(\sqrt{13})$ & $\fo_K$ & $3$ & $2$ & yes & yes \\
        $\bfQ(\sqrt{33})$ & $\fo_K$ & $2$ & $3$ & yes & yes \\
        $\bfQ(\sqrt{6})$ & $\fp_2\fp_3$ & $5$ & $4$ & yes & yes \\
        $\bfQ(\sqrt{5})$ & $\fo_K$ & $p \equiv 1 \bmod 5$ & $(p-1)^2/4$ & yes & yes \\
        \bottomrule
    \end{tabular}
    \medskip
    \caption{The lattices of Theorem \ref{thm:B}. The torsion-free and linear columns indicate whether or not the commensurability class determined by $K$ and $B$ contains a lattice with the listed Euler characteristic that acts freely on $\Omega_\fs(\bfC_p)$ and has this property.}
    \label{tab:thm-b}
\end{table}

We say more about the infinite family over $K = \bfQ(\sqrt{5})$ in Theorem \ref{thm:B}. It is analogous to the infinite family constructed by Stix and the second author. However, their family comes from a uniform construction of quaternion algebras over $\bfF_q(t)$ with $q$ varying. On the other hand, our new family of lattices is constructed from $\fO[1/p]^\times$ for varying $p$ and $\fO$ a fixed maximal order of $(-1,-1)_{K}$. This construction is based on recent work of the first author with Jiahui Yu \cite{mendel-and-yu-2026-s-arithmetic-groups-acting-simply-transitively-on-products-of-bruhat-tits-trees}, where they construct many different families of simply transitive actions on products of Bruhat--Tits trees from class number one quaternion algebras. Other infinite families can be constructed from \cite{mendel-and-yu-2026-s-arithmetic-groups-acting-simply-transitively-on-products-of-bruhat-tits-trees}. For instance, see Sections \ref{sec:cubic} and \ref{sec:family-cubic} for 3-dimensional constructions. Independently, Vast also obtained similar results about infinite families of simply transitive lattices \cite{vast-2026-new-families-of-vertex-transitive-quaternionic-lattices-in-products-of-trees}. Lastly, these families give a higher-dimensional analogue of \cite{amoros-and-milione-2019-mumford-curves-covering-p-adic-shimura-curves-and-their-fundamental-domains}, where Amor\'os and Milione use similar methods to those of \cite{mendel-and-yu-2026-s-arithmetic-groups-acting-simply-transitively-on-products-of-bruhat-tits-trees} to construct infinite families of Mumford curves from class number one quaternion algebras.

Notice that it follows from Theorem \ref{thm:A} that there is no irreducible, torsion-free lattice $\Gamma \leq H_\fs$ with $\chi(\Gamma) = 1$ and $\omega(\fs) = 2$ for eligible $\fs$. However, we can construct such lattices when $\fs$ is not eligible. In fact, we classify all irreducible, torsion-free, Euler characteristic one lattices in $G_\fs$ when $\omega(\fs) = 2$.

\begin{thmABC}
    \label{thm:C}
    Let $K$, $B$, and $\fs$ be as above with $\omega(\fs) = 2$. Let $\Gamma \leq G_\fs$ be a torsion-free lattice with $\chi(\Gamma) = 1$ in the commensurability class of $B$. Then $\Gamma$ is in one of three commensurability classes described below.
    \begin{center}
        \setlength{\tabcolsep}{16pt}
        \begin{tabular}{ccc}
            \toprule
            $K$ & $\fd(B)$ & $\fs$ \\
            \midrule
            $\bfQ$ & $5$ & $6$ \\
            $\bfQ$ & $13$ & $6$ \\
            $\bfQ(\sqrt{3})$ & $\fo_K$ & $\fp_2\fp_3$ \\
            \bottomrule
        \end{tabular}
    \end{center}
    In each case, there is exactly one such lattice up to conjugacy, and it is linear. In particular, $\fs$ is not eligible.
\end{thmABC}

The paper is organized as follows. Section \ref{sec:quaternion-algebras} gives background on quaternion algebras, Bruhat--Tits trees, and $p$-adic semilinear groups. In Section \ref{sec:arithmetic-groups}, we construct arithmetic groups from quaternion algebras and generalize this construction to semilinear arithmetic groups. We also explain Prasad's volume formula and the construction of maximal semilinear arithmetic groups $\Psi$. Section \ref{sec:fundamental-domains} gives an algorithm for computing a fundamental domain and a presentation for $\Psi$. In Section \ref{sec:uniformization}, we derive Proposition \ref{prop:uniformization} from the work of Varshavsky. Section \ref{sec:examples} constructs the lattices discussed in this introduction, together with other lattices of interest. These constructions will prove Theorems \ref{thm:A}, \ref{thm:B}, and \ref{thm:C}. Finally, Appendix \ref{app:square-complexes} displays the square complexes associated with the lattices of Section \ref{sec:examples}.

\begin{nota}
    We establish the following notation.
    \begin{itemize}
        \item If $K$ is a number field or a non-archimedean local field, then $\fo_K$ is its ring of integers. 
        \item If $\fp$ is a prime of a number field $K$, then $K_\fp$ is the completion of $K$ at $\fp$ and $N\fp = |\fo_K/\fp|$ is its absolute norm.
        \item When there is no ambiguity, $\fp_q$ and $\fq_q$ denote primes of $\fo_K$ of norm $q$.
        \item For an ideal $\fs \subseteq \fo_K$, $\omega(\fs)$ is the number of distinct prime divisors of $\fs$.
        \item For a CW-complex $X$, $X^{(k)}$ is its $k$-skeleton.
        \item For a lattice $\Gamma$, $\chi(\Gamma)$ is its orbifold Euler characteristic.
        \item $\bfH$ is the upper half plane.
    \end{itemize}
\end{nota}

\subsection*{Acknowledgments}

The authors are very grateful to Jakob Stix for many discussions about this paper, especially on $p$-adic uniformization and the Prasad volume formula, and for his feedback on earlier drafts. They thank Alan Reid for his feedback on the present manuscript, and the organizers of the 2025 UT Austin Analytic Group Theory Mini-school, where this research began.

The first author thanks his advisor, Alan Reid, for his interest in this project and for many conversations about arithmetic groups. He also thanks Matthew Stover and John Voight for helpful conversations.

\subsection*{AI Disclosure}

We used the large language models Claude (Anthropic) and ChatGPT (OpenAI) during this project, which began in March 2025. We used both models over overlapping periods. Below, we refer to them together as AI.

The authors obtained the following constructions with only light AI assistance in writing code: the lattices in Sections~\ref{sec:chi1-lattices}, \ref{sec:qsqrt13}, \ref{sec:qsqrt13-chi3}, \ref{sec:qsqrt17}, \ref{sec:cubic}, and~\ref{sec:families}, and the linear lattices in Sections~\ref{sec:qsqrt33} and~\ref{sec:qsqrt6}. In particular, we wrote a first version of the algorithm of Section~\ref{sec:fundamental-domains} with minimal AI assistance. We later used AI to optimize this algorithm substantially and to extend it to the semilinear setting. The proofs of Theorems~\ref{thm:C} and~\ref{thm:chi2-classification} rely on an enumeration of certain square complexes. AI wrote the Magma code for this enumeration, following explicit instructions from the authors.

It was known to the authors that the approach of Stix and Vdovina \cite{stix-and-vdovina-2017-simply-transitive-quaternionic-lattices-of-rank-2-over-f-q-t-and-a-non-classical-fake-quadric} generalizes to $p$-adic lattices in characteristic zero, such as those in the present paper. A literature search carried out by AI identified the papers of Varshavsky \cite{varshavsky-1998-p-adic-uniformization-of-unitary-shimura-varieties,varshavsky-1998-p-adic-uniformization-of-unitary-shimura-varieties-ii} as also relevant. AI also proposed Lemma~\ref{lem:free-action} and its proof.

During writing, we used AI to proofread the text, suggest edits, generate figures and tables, simplify proofs, and improve the exposition. Discussions with AI substantially improved the background material in Sections~\ref{sec:quaternion-algebras} and~\ref{sec:arithmetic-groups} over earlier drafts. Our account of Varshavsky's work in Section~\ref{sec:uniformization} also draws heavily on such discussions.

The authors take full responsibility for all mathematical content and code in this paper.

\section{Quaternion algebras}
\label{sec:quaternion-algebras}

In this section, we give the necessary background on quaternion algebras and orders over local and global fields. We then define the Bruhat--Tits tree in terms of orders and ideals of $M_2(K_\fp)$. We also define the action of $\PGammaL_2(K_\fp)$ on the tree. Standard references for this material are Maclachlan--Reid \cite{maclachlan-and-reid-2003-the-arithmetic-of-hyperbolic-3-manifolds}, Vign\'eras \cite{vigneras-1980-arithmetique-des-algebres-de-quaternions}, and Voight \cite{voight-2021-quaternion-algebras}.

\subsection{Quaternion algebras}
\label{sec:quaternion-definitions}

Let $K$ be a field. A \emph{quaternion algebra} $B$ over $K$ is a central simple $K$-algebra of dimension $4$. The unique standard involution, also called \emph{conjugation}, is denoted $\alpha\mapsto\overline{\alpha}$, and the \emph{reduced norm} and the \emph{reduced trace}
$$
    \Nrd(\alpha)=\alpha\overline{\alpha},\qquad
    \trd(\alpha)=\alpha+\overline{\alpha},
$$
lie in $K$ for all $\alpha\in B$. If $\Char(K)\neq2$, then $B\cong(a,b)_K$ for some $a,b\in K^\times$, where
$$
    (a,b)_K=K+ Ki+ Kj+ Kij,
    \qquad i^2=a,\quad j^2=b,\quad ij=-ji,
$$
and conjugation is
$$
    \overline{x+yi+zj+wij}=x-yi-zj-wij.
$$
By Wedderburn's theorem, either $B\cong M_2(K)$, in which case we say that $B$ \emph{splits}, or $B$ is a division algebra.

\subsection{Orders and ideals}
\label{sec:orders-ideals}

Let $R$ be a Dedekind domain with field of fractions $K$. An \emph{$R$-order} in $B$ is a subring $\fO\subseteq B$ that is a finitely generated $R$-module and contains a $K$-basis of $B$. An order is \emph{maximal} if it is not properly contained in another order, and \emph{Eichler} if it is the intersection of two maximal orders.

Let $\fO$ be an order. For a prime $\fp$ of $R$, the completions of $B$ and $\fO$ at $\fp$ are
$$
    B_\fp=B\otimes_KK_\fp,\qquad\fO_\fp=\fO\otimes_RR_\fp.
$$
Then $B_\fp$ is a quaternion algebra over $K_\fp$ and $\fO_\fp$ is an $R_\fp$-order of $B_\fp$.

A \emph{right $\fO$-ideal} is a finitely generated $R$-submodule $\fI\subseteq B$ that is \emph{locally principal}. That is, for every prime $\fp$ of $R$,
$$
    \fI_\fp=\fI\otimes_RR_\fp=x_\fp\fO_\fp
    \qquad\text{for some }x_\fp\in B_\fp^\times.
$$
Its \emph{left} and \emph{right orders} are
$$
    \fO_L(\fI)=\{\alpha\in B:\alpha\fI\subseteq\fI\},\qquad
    \fO_R(\fI)=\{\alpha\in B:\fI\alpha\subseteq\fI\}.
$$
The ideal $\fI$ is \emph{integral} if $\fI\subseteq\fO$, and its \emph{reduced norm} is the fractional ideal
$$
    \Nrd(\fI)=\sum_{\alpha\in\fI}\Nrd(\alpha)R,
    \qquad
    \Nrd(\fI)R_\fp=\Nrd(x_\fp)R_\fp.
$$
For $\alpha \in B^\times$ and an ideal (order) $\fA$, let
$$
    {}^\alpha \fA = \alpha \fA \alpha^{-1},
$$
which is still an ideal (order) of $B$. Thus, $B^\times$ acts on the ideals (orders) of $B$.

\subsection{Quaternion algebras over local fields}
\label{sec:local-quaternions}

Let $K$ be a local field of characteristic zero. If $K$ is archimedean, then $B$ is a division algebra if and only if $K=\bfR$ and $B\cong(-1,-1)_\bfR$, and then
$$
    B^\times/\bfR^\times\cong\SO(3)
$$
is compact. For the rest of this section, suppose that $K$ is non-archimedean, with valuation ring $\fo_K$, and maximal ideal $\fp$. There is a unique division quaternion algebra $A$ over $K$. It has a unique maximal order
$$
    \fM=\{\alpha\in A: \Nrd(\alpha) \in \fo_K\},
$$
and $A^\times/K^\times$ is compact. We say that $\fM$ has \emph{level} $\fn(\fM) = \fo_K$. Note that $\fM$ is the only Eichler order in $A$.

Now suppose that $B=M_2(K)$. Every Eichler order $\fO$ of $B$ is conjugate to
\begin{equation}
    \label{eq:standard-eichler-order}
    \begin{pmatrix}
        \fo_K & \fo_K \\
        \fp^e & \fo_K
    \end{pmatrix}
    =M_2(\fo_K)\cap
    \begin{pmatrix}
        \fo_K & \fp^{-e} \\
        \fp^e & \fo_K
    \end{pmatrix}
\end{equation}
for a unique $e\geq0$, and we say that $\fO$ has \emph{level} $\fn(\fO)=\fp^e$. If $e = 0$, then $\fO$ is maximal. The \emph{Bruhat--Tits tree} $T_\fp$ of $K$ is the graph with vertices and edges respectively,
$$
    \{\fM\subseteq M_2(K) : \fn(\fM) = \fo_K\},
    \qquad
    \{\{\fM,\fM'\} : \fn(\fM\cap\fM')=\fp\}.
$$
It is an $(N\fp+1)$-regular tree, and its path metric is
$$
    d(\fM,\fM')=e\iff\fn(\fM\cap\fM')=\fp^e.
$$
Note that the unoriented non-backtracking paths of length $e$ in $T_\fp$ are in bijection with the Eichler orders in $B$ of level $\fp^e$. We can describe neighbors of an order via the integral ideals of norm $\fp$. For details, see \cite[Section 26.4]{voight-2021-quaternion-algebras}.

\begin{lem}\label{lem:tree-neighbors}
    Let $\fM\subseteq B$ be a maximal order. Then there is a bijection
    $$
        \bigl\{\text{integral right }\fM\text{-ideals }\fP:
        \Nrd(\fP)=\fp\bigr\}
        \longrightarrow
        \bigl\{\text{neighbors of }\fM\text{ in }T_\fp\bigr\},
        \qquad
        \fP\longmapsto\fO_L(\fP).
    $$
\end{lem}

The group $\PGL_2(K)=B^\times/K^\times$ acts transitively on the vertices by conjugation. Furthermore, given an Eichler order $\fO$, conjugation preserves the level of $\fO$, so $\PGL_2(K)$ acts via automorphisms on $T_\fp$ with stabilizers
\begin{equation}\label{eq:vertex-stabilizer}
    K^\times\fM^\times/K^\times\cong\fM^\times/\fo_K^\times.
\end{equation}
For $g\in\GL_2(K)$ and every vertex $\fM$, we have
\begin{equation}\label{eq:distance-parity}
    d(\fM,{}^g\fM)\equiv v_\fp(\det g)\pmod 2.
\end{equation}
For more details on the Bruhat--Tits tree, see \cite{voight-2021-quaternion-algebras,maclachlan-and-reid-2003-the-arithmetic-of-hyperbolic-3-manifolds,serre-1980-trees-2003-printing}.

We will further consider the projective general semilinear group, $\PGammaL_2(K)$, which also acts on $T_\fp$. Namely, let $K_0 \subseteq K$ be the fixed field under the group of continuous automorphisms $\Aut(K) = \Aut_{\bfQ_p}(K)$. Then the projective general semilinear group is the semidirect product
$$
    \PGammaL_2(K) = \PGL_2(K) \rtimes \Gal(K/K_0)
$$
with the Galois action given by functoriality of $\PGL_2(-)$. It contains $\PGL_2(K)$ as a subgroup of finite index equal to the degree of $K/K_0$. Given $g = (\alpha,\sigma) \in \PGammaL_2(K)$ and a maximal order $\fM$, we set 
$$
    {}^g\fM = \alpha \sigma(\fM) \alpha^{-1},
$$
which is a vertex of $T_\fp$ since all maximal orders in $B$ are conjugate.

Let $\fO$ be an Eichler order in a quaternion algebra over $K$. We call $\fO$ \emph{hereditary} if $\fn(\fO)$ is square-free.

\begin{lem}\label{lem:local-units}
    Let $\fO$ be a hereditary order of a quaternion algebra over $K$. If $\alpha$ normalizes $\fO$ and $\Nrd(\alpha)\in\fo_K^\times$, then $\alpha\in\fO^\times$.
\end{lem}

\begin{proof}
    For the division algebra $A$, the claim follows immediately from the definition of $\fM$.

    We now consider the case where the quaternion algebra splits. Suppose that $\fM$ is maximal. Then by \eqref{eq:vertex-stabilizer}, we have
    $$
        \alpha=\lambda\beta, \qquad \lambda\in K^\times,
        \quad \beta\in\fM^\times.
    $$
    Observe that $\lambda \in \fo_K^\times$ since $\Nrd(\beta) \in \fo_K^\times$ and $\lambda^2\Nrd(\beta)=\Nrd(\alpha) \in \fo_K^\times$. Thus, $\alpha \in \fM^\times$.

    Now suppose that $\fO=\fM\cap\fM'$ has level $\fp$. Then $\alpha$ permutes $\fM$ and $\fM'$. Since $d(\fM,{}^\alpha\fM)$ is even by \eqref{eq:distance-parity}, $\alpha$ fixes both, so $\alpha\in\fM^\times\cap\fM'^\times=\fO^\times$ by the previous case.
\end{proof}

\subsection{Quaternion algebras over number fields}
\label{sec:global-quaternions}

Let $K$ be a number field and $B/K$ a quaternion algebra. For a place $v$ of $K$, we say that $B$ \emph{splits} at $v$ if $B_v=B\otimes_KK_v\cong M_2(K_v)$, and \emph{ramifies} at $v$ otherwise. Let $\Ram(B)$ be the set of places at which $B$ ramifies. By class field theory, the map $B \mapsto \Ram(B)$ induces a bijection
\begin{equation}
    \label{eq:quaternion-classification}
    \left\{
    \begin{array}{c}
    \text{isomorphism classes of} \\
    \text{quaternion algebras } B/K
    \end{array}
    \right\}
    \longrightarrow
    \left\{
        \begin{array}{c}
        \text{finite sets of non-complex places of } K \\
        \text{of even cardinality}
        \end{array}
    \right\}.
\end{equation}
We say that $B$ is \emph{totally definite} if $\Ram(B)$ contains all archimedean places of $K$. If $\Ram_f(B)$ denotes the set of finite places at which $B$ ramifies, then the \emph{discriminant} of $B$ is
$$
    \fd=\fd(B)=\prod_{\fp\in\Ram_f(B)}\fp.
$$
If $B$ is totally definite, then $K$ is totally real, and a totally definite quaternion algebra is determined by its discriminant.

For a square-free (possibly trivial) ideal $\fs\subseteq\fo_K$, the ring of \emph{$\fs$-integers} of $K$ is
$$
    \fo_K[1/\fs]=\{x\in K:v_\fp(x)\geq0
    \text{ for all }\fp\nmid\fs\}.
$$
We identify the ideals of $\fo_K[1/\fs]$ with the ideals of $\fo_K$ coprime to $\fs$. Let $R=\fo_K[1/\fs]$. The next lemma follows from the local-global principle for lattices \cite[Theorem 9.4.9]{voight-2021-quaternion-algebras}.

\begin{lem}
    \label{lem:local-global}
    Let $\fO$ be an $R$-order and let $S$ be a finite set of primes of $R$. Then we have a bijection
    $$
        \bigl\{R\text{-orders }\fA\subseteq B:
        \fA_\fq=\fO_\fq\text{ for all }\fq\notin S\bigr\}
        \longrightarrow
        \prod_{\fp\in S}
        \bigl\{R_\fp\text{-orders in }B_\fp\bigr\},
        \qquad
        \fA\longmapsto(\fA_\fp)_{\fp\in S}.
    $$
    Furthermore, such an order $\fA\subseteq B$ is maximal or Eichler if and only if $\fO_\fq$ has this property for all $\fq\notin S$ and $\fA_\fp$ has it for all $\fp\in S$.
\end{lem}

An $R$-order $\fO$ is \emph{hereditary} if $\fO_\fp$ is hereditary for every prime $\fp$ of $R$. Its \emph{level} is the $R$-ideal $\fn(\fO)$ with
$$
    \fn(\fO)R_\fp=\fn(\fO_\fp)
    \qquad\text{for every prime $\fp$ of $R$}.
$$
We now extend Lemma \ref{lem:local-units} to global orders.

\begin{lem}\label{lem:global-units}
    Let $\fO$ be a hereditary $R$-order. If $\alpha\in B^\times$ normalizes $\fO$ and $\Nrd(\alpha)\in R^\times$, then $\alpha\in\fO^\times$.
\end{lem}

\begin{proof}
    For every prime $\fp$ of $R$, the completion $\fO_\fp$ is hereditary and normalized by $\alpha$, so $\alpha\in\fO_\fp^\times$ by Lemma \ref{lem:local-units}. Thus, $\alpha\in\fO^\times$.
\end{proof}

\section{Quaternionic \texorpdfstring{$S$}{S}-arithmetic groups}
\label{sec:arithmetic-groups}

In this section, we establish the notation and results that we will need about quaternionic $S$-arithmetic groups. Most of this material is standard, the exception being our treatment of semilinear arithmetic groups, which are known to experts but for which we could not find an account sufficient for our purposes. Our treatment is inspired by \cite{linowitz-stover-and-voight-2015-commensurability-classes-of-fake-quadrics}.

\subsection{The ambient group}
\label{sec:ambient-group}

Let $K$ be a number field of degree $n$, and let $\fs \subseteq \fo_K$ be a non-trivial, square-free ideal. Let
$$
    \Sym(\fs) = \{\sigma \in \Sym(\{\fp \mid \fs\}) : K_\fp \cong K_{\sigma\fp} \text{ for all } \fp \mid \fs\}.
$$
We will work with the locally compact, Hausdorff, topological groups
$$
    G_\fs = \left(\prod_{\fp \mid \fs} \PGammaL_2(K_\fp)\right) \rtimes \Sym(\fs), \qquad G_\fs^* = \prod_{\fp \mid \fs} \PGL_2(K_\fp).
$$
Observe that $G_\fs^*$ has finite index in $G_\fs$. We call $G_\fs^*$ the \emph{linear subgroup} of $G_\fs$. There is a proper, vertex-transitive action of $G_\fs$ on
$$
    T_\fs = \prod_{\fp \mid \fs} T_\fp.
$$
We call a discrete group $\Gamma \leq G_\fs$ a \emph{lattice} if the quotient
$$
    X_\Gamma = \Gamma\backslash T_\fs
$$
is compact. Note that if $\Gamma$ is torsion-free, then $X_\Gamma$ is a non-positively curved cube complex. Let
$$
    \Gamma^* = \Gamma \cap G_\fs^*,
$$
which is normal and of finite index in $\Gamma$. We call $\Gamma$ \emph{reducible} if $\fs = \fs_0\fs_1$ with $\fs_i$ non-trivial for each $i$, and there exist lattices $\Lambda_i \leq G_{\fs_i}^*$ and a finite-index subgroup $\Lambda \leq \Gamma$ such that
$$
    \Lambda = \Lambda_0 \times \Lambda_1.
$$
Otherwise, we call $\Gamma$ \emph{irreducible}. Observe that $\Gamma$ is reducible if and only if $X_\Gamma$ is finitely covered by a product of non-positively curved cube complexes of positive dimension.

\subsection{Linear arithmetic groups}
\label{sec:linear-groups}

Let $B/K$ be a totally definite quaternion algebra with discriminant $\fd$ coprime to $\fs$, so $B$ splits at each $\fp \mid \fs$. Completion at each $\fp \mid \fs$ induces an embedding 
$$
    \rho_{B,\fs}:B^\times/K^\times \longrightarrow G_\fs^*. 
$$
Let $\fO$ be an order of $B$ and let 
$$
    \Phi_{\fO,\fs}^1 = \fO[1/\fs]^1K^\times/K^\times \cong \fO[1/\fs]^1 / \{\pm 1\}.
$$
It follows from Borel and Harish-Chandra \cite{borel-and-harish-chandra-1962-arithmetic-subgroups-of-algebraic-groups} that $\rho_{B,\fs}(\Phi_{\fO,\fs}^1)$ is an irreducible lattice in $G_\fs^*$. If $\fA \subseteq \fB$  are orders of $B$, then $\Phi_{\fA,\fs}^1 \leq \Phi_{\fB,\fs}^1$ with finite index. In particular, $B$ determines a commensurability class of lattices $\dC^*(B,\fs)$ in $G_\fs^*$. 

We now apply the \emph{Prasad volume formula} to find the \emph{orbifold Euler characteristic} of $\Phi_{\fO,\fs}^1$ for hereditary $\fO$. We sketch a proof for the formula specialized to our setting. Our proof follows Prasad \cite{prasad-1989-volumes-of-s-arithmetic-quotients-of-semi-simple-groups} and Vign\'eras \cite{vigneras-1980-arithmetique-des-algebres-de-quaternions}, and we recommend both sources for a more complete treatment. See also Stix--Vdovina \cite{stix-and-vdovina-2017-simply-transitive-quaternionic-lattices-of-rank-2-over-f-q-t-and-a-non-classical-fake-quadric}. First, we recall the standard fact that if $\Gamma$ is a lattice and $\Lambda \leq \Gamma$ is a finite-index subgroup then $\chi(\Lambda) = [\Gamma:\Lambda]\cdot\chi(\Gamma)$.

\begin{thm}\label{thm:prasad-volume} \cite{prasad-1989-volumes-of-s-arithmetic-quotients-of-semi-simple-groups}
    Let $\fO$ be a hereditary $\fo_K$-order of level $\fn$ coprime to $\fs$. Then
    $$
        \chi(\Phi^1_{\fO,\fs})
        =(-1)^{\omega(\fs)}2^{1-n}|\zeta_K(-1)|
        \prod_{\fp\mid\fs}(N\fp-1)
        \prod_{\fq\mid\fd}(N\fq-1)
        \prod_{\fl\mid\fn}(N\fl+1),
    $$
    where $\zeta_K$ is the Dedekind zeta function of $K$.
\end{thm}

\begin{proof}
    Let $\Delta = \fO[1/\fs]^1$, which is a \emph{principal arithmetic subgroup} of 
    $$
        G_\fs^1 = \prod_{\fp \mid \fs} \SL_2(K_\fp).
    $$
    For the measure $\mu$ normalized as in \cite[Section 3.6]{prasad-1989-volumes-of-s-arithmetic-quotients-of-semi-simple-groups}, it follows from \cite[Theorem 3.7]{prasad-1989-volumes-of-s-arithmetic-quotients-of-semi-simple-groups} that 
    $$
        \mu(G_\fs^1/\Delta) = 2^{-n}|\zeta_K(-1)|\prod_{\fp \mid \fs}(N\fp +1)\prod_{\fq \mid \fd}(N\fq - 1)\prod_{\fl\mid\fn}(N\fl + 1).
    $$
    Passing to the Euler--Poincar\'e measure as in \cite[Section 4]{borel-and-prasad-1989-finiteness-theorems-for-discrete-subgroups-of-bounded-covolume-in-semi-simple-groups}, we obtain 
    $$
        \chi(\Delta) = (-1)^{\omega(\fs)}\mu(G_\fs^1/\Delta)\prod_{\fp \mid \fs} \frac{N\fp - 1}{N\fp+1}.
    $$
    We obtain our desired formula since $\chi(\Phi_{\fO,\fs}^1) = 2\chi(\Delta)$.
\end{proof}

For the rest of this section, we assume that $\fO \subseteq B$ is a hereditary order of level $\fn$, and let
$$
    \Phi_{\fO,\fs} = \{\alpha \in B^\times : {}^\alpha\fO[1/\fs] = \fO[1/\fs]\}/K^\times. 
$$
Let $K_{>0}$ be the totally positive elements of $K$. The reduced norm induces a homomorphism 
$$
    \Nrd: \Phi_{\fO,\fs} \longrightarrow K_{>0}^\times/K^{\times 2}.
$$
The next lemma shows that $\Phi_{\fO,\fs}^1$ is the kernel of this map, so $\Phi_{\fO,\fs}^1$ has finite index in $\Phi_{\fO,\fs}$ and $\rho_{B,\fs}(\Phi_{\fO,\fs}) \in \dC^*(B,\fs)$.

\begin{lem}\label{lem:norm-kernel}
    $\Phi_{\fO,\fs}^1$ is the kernel of $\Nrd$. In particular,
    $$
        \chi(\Phi_{\fO,\fs}) = \chi(\Phi_{\fO,\fs}^1)/|\Nrd(\Phi_{\fO,\fs})|.
    $$
\end{lem}

\begin{proof}
    It is clear that $\Phi_{\fO,\fs}^1$ is contained in the kernel of $\Nrd$. Conversely, an element of the kernel has a representative $\alpha \in B^\times$ with $\Nrd(\alpha) = 1$, so $\alpha \in \fO[1/\fs]^1$ by Lemma \ref{lem:global-units}.
\end{proof}

\begin{thm}\cite[Theorem 1.4]{borel-and-prasad-1989-finiteness-theorems-for-discrete-subgroups-of-bounded-covolume-in-semi-simple-groups}
    \label{thm:phi-maximal}
    Let $\Gamma \in \dC^*(B,\fs)$. Then $\Gamma \leq \rho_{B,\fs}(\Phi_{\fO,\fs})$ for some hereditary order $\fO \subseteq B$ of level $\fn$ coprime to $\fs$.
\end{thm}

We now give a precise formula for $|\Nrd(\Phi_{\fO,\fs})|$. To this end, we introduce the \emph{$\ft$-unit groups} 
$$
    \Phi_{\fO,\ft}^\times = 
    \fO[1/\ft]^\times K^\times/K^\times \cong
    \fO[1/\ft]^\times/\fo_K[1/\ft]^\times, \qquad \fs \mid \ft \mid \fs\fd.
$$

\begin{lem}\label{lem:t-units}
    For an ideal $\ft = \fs\ft'$ with $\ft' \mid \fd$, we have
    $$
        \Phi_{\fO,\ft}^\times 
        = \ker(\Phi_{\fO,\fs} \longrightarrow K_{>0}^\times/K^{\times 2}\fo_K[1/\ft]^\times_{>0}), \qquad \fo_K[1/\ft]^\times_{>0} = K^\times_{>0} \cap \fo_K[1/\ft]^\times.
    $$
\end{lem}

\begin{proof}
    We first show that $\Phi_{\fO,\ft}^\times \subseteq \Phi_{\fO,\fs}$. Let $\alpha$ represent an element of $\Phi_{\fO,\ft}^\times$. Then $\alpha \in \fO_\fp^\times$ for all $\fp \nmid\ft$. Furthermore, $\alpha$ normalizes the unique maximal order $\fO_\fp$ of $B_\fp$ for each $\fp \mid \fd$. Thus, $\alpha$ normalizes $\fO[1/\fs]$ by the local-global principle (Lemma \ref{lem:local-global}). We get the desired inclusion, and it is straightforward to see that $[\alpha]$ is in the given kernel.

    Now let $\alpha$ represent an element of the kernel, and $R = \fo_K[1/\ft]$. Then we can pick $\alpha$ such that $\Nrd(\alpha) \in R^\times$. Since $\alpha$ normalizes $\fO[1/\ft]$, we are done by Lemma \ref{lem:global-units}.
\end{proof}

We now compute the indices of 
$$
    \Phi_{\fO,\fs}^1 \subseteq \Phi_{\fO,\fs}^\times \subseteq \Phi_{\fO,\fs}.
$$
Let $R = \fo_K[1/\fs]$ and let $h(R)$ and $h^+(R)$ be the orders of the class group $\Cl(R)$ and the narrow class group $\Cl^+(R)$ respectively. First, it follows from \cite[Corollary 28.9.21]{voight-2021-quaternion-algebras} that 
\begin{equation}
    \label{eq:normalizer-unit-index}
    [\Phi_{\fO,\fs}:\Phi^\times_{\fO,\fs}]=|\bigl\{\fa\mid\fd\fn:[\fa]\in\Cl^+(R)^2\bigr\}|\cdot|\im\bigl(\Cl^+(R)[2]\to\Cl(R)\bigr)|,
\end{equation}
where $\Cl^+(R)[2]$ is the 2-torsion subgroup of $\Cl^+(R)$. Thus, to compute $|\Nrd(\Phi_{\fO,\fs})|$, we need only find $[\Phi^\times_{\fO,\fs}:\Phi^1_{\fO,\fs}]$.

\begin{prop}\label{prop:norm-index}
    We have
    \begin{equation}
        \label{eq:norm-one-index}
        [\Phi^\times_{\fO,\fs}:\Phi^1_{\fO,\fs}]=2^{\omega(\fs)}\frac{h^+(R)}{h(R)}.
    \end{equation}
    In particular, $|\Nrd(\Phi_{\fO,\fs})| = 2^{\omega(\fd\fn\fs)}$ if $h^+(R)=1$.
\end{prop}

\begin{proof}
    Since $B$ splits at each $\fp \mid \fs$, it follows from \cite[Corollary 31.1.11]{voight-2021-quaternion-algebras} that 
    $$
        \Nrd(\fO[1/\fs]^\times) = R_{>0}^\times.
    $$
    Therefore, $\Phi_{\fO,\fs}^\times/\Phi_{\fO,\fs}^1 \cong R_{>0}^\times/R^{\times 2}$. By Dirichlet's $S$-unit theorem \cite[Chapter I, Corollary 11.7]{neukirch-1999-algebraic-number-theory}, we have
    \begin{equation}
        \label{eq:unit-squares-index}
        [R^\times:R^{\times 2}] = 2^{n+\omega(\fs)}.
    \end{equation}
    By \cite[Chapter VI, Proposition 1.11]{neukirch-1999-algebraic-number-theory}, we also get 
    \begin{equation}
        \label{eq:narrow-class-ratio}
        \frac{h^+(R)}{h(R)} = \frac{2^n}{[R^\times:R^\times_{>0}]}.
    \end{equation}
    Combining \eqref{eq:unit-squares-index} and \eqref{eq:narrow-class-ratio}, we get 
    $$
        [R^\times_{>0}:R^{\times 2}] = \frac{[R^\times : R^{\times 2}]}{[R^\times : R^{\times}_{>0}]} = \frac{2^{n+\omega(\fs)}}{2^n h(R)/h^+(R)} = 2^{\omega(\fs)}\frac{h^+(R)}{h(R)}.
    $$
    We thus get \eqref{eq:norm-one-index}. The order of $\Nrd(\Phi_{\fO,\fs})$ when $h^+(R) = 1$ now follows from \eqref{eq:normalizer-unit-index} and \eqref{eq:norm-one-index}.
\end{proof}

\subsection{Semilinear arithmetic groups}
\label{sec:semilinear-groups}

We now define lattices inside of $G_\fs$ coming from quaternions. Let $K_0 \subseteq K$ be the fixed field of $\Aut(K) = \Aut_\bfQ(K)$. If $\sigma \in \Gal(K/K_0)$, then $\sigma$ induces an isomorphism of $K_0$-algebras,
$$
    \sigma_*: B \to B', \qquad \fd(B') = \sigma\fd(B).
$$
If $\fd(B') = \fd(B)$, then there is a $K$-isomorphism $B' \cong B$. Thus, for $\sigma\fd(B) = \fd(B)$, we may assume that $\sigma_*$ is a $K_0$-automorphism of $B$. It is unique up to $\Aut_K(B)$. We often drop the subscript and denote the induced automorphism by $\sigma$. For ideals $\fa_1,\dots,\fa_r$ of $\fo_K$, let 
$$
    \Sigma(\fa_1,\dots,\fa_r) = \{\sigma \in \Gal(K/K_0) : \sigma\fa_i = \fa_i\}.
$$
For $\alpha \in \Aut_{K_0}(B)$, the restriction of $\alpha$ to the center $K$ of $B$ yields an element of $\Sigma(\fd)$, and every $\sigma \in \Sigma(\fd)$ arises this way. In particular, we have a short exact sequence
$$
     1 \longrightarrow \Aut_K(B) \longrightarrow \Aut_{K_0}(B) \longrightarrow \Sigma(\fd) \longrightarrow 1.
$$
For $\sigma \in \Sigma(\fd)$ and $\alpha \in \Aut_K(B)$, we write $(\alpha,\sigma)$ for $\alpha \circ \sigma$. Every element of $\Aut_{K_0}(B)$ is uniquely of this form. Set
$$
    \Aut_{K_0}(B,\fs) = 
    \{g \in \Aut_{K_0}(B) : g|_{K} \in \Sigma(\fd,\fs)\}.
$$
Then we get another short exact sequence 
$$
    1 \longrightarrow \Aut_K(B) \longrightarrow \Aut_{K_0}(B,\fs) \longrightarrow \Sigma(\fd,\fs) \longrightarrow 1.
$$
By the Skolem--Noether theorem \cite[Theorem 7.7.1]{voight-2021-quaternion-algebras}, the conjugation action of $B^\times$ on $B$ induces an isomorphism $B^\times/K^\times \cong \Aut_K(B)$. Under this identification, we extend $\rho_{B,\fs}$ to the embedding 
$$
    \rho_{B,\fs} : \Aut_{K_0}(B,\fs) \longrightarrow G_\fs.
$$
We have 
$$
    \rho_{B,\fs}^{-1}(G_\fs^*) = \Aut_K(B).
$$
In particular, if $\Gamma \leq \rho_{B,\fs}(\Aut_{K_0}(B,\fs))$, then restriction to $K$ gives a homomorphism 
$$
    \Sigma: \Gamma \to \Sigma(\fd,\fs)
$$
with kernel $\Gamma^*$.

Let $\dC(B,\fs) \supseteq \dC^*(B,\fs)$ be the commensurability class of lattices in $G_\fs$ determined by $B$. It follows from Margulis' \emph{arithmeticity theorem} \cite{margulis-1991-discrete-subgroups-of-semisimple-lie-groups} that if $\omega(\fs) \geq 2$ and $\Gamma \leq G_\fs$ is an irreducible lattice, then there exists a totally real number field $F$, an ideal $\ft$ of $\fo_F$, and a totally definite quaternion algebra $B/F$ with discriminant coprime to $\ft$ such that $G_\fs \cong G_\ft$ and $\Gamma \in \dC(B,\ft)$.

We now identify the commensurator of lattices in $\dC(B,\fs)$. The next proposition is a $p$-adic analogue of \cite[Theorem 3.3]{linowitz-stover-and-voight-2015-commensurability-classes-of-fake-quadrics}, and its proof follows theirs.

\begin{prop}
    \label{prop:commensurator}
    Let $\Gamma \in \dC(B,\fs)$. Then
    $$
        \Comm_{G_\fs}(\Gamma) = \rho_{B,\fs}(\Aut_{K_0}(B,\fs)).
    $$
\end{prop}

\begin{proof}
    Let $B_\fs = \prod_{\fp \mid \fs} B_\fp$ and identify $K$ and $B$ in $B_\fs$ under the diagonal embedding. Let $\Aut(B_\fs)$ be the continuous, ring automorphisms of $B_\fs$. If $g \in \Aut_{K_0}(B,\fs)$, then there is a unique, continuous extension of $g$ to $B_\fs$ since $B$ is dense in $B_\fs$ \cite[Proposition 28.7.3]{voight-2021-quaternion-algebras}. Furthermore, if $g \in \Aut(B_\fs)$ with $g(B) = B$, then $g|_K$ permutes the factors of $B_\fs$, so $g|_K \in \Sigma(\fd(B),\fs)$.
    
    There is an isomorphism $\Aut(B_\fs) \cong G_\fs$ which extends $\rho_{B,\fs}$. We will work with $\Aut(B_\fs)$.
    
    It suffices to find the commensurator for $\Phi_{\fO,\fs}^1$ in $\Aut(B_\fs)$, where $\fO$ is a maximal order of $B$. It is clear that $\Aut_{K_0}(B,\fs)$ commensurates $\Phi_{\fO,\fs}^1$. We show the reverse inclusion. Let $g$ commensurate $\Phi_{\fO,\fs}^1$ in $\Aut(B_\fs)$, and let $\Lambda$ be the preimage of 
    $$
        \Phi_{\fO,\fs}^1 \cap g^{-1}\Phi_{\fO,\fs}^1 g
    $$
    in $\fO[1/\fs]^1 \subseteq B_\fs$. Then $\Lambda$ is finitely generated, and $\Lambda^{(2)} = \lag \lambda^2 \mid \lambda \in \Lambda\rag$ has finite index in $\Lambda$. If $\lambda \in \Lambda$, then $g(\lambda^2) \in \fO[1/\fs]^1$. It follows from \cite[Section 3]{maclachlan-and-reid-2003-the-arithmetic-of-hyperbolic-3-manifolds} that $g(B) = B$. Thus, $g \in \Aut_{K_0}(B,\fs)$.
\end{proof}

Going forward, we drop the map $\rho_{B,\fs}$ and identify the relevant lattices with their image in $G_\fs$. For a hereditary order $\fO \subseteq B$ of level $\fn$ coprime to $\fs$, let
$$
    \Psi_{\fO,\fs} = \{g \in \Aut_{K_0}(B) : g(\fO[1/\fs]) = \fO[1/\fs]\}.
$$
Then $\Psi_{\fO,\fs} \cap G_\fs^* = \Phi_{\fO,\fs}$ and
$$
    [\Psi_{\fO,\fs}:\Phi_{\fO,\fs}] = |\{\sigma \in \Sigma(\fd,\fs) : \sigma(\fO[1/\fs]) \text{ is conjugate to } \fO[1/\fs]\}| < \infty.
$$ 
In particular, $\Psi_{\fO,\fs} \in \dC(B,\fs)$.  We now extend Theorem \ref{thm:phi-maximal} to the semilinear setting, adapting the proof of Borel and Prasad.

\begin{thm}
    \label{thm:psi-maximal}
    Let $\Gamma \in \dC(B,\fs)$. Then $\Gamma \leq \Psi_{\fO,\fs}$ for some hereditary order $\fO$ of level $\fn$ coprime to $\fs$.
\end{thm}

\begin{proof}
    Let $R=\fo_K[1/\fs]$ and let $\fM$ be a maximal $R$-order. By Proposition \ref{prop:commensurator}, $\Gamma \leq \Aut_{K_0}(B,\fs)$, so $\Gamma$ acts on the set of $R$-orders in $B$. Since $\Gamma \cap \Phi_{\fM,\fs}$ has finite index in $\Gamma$ and fixes $\fM$, the orbit
    $$
        \dM = \Gamma \cdot \fM = \{\fM_1,\dots,\fM_m\}
    $$
    is finite. For $\fp \nmid \fd\fs$, the vertices $\dM_\fp = \{\fM_{1,\fp},\dots,\fM_{m,\fp}\}$ of $T_\fp$ have a unique circumcenter, which is a vertex or the midpoint of an edge. Let $\fO_\fp$ be the corresponding maximal order or Eichler order of level $\fp$. For $\fp \mid \fd$, let $\fO_\fp$ be the unique maximal order of $B_\fp$. For all but finitely many $\fp$, $\dM_\fp$ consists of exactly one order $\fM_\fp$. Thus, by Lemma \ref{lem:local-global}, there exists a hereditary $R$-order $\fO$ with localizations $\fO_\fp$.

    Let $g = (\alpha,\sigma) \in \Gamma$. Recall from Proposition \ref{prop:commensurator} that $\Gamma \subseteq \Aut_{K_0}(B,\fs)$. We now show that $g(\fO_\fp) = \fO_{\sigma\fp}$ for all $\fp \nmid \fs$. If $\fp \mid \fd$, then $\sigma\fp \mid \fd$, and $g(\fO_\fp) = \fO_{\sigma\fp}$ since $B_{\sigma\fp}$ has a unique maximal order. Now suppose that $\fp \nmid \fs\fd$. Since $g$ permutes $\dM$, $g$ induces a tree isomorphism $T_\fp \to T_{\sigma\fp}$ carrying $\dM_\fp$ to $\dM_{\sigma\fp}$. Hence it takes the circumcenter, $\fO_\fp$, of $\dM_\fp$ to the circumcenter, $\fO_{\sigma\fp}$, of $\dM_{\sigma\fp}$. Thus $g(\fO)$ and $\fO$ have the same localizations, so $g(\fO) = \fO$ by Lemma \ref{lem:local-global}, and $\Gamma \leq \Psi_{\fO,\fs}$.
\end{proof}

We will also use the following lemma.

\begin{lem}
    \label{lem:conjugacy-in-psi}
    Let $\fM \subseteq B$ be a maximal order, and let $\Gamma, \Gamma' \leq \Psi_{\fM,\fs}$ have index $k$, where $k < N\fl + 1$ for every prime $\fl \nmid \fd\fs$. If $g\Gamma g^{-1} = \Gamma'$ for $g \in G_\fs$, then $g \in \Psi_{\fM,\fs}$.
\end{lem}

\begin{proof}
    First observe that $g \in \Comm_{G_\fs}(\Psi_{\fM,\fs}) = \Aut_{K_0}(B,\fs)$. Thus, $\Gamma'$ fixes the order $g(\fM)[1/\fs]$. Therefore, the $\Psi_{\fM,\fs}$-orbit of $g(\fM)[1/\fs]$ has size at most $k$. Suppose that $g(\fM)_\fl \neq \fM_\fl$ for some $\fl \nmid \fd\fs$, so 
    $$
        d = d(g(\fM)_\fl,\fM_\fl) > 0.
    $$
    By strong approximation, $\fM[1/\fs]^1$ is dense in $\fM^1_\fl$, so $\fM[1/\fs]^1$ acts transitively on the orders of distance $d$ from $\fM_\fl$. Therefore, the $\Psi_{\fM,\fs}$-orbit of $g(\fM)[1/\fs]$ has size at least $N\fl + 1 > k$. This is a contradiction. Thus, $\fM_\fl = g(\fM)_\fl$ for all $\fl \nmid \fd\fs$, so $\fM[1/\fs] = g(\fM)[1/\fs]$.
\end{proof}

\subsection{Congruence subgroups}
\label{sec:congruence-subgroups}

Let $B$, $K$, and $\fs$ be as before, and fix an ideal $\fa \subseteq \fo_K$ coprime to $\fs$. Let $\fO \subseteq B$ be a hereditary order of level $\fn$ coprime to $\fs$. Observe that 
$$
    \fO/\fa\fO \cong \fO[1/\fs]/\fa\fO[1/\fs].
$$
Let $\Gamma \leq \Psi_{\fO,\fs}$, and let $\dA = \Sigma(\Gamma) \cdot \fa$. Observe that every ideal in $\dA$ is coprime to $\fs$, and $\Gamma$ acts on the finite set
$$
    \coprod_{\fa' \in \dA} \fO/\fa'\fO.
$$
Let $\Gamma_\fa$ be the stabilizer of $\fO/\fa\fO$, which has finite index in $\Gamma$, and set
$$
    \pi_\fa = \pi_{\Gamma,\fa}: \Gamma_\fa \longrightarrow \Aut(\fO/\fa\fO).
$$
We call $\Lambda \leq \Gamma$ a \emph{congruence subgroup of level} $\fa$ of $\Gamma$ if $\Lambda$ contains the kernel $\Gamma(\fa) = \ker(\pi_\fa)$. For example, suppose that $\Gamma = \Phi_{\fO,\ft}^\times$ for some ideal $\fs \mid \ft \mid \fs\fd$ coprime to $\fa$. Then $\Gamma_\fa = \Gamma$ and $\pi_\fa$ is induced by reduction modulo $\fa$.

\section{Computing fundamental domains}

\label{sec:fundamental-domains}

Throughout this section, $K$ is a totally real number field, $B$ is a totally definite quaternion algebra of discriminant $\fd$, $\fs$ is a non-trivial square-free ideal coprime to $\fd$, and $\fO \subseteq B$ is a hereditary order of level $\fn$ coprime to $\fs$. Furthermore, $\fa$ is always an ideal dividing $\fs$. In this section, we describe a general procedure for computing fundamental domains of $\Phi_{\fO,\fs}$ and $\Psi_{\fO,\fs}$ in $T_\fs$. We implement this algorithm in Magma \cite{bosma-cannon-and-playoust-1997-the-magma-algebra-system-i-the-user-language}, which supports a wide array of computations with quaternion algebras \cite{voight-2005-quadratic-forms-and-quaternion-algebras-algorithms-and-arithmetic,kirschmer-and-voight-2010-algorithmic-enumeration-of-ideal-classes-for-quaternion-orders,cannon-bosma-fieker-and-steel-2026-handbook-of-magma-functions}. Our Magma implementation of this algorithm, together with the other code used in this paper, will be made available at \url{https://github.com/jonah-mendel}.

Below, we point out the key places where we use Magma. Our method is similar to that of Masdeu--Torrents \cite{masdeu-and-torrents-2025-fundamental-domains-for-quaternionic-s-arithmetic-groups-over-totally-real-fields}, except that we treat semilinear lattices and allow several primes to be inverted.

Observe that $T_\fs$ is a CAT(0) cube complex. We denote a cube of $T_\fs$ by tuples of orders 
$$
    \fc = (\fA_\fp)_{\fp \mid \fs},
$$
where $\fA_\fp \subseteq M_2(K_\fp)$ is hereditary, and the dimension of $\fc$ is 
$$
    \dim(\fc) = |\{\fp \mid \fs : \fn(\fA_\fp) = \fp\}|.
$$
Observe that $\fc = (\fA_\fp)$ is a \emph{face} of $\fc' = (\fA'_\fp)$ if and only if $\fA'_\fp \subseteq \fA_\fp$ for each $\fp \mid \fs$, and the \emph{codimension} of $\fc$ in $\fc'$ is $|\{\fp \mid \fs : \fn(\fA_\fp')/\fn(\fA_\fp) = \fp\}|$. For $\fa \mid \fs$, let 
$$
    T_\fs(\fa) = \{(\fA_\fp)_{\fp \mid \fs} : \fn(\fA_\fp) = \fp \text{ if and only if } \fp \mid \fa\},
$$
so the $k$-cubes of $T_\fs$ are the union of all $T_\fs(\fa)$ with $\omega(\fa) = k$.

We now parametrize the cubes of $T_\fs$ by global orders. Suppose that $\fA \subseteq B$ is hereditary of level $\fn\fa$ for $\fa \mid \fs$ and satisfies $\fA[1/\fs] = \fO[1/\fs]$. Then $\fA_\fp = \fO_\fp$ for all $\fp \nmid \fs$. By the local-global principle (Lemma \ref{lem:local-global}), there is a $\Phi_{\fO,\fs}$-equivariant bijection
\begin{equation}
    \label{eq:cube-parametrization}
    \{\fA \subseteq B \text{ hereditary of level } \fn\fa : \fA[1/\fs] = \fO[1/\fs]\} \to T_\fs(\fa), \quad \fA \mapsto (\fA_\fp)_{\fp \mid \fs},
\end{equation}
where $\Phi_{\fO,\fs}$ acts on the left-hand set via conjugation. From now on, we refer to elements of $T_\fs(\fa)$ by the corresponding global order. Note that $\fA$ represents a face of $\fB$ if and only if $\fB \subseteq \fA$, and we say that $\fA$ is a codimension $k$ \emph{superorder} of $\fB$, where $k = \omega(\fn(\fB)/\fn(\fA))$.

The action of $\Phi_{\fO,\fs}$ on $T_\fs$ extends to an action of $\Psi_{\fO,\fs}$ which no longer preserves the sets $T_\fs(\fa)$, but still sends $k$-cubes to $k$-cubes. Let $\fA \in T_\fs(\fa)$. We want to compute the $\Psi_{\fO,\fs}$-stabilizer $\Psi_\fA$ of $\fA$. First, the $\Phi_{\fO,\fs}$-stabilizer of $\fA$ is 
$$
    \Phi_{\fA} = N_{B^\times}(\fA)/K^\times.
$$
One can compute the finite group $\Phi_\fA$ directly in Magma. Then, we run through the elements $\sigma \in \Sigma(\fd,\fs,\fa)$ and try to find an element $\alpha_\sigma\in B^\times$ such that ${}^{\alpha_\sigma}\sigma(\fA) = \fA$. Such an element need not exist. Let $\Sigma(\fA,\fs)$ be the collection of $\sigma$ such that $\alpha_\sigma$ exists. Observe that for each $\sigma, \tau \in \Sigma(\fA,\fs)$, there exists $\beta_{\sigma,\tau} \in \Phi_\fA$ such that 
$$
    (\alpha_\sigma,\sigma)(\alpha_\tau,\tau) = \beta_{\sigma,\tau}(\alpha_{\sigma\tau},\sigma\tau).
$$
We get an isomorphism
\begin{equation}
    \label{eq:cube-stabilizer}
    \Psi_\fA \cong \bigl\langle \Phi_\fA,\ t_\sigma \bigm| t_\sigma t_\tau = \beta_{\sigma,\tau} t_{\sigma\tau},\ t_\sigma \phi t_\sigma^{-1} = \alpha_\sigma \sigma(\phi) \alpha_\sigma^{-1} \bigr\rangle,
\end{equation}
where $\sigma,\tau$ range over $\Sigma(\fA,\fs)$.

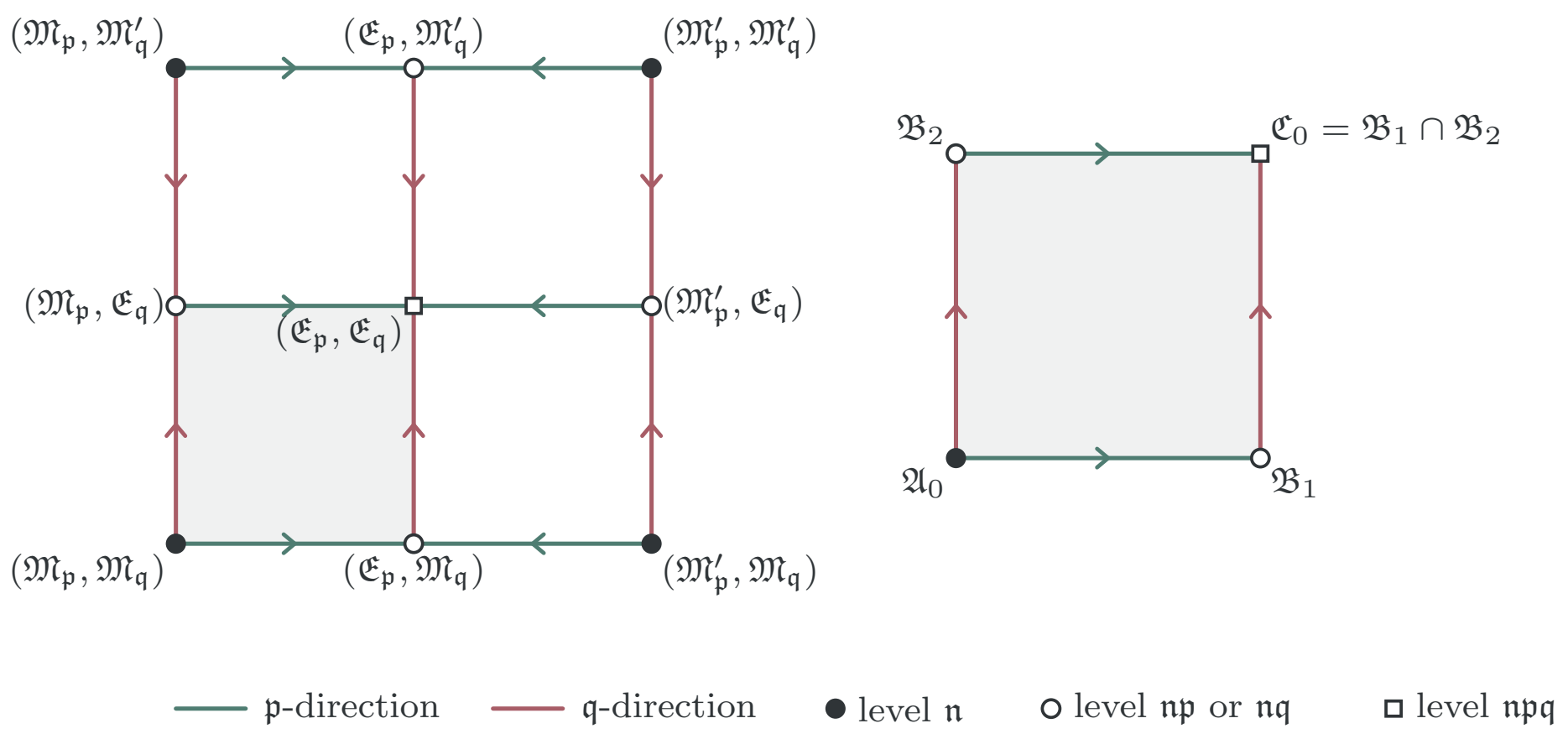
\begin{figure}[ht!]
     \centering
     \definecolor{SoftInk}{HTML}{2F3437}
\definecolor{PrimeOne}{HTML}{4F7D72}
\definecolor{PrimeTwo}{HTML}{A85D67}
\begin{tikzpicture}[x=0.75cm,y=0.75cm,
  every node/.style={font=\scriptsize,text=SoftInk},
  side/.style={line width=0.95pt,line cap=round},
  tip/.style={line width=0.9pt,line cap=round,line join=round},
  cell/.style={fill=SoftInk!7},
  label/.style={inner sep=2pt},
  vzero/.style={circle,fill=SoftInk,inner sep=1.6pt},
  vone/.style={circle,fill=white,draw=SoftInk,line width=0.7pt,inner sep=1.4pt},
  vtwo/.style={rectangle,fill=white,draw=SoftInk,line width=0.7pt,inner sep=1.7pt}]
  \fill[cell] (0,0) rectangle (2.5,2.5);
  \draw[side,draw=PrimeOne] (0,0) -- (2.5,0);
  \draw[tip,PrimeOne,shift={(1.25,0)},rotate=0] (-2.5pt,2.1pt) -- (0,0) -- (-2.5pt,-2.1pt);
  \draw[side,draw=PrimeOne] (5,0) -- (2.5,0);
  \draw[tip,PrimeOne,shift={(3.75,0)},rotate=180] (-2.5pt,2.1pt) -- (0,0) -- (-2.5pt,-2.1pt);
  \draw[side,draw=PrimeOne] (0,2.5) -- (2.5,2.5);
  \draw[tip,PrimeOne,shift={(1.25,2.5)},rotate=0] (-2.5pt,2.1pt) -- (0,0) -- (-2.5pt,-2.1pt);
  \draw[side,draw=PrimeOne] (5,2.5) -- (2.5,2.5);
  \draw[tip,PrimeOne,shift={(3.75,2.5)},rotate=180] (-2.5pt,2.1pt) -- (0,0) -- (-2.5pt,-2.1pt);
  \draw[side,draw=PrimeOne] (0,5) -- (2.5,5);
  \draw[tip,PrimeOne,shift={(1.25,5)},rotate=0] (-2.5pt,2.1pt) -- (0,0) -- (-2.5pt,-2.1pt);
  \draw[side,draw=PrimeOne] (5,5) -- (2.5,5);
  \draw[tip,PrimeOne,shift={(3.75,5)},rotate=180] (-2.5pt,2.1pt) -- (0,0) -- (-2.5pt,-2.1pt);
  \draw[side,draw=PrimeTwo] (0,0) -- (0,2.5);
  \draw[tip,PrimeTwo,shift={(0,1.25)},rotate=90] (-2.5pt,2.1pt) -- (0,0) -- (-2.5pt,-2.1pt);
  \draw[side,draw=PrimeTwo] (0,5) -- (0,2.5);
  \draw[tip,PrimeTwo,shift={(0,3.75)},rotate=-90] (-2.5pt,2.1pt) -- (0,0) -- (-2.5pt,-2.1pt);
  \draw[side,draw=PrimeTwo] (2.5,0) -- (2.5,2.5);
  \draw[tip,PrimeTwo,shift={(2.5,1.25)},rotate=90] (-2.5pt,2.1pt) -- (0,0) -- (-2.5pt,-2.1pt);
  \draw[side,draw=PrimeTwo] (2.5,5) -- (2.5,2.5);
  \draw[tip,PrimeTwo,shift={(2.5,3.75)},rotate=-90] (-2.5pt,2.1pt) -- (0,0) -- (-2.5pt,-2.1pt);
  \draw[side,draw=PrimeTwo] (5,0) -- (5,2.5);
  \draw[tip,PrimeTwo,shift={(5,1.25)},rotate=90] (-2.5pt,2.1pt) -- (0,0) -- (-2.5pt,-2.1pt);
  \draw[side,draw=PrimeTwo] (5,5) -- (5,2.5);
  \draw[tip,PrimeTwo,shift={(5,3.75)},rotate=-90] (-2.5pt,2.1pt) -- (0,0) -- (-2.5pt,-2.1pt);
  \node[vzero] at (0,0) {};
  \node[label,anchor=north east] at (0,0) {$(\fM_\fp,\fM_\fq)$};
  \node[vzero] at (5,0) {};
  \node[label,anchor=north west] at (5,0) {$(\fM'_\fp,\fM_\fq)$};
  \node[vzero] at (0,5) {};
  \node[label,anchor=south east] at (0,5) {$(\fM_\fp,\fM'_\fq)$};
  \node[vzero] at (5,5) {};
  \node[label,anchor=south west] at (5,5) {$(\fM'_\fp,\fM'_\fq)$};
  \node[vone] at (2.5,0) {};
  \node[label,anchor=north] at (2.5,0) {$(\fE_\fp,\fM_\fq)$};
  \node[vone] at (2.5,5) {};
  \node[label,anchor=south] at (2.5,5) {$(\fE_\fp,\fM'_\fq)$};
  \node[vone] at (0,2.5) {};
  \node[label,anchor=east] at (0,2.5) {$(\fM_\fp,\fE_\fq)$};
  \node[vone] at (5,2.5) {};
  \node[label,anchor=west] at (5,2.5) {$(\fM'_\fp,\fE_\fq)$};
  \node[vtwo] at (2.5,2.5) {};
  \node[label,anchor=north east] at (2.5,2.5) {$(\fE_\fp,\fE_\fq)$};
  \fill[cell] (8.2,0.9) rectangle (11.4,4.1);
  \draw[side,draw=PrimeOne] (8.2,0.9) -- (11.4,0.9);
  \draw[tip,PrimeOne,shift={(9.8,0.9)},rotate=0] (-2.5pt,2.1pt) -- (0,0) -- (-2.5pt,-2.1pt);
  \draw[side,draw=PrimeOne] (8.2,4.1) -- (11.4,4.1);
  \draw[tip,PrimeOne,shift={(9.8,4.1)},rotate=0] (-2.5pt,2.1pt) -- (0,0) -- (-2.5pt,-2.1pt);
  \draw[side,draw=PrimeTwo] (8.2,0.9) -- (8.2,4.1);
  \draw[tip,PrimeTwo,shift={(8.2,2.5)},rotate=90] (-2.5pt,2.1pt) -- (0,0) -- (-2.5pt,-2.1pt);
  \draw[side,draw=PrimeTwo] (11.4,0.9) -- (11.4,4.1);
  \draw[tip,PrimeTwo,shift={(11.4,2.5)},rotate=90] (-2.5pt,2.1pt) -- (0,0) -- (-2.5pt,-2.1pt);
  \node[vzero] at (8.2,0.9) {};
  \node[label,anchor=north east] at (8.2,0.9) {$\fA_0$};
  \node[vone] at (11.4,0.9) {};
  \node[label,anchor=north west] at (11.4,0.9) {$\fB_1$};
  \node[vone] at (8.2,4.1) {};
  \node[label,anchor=south east] at (8.2,4.1) {$\fB_2$};
  \node[vtwo] at (11.4,4.1) {};
  \node[label,anchor=south west] at (11.4,4.1) {$\fC_0 = \fB_1 \cap \fB_2$};
  \begin{scope}[yshift=-1.3cm,x=1cm]
  \draw[side,draw=PrimeOne] (0,0) -- (0.55,0);
  \node[label,anchor=west] at (0.62,0) {$\fp$-direction};
  \draw[side,draw=PrimeTwo] (2.5,0) -- (3.05,0);
  \node[label,anchor=west] at (3.12,0) {$\fq$-direction};
  \node[vzero] at (5.2,0) {}; \node[label,anchor=west] at (5.3,0) {level $\fn$};
  \node[vone] at (6.9,0) {}; \node[label,anchor=west] at (7.0,0) {level $\fn\fp$ or $\fn\fq$};
  \node[vtwo] at (9.6,0) {}; \node[label,anchor=west] at (9.7,0) {level $\fn\fp\fq$};
  \end{scope}
\end{tikzpicture}
     \caption{The figure on the left shows a square of $T_\fp \times T_\fq$, and its subdivision in $C_\fs$, for $\fs = \fp\fq$. Here $\fM_\fp, \fM'_\fp$ are adjacent maximal orders of $M_2(K_\fp)$, $\fE_\fp = \fM_\fp \cap \fM'_\fp$, and $\fE_\fq = \fM_\fq \cap \fM'_\fq$. The figure on the right shows the global parametrization of the bottom-left cube.}
     \label{fig:cubical-subdivision}
\end{figure}

Let $C_\fs$ be the cube subdivision of $T_\fs$, so the vertices of $C_\fs$ are in bijection with the union of all $T_\fs(\fa)$ for $\fa \mid \fs$. There is an (oriented) edge from $\fA$ to $\fB$ if $\fA$ is a codimension 1 superorder of $\fB$, and there is a square with vertices $\fA_0$, $\fB_1$, $\fB_2$, and $\fC_0$ if $\fA_0$ is a codimension 1 superorder of $\fB_1$ and $\fB_2$, the $\fB_i$ have different levels, and $\fC_0 = \fB_1 \cap \fB_2$. See Figure \ref{fig:cubical-subdivision}.

The action of $\Psi_{\fO,\fs}$ on $T_\fs$ induces an orientation-preserving action on $C_\fs$. We now find representative 0-cells $\dV$, 1-cells $\dE$, and 2-cells $\dF$ for this action so that their union is a connected fundamental domain for $C_\fs^{(2)}$.

Let $\dV(k)$ be the representative $k$-cubes of $\dV$. We first construct $\dV(0)$ via a breadth-first search. We begin the search with $\dV(0) = \{\fO\}$. Suppose we are at an order $\fV \in \dV(0)$. We compute its neighbors in $T_\fs$,
$$
    \dN(\fV) = \{\fO_L(\fI) : \text{$\fI$ is an integral right $\fV$-ideal of norm $\fp$ for some $\fp \mid \fs$}\}.
$$
We then go through each $\fN \in \dN(\fV)$, and check whether $\fN$ is $\Psi_{\fO,\fs}$-equivalent to an order in $\dV(0)$. In practice, this means testing whether $\sigma(\fN)$ is conjugate to an order in $\dV(0)$ for some $\sigma \in \Sigma(\fd,\fs)$. If $\fN$ is not equivalent to any order in $\dV(0)$, then it is added to $\dV(0)$. The breadth-first search eventually terminates because the quotient is compact, and we obtain the set $\dV(0)$. Note that we can also keep track of edges, thus giving us a connected spanning tree $\dS \subseteq \dV(1)$.

For $k \geq 1$, every $k$-cube is in the $\Psi_{\fO,\fs}$-orbit of a cube containing some vertex of $\dV(0)$. Thus, we can go through the finitely many $k$-cubes containing a vertex in $\dV(0)$ and find representatives $\dV(k)$ for $1 \leq k \leq \omega(\fs)$ such that the union of the cubes, $\dV = \bigcup_k \dV(k)$, is a connected subset of $T_\fs$. Note that
$$
    \chi(\Psi_{\fO,\fs}) = \sum_{k=0}^{\omega(\fs)}(-1)^k \sum_{\fA \in \dV(k)} \frac{1}{|\Psi_\fA|},
$$
which can be compared with the value of $\chi(\Psi_{\fO,\fs})$ one obtains with the Prasad volume formula.

Now that $\dV$ is constructed, one goes through each order in $\dV$ and constructs the edge representatives $\dE$ of $C_\fs$. Another sweep constructs the face representatives $\dF$. Then, by applying a general procedure of Brown \cite{brown-1984-presentations-for-groups-acting-on-simply-connected-complexes}, we obtain a finite presentation for $\Psi_{\fO,\fs}$ and representatives in $\Aut_{K_0}(B,\fs)$ for each generator. Because we find the presentation from the cubical subdivision, we also have the image of $\Psi_\fA$ in the presentation for each cube $\fA \in \dV$, which we will need for constructing interesting examples in later sections.

\begin{lem}
    \label{lem:torsion-stabilizes-cubes}
    Let $g \in \Psi_{\fO,\fs}$ be finite order. Then $g$ lies in a conjugate of $\Psi_\fA$ for some cube $\fA \in \dV$.
\end{lem}

\begin{proof}
    By the Bruhat--Tits fixed-point lemma \cite[Corollary II.2.8]{bridson-and-haefliger-1999-metric-spaces-of-non-positive-curvature}, $g$ fixes a point in $T_\fs$. Since $\Psi_{\fO,\fs}$ acts cubically, $g$ necessarily fixes a cube $\fA'$ in $T_\fs$. By our assumptions, there exists $h \in \Psi_{\fO,\fs}$ such that $\fA = {}^h \fA' \in \dV$, so $h^{-1}\Psi_\fA h$ contains $g$.
\end{proof}

Now, let $\Gamma \leq \Psi_{\fO,\fs}$ be a finite-index subgroup.  By Lemma \ref{lem:torsion-stabilizes-cubes}, any torsion element of $\Gamma$ lies in the stabilizer of a translate of one of the cubes in $\dV$.  Thus, $\Gamma$ is torsion-free if and only if $\Gamma \cap g\Psi_\fA g^{-1}$ is trivial for every $\fA \in \dV$ and every $g \in \Psi_{\fO,\fs}$. In Magma, this is a finite check: for each $\fA \in \dV$ and each non-trivial $h \in \Psi_\fA$, we compute the action of $h$ on $\Gamma \backslash \Psi_{\fO,\fs}$ and check that it fixes no coset, since $h$ fixes $\Gamma g$ if and only if $ghg^{-1} \in \Gamma$.

\section{Non-archimedean uniformization and fake quadrics}

\label{sec:uniformization}

In this section, we derive Proposition \ref{prop:uniformization} from the work of Varshavsky \cite{varshavsky-1998-p-adic-uniformization-of-unitary-shimura-varieties,varshavsky-1998-p-adic-uniformization-of-unitary-shimura-varieties-ii}. Furthermore, we identify the commensurability class of the arithmetic group $\Lambda_\Gamma$.

Let $K$ be a totally real number field, let $p$ be a rational prime, and let $\fs$ be an eligible square-free ideal of $\fo_K$ over $p$. Fix a finite extension $L/\bfQ_p$ and embeddings 
$$
    i_\fp : K_\fp \hookrightarrow L, \qquad i_\fp(K_\fp) = i_\fq(K_\fq) \quad\text{when $K_\fp \cong K_\fq$ for $\fp,\fq \mid \fs$}.
$$
We identify each $p$-adic field $K_\fp$ with its image in $L$. Fix a completion $\bfC_p$ of the algebraic closure of $L$.

Fix a prime $\fp \mid \fs$. The \emph{Drinfeld upper half plane} is
$$
    \Omega_\fp = (\bfP^{1}_{K_\fp})^{\an} \smallsetminus \bfP^1(K_\fp).
$$
Following \cite{varshavsky-1998-p-adic-uniformization-of-unitary-shimura-varieties-ii}, we take $\Omega_\fp$ to be a Berkovich space \cite{berkovich-1990-spectral-theory-and-analytic-geometry-over-non-archimedean-fields}.  Berkovich shows that $\PGL_2(K_\fp)$ is the group of all $K_\fp$-analytic automorphisms of $\Omega_\fp$ \cite[Theorem 1]{berkovich-1995-the-automorphism-group-of-the-drinfeld-half-plane}. Now consider the Berkovich space defined over $L$,
$$
    \Omega_\fs = \prod_{\fp \mid \fs}(\Omega_\fp \otimes_{K_\fp} L).
$$
Permuting factors of $\Omega_\fs$ is $L$-analytic, so
$$
    H_\fs = \left(\prod_{\fp\mid\fs} \PGL_2(K_\fp)\right) \rtimes \Sym(\fs)
$$
acts on $\Omega_\fs$ via $L$-analytic automorphisms.

We are most interested in lattices $\Gamma \leq H_\fs$ that act freely on $\Omega_\fs(\bfC_p)$. Let $B/K$ be a totally definite quaternion algebra. We give a criterion to check whether an element $g \in \Aut_{K_0}(B,\fs)$ fixes a point. Recall that, if $\beta \in B$ and $\beta \not\in K$, then $K(\beta)$ is a quadratic extension of $K$.

\begin{lem}
    \label{lem:free-action}
    Let $g = (\alpha,\sigma) \in \Aut_{K_0}(B,\fs)$ such that $\rho_{B,\fs}(g) \in H_\fs$. Let $m$ be the order of $\sigma$ and $\beta \in B^\times$ represent $g^m$. Then $g$ fixes a point of $\Omega_\fs(\bfC_p)$ if and only if $\beta \in K^\times$ or no prime $\fp \mid \fs$ splits in $K(\beta)$.
\end{lem}

\begin{proof}
    The $\bfC_p$ points of $\Omega_\fp$ and the product $\Omega_\fs$ are
    $$
        \Omega_\fp(\bfC_p) = \bfP^1(\bfC_p) \smallsetminus \bfP^1(K_\fp), \qquad 
        \Omega_\fs(\bfC_p) = \prod_{\fp \mid \fs} \Omega_\fp(\bfC_p).
    $$
    Suppose that $\beta \not\in K^\times$, and let $\fp \mid \fs$. Then $\beta$ acts on $\bfP^1(\bfC_p)$ with two fixed points. These points lie in $\bfP^1(K_\fp)$ if and only if $K(\beta) \otimes_K K_\fp \cong K_\fp \times K_\fp$. In particular, $\beta$ fixes a point of $\Omega_\fp(\bfC_p)$ if and only if $\beta \in K^\times$ or $\fp$ does not split in $K(\beta)$. 

    First suppose that $g$ fixes a point of $\Omega_\fs(\bfC_p)$. Then $\beta = g^m$ fixes a point of $\Omega_\fp(\bfC_p)$ for each $\fp \mid \fs$, so either $\beta \in K^\times$ or no prime $\fp \mid \fs$ splits in $K(\beta)$. Now suppose the converse. Pick representatives $\fp_i$ for the $\sigma$-orbits of $\{\fp \mid \fs\}$, and for each $i$, pick $x_{\fp_i} \in \Omega_{\fp_i}(\bfC_p)$ that is fixed by $\beta$. Since $\rho_{B,\fs}(g) \in H_\fs$, the size of each $\lag\sigma\rag$-orbit of $\fp_i$ is $m$. Indeed, if $\sigma^k \fp_i = \fp_i$ for $0 < k < m$, then $\sigma^k$ induces a non-trivial automorphism of $K_{\fp_i}$, which is not allowed since $\rho_{B,\fs}(g^k) \in H_\fs$. Let
    $$
        x_{\sigma^k\fp_i} = g^k x_{\fp_i}, \quad 0 \leq k < m.
    $$
    Then $x = (x_\fp)_{\fp \mid \fs}$ is fixed by $g$.
\end{proof}

Note that the $\Gamma$-stabilizer of a point in $\Omega_\fs(\bfC_p)$ is finite as it acts discretely on $\Omega_\fs(\bfC_p)$. To see this, apply \cite[Lemma 6]{berkovich-1995-the-automorphism-group-of-the-drinfeld-half-plane} to each factor. Recall from Section \ref{sec:fundamental-domains} that, given a hereditary order $\fO$, we identify all the maximal finite subgroups of $\Psi_{\fO,\fs}$ up to conjugation. Therefore, given $\Gamma \leq \Psi_{\fO,\fs}$, we can check whether $\Gamma$ acts freely by applying Lemma \ref{lem:free-action} to all the torsion elements in $\Gamma$. 

We now prove Proposition \ref{prop:uniformization}. We actually prove a stronger statement that gives information on the commensurability class of $\Lambda_\Gamma$. First, we need some notation. Fix an embedding $\iota: L \hookrightarrow \bfC$. For each $\fp \mid \fs$, the composite
$$
    K \to K_\fp \xrightarrow{\ i_\fp\ } L \xrightarrow{\ \iota\ } \bfC
$$
is a real place $\infty_{\fp}$ of $K$, and these $r$ places are distinct. Let $B^\dag/K$ be the quaternion algebra that splits at each $\infty_{\fp}$, ramifies at the other real places and at each $\fp \mid \fs$, and is locally isomorphic to $B$ at every other finite prime. In particular, $\fd(B^\dag) = \fd(B)\fs$. Note that if $r = [K : \bfQ]$, then $B^\dag$ does not depend on $\iota$.

\begin{prop}\cite{varshavsky-1998-p-adic-uniformization-of-unitary-shimura-varieties-ii}
    \label{prop:uniformization-commensurability}
    Let $K$, $B$, $\fs$, $L$, and $\iota$ be as above, and assume that $\fs$ is eligible with $r = \omega(\fs) \geq 2$. Let $\Gamma \in \dC(B,\fs)$ such that $\Gamma \leq H_\fs$ and $\Gamma$ acts freely on $\Omega_\fs(\bfC_p)$. Then there is a unique (up to isomorphism) smooth projective variety $Y_\Gamma/L$ with ample canonical bundle, and an isomorphism of $L$-analytic spaces,
    $$
        (Y_\Gamma)^\an \cong \Gamma \backslash \Omega_\fs,
    $$
    such that the structure sheaf $\dO_{Y_\Gamma}$ has Euler characteristic $\chi(\dO_{Y_\Gamma}) = \chi(\Gamma)$ and $h^1(\dO_{Y_\Gamma}) = 0$, and for every partition $I$ of $r$, the Chern numbers of $Y_\Gamma$ are
    $$
        c_I(Y_\Gamma) = \chi(\Gamma) c_I((\bfP^1)^r).
    $$
    Furthermore, for every embedding $L \hookrightarrow \bfC$, the resulting complex variety $Y_{\Gamma,\bfC}$ satisfies $(Y_{\Gamma,\bfC})^\an \cong \Lambda_\Gamma\backslash \bfH^r$ for an irreducible, cocompact, torsion-free lattice $\Lambda_\Gamma$ in the group of holomorphic isometries of $\bfH^r$, which is in the commensurability class determined by $B^\dag$.
\end{prop}

\begin{proof}
    First suppose that $\Gamma$ is torsion-free and linear. Then the existence of $Y = Y_\Gamma$ and its claimed properties follows from \cite[Lemma 2.9]{varshavsky-1998-p-adic-uniformization-of-unitary-shimura-varieties-ii}. The claim about its Chern numbers follows from \cite[Proposition A.1]{varshavsky-1998-p-adic-uniformization-of-unitary-shimura-varieties-ii}, and then 
    $$
        \chi(\dO_Y) = \chi(\Gamma)\chi(\dO_{(\bfP^1)^r}) = \chi(\Gamma)
    $$
    by the Hirzebruch--Riemann--Roch theorem \cite[Corollary 15.2.2]{fulton-1998-intersection-theory}. Let $Y_\bfC  = Y \otimes_L \bfC$ be the complex variety obtained from base change. Varshavsky shows in \cite[Corollary 3.14]{varshavsky-1998-p-adic-uniformization-of-unitary-shimura-varieties-ii} that $\bfH^r$ is the universal cover of $(Y_\bfC)^\an$, and there exists an irreducible lattice $\Lambda_\Gamma \leq \PSL_2(\bfR)^r$ such that $(Y_\bfC)^\an \cong \Lambda_\Gamma\backslash\bfH^r$. It further follows from \cite[Section 5]{varshavsky-1998-p-adic-uniformization-of-unitary-shimura-varieties-ii} that $\Lambda_\Gamma$ is commensurable with $(\fO^\dag)^1/\{\pm 1\}$ for a maximal order $\fO^\dag$ of $B^\dag$. By Margulis' normal subgroup theorem, $\Lambda_\Gamma$ has vanishing first Betti number, so $h^1(\dO_{Y_\bfC}) = 0$. Thus, $h^1(\dO_Y) = 0$ \cite[Tag 02KH]{the-stacks-project-authors-2026-the-stacks-project}. We have thus proven our claim for torsion-free $\Gamma$.

    We now consider the general case. By Selberg's lemma, there is a finite-index, normal, torsion-free, linear subgroup $\Gamma_0 \triangleleft \Gamma$. Let $Y_0 = Y_{\Gamma_0}$ and $\Delta = \Gamma/\Gamma_0$. The group $\Delta$ acts on $(Y_0)^\an \cong \Gamma_0\backslash\Omega_\fs$, and since analytification is fully faithful on proper varieties \cite[Corollary 3.4.12]{berkovich-1990-spectral-theory-and-analytic-geometry-over-non-archimedean-fields}, this action comes from a unique action of $\Delta$ on $Y_0$. Since $Y_0$ is projective, the quotient $Y = \Delta\backslash Y_0$ exists and is separated and of finite type over $L$, and the quotient map $\pi : Y_0 \to Y$ is finite and surjective \cite[Expos\'e V, Proposition 1.8 and Corollary 1.5]{grothendieck-and-raynaud-2004-revetements-etales-et-groupe-fondamental-sga-1}. In particular, $Y$ is proper \cite[Tag 03GN]{the-stacks-project-authors-2026-the-stacks-project}. Then $Y^\an \cong \Delta\backslash(Y_0)^\an \cong \Gamma\backslash\Omega_\fs$ \cite[Lemma 1.2]{allcock-and-kato-2017-a-fake-projective-plane-via-2-adic-uniformization-with-torsion}, and $Y$ is unique by \cite[Corollary 3.4.12]{berkovich-1990-spectral-theory-and-analytic-geometry-over-non-archimedean-fields}.

    Since $\Gamma$ acts freely on $\Omega_\fs(\bfC_p)$, the group $\Delta$ acts freely on $Y_0(\bfC_p) = \Gamma_0\backslash\Omega_\fs(\bfC_p)$, so $\pi$ is finite and \'etale of degree $|\Delta|$ \cite[Expos\'e V, \S2 and Corollary 2.3]{grothendieck-and-raynaud-2004-revetements-etales-et-groupe-fondamental-sga-1}. Hence $Y$ is smooth, and $\pi^*\Omega^1_Y \cong \Omega^1_{Y_0}$ \cite[Expos\'e II, Th\'eor\`eme 4.3]{grothendieck-and-raynaud-2004-revetements-etales-et-groupe-fondamental-sga-1}. In particular, the canonical bundle of $Y_0$ is the pullback of that of $Y$, so the canonical bundle of $Y$ is ample \cite[Tag 0B5V]{the-stacks-project-authors-2026-the-stacks-project}, and hence $Y$ is projective. Likewise, the tangent bundle of $Y_0$ is the pullback of that of $Y$, so by multiplicativity, we have
    $$
        c_I(Y_0) = |\Delta|c_I(Y), \qquad \chi(\dO_{Y_0}) = |\Delta|\chi(\dO_Y).
    $$
    Since $\chi(\Gamma_0) = |\Delta|\chi(\Gamma)$, the claims about $\chi(\dO_Y)$ and the Chern numbers follow from the torsion-free case.

    Now let $Y_{0,\bfC} = Y_0 \otimes_L \bfC$ and $Y_\bfC = Y \otimes_L \bfC$ be the complex varieties obtained by base change. Then $Y_{0,\bfC} \to Y_\bfC$ is finite and \'etale by \cite[Tags 01WL and 02GO]{the-stacks-project-authors-2026-the-stacks-project}. Thus, the induced map 
    \begin{equation}
        \label{eq:complex-covering}
        (Y_{0,\bfC})^\an \to (Y_\bfC)^\an
    \end{equation}
    is finite and \'etale \cite[Expos\'e XII, Propositions 3.1 and 3.2]{grothendieck-and-raynaud-2004-revetements-etales-et-groupe-fondamental-sga-1}. Since $(Y_{0,\bfC})^\an$ is connected, with universal cover $\bfH^r$, and $\Delta$ still acts freely after base change, \eqref{eq:complex-covering} is a finite Galois covering map with deck group $\Delta$. Hence, $(Y_\bfC)^\an$ is connected and $\bfH^r$ is the universal cover. Let $\Lambda = \Lambda_\Gamma$ be the deck transformations of this covering, and let $\Lambda_0 \leq \Lambda$ be the subgroup associated to \eqref{eq:complex-covering}. Then $\Lambda_0$ is conjugate to $\Lambda_{\Gamma_0}$, so $\Lambda$ is in the commensurability class determined by $B^\dag$. Finally, $\Lambda$ is torsion-free since it acts freely on the contractible space $\bfH^r$ and cocompact and irreducible since it is commensurable with $\Lambda_0$. We can again apply Margulis' normal subgroup theorem \cite{margulis-1991-discrete-subgroups-of-semisimple-lie-groups}  and base change \cite[Tag 02KH]{the-stacks-project-authors-2026-the-stacks-project} to show that $h^1(\dO_Y) = 0$.
\end{proof}

\begin{proof}[Proof of Theorem \ref{thm:A}]
    Let $\Gamma$ be as in the statement of the theorem and let $n = [K:\bfQ]$.

    In \cite{linowitz-stover-and-voight-2015-commensurability-classes-of-fake-quadrics}, Linowitz--Stover--Voight enumerate all number fields $F/\bfQ$ and all quaternion algebras $A/F$ such that there can exist a torsion-free lattice in the group of holomorphic isometries of $\bfH^2$ in the commensurability class determined by $A$ with Euler characteristic 4. Let $\dL\dS\dV$ be the list of pairs $(F,A)$.

    Fix an embedding $L \hookrightarrow \bfC$. By Proposition \ref{prop:uniformization-commensurability}, we have $(Y_{\Gamma,\bfC})^\an \cong \Lambda_\Gamma \backslash \bfH^2$ for an irreducible, cocompact, torsion-free lattice $\Lambda_\Gamma$ with $\chi(\Lambda_\Gamma) = 4$. Furthermore, $\Lambda_\Gamma$ lies in the commensurability class of lattices in the group of holomorphic isometries of $\bfH^2$ determined by $B^\dag/K$ with $\fd(B^\dag) = \fd(B)\fs$. Therefore, $(K,B^\dag) \in \dL\dS\dV$.

    By Theorem \ref{thm:psi-maximal}, there exists a hereditary order $\fO \subseteq B$ of level $\fn$ such that $\Gamma \leq \Psi_{\fO,\fs}$ with finite index. Note that 
    $$
        [\Psi_{\fO,\fs}:\Phi_{\fO,\fs}] \le |\Sigma(\fd,\fs)| \le |\Aut(K)| \le n.
    $$
    Let $\fM$ be a maximal order of $B$. By the Prasad volume formula (Theorem \ref{thm:prasad-volume}), together with \eqref{eq:normalizer-unit-index} and Proposition \ref{prop:norm-index}, we have
    \begin{equation}
        \label{eq:psi-euler-bound}
        1 \geq \chi(\Psi_{\fO,\fs}) \geq \frac{1}{n}\cdot\chi(\Phi_{\fO,\fs}) 
        \geq \frac{1}{n} \cdot \chi(\Phi_{\fM,\fs})\cdot \prod_{\fl \mid \fn} \frac{N\fl + 1}{2} \geq \frac{1}{n} \cdot \chi(\Phi_{\fM,\fs}).
    \end{equation}
    
    Suppose that $(F,A) \in \dL\dS\dV$, and let $\ft \mid \fd(A)$ be an eligible ideal that is a product of two primes. Let $A^\ddag/F$ be the totally definite quaternion algebra with $\fd(A^\ddag) = \fd(A)/\ft$. Let $\fM$ be a maximal order of $A^\ddag$. By \eqref{eq:psi-euler-bound}, if the commensurability class of $A^\ddag$ in $G_\ft$ contains an Euler characteristic one lattice, then $\chi(\Phi_{\fM,\ft}) \leq [F:\bfQ]$. We enumerate all such pairs $(A^\ddag,\ft)$. For each pair, we compute $\Psi_{\fO,\ft}$ for each of the finitely many hereditary orders $\fO$ of $A^\ddag$ with $\chi(\Psi_{\fO,\ft}) \leq 1$. We can then enumerate all Euler characteristic one $\Gamma$ acting freely on $\Omega_\ft(\bfC_p)$ by applying Lemma \ref{lem:free-action}.
    
    We obtain exactly two commensurability classes which are described in the next section. In each example, $K$ is real quadratic, so the commensurability class of $\Lambda_\Gamma$ does not depend on the embedding $L \hookrightarrow \bfC$. See Sections \ref{sec:qsqrt7} and \ref{sec:qsqrt6} for an explicit description of these examples.
\end{proof}

\section{Examples of small lattices}

\label{sec:examples}

In this section, we study lattices defined over number fields with small Euler characteristics. We are most interested in torsion-free lattices, lattices that act freely on $\Omega_\fs$, and lattices that act freely on $T_\fs^{(0)}$. These examples are all computed in Magma using the algorithm described in Section \ref{sec:fundamental-domains}. The code accompanying these examples will be provided in the supplementary code (see Section \ref{sec:fundamental-domains}).

We will use the following combinatorial lemma.

\begin{lem} \label{lem:euler-vertex-count}
    Let $\Gamma \leq G_\fs$ be a lattice that acts freely on $T_\fs^{(0)}$. Then
    $$
        \chi(\Gamma) = 
        V_\Gamma \cdot \prod_{\fp \mid \fs} \frac{1 - N\fp}{2}, 
        \qquad V_\Gamma = |X_\Gamma^{(0)}|.
    $$
    If $\Gamma$ is torsion-free, then $X_\Gamma$ is a cube complex with $\chi(X_\Gamma) = \chi(\Gamma)$.
\end{lem}

\begin{proof}
    Let $v \in T_\fs^{(0)}$. Then, for each $\fa \mid \fs$, there are $\prod_{\fp \mid \fa} (N\fp + 1)$ cubes in $T_\fs(\fa)$ containing $v$. Therefore, the number of $k$-cubes containing $v$ is
    $$
        e_k = \sum_{\substack{\fa \mid \fs \\ \omega(\fa) = k}} \prod_{\fp \mid \fa} (N\fp + 1).
    $$
    Let $T_\fs[k]$ be the $k$-cubes of $T_\fs$, and let
    $$
        P_k = \{(v,c) \in T_\fs[0] \times T_\fs[k] : v \in c\}.
    $$
    We count the number of elements in $\Gamma\backslash P_k$ in two different ways. Since $\Gamma$ acts freely on $T_\fs[0]$, we have 
    \begin{equation}
        \label{eq:incidences-by-vertex}
        |\Gamma \backslash P_k| = V_\Gamma \cdot e_k.
    \end{equation}
    If $c$ is a $k$-cube, then $c$ has $2^k$ vertices. Since $\Gamma$ acts freely on the vertices, there are $2^k/|\Gamma_c|$ orbits in $\Gamma\backslash P_k$ lifting to a pair $(v,c)$, where $\Gamma_c$ is the stabilizer of $c$. Therefore,
    \begin{equation}
        \label{eq:incidences-by-cube}
        |\Gamma \backslash P_k| = \sum_c \frac{2^k}{|\Gamma_c|},
    \end{equation}
    where the sum ranges over representatives of $\Gamma\backslash T_\fs[k]$. Combining \eqref{eq:incidences-by-vertex} and \eqref{eq:incidences-by-cube}, we get
    $$
        \chi(\Gamma) = 
        \sum_{k = 0}^{\omega(\fs)} (-1)^k \sum_{c \in \Gamma\backslash T_\fs[k]} \frac{1}{|\Gamma_c|} = 
        V_\Gamma\cdot\sum_{k=0}^{\omega(\fs)} (-1)^k \frac{e_k}{2^k} = 
        V_\Gamma\cdot\prod_{\fp \mid \fs} \left(1 - \frac{N\fp + 1}{2}\right) = 
        V_\Gamma\cdot\prod_{\fp \mid \fs} \frac{1 - N\fp}{2}.
    $$
    Finally, when $\Gamma$ is torsion-free, we get $\chi(X_\Gamma) = \chi(\Gamma)$ since $T_\fs$ is contractible.
\end{proof}

\subsection{Euler characteristic one lattices}
\label{sec:chi1-lattices}

In this section, we enumerate all irreducible, torsion-free lattices $\Gamma \leq G_\fs$ with $\omega(\fs) = 2$ such that $\chi(\Gamma) = 1$, proving Theorem \ref{thm:C}.

\begin{proof}[Proof of Theorem \ref{thm:C}]
    Let $\fs = \fp\fq$. By Lemma \ref{lem:euler-vertex-count}, $(N\fp,N\fq,V_\Gamma)$ must be a solution to
    $$
        V_\Gamma(N\fp - 1)(N\fq - 1)/4 = 1
    $$
    which has exactly 4 solutions:
    $$
        (N\fp,N\fq,V_\Gamma) = \quad (3,3,1), \quad (2,3,2), \quad (2,2,4), \quad (2,5,1).
    $$
    If $N\fp \neq N\fq$, then $\Gamma$ cannot exchange $T_\fp$ and $T_\fq$. Thus, $\Gamma$ acts freely on the edges in the $T_\fp$ factor, so there are $V_\Gamma(N\fp+1)/2$ edges in $X_\Gamma$ that lift to an edge in $T_\fp$. If $(N\fp,N\fq,V_\Gamma) = (2,5,1)$, then $\Gamma$ must invert an edge of $T_\fp$, so $\Gamma$ is not torsion-free. This eliminates one of the possible solutions.
    
    In Magma, we compute all square complexes $X$ covered by $T_\fs$ with $(N\fp,N\fq,V_\Gamma)$ as above. This computation is feasible because $X$ has so few cells. Fix such a square complex $X$. If $\Gamma = \pi_1(X)$ is isomorphic to an arithmetic lattice in $G_\fs$, then both $\Gamma^\ab$ and $[\Gamma,\Gamma]^\ab$ are finite by Margulis' normal subgroup theorem \cite{margulis-1991-discrete-subgroups-of-semisimple-lie-groups}. We enumerate all such $X$ below (up to isometry).
    \begin{center}
        \setlength{\tabcolsep}{16pt}
        \begin{tabular}{ccc}
            \toprule
            $(N\fp,N\fq,V_\Gamma)$ & square complexes $X$ up to isometry & $\Gamma^\ab$ and $[\Gamma,\Gamma]^\ab$ finite \\
            \midrule
            $(3,3,1)$ & $66$ & $3$ \\
            $(2,3,2)$ & $179$ & $3$ \\
            $(2,2,4)$ & $166$ & $2$ \\
            \bottomrule
        \end{tabular}
    \end{center}
    Let $\Gamma$ be one of the $(3,3,1)$ or $(2,2,4)$ examples. Using the \texttt{L2Quotients} function in Magma, we show that there is a prime $p$ such that $\Gamma$ contains an index 2 subgroup $\Gamma_0$ that admits infinitely many surjections to groups of the form $\PGL_2(\bfF)$ where $\bfF$ is a finite field of characteristic $p$, which is impossible for an irreducible lattice in $G_\fs$ by Margulis' superrigidity theorem \cite{margulis-1991-discrete-subgroups-of-semisimple-lie-groups}. Therefore, $\Gamma$ cannot be commensurable with an irreducible lattice in $G_\fs$. See the Magma handbook for details on \texttt{L2Quotients} \cite{cannon-bosma-fieker-and-steel-2026-handbook-of-magma-functions}. This leaves the three groups $\Gamma$ associated to the solution $(2,3,2)$. We identify each of these groups with one of the cases in Sections \ref{sec:q-disc5}, \ref{sec:q-disc13}, and \ref{sec:qsqrt3}.
\end{proof}

\begin{rmk}\label{rmk:chi1-positive-characteristic}
    Recall that in \cite{stix-and-vdovina-2017-simply-transitive-quaternionic-lattices-of-rank-2-over-f-q-t-and-a-non-classical-fake-quadric} and \cite{rungtanapirom-2018-quaternionic-arithmetic-lattices-of-rank-2-and-a-fake-quadric-in-characteristic-2}, Stix--Vdovina and Rungtanapirom respectively construct torsion-free lattices,
    $$
        \Gamma_{\text{SV}} \leq \PGL_2\left(\bfF_3((t))\right) \times \PGL_2\left(\bfF_3((t))\right), \qquad \Gamma_{\text{R}} \leq \PGL_2\left(\bfF_2((t))\right) \times \PGL_2\left(\bfF_2((t))\right),
    $$
    with Euler characteristic 1. These give examples for the $(3,3,1)$ and $(2,2,4)$ cases respectively. We give their associated square complexes in Figures \ref{fig:chi1-331-stix-vdovina} and \ref{fig:chi1-224-rungtanapirom}. With Magma, we can show that the remaining lattices are also arithmetic.
    \begin{enumerate}
        \item The $(3,3,1)$ lattice of Figure \ref{fig:chi1-331-swap} is commensurable with $\Gamma_{\text{SV}}$.
        \item The $(2,2,4)$ lattice of Figure \ref{fig:chi1-224-swap} is commensurable with $\Gamma_{\text{R}}$.
        \item The remaining $(3,3,1)$ lattice $\Gamma_{\text{JW2}}$ in Figure \ref{fig:chi1-331-janzen-wise} was constructed by Janzen and Wise \cite{janzen-and-wise-2009-a-smallest-irreducible-lattice-in-the-product-of-trees}, where they conjectured it is irreducible. It was stated in \cite{caprace-and-vast-2026-fractal-anti-tori} that the group $\Gamma_{\text{JW2}}$ has a surjective map to an arithmetic lattice over $\bfF_3(t)$, thus proving irreducibility. We can show in Magma that $\Gamma_{\text{JW2}}$ is arithmetic, but it is not commensurable with $\Gamma_{\text{SV}}$.
    \end{enumerate}   
    Given that the focus of this paper is characteristic 0 lattices, we leave the details of these computations to the interested reader.
\end{rmk}

\subsubsection{Rational example of discriminant 5}
\label{sec:q-disc5}

Let 
$$
    B = (-2,-5)_\bfQ, \qquad \fd(B) = 5,
$$
which contains a unique maximal order $\fO \subseteq B$ up to conjugation. If $\fs = 6$, then $\chi(\Phi_{\fO,\fs}) = 1/12$, and
$$
    \Phi_{\fO,\fs} = \left\langle a, b, c \;\middle|\;
    \begin{array}{l}
        a^2,\ b^2,\ c^2,\ (ab)^6, \\
        (cababa)^2,\ (cbcaca)^2
    \end{array}
    \right\rangle.
$$
The generators are represented in $B^\times$ by
$$
    a = j, \quad
    b = 5i + 6j + ij, \quad
    c = i + 2j + ij.
$$
There exist congruence subgroups of level 5,
$$
    \Lambda = \langle (ca)^2,\ abcabab \rangle \leq \Gamma = \langle c,\ aca,\ (ab)^3 \rangle \leq \Phi_{\fO,\fs}^\times \leq \Phi_{\fO,\fs}, \qquad [\Gamma:\Lambda] = 2,
$$
such that $\Gamma$ acts simply transitively on the vertices of $T_\fs$, so $\chi(\Gamma) = 1/2$ and $\chi(\Lambda) = 1$ by Lemma \ref{lem:euler-vertex-count}. Also, $\Lambda$ is torsion-free. We display the square complex $X_\Lambda$ and the covering $X_\Lambda \to X_\Gamma$ in Figure \ref{fig:q-disc5-covering}.

These groups $\Gamma$ and $\Lambda$ were also constructed by hand in unpublished work of the second author with Jakob Stix.

\begin{rmk}\label{rmk:janzen-wise}
    Let
    $$
        \Gamma_{\text{JW}} = \lag u, v, x, y \;\middle|\; uxuy,\ ux^{-1}vy^{-1},\ uy^{-1}v^{-1}x^{-1},\ vxv^{-1}y^{-1} \rag
    $$
    be the group of Janzen--Wise \cite{janzen-and-wise-2009-a-smallest-irreducible-lattice-in-the-product-of-trees}, with the presentation as in \cite[Section 2]{caprace-and-vast-2026-fractal-anti-tori}. Janzen--Wise show that $\Gamma_{\text{JW}}$ is irreducible, and it was later observed that $\Gamma_{\text{JW}}$ is not residually finite \cite{bondarenko-and-kivva-2022-automaton-groups-and-complete-square-complexes,caprace-and-wesolek-2018-indicability-residual-finiteness-and-simple-subquotients-of-groups-acting-on-trees}. We get
    $$
        \Gamma_{\text{JW}}/\ncl{x^3} \cong \lag v, x \mid x^3,\ z \rag, \qquad
        z = xv^2xv^{-1}x^{-1}v^2x^{-1}v^{-1}.
    $$
    The assignment $v \mapsto c(ab)^3$, $x \mapsto (ba)^2$ defines an isomorphism from this group onto
    $$
        \Delta_{\text{JW}} = \lag c(ab)^3,\ (ba)^2 \rag \leq \Phi_{\fO,\fs}, \qquad [\Phi_{\fO,\fs}:\Delta_{\text{JW}}] = 4.
    $$
    Conjugation by $x$ permutes $v_i = x^ivx^{-i}$ for $i \in \bfZ/3\bfZ$, and the kernel of the map $\Delta_{\text{JW}} \to \bfZ/3\bfZ$ is the group $\Lambda$ described above. We get
    $$
        \Lambda = \lag v_0, v_1, v_2 \mid v_{i+1}^2 v_{i-1}^{-1} v_{i+1}^2 v_i^{-1},\ i \in \bfZ/3\bfZ \rag.
    $$
    See also \cite[Remark 2.5]{caprace-and-vast-2026-fractal-anti-tori}, where a surjection like this was referenced by Caprace and Vast.
\end{rmk}

\subsubsection{Rational example of discriminant 13}
\label{sec:q-disc13}

Let
$$
    B = (-2,-13)_\bfQ, \qquad \fd(B) = 13.
$$
Then $B$ contains a unique maximal order $\fO$ up to conjugation. If $\fs = 6$, then $\chi(\Phi_{\fO,\fs}) = 1/4$, and
$$
    \Phi_{\fO,\fs} = \left\langle a, b, c \;\middle|\;
    \begin{array}{l}
        a^2,\ b^2,\ c^2,\ (baca)^2, \\
        (acbcacacbc)^2,\ (bcbcacacbca)^2
    \end{array}
    \right\rangle.
$$
The unit subgroup $\Gamma = \Phi_{\fO,\fs}^\times$ acts simply transitively on the vertices of $T_\fs$, so $\chi(\Gamma) = 1/2$ by Lemma \ref{lem:euler-vertex-count}. Furthermore, there exists a torsion-free, index 2 subgroup $\Lambda \leq \Gamma$, so $\chi(\Lambda) = 1$. We display the covering $X_\Lambda \to X_\Gamma$ in Figure \ref{fig:q-disc13-covering}.

We remark that Chiu first pointed out that $\Phi_{\fO,p}^\times$ acts simply transitively on the associated Bruhat--Tits tree $T_p$ for all primes $p \neq 13$ \cite{chiu-1992-cubic-ramanujan-graphs}. He was motivated by the construction of Cayley Ramanujan graphs.

\subsubsection{An example defined over \texorpdfstring{$\bfQ(\sqrt{3})$}{Q(√3)}}
\label{sec:qsqrt3}

Let
$$
    K = \bfQ(\sqrt{3}), \qquad B = (-1,-1)_K, \qquad \fd(B) = \fo_K.
$$
Then $B$ contains two maximal orders $\fO$ and $\fO'$ up to conjugation. The orders are distinguished by their unit groups. If $\fs = \fp_2\fp_3$, then $\chi(\Psi_{\fO,\fs}) = 1/48$, and it does not act transitively on $T_\fs^{(0)}$. The presentation of $\Psi_{\fO,\fs}$ is quite complicated, but we will only consider subgroups of 
$$
    \Phi_{\fO,\fs}^\times = \left\langle a, b, c \;\middle|\;
    \begin{array}{l}
        a^2,\ b^2,\ c^2, \\
        (babacacbc)^2,\ c(ba)^2bcbcac(ba)^4
    \end{array}
    \right\rangle.
$$
There exists a torsion-free, congruence subgroup of level $\fp_{11}$,
$$
   \Gamma = \langle bcabab,\ cabacbcac \rangle \leq \Phi_{\fO,\fs}^\times,
$$
with $V_\Gamma = 2$, so $\chi(\Gamma) = 1$ by Lemma \ref{lem:euler-vertex-count}. It is noteworthy that there exists a torsion-free subgroup $\Gamma$ with the same number of vertex orbits as $\Psi_{\fO,\fs}$. This appears to be very rare when $\Psi_{\fO,\fs}$ has more than one vertex orbit. However, there are many examples when $\Psi_{\fO,\fs}$ (or $\Phi_{\fO,\fs}$) has one vertex orbit; see \cite{mendel-and-yu-2026-s-arithmetic-groups-acting-simply-transitively-on-products-of-bruhat-tits-trees}. We display the square complex $X_\Gamma$ in Figure \ref{fig:qsqrt3-complex}.

\subsection{Small 3-adic lattices}
\label{sec:3-adic}

In this section, we give examples of interesting lattices $\Gamma \leq G_\fs$ with $\fs = \fp\fq$ and $N\fp = N\fq = 3$.

\begin{thm}
    \label{thm:chi2-classification}
    Let $K$ be a totally real number field, let $\fs$ be an eligible square-free ideal of $\fo_K$ with $\omega(\fs) = 2$, and let $B/K$ be a totally definite quaternion algebra with discriminant coprime to $\fs$. Let $\Lambda \in \dC^*(B,\fs)$ be torsion-free with $\chi(\Lambda) = 2$. Then $K = \bfQ(\sqrt{13})$, $\fd(B) = \fo_K$, and $\fs$ is over $3$. Moreover, $\Lambda$ is the lattice described in Section \ref{sec:qsqrt13}.
\end{thm}

\begin{proof}
    By Lemma \ref{lem:euler-vertex-count}, $\fs = \fp_3\fq_3$ or $\fp_2\fq_2$ with $V_\Lambda = 2$ or $8$ respectively. Since $\Lambda \le G_\fs^*$, $\Lambda$ does not exchange the two trees (so $X_\Lambda$ is a \emph{VH-complex}). We enumerate all such complexes $X$ and keep those whose fundamental group $\Lambda$ has $\Lambda^\ab$ and $[\Lambda,\Lambda]^\ab$ finite, as required by Margulis' normal subgroup theorem. Since $\Lambda$ is linear, strong approximation shows that $\Lambda$ surjects onto $\PSL_2(\bfF)$ or $\PGL_2(\bfF)$ for finite fields $\bfF$ of almost every characteristic. Using \texttt{L2Quotients}, we find exactly six $\Lambda$ with such an infinite family of quotients. Five of these complexes have a subgroup of finite index with infinite abelianization, which contradicts the normal subgroup theorem. The sixth is the lattice $\Lambda$ in Section \ref{sec:qsqrt13}.
\end{proof}

\subsubsection{An example defined over \texorpdfstring{$\bfQ(\sqrt{13})$}{Q(√13)}}
\label{sec:qsqrt13}

Let
$$
    K = \bfQ(\sqrt{13}), \qquad B = (-1,-1)_K, \qquad \fd(B) = \fo_K.
$$
Then $B$ has a unique maximal order $\fO$ up to conjugation. If $\fs = \fp_3\fq_3$, then $\chi(\Phi_{\fO,\fs}^\times) = 1/12$, and
$$
    \Phi_{\fO,\fs}^\times = \left\langle a, b, c \;\middle|\;
    \begin{array}{l}
        a^3,\ b^2,\ c^2,\ (bc)^6, \\
        (a^{-1}cab)^2,\ (bcbca^{-1})^2
    \end{array}
    \right\rangle.
$$
There exist subgroups
$$
    \Lambda = \langle (cbcba^{-1}b)^2,\ bacb \rangle \leq \Gamma = \langle cbcba^{-1}b,\ bacb \rangle \leq \Phi_{\fO,\fs}^\times, \qquad [\Gamma:\Lambda] = 2,
$$
such that $\Gamma$ has torsion but acts simply transitively on the vertices of $T_\fs$, so $\chi(\Gamma) = 1$ by Lemma \ref{lem:euler-vertex-count}. Furthermore, $\Lambda$ is torsion-free with $\chi(\Lambda) = 2$. We remark that there are two non-conjugate copies of $\Gamma$ in $\Phi_{\fO,\fs}$, but $\Gamma$ is unique up to conjugation in $\Psi_{\fO,\fs}$. We display the covering $X_\Lambda \to X_\Gamma$ in Figure \ref{fig:qsqrt13-covering}.

\subsubsection{An example defined over \texorpdfstring{$\bfQ(\sqrt{7})$}{Q(√7)}}
\label{sec:qsqrt7}

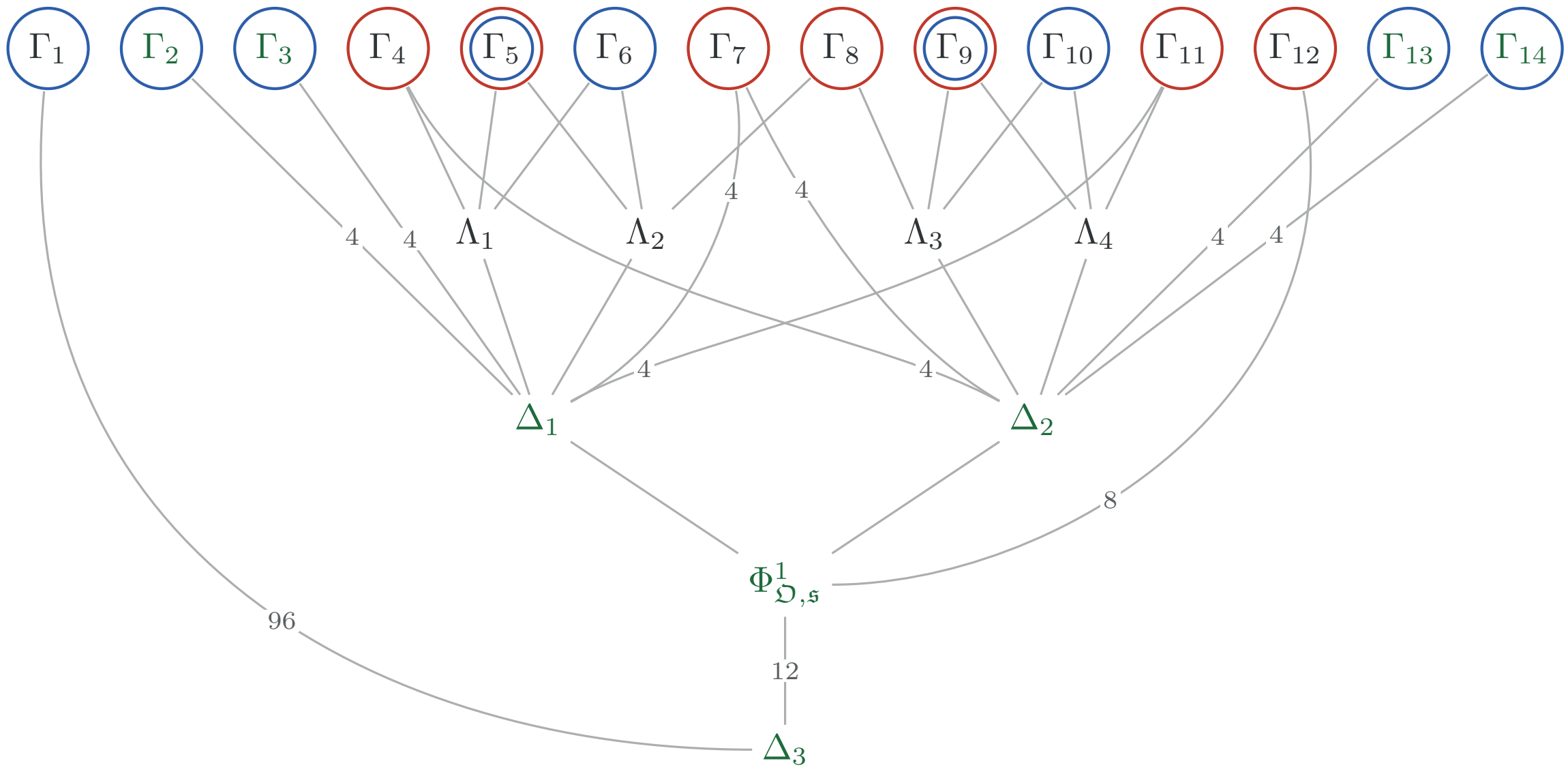
\begin{figure}[htbp]
      \centering
      \definecolor{SoftInk}{HTML}{2F3437}
\definecolor{LinearGreen}{HTML}{1E6B3C}
\definecolor{FreeOmega}{HTML}{C0392B}
\definecolor{FreeVertex}{HTML}{2E5EAA}
\begin{tikzpicture}[x=1cm,y=1cm,
  every node/.style={text=SoftInk},
  lattice/.style={inner sep=2.5pt,font=\small},
  incl/.style={draw=SoftInk!40,line width=0.6pt},
  gamma/.style={lattice,circle,inner sep=0pt,minimum size=0.82cm},
  omegaring/.style={circle,draw=FreeOmega,line width=0.8pt,minimum size=0.78cm,inner sep=0pt},
  vertexring/.style={circle,draw=FreeVertex,line width=0.8pt,minimum size=0.78cm,inner sep=0pt},
  inner/.style={minimum size=0.6cm},
  degree/.style={font=\tiny,fill=white,inner sep=0.8pt,text=SoftInk!80}]
  \node[gamma] (Gamma10) at (-7.15,7.0) {$\Gamma_{1}$};
  \node[gamma,text=LinearGreen] (Gamma9) at (-6.05,7.0) {$\Gamma_{2}$};
  \node[gamma,text=LinearGreen] (Gamma11) at (-4.95,7.0) {$\Gamma_{3}$};
  \node[gamma] (Gamma5) at (-3.85,7.0) {$\Gamma_{4}$};
  \node[gamma] (Gamma2) at (-2.75,7.0) {$\Gamma_{5}$};
  \node[gamma] (Gamma12) at (-1.65,7.0) {$\Gamma_{6}$};
  \node[gamma] (Gamma4) at (-0.55,7.0) {$\Gamma_{7}$};
  \node[gamma] (Gamma7) at (0.55,7.0) {$\Gamma_{8}$};
  \node[gamma] (Gamma1) at (1.65,7.0) {$\Gamma_{9}$};
  \node[gamma] (Gamma14) at (2.75,7.0) {$\Gamma_{10}$};
  \node[gamma] (Gamma6) at (3.85,7.0) {$\Gamma_{11}$};
  \node[gamma] (Gamma3) at (4.95,7.0) {$\Gamma_{12}$};
  \node[gamma,text=LinearGreen] (Gamma8) at (6.05,7.0) {$\Gamma_{13}$};
  \node[gamma,text=LinearGreen] (Gamma13) at (7.15,7.0) {$\Gamma_{14}$};
  \node[lattice] (Lambda1) at (-3.0,5.2) {$\Lambda_{1}$};
  \node[lattice] (Lambda2) at (-1.35,5.2) {$\Lambda_{2}$};
  \node[lattice] (Lambda3) at (1.35,5.2) {$\Lambda_{3}$};
  \node[lattice] (Lambda4) at (3.0,5.2) {$\Lambda_{4}$};
  \node[lattice,text=LinearGreen] (DeltaA) at (-2.4,3.4) {$\Delta_1$};
  \node[lattice,text=LinearGreen] (DeltaB) at (2.4,3.4) {$\Delta_2$};
  \node[lattice,text=LinearGreen] (Phi1) at (0.0,1.8) {$\Phi^1_{\fO,\fs}$};
  \node[lattice,text=LinearGreen] (DeltaC) at (0.0,0.2) {$\Delta_3$};
  \draw[incl] (Phi1) -- node[degree,pos=0.5] {$12$} (DeltaC);
  \draw[incl] (DeltaA) -- (Phi1);
  \draw[incl] (DeltaB) -- (Phi1);
  \draw[incl] (Lambda1) -- (DeltaA);
  \draw[incl] (Lambda2) -- (DeltaA);
  \draw[incl] (Lambda3) -- (DeltaB);
  \draw[incl] (Lambda4) -- (DeltaB);
  \draw[incl] (Gamma10) .. controls (-7.6,2.2) and (-4.0,0.2) .. node[degree,pos=0.55] {$96$} (DeltaC);
  \draw[incl] (Gamma9) -- node[degree,pos=0.5] {$4$} (DeltaA);
  \draw[incl] (Gamma11) -- node[degree,pos=0.5] {$4$} (DeltaA);
  \draw[incl] (Gamma5) -- (Lambda1);
  \draw[incl] (Gamma5) .. controls (-2.6,4.6) and (0.6,4.4) .. node[degree,pos=0.86] {$4$} (DeltaB);
  \draw[incl] (Gamma2) -- (Lambda1);
  \draw[incl] (Gamma2) -- (Lambda2);
  \draw[incl] (Gamma12) -- (Lambda1);
  \draw[incl] (Gamma12) -- (Lambda2);
  \draw[incl] (Gamma4) .. controls (-0.3,5.6) and (-0.9,4.2) .. node[degree,pos=0.3] {$4$} (DeltaA);
  \draw[incl] (Gamma4) .. controls (0.1,5.6) and (1.0,4.2) .. node[degree,pos=0.3] {$4$} (DeltaB);
  \draw[incl] (Gamma7) -- (Lambda2);
  \draw[incl] (Gamma7) -- (Lambda3);
  \draw[incl] (Gamma1) -- (Lambda3);
  \draw[incl] (Gamma1) -- (Lambda4);
  \draw[incl] (Gamma14) -- (Lambda3);
  \draw[incl] (Gamma14) -- (Lambda4);
  \draw[incl] (Gamma6) -- (Lambda4);
  \draw[incl] (Gamma6) .. controls (2.6,4.6) and (-0.6,4.4) .. node[degree,pos=0.86] {$4$} (DeltaA);
  \draw[incl] (Gamma3) .. controls (5.6,3.6) and (2.6,1.8) .. node[degree,pos=0.6] {$8$} (Phi1);
  \draw[incl] (Gamma8) -- node[degree,pos=0.5] {$4$} (DeltaB);
  \draw[incl] (Gamma13) -- node[degree,pos=0.5] {$4$} (DeltaB);
  \node[vertexring] at (Gamma10) {};
  \node[vertexring] at (Gamma9) {};
  \node[vertexring] at (Gamma11) {};
  \node[omegaring] at (Gamma5) {};
  \node[omegaring] at (Gamma2) {};
  \node[vertexring,inner] at (Gamma2) {};
  \node[vertexring] at (Gamma12) {};
  \node[omegaring] at (Gamma4) {};
  \node[omegaring] at (Gamma7) {};
  \node[omegaring] at (Gamma1) {};
  \node[vertexring,inner] at (Gamma1) {};
  \node[vertexring] at (Gamma14) {};
  \node[omegaring] at (Gamma6) {};
  \node[omegaring] at (Gamma3) {};
  \node[vertexring] at (Gamma8) {};
  \node[vertexring] at (Gamma13) {};
\end{tikzpicture}
      \caption{Inclusions up to conjugacy among the subgroups of $\Psi_{\fO,\fs}$ for $K = \bfQ(\sqrt{7})$. Unlabelled edges have index $2$. Green labels denote linear subgroups. A red circle marks a $\Gamma_i$ acting freely on $\Omega_\fs$, and a blue circle one acting simply transitively on the vertices of $T_\fs$.}
      \label{fig:qsqrt7-subgroups}
\end{figure}

This next family of examples is quite spectacular. Let 
$$
    K = \bfQ(\sqrt{7}), \qquad B = (2\sqrt{7}-7,-2)_K, \qquad \fd(B) = \fp_2\fp_7.
$$
Then $B$ contains a unique maximal order $\fO$ up to conjugation. If $\fs = 3\fo_K = \fp_3\fq_3$, then $\chi(\Psi_{\fO,\fs}) = 1/4$, and
$$
    \Psi_{\fO,\fs} = \left\langle a, b, c \;\middle|\;
    \begin{array}{l}
        a^4,\ b^2,\ c^2,\ (a^2c)^2,\ (baca)^2,\ (baba^{-1}ba^2)^2, \\
        baba^{-1}bcacbaba^{-1}ba
    \end{array}
    \right\rangle.
$$
There are twenty conjugacy classes of subgroups with index $4$, none of which are torsion-free. It follows from Lemma \ref{lem:conjugacy-in-psi} that these groups are mutually non-conjugate in $G_\fs$. Consider fourteen of these subgroups:
$$
\begin{array}{lll}
    \Gamma_1 = \langle a^2c,\ aba^{-1},\ cba \rangle, &
    \Gamma_2 = \langle b,\ aba,\ aca^{-1} \rangle, &
    \Gamma_3 = \langle b,\ a^2c,\ aba \rangle, \\
    \Gamma_4 = \langle ba^{-1},\ a^2 \rangle, &
    \Gamma_5 = \langle ba^{-1},\ a^2c,\ cba \rangle, &
    \Gamma_6 = \langle ba^{-1},\ c,\ a^{-1}ba^{-2} \rangle, \\
    \Gamma_7 = \langle a,\ bab \rangle, &
    \Gamma_8 = \langle a^2,\ abc,\ bac \rangle, &
    \Gamma_9 = \langle ca^{-1},\ a^2ba^{-2} \rangle, \\
    \Gamma_{10} = \langle ca^{-1},\ a^2b \rangle, &
    \Gamma_{11} = \langle ca^{-1},\ a^2,\ baba^{-1} \rangle, &
    \Gamma_{12} = \langle a,\ cb \rangle, \\
    \Gamma_{13} = \langle b,\ c,\ a^{-1}ba \rangle, &
    \Gamma_{14} = \langle c,\ a^2b,\ aba \rangle. &
\end{array}
$$
Also,
\begin{itemize}
    \item $\Gamma_4, \Gamma_5, \Gamma_7, \Gamma_8, \Gamma_9, \Gamma_{11}, \Gamma_{12}$ act freely on $\Omega_\fs$;
    \item $\Gamma_1, \Gamma_2, \Gamma_3, \Gamma_5, \Gamma_6, \Gamma_9, \Gamma_{10}, \Gamma_{13}, \Gamma_{14}$ act simply transitively on the vertices of $T_\fs$;
    \item $\Gamma_2, \Gamma_3, \Gamma_{13}, \Gamma_{14}$ are linear.
\end{itemize}
We display the complexes $X_{\Gamma_i}$ of the nine lattices that act simply transitively on the vertices of $T_\fs$ in Figures \ref{fig:qsqrt7-gamma1}--\ref{fig:qsqrt7-gamma14}. There are four torsion-free groups of index $8$ in $\Psi_{\fO,\fs}$, up to conjugacy. These groups are
$$
\begin{array}{ll}
    \Lambda_1 = \langle cbca^{-1},\ ba^{-1},\ caca^{-1} \rangle, &
    \Lambda_2 = \langle ba^{-1}c,\ abc,\ acb \rangle, \\
    \Lambda_3 = \langle bac,\ abc,\ acb \rangle, &
    \Lambda_4 = \langle ca^{-1},\ baba^{-1} \rangle.
\end{array}
$$
We display the square complexes $X_{\Lambda_i}$ in Figures \ref{fig:qsqrt7-lambda1}--\ref{fig:qsqrt7-lambda4}. None of them are linear, and their intersections with $\Phi_{\fO,\fs}$ are
\begin{align*}
    \Delta_1 &= \Lambda_1 \cap \Phi_{\fO,\fs} = \Lambda_2 \cap \Phi_{\fO,\fs} = \langle (a^{-1}b)^2,\ (ab)^2 \rangle, \\
    \Delta_2 &= \Lambda_3 \cap \Phi_{\fO,\fs} = \Lambda_4 \cap \Phi_{\fO,\fs} = \langle baba^{-1},\ ababa^{-2} \rangle.
\end{align*}
The subgroups $\Delta_1$ and $\Delta_2$ are normal in $\Psi_{\fO,\fs}$, of index $16$, and 
$$
    \Delta_1 \cap \Delta_2 = \Phi^1_{\fO,\fs}.
$$
Let $\Delta_3$ be the intersection of the normal cores of $\Gamma_1, \dots, \Gamma_{14}$. Then $\Delta_3$ is a torsion-free normal congruence subgroup of level $\fp_2^2\fp_7$ and index $384$ in $\Psi_{\fO,\fs}$. Figure \ref{fig:qsqrt7-subgroups} shows the inclusions among these subgroups.

\subsubsection{Another example defined over \texorpdfstring{$\bfQ(\sqrt{13})$}{Q(√13)}}
\label{sec:qsqrt13-chi3}

Let
$$
    K = \bfQ(\sqrt{13}), \qquad \fd(B) = \fp_4\fp_{13}.
$$
Then $B$ contains two maximal orders $\fO$ and $\fO'$ up to conjugation. Both orders have the same unit groups but are distinguished by their normalizers in $B^\times$. If $\fs = \fp_3\fq_3$, then $\chi(\Psi_{\fO,\fs}) = 3/8$. We find that $\Gamma = \Phi_{\fO,\fs}^\times$ is torsion-free with $V_\Gamma = 3$, so $\chi(\Gamma) = 3$ by Lemma \ref{lem:euler-vertex-count}. Up to conjugacy, $\Gamma$ is the unique torsion-free subgroup of index $8$ in $\Psi_{\fO,\fs}$.

\subsection{Small 2-adic lattices}
\label{sec:2-adic}

In this section, we give the $2$-adic lattices described in Theorem \ref{thm:B}. We also prove the following theorem.

\begin{thm}\label{thm:2-adic-torsion-free}
    There exist irreducible, torsion-free lattices $\Lambda \leq \PGL_2(\bfQ_2) \times \PGL_2(\bfQ_2)$ with $\chi(\Lambda) = 3$ and $4$.
\end{thm}

\subsubsection{An example defined over \texorpdfstring{$\bfQ(\sqrt{17})$}{Q(√17)}}
\label{sec:qsqrt17}

Let
$$
    K = \bfQ(\sqrt{17}), \qquad B = (-1,-3)_K, \qquad \fd(B) = \fo_K.
$$
Then $B$ contains a unique maximal order $\fO$ up to conjugation. If $\fs = 2\fo_K = \fp_2\fq_2$, then $\chi(\Psi_{\fO,\fs}) = 1/48$, and
$$
    \Psi_{\fO,\fs} = \left\langle a, b, c \;\middle|\;
    \begin{array}{l}
        a^8,\ b^4,\ c^2,\ (a^{-1}b)^2,\ (ac)^2,\ cb^{-1}ca^3ba^3, \\
        aca^{-2}(cb)^3,\ ab^{-1}a^{-3}b^{-1}ca^{-1}bc
    \end{array}
    \right\rangle.
$$
There exist congruence subgroups of level $\fp_{47}$,
$$
    \Lambda \leq \Gamma \leq \Phi_{\fO,\fs}^\times \leq \Psi_{\fO,\fs}, \qquad [\Gamma:\Lambda] = 2,
$$
with $\chi(\Gamma) = 2$ and $\chi(\Lambda) = 4$. The lattice $\Gamma$ is not torsion-free, but it acts freely on $\Omega_\fs$ and freely on the vertices of $T_\fs$. The lattice $\Lambda$ is torsion-free. We display the one-skeleton of $X_\Lambda$ and the covering $X_\Lambda \to X_\Gamma$ in Figure \ref{fig:qsqrt17-covering}.

\subsubsection{An example defined over \texorpdfstring{$\bfQ(\sqrt{33})$}{Q(√33)}}
\label{sec:qsqrt33}

Let
$$
    K = \bfQ(w), \qquad w = \frac{1+\sqrt{33}}{2}, \qquad B = (-1,w-4)_K, \qquad \fd(B) = \fo_K, \qquad \fs = 2\fo_K = \fp_2\fq_2.
$$
We will consider hereditary orders of different levels, which we describe in the table below.
\begin{center}
    \setlength{\tabcolsep}{12pt}
    \begin{tabular}{cccc}
        \toprule
        level of $\fO$ & orders up to conjugacy & $\Psi_{\fO,\fs}$ up to conjugacy & $\chi(\Psi_{\fO,\fs})$ \\
        \midrule
        $\fo_K$ & $3$ & $1$ & $1/16$ \\
        $\fp_3$ & $5$ & $1$ & $1/8$ \\
        $\fp_{11}$ & $3$ & $1$ & $3/8$ \\
        $\fp_3\fp_{11}$ & $8$ & $1$ & $3/4$ \\
        \bottomrule
    \end{tabular}
\end{center}
Notice that $\Psi_{\fO,\fs}$ is well-defined for each level we consider. Fix a maximal order $\fM$, and let $\fM_3$ and $\fM_{11}$ be maximal orders such that 
$$
    \fO(3) = \fM \cap \fM_3, \qquad \fO(11) = \fM \cap \fM_{11}, \qquad \fO(33) = \fO(3) \cap \fO(11)
$$
are hereditary of level $\fp_3$, $\fp_{11}$, and $\fp_3\fp_{11}$ respectively. To ease notation, let
$$
    \Psi = \Psi_{\fM,\fs}, \qquad \Psi_3 = \Psi_{\fO(3),\fs}, \qquad \Psi_{11} = \Psi_{\fO(11),\fs}, \qquad \Psi_{33} = \Psi_{\fO(33),\fs},
$$
and we use similar notation for the linear groups $\Phi$.

\begin{figure}[htbp]
      \centering
      \definecolor{SoftInk}{HTML}{2F3437}
\definecolor{LinearGreen}{HTML}{1E6B3C}
\definecolor{Commensurator}{HTML}{8B5FBF}
\definecolor{TorsionFree}{HTML}{C0392B}
\definecolor{NormalizerBlue}{HTML}{2F6DB5}
\begin{tikzpicture}[x=1cm,y=1cm,
  every node/.style={text=SoftInk},
  lattice/.style={inner sep=2.5pt,font=\small},
  small/.style={lattice},
  incl/.style={draw=SoftInk!40,line width=0.6pt},
  degree/.style={font=\tiny,fill=white,inner sep=0.8pt,text=SoftInk!80},
  comm/.style={draw=Commensurator,line width=0.8pt,densely dotted,<->,>=stealth},
  commlabel/.style={font=\scriptsize,text=Commensurator,fill=white,inner sep=0.6pt},
  tf/.style={lattice,circle,draw=TorsionFree,line width=0.7pt,minimum size=0.56cm,inner sep=1pt},
  of/.style={lattice,circle,minimum size=0.56cm,inner sep=1pt}]
  \node[lattice] (PsiM) at (-5.19,8.6) {$\Psi$};
  \node[lattice] (Psi3) at (4.925,7.2) {$\Psi_3$};
  \node[lattice] (Psi11) at (-2.6,5.8) {$\Psi_{11}$};
  \node[lattice] (E11) at (-5.19,4.2) {$\Psi\cap\Psi_{11}$};
  \node[lattice,text=LinearGreen] (PhiO) at (-2.6,4.2) {$\Phi_{11}$};
  \node[lattice] (Psi33) at (2.625,4.2) {$\Psi_{33}$};
  \node[lattice,text=LinearGreen] (PhiMO) at (-4.35,2.4) {$\Phi\cap\Phi_{11}$};
  \node[lattice] (E3) at (0.325,2.4) {$\Psi_{11}\cap\Psi_{33}$};
  \node[lattice] (E33) at (4.925,2.4) {$\Psi_3\cap\Psi_{33}$};
  \node[tf,text=LinearGreen] (Lambda) at (-4.35,0.6) {$\Lambda$};
  \node[tf] (D1) at (-6.6,0.6) {$\Delta_1$};
  \node[tf] (D1p) at (-5.45,0.6) {$\Delta_1'$};
  \node[of,text=LinearGreen] (Gamma) at (-2.6,0.6) {$\Gamma$};
  \node[tf] (D2) at (-1.4,0.6) {$\Delta_2$};
  \node[tf] (D2p) at (-0.25,0.6) {$\Delta_2'$};
  \node[of] (P1) at (0.9,0.6) {$\Pi_1$};
  \node[of] (P1p) at (2.05,0.6) {$\Pi_1'$};
  \node[tf] (D3) at (3.2,0.6) {$\Delta_3$};
  \node[tf] (D3p) at (4.35,0.6) {$\Delta_3'$};
  \node[of] (P2) at (5.5,0.6) {$\Pi_2$};
  \node[of] (P2p) at (6.65,0.6) {$\Pi_2'$};
  \draw[incl] (PsiM) -- node[degree,pos=0.5] {$12$} (E11);
  \draw[incl] (Psi11) -- (E11);
  \draw[incl] (Psi11) -- node[degree,pos=0.5] {$4$} (E3);
  \draw[incl] (Psi11) -- (PhiO);
  \draw[incl] (PhiO) -- node[degree,pos=0.5] {$4$} (Gamma);
  \draw[incl] (PhiO) -- (PhiMO);
  \draw[incl] (Psi33) -- (E3);
  \draw[incl] (Psi33) -- (E33);
  \draw[incl] (Psi3) -- node[degree,pos=0.5] {$12$} (E33);
  \draw[incl] (E11) -- (PhiMO);
  \draw[incl] (PhiMO) -- (Lambda);
  \draw[incl] (E11) -- node[degree,pos=0.5] {$4$} (D1);
  \draw[incl] (E11) -- node[degree,pos=0.5] {$4$} (D1p);
  \foreach \x in {D2, D2p, P1, P1p} \draw[incl] (E3) -- (\x);
  \foreach \x in {D3, D3p, P2, P2p} \draw[incl] (E33) -- (\x);
  \draw[comm] (D1.south) to[bend right=45] node[commlabel,below] {$g_{11}$} (D1p.south);
  \draw[comm] (D2.south) to[bend right=45] node[commlabel,below] {$g_3$} (D2p.south);
  \draw[comm] (P1.south) to[bend right=45] node[commlabel,below] {$g_3$} (P1p.south);
  \draw[comm] (D3.south) to[bend right=45] node[commlabel,below] {$g_{11}$} (D3p.south);
  \draw[comm] (P2.south) to[bend right=45] node[commlabel,below] {$g_{11}$} (P2p.south);
  \draw[draw=NormalizerBlue,line width=0.8pt,densely dotted,line cap=round]
    (Lambda.north east) .. controls (-3.7,1.3) and (-3.05,2.2) .. (-3.45,3.25)
    .. controls (-3.8,4.15) and (-3.4,5.0) .. (-2.75,5.5);
\end{tikzpicture}
      \caption{The lattices described in Section \ref{sec:qsqrt33}. The bottom row consists of lattices with Euler characteristic $3$ acting freely on $\Omega_\fs$, and the circled ones are torsion-free. Edges are inclusions up to conjugation in the commensurator, unlabelled edges have index $2$, and green labels denote linear subgroups. The dotted purple arrows are conjugation by $g_3$ and $g_{11}$. The normalizer in the commensurator of each lattice in the bottom row is the group directly above it, except for $\Lambda$, whose normalizer is $\Psi_{11}$.}
      \label{fig:qsqrt33-subgroups}
\end{figure}
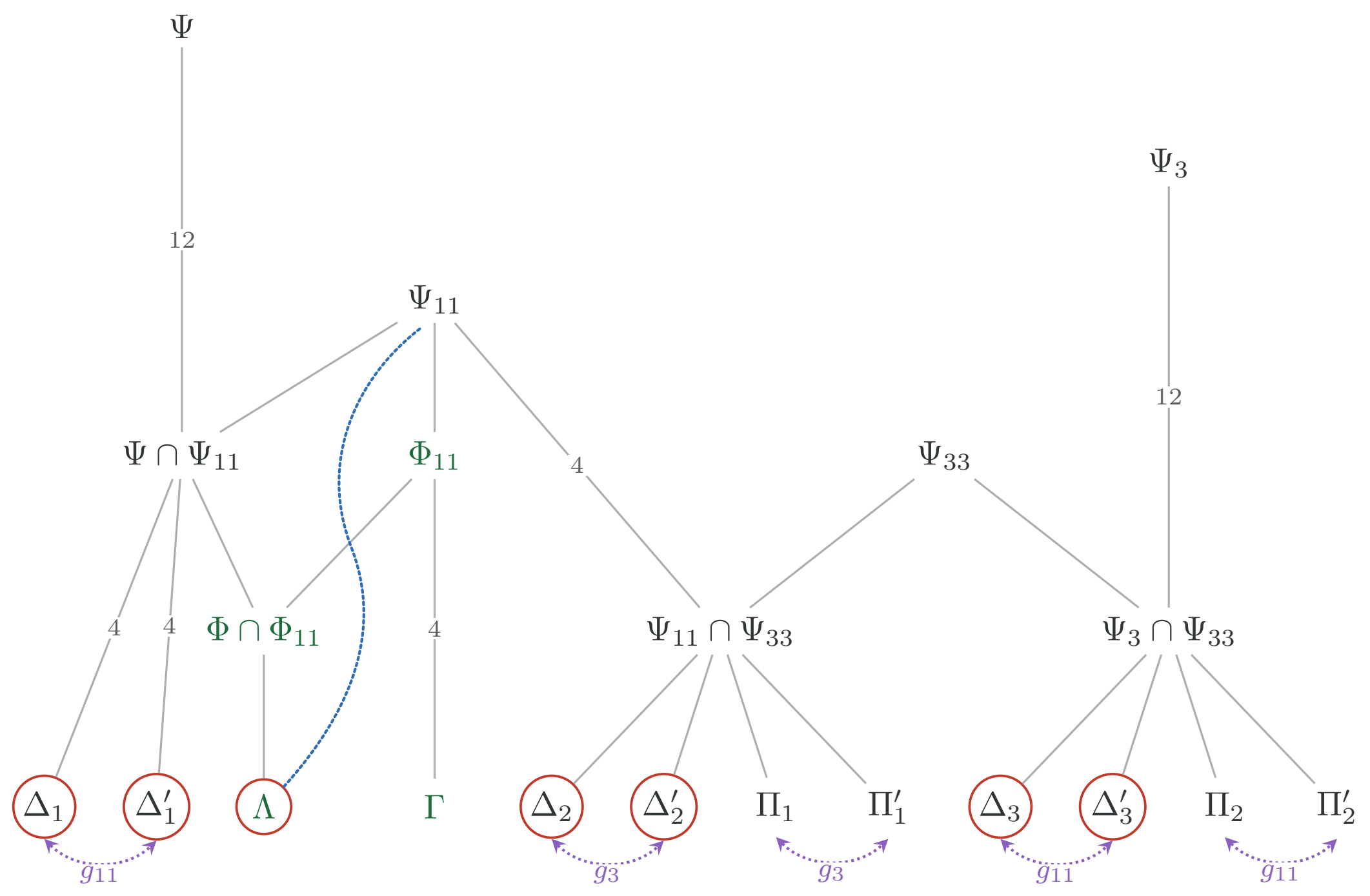

Figure \ref{fig:qsqrt33-subgroups} shows how seven isomorphism classes of lattices with Euler characteristic $3$ lie inside $\Psi$, $\Psi_3$, $\Psi_{11}$, and $\Psi_{33}$. Each of these lattices but $\Pi_2$ acts freely on the vertices of $T_\fs$. We display the one-skeleta of these lattices in Figure \ref{fig:qsqrt33-one-skeleta}.

We describe how the lattices sit inside the maximal groups. Every element of $\Psi_{11}$ either fixes both $\fM$ and $\fM_{11}$ or interchanges them, so $\Psi \cap \Psi_{11}$ has index $2$ in $\Psi_{11}$. Similarly, $\Psi_{11} \cap \Psi_{33}$ and $\Psi_3 \cap \Psi_{33}$ have index $2$ in $\Psi_{33}$, since at $\fp_3$ (respectively $\fp_{11}$) every element of $\Psi_{33}$ fixes or interchanges the two maximal orders containing $\fO(33)$. Fix involutions
$$
    g_{11} \in \Psi_{11} \cap \Psi_{33}, \qquad g_3 \in \Psi_3 \cap \Psi_{33}
$$
interchanging $\fM$ with $\fM_{11}$ and with $\fM_3$ respectively, so that $g_{11} \notin \Psi \cup \Psi_3$ and $g_3 \notin \Psi \cup \Psi_{11}$. Conjugation by $g_{11}$ preserves the normal subgroup $\Psi \cap \Psi_{11}$ of $\Psi_{11}$, but since $g_{11} \notin \Psi$, it need not preserve the $\Psi$-conjugacy class of a subgroup of $\Psi \cap \Psi_{11}$. In the same way, $g_3$ acts on $\Psi_{11} \cap \Psi_{33}$ and $g_{11}$ acts on $\Psi_3 \cap \Psi_{33}$. 

In Figure \ref{fig:qsqrt33-subgroups}, notice that we sometimes list a lattice twice. Each such pair consists of a lattice $\Xi$ and its conjugate $\Xi' = g\Xi g^{-1}$, where $g$ is $g_3$ or $g_{11}$. Both $\Xi$ and $\Xi'$ lie in one of the intersections above, so they are conjugate in the maximal lattice containing $g$. However, they are not conjugate in the other maximal lattice containing the intersection. The following table records, for each $\Xi$, the intersection containing $\Xi$ and $\Xi'$, the element $g$, and the lattice in which $\Xi$ and $\Xi'$ are not conjugate.
$$
    \begin{array}{llll}
        \toprule
        \Xi & \Xi, \Xi' \leq & g & \text{not conjugate in} \\
        \midrule
        \Delta_1 & \Psi \cap \Psi_{11} & g_{11} & \Psi \\
        \Delta_2,\ \Pi_1 & \Psi_{11} \cap \Psi_{33} & g_3 & \Psi_{11} \\
        \Delta_3,\ \Pi_2 & \Psi_3 \cap \Psi_{33} & g_{11} & \Psi_3 \\
        \bottomrule
    \end{array}
$$

\subsection{Small 5-adic lattices}
\label{sec:5-adic}

\subsubsection{An example defined over \texorpdfstring{$\bfQ(\sqrt{6})$}{Q(√6)}}
\label{sec:qsqrt6}

Let
$$
    K = \bfQ(\sqrt{6}), \qquad B = (-1,-3-\sqrt{6})_K, \qquad \fd(B) = \fp_2\fp_3,
$$
which has a unique maximal order $\fO$ up to conjugation. Let $\fs = 5\fo_K = \fp_5\fq_5$. Then $\chi(\Psi_{\fO,\fs}) = 1/4$, and
$$
    \Psi_{\fO,\fs} = \left\langle a, b, c, d \;\middle|\;
    \begin{array}{l}
        a^2,\ c^2,\ d^2,\ (ab^2)^2,\ (cb^2)^2,\ ab^3ab^{-1}, \\
        (cb^{-1}cb^{-1}a)^2,\ (cb^{-1}cba)^2,\ cbcdcdcb^{-1}cd,\ bcdcb^{-1}cbadab^{-1}c
    \end{array}
    \right\rangle.
$$
Up to conjugacy, there are 7 subgroups of index $4$ in $\Psi_{\fO,\fs}$ that act freely on $\Omega_\fs$,
$$
\begin{array}{lll}
    \Gamma_1 = \langle ba,\ c,\ d \rangle, &
    \Gamma_2 = \langle dca,\ b \rangle, &
    \Gamma_3 = \langle ba,\ ab,\ cd \rangle, \\
    \Gamma_4 = \langle ba,\ ad \rangle, &
    \Gamma_5 = \langle cd,\ b \rangle, &
    \Gamma_6 = \langle bc,\ bd \rangle, \\
    \Gamma_7 = \langle cba,\ b^2,\ cd \rangle. & &
\end{array}
$$
None of the $\Gamma_i$ are torsion-free or linear. It follows from Lemma \ref{lem:conjugacy-in-psi} that these groups are not conjugate in $G_\fs$. Up to conjugacy, the torsion-free subgroups of index $16$ in $\Psi_{\fO,\fs}$ are exactly
$$
    \Lambda_1 = \langle bc,\ db \rangle, \qquad
    \Lambda_2 = \langle dab^{-2},\ b^2dc,\ bcba \rangle.
$$
Both act simply transitively on the vertices of $T_\fs$. Each $\Gamma_i$ and $\Lambda_i$ is a congruence subgroup of $\Psi_{\fO,\fs}$. In particular, $\Lambda_2$ is a torsion-free, normal, linear congruence subgroup of level $\fp_3$. The lattice $\Lambda_1$ is not linear. We display the square complexes $X_{\Lambda_1}$ and $X_{\Lambda_2}$ in Figures \ref{fig:qsqrt6-lambda1} and \ref{fig:qsqrt6-lambda2}.

\subsubsection{A cubic example}
\label{sec:cubic}

Let $\zeta_{13}$ be a primitive 13th root of unity, and let $K$ be the unique cubic subfield of $\bfQ(\zeta_{13})$. Then $K$ is totally real and Galois over $\bfQ$. Furthermore, $13$ is the only ramified prime and $5$ splits completely. Let
$$
    B = (-2,-13)_K, \qquad \fd(B) = \fp_{13}, \qquad \fs = 5\fo_K = \fp\fq\fr.
$$
Then $B$ has a unique maximal order $\fO$ up to conjugation, and $\fO^\times/\fo_K^\times$ is trivial. Furthermore, $\Phi_{\fO,\fs}^\times$ acts simply transitively on the vertices of $T_\fs$. It is not torsion-free, but it contains a unique torsion-free subgroup $\Lambda$ of index $2$. Hence $X_\Lambda = \Lambda \backslash T_\fs$ has two vertices, and $\chi(\Lambda) = -16$ by Lemma \ref{lem:euler-vertex-count}. Compare with \cite[Example 2.36]{rungtanapirom-stix-and-vdovina-2019-infinite-series-of-quaternionic-1-vertex-cube-complexes-the-doubling-construction-and-explicit-cubical-ramanujan-complexes}, where the authors construct a torsion-free lattice in
$$
    \PGL_2(\bfF_5((t))) \times \PGL_2(\bfF_5((t))) \times \PGL_2(\bfF_5((t))) 
$$
acting simply transitively on the associated product of 6-valent trees.

\subsection{Infinite families}
\label{sec:families}

Let $K$ be a number field, let $B/K$ be a totally definite quaternion algebra, and let $\fM \subseteq B$ be a maximal order. We say that $B$ has \emph{class number one} if every right $\fM$-ideal is principal.

In \cite{mendel-and-yu-2026-s-arithmetic-groups-acting-simply-transitively-on-products-of-bruhat-tits-trees}, the first author and Yu show that if $B$ has class number one, then there exists a finite set of primes $\kappa$ such that, for every square-free ideal $\fs$ with no prime divisors in $\kappa$, there exists a congruence subgroup $\Gamma \leq \Phi_{\fM,\fs}^\times$ of level $\fa \subseteq \fo_K$ that acts simply transitively on the vertices of $T_{\fs}$. We use their construction to find families of torsion-free lattices. We give two examples below.

\subsubsection{A family defined over \texorpdfstring{$\bfQ(\sqrt{5})$}{Q(√5)}}

\label{sec:family-qsqrt5}

\begin{thm}\label{thm:family-qsqrt5}
    Let $p$ be a prime with $p \equiv 1 \bmod 5$. Then there exists a torsion-free, irreducible lattice
    $$
        \Gamma \leq \PGL_2(\bfQ_p) \times \PGL_2(\bfQ_p)
    $$
    with $V_\Gamma = 1$.
\end{thm}

\begin{proof}
    Let $K = \bfQ(\sqrt{5})$ and $\fO$ be a maximal order of the totally definite quaternion algebra $B/K$ with $\fd(B) = \fo_K$. Since $p \equiv 1 \bmod 5$, the prime $p$ totally splits in $K$. Let $\fs = p\fo_K$. Then $\Phi_{\fO,\fs}^\times$ acts transitively on the vertices of $T_{\fs}$. 
    
    Let $\fp_5$ be the unique prime of norm $5$. In Magma \cite{bosma-cannon-and-playoust-1997-the-magma-algebra-system-i-the-user-language}, we find that $\pi_{\fp_5}$ is injective on $\fO^\times/\fo_K^\times$, and there exists a subgroup
    $$
        \text{$H \leq \PGL_2(\bfF_5)$ complementary to $\pi_{\fp_5}(\fO^\times/\fo_K^\times)$}.
    $$
    It follows from \cite{mendel-and-yu-2026-s-arithmetic-groups-acting-simply-transitively-on-products-of-bruhat-tits-trees} that $\Gamma = \pi_{\fp_5}^{-1}(H)$ acts simply transitively on the vertices of $T_{\fs}$. Let $\Theta_{\fp_5}(H)$ be the torsion-obstruction set defined in \cite{mendel-and-yu-2026-s-arithmetic-groups-acting-simply-transitively-on-products-of-bruhat-tits-trees}. Then it follows from \cite[Appendix A]{mendel-and-yu-2026-s-arithmetic-groups-acting-simply-transitively-on-products-of-bruhat-tits-trees} that $\Theta_{\fp_5}(H)$ does not contain the trivial class in $\bfF_5^\times/\bfF_5^{\times 2}$. If $a \in \fo_K$ is a totally positive generator of a prime $\fp \mif \fs$, then $a \bmod \fp_5$ is a square. Therefore, by \cite[Proposition 43]{mendel-and-yu-2026-s-arithmetic-groups-acting-simply-transitively-on-products-of-bruhat-tits-trees}, $\Gamma$ is torsion-free.
\end{proof}

\begin{cor}\label{cor:family-qsqrt5}
    There exist irreducible, arithmetic complex surfaces uniformized by $\bfH \times \bfH$ with holomorphic Euler characteristic $(p-1)^2/4$ for each prime $p \equiv 1 \bmod 5$. These surfaces are pairwise non-commensurable.
\end{cor}

\subsubsection{A cubic family}

\label{sec:family-cubic}

\begin{thm}\label{thm:family-cubic}
    Let $p$ be a prime with $p \equiv \pm 1 \bmod 13$. Then there exists a torsion-free, irreducible lattice
    $$
        \Gamma \leq \PGL_2(\bfQ_p) \times \PGL_2(\bfQ_p) \times \PGL_2(\bfQ_p)
    $$
    with $V_\Gamma = 1$.
\end{thm}

\begin{proof}
    The proof is similar to that of the previous theorem. Let $K$, $B$, and $\fO$ be as in Section \ref{sec:cubic}. It follows from class field theory that $p$ totally splits in $K$ if and only if
    $$
        p \equiv \pm 1, \pm 5 \bmod 13.
    $$
    For any such prime $p$, let $\fs = p\fo_K$. Then $\Gamma = \Phi_{\fO,\fs}^\times$ acts simply transitively on the vertices of $T_{\fs}$. We now assume that $p \equiv \pm 1 \bmod 13$. Let $\Theta_{\fp_{13}}$ be the torsion-obstruction set defined in \cite{mendel-and-yu-2026-s-arithmetic-groups-acting-simply-transitively-on-products-of-bruhat-tits-trees}. Then $1 \not\in \Theta_{\fp_{13}}$ by \cite[Appendix A]{mendel-and-yu-2026-s-arithmetic-groups-acting-simply-transitively-on-products-of-bruhat-tits-trees}. If $a \in \fo_K$ is a totally positive generator of a prime $\fp \mif \fs$, then $a \bmod \fp_{13}$ is a square. Therefore, by \cite[Proposition 43]{mendel-and-yu-2026-s-arithmetic-groups-acting-simply-transitively-on-products-of-bruhat-tits-trees}, $\Gamma$ is torsion-free.
\end{proof}

\begin{cor}\label{cor:family-cubic}
    There exist irreducible, arithmetic complex threefolds uniformized by $\bfH \times \bfH \times \bfH$ with holomorphic Euler characteristic $-(p-1)^3/8$ for each prime $p \equiv \pm 1 \bmod 13$. These threefolds are pairwise non-commensurable.
\end{cor}

\appendix

\section{Square complexes}
\label{app:square-complexes}

In this appendix, we display the square complexes described in Section \ref{sec:examples}.

\begin{figure}[htbp]
    \centering
    \definecolor{SoftInk}{HTML}{2F3437}
\definecolor{EdgeColour}{HTML}{4F7D72}
\definecolor{PrimeOne}{HTML}{4F7D72}
\definecolor{PrimeTwo}{HTML}{A85D67}
\begin{tikzpicture}[x=1cm,y=1cm,
  every node/.style={font=\scriptsize,text=SoftInk},
  edge/.style={line width=1.05pt,line cap=round},
  side/.style={line width=0.95pt,line cap=round},
  tip/.style={line width=0.9pt,line cap=round,line join=round},
  cell/.style={fill=SoftInk!5},
  label/.style={inner sep=1.2pt},
  vert/.style={circle,fill=SoftInk,inner sep=1.8pt},
  corner/.style={circle,fill=SoftInk,inner sep=1.2pt}]
  \draw[edge,draw=PrimeOne] (0,0.5) .. controls (2.231,1.098) and (0.598,2.731) .. (0,0.5);
  \draw[tip,draw=PrimeOne,shift={(1.061,1.561)},rotate=135.0] (-2.5pt,2.1pt) -- (0,0) -- (-2.5pt,-2.1pt);
  \node[label,anchor=225] at (1.146,1.646) {$E_{1}$};
  \draw[edge,draw=PrimeOne] (0,0.5) .. controls (-0.598,2.731) and (-2.231,1.098) .. (0,0.5);
  \draw[tip,draw=PrimeOne,shift={(-1.061,1.561)},rotate=225.0] (-2.5pt,2.1pt) -- (0,0) -- (-2.5pt,-2.1pt);
  \node[label,anchor=315] at (-1.146,1.646) {$E_{2}$};
  \draw[edge,draw=PrimeTwo] (0,0.5) .. controls (-2.231,-0.098) and (-0.598,-1.731) .. (0,0.5);
  \draw[tip,draw=PrimeTwo,shift={(-1.061,-0.561)},rotate=315.0] (-2.5pt,2.1pt) -- (0,0) -- (-2.5pt,-2.1pt);
  \node[label,anchor=405] at (-1.146,-0.646) {$E_{3}$};
  \draw[edge,draw=PrimeTwo] (0,0.5) .. controls (0.598,-1.731) and (2.231,-0.098) .. (0,0.5);
  \draw[tip,draw=PrimeTwo,shift={(1.061,-0.561)},rotate=405.0] (-2.5pt,2.1pt) -- (0,0) -- (-2.5pt,-2.1pt);
  \node[label,anchor=495] at (1.146,-0.646) {$E_{4}$};
  \node[vert] at (0,0.5) {};
  \fill[cell] (2.900,0.550) rectangle (4.200,1.850);
  \draw[side,draw=PrimeOne] (2.900,0.550) -- (2.900,1.850);
  \draw[tip,draw=PrimeOne,shift={(2.900,1.200)},rotate=90.0] (-2.5pt,2.1pt) -- (0,0) -- (-2.5pt,-2.1pt);
  \node[label,anchor=east] at (2.810,1.200) {$E_{1}$};
  \draw[side,draw=PrimeTwo] (2.900,1.850) -- (4.200,1.850);
  \draw[tip,draw=PrimeTwo,shift={(3.550,1.850)},rotate=180.0] (-2.5pt,2.1pt) -- (0,0) -- (-2.5pt,-2.1pt);
  \node[label,anchor=south] at (3.550,1.940) {$E_{3}$};
  \draw[side,draw=PrimeOne] (4.200,1.850) -- (4.200,0.550);
  \draw[tip,draw=PrimeOne,shift={(4.200,1.200)},rotate=-90.0] (-2.5pt,2.1pt) -- (0,0) -- (-2.5pt,-2.1pt);
  \node[label,anchor=west] at (4.290,1.200) {$E_{1}$};
  \draw[side,draw=PrimeTwo] (4.200,0.550) -- (2.900,0.550);
  \draw[tip,draw=PrimeTwo,shift={(3.550,0.550)},rotate=180.0] (-2.5pt,2.1pt) -- (0,0) -- (-2.5pt,-2.1pt);
  \node[label,anchor=north] at (3.550,0.460) {$E_{4}$};
  \node[corner] at (2.900,0.550) {};
  \node[corner] at (2.900,1.850) {};
  \node[corner] at (4.200,1.850) {};
  \node[corner] at (4.200,0.550) {};
  \node[inner sep=1.4pt,fill=SoftInk!5,rounded corners=1pt] at (3.550,1.200) {$F_{1}$};
  \fill[cell] (5.650,0.550) rectangle (6.950,1.850);
  \draw[side,draw=PrimeOne] (5.650,0.550) -- (5.650,1.850);
  \draw[tip,draw=PrimeOne,shift={(5.650,1.200)},rotate=90.0] (-2.5pt,2.1pt) -- (0,0) -- (-2.5pt,-2.1pt);
  \node[label,anchor=east] at (5.560,1.200) {$E_{1}$};
  \draw[side,draw=PrimeTwo] (5.650,1.850) -- (6.950,1.850);
  \draw[tip,draw=PrimeTwo,shift={(6.300,1.850)},rotate=0.0] (-2.5pt,2.1pt) -- (0,0) -- (-2.5pt,-2.1pt);
  \node[label,anchor=south] at (6.300,1.940) {$E_{3}$};
  \draw[side,draw=PrimeOne] (6.950,1.850) -- (6.950,0.550);
  \draw[tip,draw=PrimeOne,shift={(6.950,1.200)},rotate=-90.0] (-2.5pt,2.1pt) -- (0,0) -- (-2.5pt,-2.1pt);
  \node[label,anchor=west] at (7.040,1.200) {$E_{2}$};
  \draw[side,draw=PrimeTwo] (6.950,0.550) -- (5.650,0.550);
  \draw[tip,draw=PrimeTwo,shift={(6.300,0.550)},rotate=180.0] (-2.5pt,2.1pt) -- (0,0) -- (-2.5pt,-2.1pt);
  \node[label,anchor=north] at (6.300,0.460) {$E_{3}$};
  \node[corner] at (5.650,0.550) {};
  \node[corner] at (5.650,1.850) {};
  \node[corner] at (6.950,1.850) {};
  \node[corner] at (6.950,0.550) {};
  \node[inner sep=1.4pt,fill=SoftInk!5,rounded corners=1pt] at (6.300,1.200) {$F_{2}$};
  \fill[cell] (8.400,0.550) rectangle (9.700,1.850);
  \draw[side,draw=PrimeOne] (8.400,0.550) -- (8.400,1.850);
  \draw[tip,draw=PrimeOne,shift={(8.400,1.200)},rotate=90.0] (-2.5pt,2.1pt) -- (0,0) -- (-2.5pt,-2.1pt);
  \node[label,anchor=east] at (8.310,1.200) {$E_{1}$};
  \draw[side,draw=PrimeTwo] (8.400,1.850) -- (9.700,1.850);
  \draw[tip,draw=PrimeTwo,shift={(9.050,1.850)},rotate=180.0] (-2.5pt,2.1pt) -- (0,0) -- (-2.5pt,-2.1pt);
  \node[label,anchor=south] at (9.050,1.940) {$E_{4}$};
  \draw[side,draw=PrimeOne] (9.700,1.850) -- (9.700,0.550);
  \draw[tip,draw=PrimeOne,shift={(9.700,1.200)},rotate=90.0] (-2.5pt,2.1pt) -- (0,0) -- (-2.5pt,-2.1pt);
  \node[label,anchor=west] at (9.790,1.200) {$E_{2}$};
  \draw[side,draw=PrimeTwo] (9.700,0.550) -- (8.400,0.550);
  \draw[tip,draw=PrimeTwo,shift={(9.050,0.550)},rotate=360.0] (-2.5pt,2.1pt) -- (0,0) -- (-2.5pt,-2.1pt);
  \node[label,anchor=north] at (9.050,0.460) {$E_{4}$};
  \node[corner] at (8.400,0.550) {};
  \node[corner] at (8.400,1.850) {};
  \node[corner] at (9.700,1.850) {};
  \node[corner] at (9.700,0.550) {};
  \node[inner sep=1.4pt,fill=SoftInk!5,rounded corners=1pt] at (9.050,1.200) {$F_{3}$};
  \fill[cell] (11.150,0.550) rectangle (12.450,1.850);
  \draw[side,draw=PrimeOne] (11.150,0.550) -- (11.150,1.850);
  \draw[tip,draw=PrimeOne,shift={(11.150,1.200)},rotate=90.0] (-2.5pt,2.1pt) -- (0,0) -- (-2.5pt,-2.1pt);
  \node[label,anchor=east] at (11.060,1.200) {$E_{2}$};
  \draw[side,draw=PrimeTwo] (11.150,1.850) -- (12.450,1.850);
  \draw[tip,draw=PrimeTwo,shift={(11.800,1.850)},rotate=180.0] (-2.5pt,2.1pt) -- (0,0) -- (-2.5pt,-2.1pt);
  \node[label,anchor=south] at (11.800,1.940) {$E_{3}$};
  \draw[side,draw=PrimeOne] (12.450,1.850) -- (12.450,0.550);
  \draw[tip,draw=PrimeOne,shift={(12.450,1.200)},rotate=-90.0] (-2.5pt,2.1pt) -- (0,0) -- (-2.5pt,-2.1pt);
  \node[label,anchor=west] at (12.540,1.200) {$E_{2}$};
  \draw[side,draw=PrimeTwo] (12.450,0.550) -- (11.150,0.550);
  \draw[tip,draw=PrimeTwo,shift={(11.800,0.550)},rotate=360.0] (-2.5pt,2.1pt) -- (0,0) -- (-2.5pt,-2.1pt);
  \node[label,anchor=north] at (11.800,0.460) {$E_{4}$};
  \node[corner] at (11.150,0.550) {};
  \node[corner] at (11.150,1.850) {};
  \node[corner] at (12.450,1.850) {};
  \node[corner] at (12.450,0.550) {};
  \node[inner sep=1.4pt,fill=SoftInk!5,rounded corners=1pt] at (11.800,1.200) {$F_{4}$};
  \draw[edge,draw=PrimeOne] (2.90,-0.75) -- (3.45,-0.75);
  \node[anchor=west] at (3.55,-0.75) {$\fp_3$-factor};
  \draw[edge,draw=PrimeTwo] (5.50,-0.75) -- (6.05,-0.75);
  \node[anchor=west] at (6.15,-0.75) {$\fq_3$-factor};
\end{tikzpicture}
    \caption{The square complex $X_{\Gamma_{\text{SV}}}$ of Stix--Vdovina \cite{stix-and-vdovina-2017-simply-transitive-quaternionic-lattices-of-rank-2-over-f-q-t-and-a-non-classical-fake-quadric}.}
    \label{fig:chi1-331-stix-vdovina}
\end{figure}
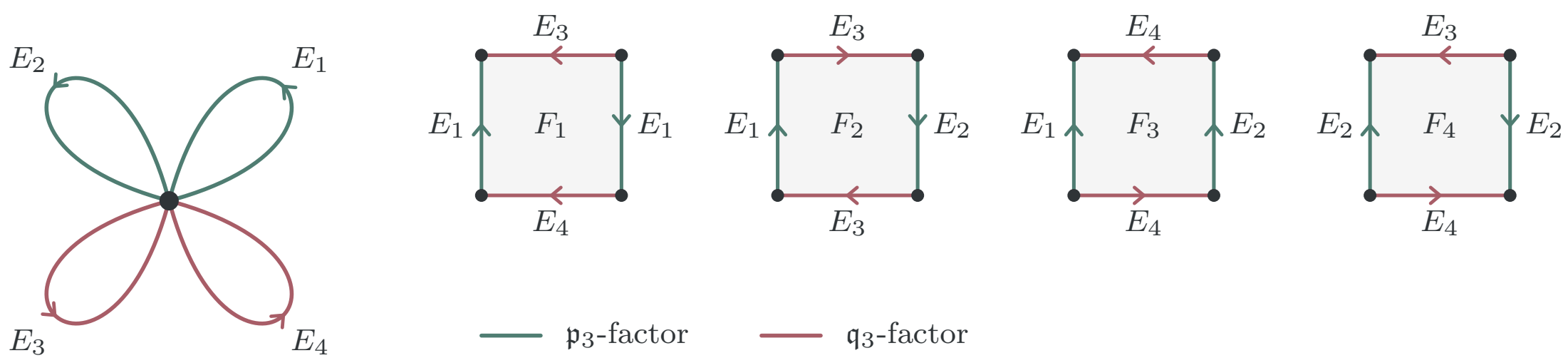

\begin{figure}[htbp]
    \centering
    \definecolor{SoftInk}{HTML}{2F3437}
\definecolor{EdgeColour}{HTML}{4F7D72}
\definecolor{PrimeOne}{HTML}{4F7D72}
\definecolor{PrimeTwo}{HTML}{A85D67}
\begin{tikzpicture}[x=1cm,y=1cm,
  every node/.style={font=\scriptsize,text=SoftInk},
  edge/.style={line width=1.05pt,line cap=round},
  side/.style={line width=0.95pt,line cap=round},
  tip/.style={line width=0.9pt,line cap=round,line join=round},
  cell/.style={fill=SoftInk!5},
  label/.style={inner sep=1.2pt},
  vert/.style={circle,fill=SoftInk,inner sep=1.8pt},
  corner/.style={circle,fill=SoftInk,inner sep=1.2pt}]
  \draw[edge,draw=EdgeColour] (0,0.5) .. controls (2.231,1.098) and (0.598,2.731) .. (0,0.5);
  \draw[tip,draw=EdgeColour,shift={(1.061,1.561)},rotate=135.0] (-2.5pt,2.1pt) -- (0,0) -- (-2.5pt,-2.1pt);
  \node[label,anchor=225] at (1.146,1.646) {$E_{1}$};
  \draw[edge,draw=EdgeColour] (0,0.5) .. controls (-0.598,2.731) and (-2.231,1.098) .. (0,0.5);
  \draw[tip,draw=EdgeColour,shift={(-1.061,1.561)},rotate=225.0] (-2.5pt,2.1pt) -- (0,0) -- (-2.5pt,-2.1pt);
  \node[label,anchor=315] at (-1.146,1.646) {$E_{2}$};
  \draw[edge,draw=EdgeColour] (0,0.5) .. controls (-2.231,-0.098) and (-0.598,-1.731) .. (0,0.5);
  \draw[tip,draw=EdgeColour,shift={(-1.061,-0.561)},rotate=315.0] (-2.5pt,2.1pt) -- (0,0) -- (-2.5pt,-2.1pt);
  \node[label,anchor=405] at (-1.146,-0.646) {$E_{3}$};
  \draw[edge,draw=EdgeColour] (0,0.5) .. controls (0.598,-1.731) and (2.231,-0.098) .. (0,0.5);
  \draw[tip,draw=EdgeColour,shift={(1.061,-0.561)},rotate=405.0] (-2.5pt,2.1pt) -- (0,0) -- (-2.5pt,-2.1pt);
  \node[label,anchor=495] at (1.146,-0.646) {$E_{4}$};
  \node[vert] at (0,0.5) {};
  \fill[cell] (2.900,0.550) rectangle (4.200,1.850);
  \draw[side,draw=EdgeColour] (2.900,0.550) -- (2.900,1.850);
  \draw[tip,draw=EdgeColour,shift={(2.900,1.200)},rotate=90.0] (-2.5pt,2.1pt) -- (0,0) -- (-2.5pt,-2.1pt);
  \node[label,anchor=east] at (2.810,1.200) {$E_{1}$};
  \draw[side,draw=EdgeColour] (2.900,1.850) -- (4.200,1.850);
  \draw[tip,draw=EdgeColour,shift={(3.550,1.850)},rotate=0.0] (-2.5pt,2.1pt) -- (0,0) -- (-2.5pt,-2.1pt);
  \node[label,anchor=south] at (3.550,1.940) {$E_{1}$};
  \draw[side,draw=EdgeColour] (4.200,1.850) -- (4.200,0.550);
  \draw[tip,draw=EdgeColour,shift={(4.200,1.200)},rotate=-90.0] (-2.5pt,2.1pt) -- (0,0) -- (-2.5pt,-2.1pt);
  \node[label,anchor=west] at (4.290,1.200) {$E_{2}$};
  \draw[side,draw=EdgeColour] (4.200,0.550) -- (2.900,0.550);
  \draw[tip,draw=EdgeColour,shift={(3.550,0.550)},rotate=180.0] (-2.5pt,2.1pt) -- (0,0) -- (-2.5pt,-2.1pt);
  \node[label,anchor=north] at (3.550,0.460) {$E_{2}$};
  \node[corner] at (2.900,0.550) {};
  \node[corner] at (2.900,1.850) {};
  \node[corner] at (4.200,1.850) {};
  \node[corner] at (4.200,0.550) {};
  \node[inner sep=1.4pt,fill=SoftInk!5,rounded corners=1pt] at (3.550,1.200) {$F_{1}$};
  \fill[cell] (5.650,0.550) rectangle (6.950,1.850);
  \draw[side,draw=EdgeColour] (5.650,0.550) -- (5.650,1.850);
  \draw[tip,draw=EdgeColour,shift={(5.650,1.200)},rotate=90.0] (-2.5pt,2.1pt) -- (0,0) -- (-2.5pt,-2.1pt);
  \node[label,anchor=east] at (5.560,1.200) {$E_{1}$};
  \draw[side,draw=EdgeColour] (5.650,1.850) -- (6.950,1.850);
  \draw[tip,draw=EdgeColour,shift={(6.300,1.850)},rotate=0.0] (-2.5pt,2.1pt) -- (0,0) -- (-2.5pt,-2.1pt);
  \node[label,anchor=south] at (6.300,1.940) {$E_{3}$};
  \draw[side,draw=EdgeColour] (6.950,1.850) -- (6.950,0.550);
  \draw[tip,draw=EdgeColour,shift={(6.950,1.200)},rotate=-90.0] (-2.5pt,2.1pt) -- (0,0) -- (-2.5pt,-2.1pt);
  \node[label,anchor=west] at (7.040,1.200) {$E_{3}$};
  \draw[side,draw=EdgeColour] (6.950,0.550) -- (5.650,0.550);
  \draw[tip,draw=EdgeColour,shift={(6.300,0.550)},rotate=180.0] (-2.5pt,2.1pt) -- (0,0) -- (-2.5pt,-2.1pt);
  \node[label,anchor=north] at (6.300,0.460) {$E_{4}$};
  \node[corner] at (5.650,0.550) {};
  \node[corner] at (5.650,1.850) {};
  \node[corner] at (6.950,1.850) {};
  \node[corner] at (6.950,0.550) {};
  \node[inner sep=1.4pt,fill=SoftInk!5,rounded corners=1pt] at (6.300,1.200) {$F_{2}$};
  \fill[cell] (8.400,0.550) rectangle (9.700,1.850);
  \draw[side,draw=EdgeColour] (8.400,0.550) -- (8.400,1.850);
  \draw[tip,draw=EdgeColour,shift={(8.400,1.200)},rotate=90.0] (-2.5pt,2.1pt) -- (0,0) -- (-2.5pt,-2.1pt);
  \node[label,anchor=east] at (8.310,1.200) {$E_{1}$};
  \draw[side,draw=EdgeColour] (8.400,1.850) -- (9.700,1.850);
  \draw[tip,draw=EdgeColour,shift={(9.050,1.850)},rotate=0.0] (-2.5pt,2.1pt) -- (0,0) -- (-2.5pt,-2.1pt);
  \node[label,anchor=south] at (9.050,1.940) {$E_{4}$};
  \draw[side,draw=EdgeColour] (9.700,1.850) -- (9.700,0.550);
  \draw[tip,draw=EdgeColour,shift={(9.700,1.200)},rotate=-90.0] (-2.5pt,2.1pt) -- (0,0) -- (-2.5pt,-2.1pt);
  \node[label,anchor=west] at (9.790,1.200) {$E_{4}$};
  \draw[side,draw=EdgeColour] (9.700,0.550) -- (8.400,0.550);
  \draw[tip,draw=EdgeColour,shift={(9.050,0.550)},rotate=180.0] (-2.5pt,2.1pt) -- (0,0) -- (-2.5pt,-2.1pt);
  \node[label,anchor=north] at (9.050,0.460) {$E_{3}$};
  \node[corner] at (8.400,0.550) {};
  \node[corner] at (8.400,1.850) {};
  \node[corner] at (9.700,1.850) {};
  \node[corner] at (9.700,0.550) {};
  \node[inner sep=1.4pt,fill=SoftInk!5,rounded corners=1pt] at (9.050,1.200) {$F_{3}$};
  \fill[cell] (11.150,0.550) rectangle (12.450,1.850);
  \draw[side,draw=EdgeColour] (11.150,0.550) -- (11.150,1.850);
  \draw[tip,draw=EdgeColour,shift={(11.150,1.200)},rotate=90.0] (-2.5pt,2.1pt) -- (0,0) -- (-2.5pt,-2.1pt);
  \node[label,anchor=east] at (11.060,1.200) {$E_{2}$};
  \draw[side,draw=EdgeColour] (11.150,1.850) -- (12.450,1.850);
  \draw[tip,draw=EdgeColour,shift={(11.800,1.850)},rotate=0.0] (-2.5pt,2.1pt) -- (0,0) -- (-2.5pt,-2.1pt);
  \node[label,anchor=south] at (11.800,1.940) {$E_{3}$};
  \draw[side,draw=EdgeColour] (12.450,1.850) -- (12.450,0.550);
  \draw[tip,draw=EdgeColour,shift={(12.450,1.200)},rotate=-90.0] (-2.5pt,2.1pt) -- (0,0) -- (-2.5pt,-2.1pt);
  \node[label,anchor=west] at (12.540,1.200) {$E_{2}$};
  \draw[side,draw=EdgeColour] (12.450,0.550) -- (11.150,0.550);
  \draw[tip,draw=EdgeColour,shift={(11.800,0.550)},rotate=180.0] (-2.5pt,2.1pt) -- (0,0) -- (-2.5pt,-2.1pt);
  \node[label,anchor=north] at (11.800,0.460) {$E_{4}$};
  \node[corner] at (11.150,0.550) {};
  \node[corner] at (11.150,1.850) {};
  \node[corner] at (12.450,1.850) {};
  \node[corner] at (12.450,0.550) {};
  \node[inner sep=1.4pt,fill=SoftInk!5,rounded corners=1pt] at (11.800,1.200) {$F_{4}$};
\end{tikzpicture}
    \caption{The square complex $X_{\Gamma_{\text{JW2}}}$ of Janzen--Wise \cite[Section 2]{janzen-and-wise-2009-a-smallest-irreducible-lattice-in-the-product-of-trees}.}
    \label{fig:chi1-331-janzen-wise}
\end{figure}
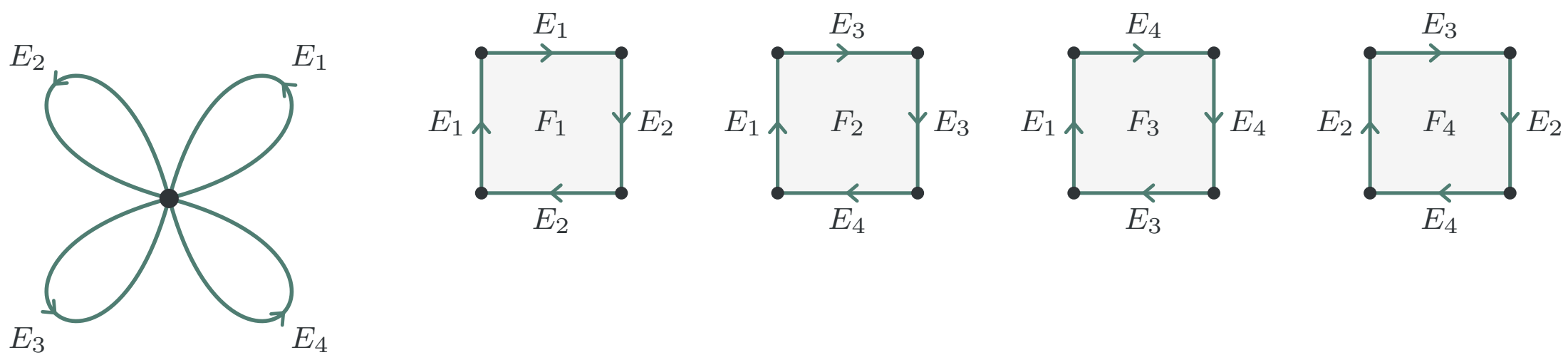

\begin{figure}[htbp]
    \centering
    \definecolor{SoftInk}{HTML}{2F3437}
\definecolor{EdgeColour}{HTML}{4F7D72}
\definecolor{PrimeOne}{HTML}{4F7D72}
\definecolor{PrimeTwo}{HTML}{A85D67}
\begin{tikzpicture}[x=1cm,y=1cm,
  every node/.style={font=\scriptsize,text=SoftInk},
  edge/.style={line width=1.05pt,line cap=round},
  side/.style={line width=0.95pt,line cap=round},
  tip/.style={line width=0.9pt,line cap=round,line join=round},
  cell/.style={fill=SoftInk!5},
  label/.style={inner sep=1.2pt},
  vert/.style={circle,fill=SoftInk,inner sep=1.8pt},
  corner/.style={circle,fill=SoftInk,inner sep=1.2pt}]
  \draw[edge,draw=EdgeColour] (0,0.5) .. controls (2.231,1.098) and (0.598,2.731) .. (0,0.5);
  \draw[tip,draw=EdgeColour,shift={(1.061,1.561)},rotate=135.0] (-2.5pt,2.1pt) -- (0,0) -- (-2.5pt,-2.1pt);
  \node[label,anchor=225] at (1.146,1.646) {$E_{1}$};
  \draw[edge,draw=EdgeColour] (0,0.5) .. controls (-0.598,2.731) and (-2.231,1.098) .. (0,0.5);
  \draw[tip,draw=EdgeColour,shift={(-1.061,1.561)},rotate=225.0] (-2.5pt,2.1pt) -- (0,0) -- (-2.5pt,-2.1pt);
  \node[label,anchor=315] at (-1.146,1.646) {$E_{2}$};
  \draw[edge,draw=EdgeColour] (0,0.5) .. controls (-2.231,-0.098) and (-0.598,-1.731) .. (0,0.5);
  \draw[tip,draw=EdgeColour,shift={(-1.061,-0.561)},rotate=315.0] (-2.5pt,2.1pt) -- (0,0) -- (-2.5pt,-2.1pt);
  \node[label,anchor=405] at (-1.146,-0.646) {$E_{3}$};
  \draw[edge,draw=EdgeColour] (0,0.5) .. controls (0.598,-1.731) and (2.231,-0.098) .. (0,0.5);
  \draw[tip,draw=EdgeColour,shift={(1.061,-0.561)},rotate=405.0] (-2.5pt,2.1pt) -- (0,0) -- (-2.5pt,-2.1pt);
  \node[label,anchor=495] at (1.146,-0.646) {$E_{4}$};
  \node[vert] at (0,0.5) {};
  \fill[cell] (2.900,0.550) rectangle (4.200,1.850);
  \draw[side,draw=EdgeColour] (2.900,0.550) -- (2.900,1.850);
  \draw[tip,draw=EdgeColour,shift={(2.900,1.200)},rotate=90.0] (-2.5pt,2.1pt) -- (0,0) -- (-2.5pt,-2.1pt);
  \node[label,anchor=east] at (2.810,1.200) {$E_{1}$};
  \draw[side,draw=EdgeColour] (2.900,1.850) -- (4.200,1.850);
  \draw[tip,draw=EdgeColour,shift={(3.550,1.850)},rotate=0.0] (-2.5pt,2.1pt) -- (0,0) -- (-2.5pt,-2.1pt);
  \node[label,anchor=south] at (3.550,1.940) {$E_{1}$};
  \draw[side,draw=EdgeColour] (4.200,1.850) -- (4.200,0.550);
  \draw[tip,draw=EdgeColour,shift={(4.200,1.200)},rotate=-90.0] (-2.5pt,2.1pt) -- (0,0) -- (-2.5pt,-2.1pt);
  \node[label,anchor=west] at (4.290,1.200) {$E_{2}$};
  \draw[side,draw=EdgeColour] (4.200,0.550) -- (2.900,0.550);
  \draw[tip,draw=EdgeColour,shift={(3.550,0.550)},rotate=180.0] (-2.5pt,2.1pt) -- (0,0) -- (-2.5pt,-2.1pt);
  \node[label,anchor=north] at (3.550,0.460) {$E_{3}$};
  \node[corner] at (2.900,0.550) {};
  \node[corner] at (2.900,1.850) {};
  \node[corner] at (4.200,1.850) {};
  \node[corner] at (4.200,0.550) {};
  \node[inner sep=1.4pt,fill=SoftInk!5,rounded corners=1pt] at (3.550,1.200) {$F_{1}$};
  \fill[cell] (5.650,0.550) rectangle (6.950,1.850);
  \draw[side,draw=EdgeColour] (5.650,0.550) -- (5.650,1.850);
  \draw[tip,draw=EdgeColour,shift={(5.650,1.200)},rotate=90.0] (-2.5pt,2.1pt) -- (0,0) -- (-2.5pt,-2.1pt);
  \node[label,anchor=east] at (5.560,1.200) {$E_{1}$};
  \draw[side,draw=EdgeColour] (5.650,1.850) -- (6.950,1.850);
  \draw[tip,draw=EdgeColour,shift={(6.300,1.850)},rotate=0.0] (-2.5pt,2.1pt) -- (0,0) -- (-2.5pt,-2.1pt);
  \node[label,anchor=south] at (6.300,1.940) {$E_{4}$};
  \draw[side,draw=EdgeColour] (6.950,1.850) -- (6.950,0.550);
  \draw[tip,draw=EdgeColour,shift={(6.950,1.200)},rotate=-90.0] (-2.5pt,2.1pt) -- (0,0) -- (-2.5pt,-2.1pt);
  \node[label,anchor=west] at (7.040,1.200) {$E_{2}$};
  \draw[side,draw=EdgeColour] (6.950,0.550) -- (5.650,0.550);
  \draw[tip,draw=EdgeColour,shift={(6.300,0.550)},rotate=180.0] (-2.5pt,2.1pt) -- (0,0) -- (-2.5pt,-2.1pt);
  \node[label,anchor=north] at (6.300,0.460) {$E_{2}$};
  \node[corner] at (5.650,0.550) {};
  \node[corner] at (5.650,1.850) {};
  \node[corner] at (6.950,1.850) {};
  \node[corner] at (6.950,0.550) {};
  \node[inner sep=1.4pt,fill=SoftInk!5,rounded corners=1pt] at (6.300,1.200) {$F_{2}$};
  \fill[cell] (8.400,0.550) rectangle (9.700,1.850);
  \draw[side,draw=EdgeColour] (8.400,0.550) -- (8.400,1.850);
  \draw[tip,draw=EdgeColour,shift={(8.400,1.200)},rotate=90.0] (-2.5pt,2.1pt) -- (0,0) -- (-2.5pt,-2.1pt);
  \node[label,anchor=east] at (8.310,1.200) {$E_{2}$};
  \draw[side,draw=EdgeColour] (8.400,1.850) -- (9.700,1.850);
  \draw[tip,draw=EdgeColour,shift={(9.050,1.850)},rotate=0.0] (-2.5pt,2.1pt) -- (0,0) -- (-2.5pt,-2.1pt);
  \node[label,anchor=south] at (9.050,1.940) {$E_{4}$};
  \draw[side,draw=EdgeColour] (9.700,1.850) -- (9.700,0.550);
  \draw[tip,draw=EdgeColour,shift={(9.700,1.200)},rotate=-90.0] (-2.5pt,2.1pt) -- (0,0) -- (-2.5pt,-2.1pt);
  \node[label,anchor=west] at (9.790,1.200) {$E_{3}$};
  \draw[side,draw=EdgeColour] (9.700,0.550) -- (8.400,0.550);
  \draw[tip,draw=EdgeColour,shift={(9.050,0.550)},rotate=180.0] (-2.5pt,2.1pt) -- (0,0) -- (-2.5pt,-2.1pt);
  \node[label,anchor=north] at (9.050,0.460) {$E_{3}$};
  \node[corner] at (8.400,0.550) {};
  \node[corner] at (8.400,1.850) {};
  \node[corner] at (9.700,1.850) {};
  \node[corner] at (9.700,0.550) {};
  \node[inner sep=1.4pt,fill=SoftInk!5,rounded corners=1pt] at (9.050,1.200) {$F_{3}$};
  \fill[cell] (11.150,0.550) rectangle (12.450,1.850);
  \draw[side,draw=EdgeColour] (11.150,0.550) -- (11.150,1.850);
  \draw[tip,draw=EdgeColour,shift={(11.150,1.200)},rotate=90.0] (-2.5pt,2.1pt) -- (0,0) -- (-2.5pt,-2.1pt);
  \node[label,anchor=east] at (11.060,1.200) {$E_{1}$};
  \draw[side,draw=EdgeColour] (11.150,1.850) -- (12.450,1.850);
  \draw[tip,draw=EdgeColour,shift={(11.800,1.850)},rotate=0.0] (-2.5pt,2.1pt) -- (0,0) -- (-2.5pt,-2.1pt);
  \node[label,anchor=south] at (11.800,1.940) {$E_{3}$};
  \draw[side,draw=EdgeColour] (12.450,1.850) -- (12.450,0.550);
  \draw[tip,draw=EdgeColour,shift={(12.450,1.200)},rotate=-90.0] (-2.5pt,2.1pt) -- (0,0) -- (-2.5pt,-2.1pt);
  \node[label,anchor=west] at (12.540,1.200) {$E_{4}$};
  \draw[side,draw=EdgeColour] (12.450,0.550) -- (11.150,0.550);
  \draw[tip,draw=EdgeColour,shift={(11.800,0.550)},rotate=180.0] (-2.5pt,2.1pt) -- (0,0) -- (-2.5pt,-2.1pt);
  \node[label,anchor=north] at (11.800,0.460) {$E_{4}$};
  \node[corner] at (11.150,0.550) {};
  \node[corner] at (11.150,1.850) {};
  \node[corner] at (12.450,1.850) {};
  \node[corner] at (12.450,0.550) {};
  \node[inner sep=1.4pt,fill=SoftInk!5,rounded corners=1pt] at (11.800,1.200) {$F_{4}$};
\end{tikzpicture}
    \caption{The third $(3,3,1)$ square complex in the proof of Theorem \ref{thm:C}.}
    \label{fig:chi1-331-swap}
\end{figure}
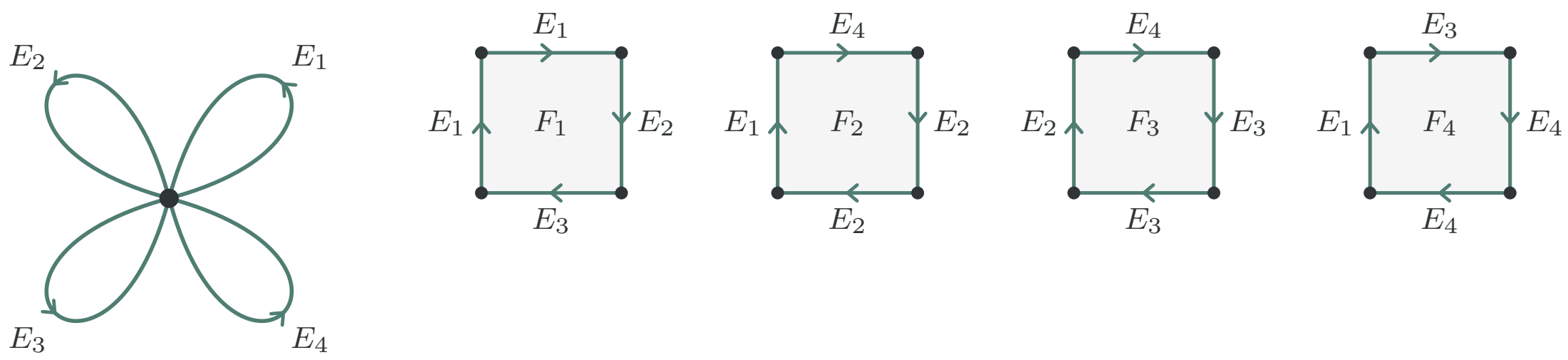

\begin{figure}[htbp]
    \centering
    \input{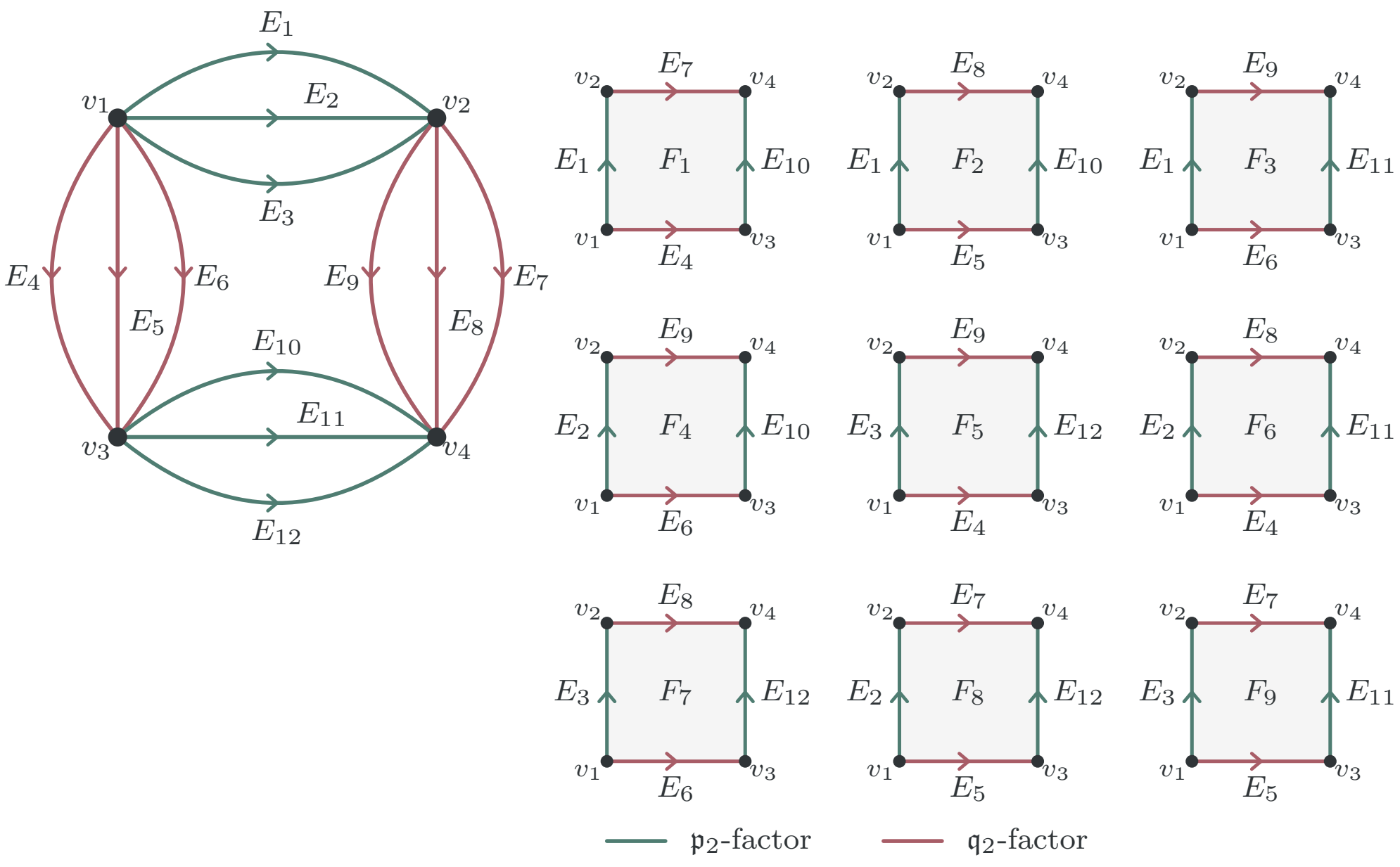}
    \caption{The square complex $X_{\Gamma_{\text{R}}}$ of Rungtanapirom \cite{rungtanapirom-2018-quaternionic-arithmetic-lattices-of-rank-2-and-a-fake-quadric-in-characteristic-2}.}
    \label{fig:chi1-224-rungtanapirom}
\end{figure}

\begin{figure}[htbp]
    \centering
    \input{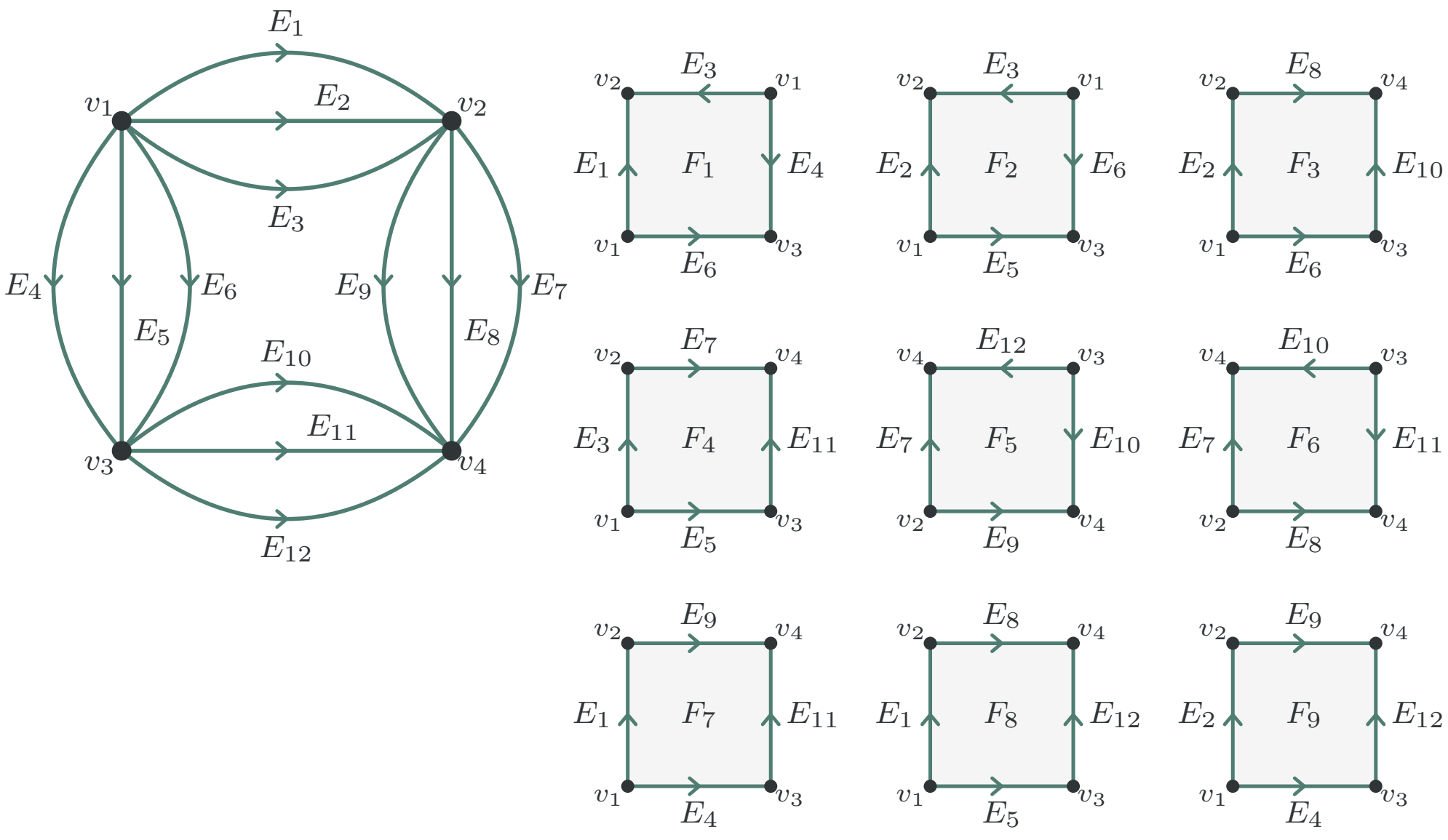}
    \caption{The second $(2,2,4)$ square complex in the proof of Theorem \ref{thm:C}.}
    \label{fig:chi1-224-swap}
\end{figure}

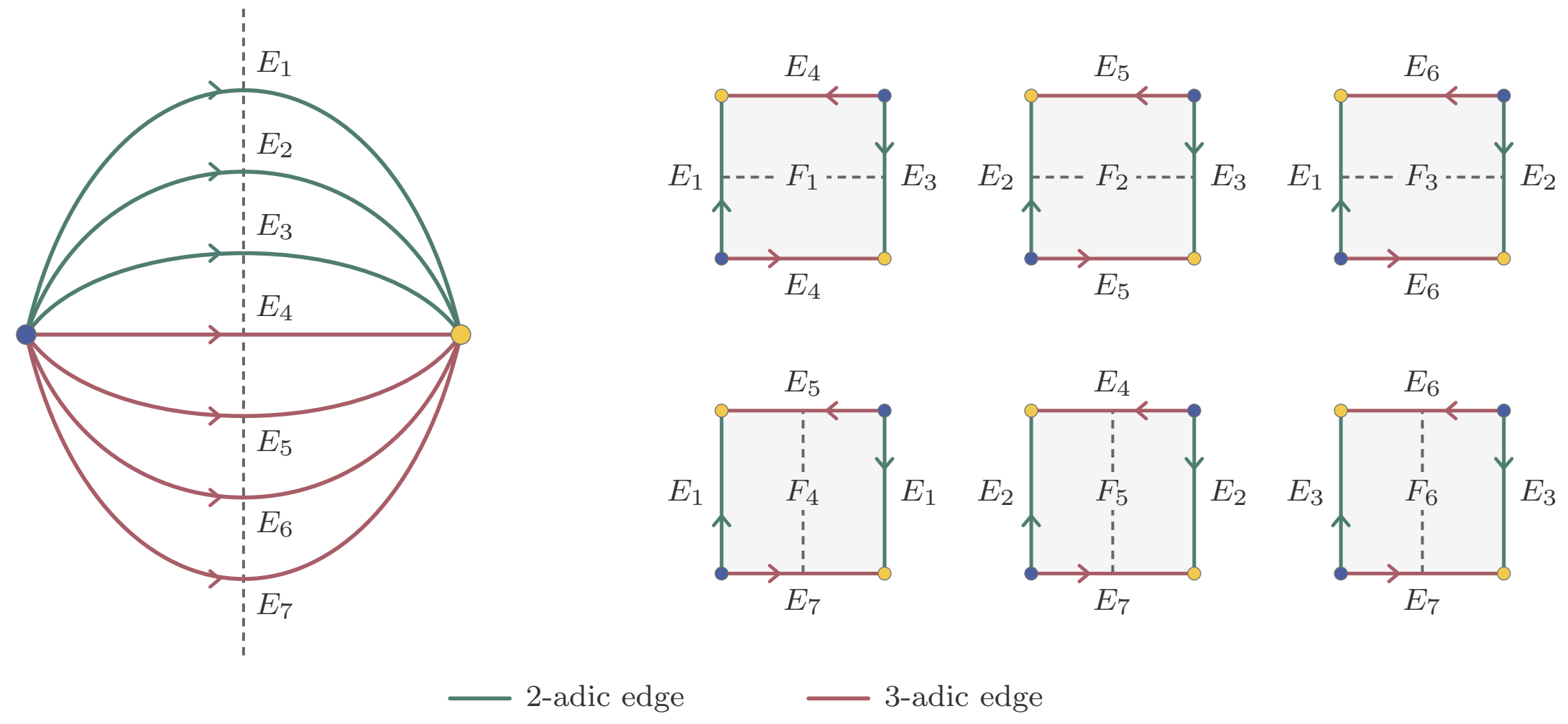
\begin{figure}[htbp]
    \centering
    \definecolor{SoftInk}{HTML}{2F3437}
\definecolor{TwoAdic}{HTML}{4F7D72}
\definecolor{ThreeAdic}{HTML}{A85D67}
\definecolor{VertexOne}{HTML}{4B5F9E}
\definecolor{VertexTwo}{HTML}{F2C94C}
\begin{tikzpicture}[x=1cm,y=1cm,
  every node/.style={font=\scriptsize,text=SoftInk},
  edge/.style={line width=1.05pt,line cap=round},
  side/.style={line width=0.95pt,line cap=round},
  tip/.style={line width=0.9pt,line cap=round,line join=round},
  cell/.style={fill=SoftInk!5},
  mirror/.style={draw=SoftInk!75,line width=0.75pt,dash pattern=on 2.2pt off 1.6pt},
  label/.style={inner sep=1.2pt},
  vone/.style={circle,fill=VertexOne,draw=SoftInk!70,line width=0.3pt,inner sep=1.8pt},
  vtwo/.style={circle,fill=VertexTwo,draw=SoftInk!70,line width=0.3pt,inner sep=1.8pt},
  cone/.style={circle,fill=VertexOne,draw=SoftInk!70,line width=0.25pt,inner sep=1.2pt},
  ctwo/.style={circle,fill=VertexTwo,draw=SoftInk!70,line width=0.25pt,inner sep=1.2pt}]
  \draw[mirror] (2,3) -- (2,-3);
  \draw[edge,draw=TwoAdic] (0,0) .. controls (0.68,3) and (3.32,3) .. (4,0);
  \draw[tip,TwoAdic,shift={(1.78,2.25)},rotate=0] (-2.5pt,2.1pt) -- (0,0) -- (-2.5pt,-2.1pt);
  \node[label,anchor=south west] at (2.06,2.33) {$E_{1}$};
  \draw[edge,draw=TwoAdic] (0,0) .. controls (0.68,2) and (3.32,2) .. (4,0);
  \draw[tip,TwoAdic,shift={(1.78,1.5)},rotate=0] (-2.5pt,2.1pt) -- (0,0) -- (-2.5pt,-2.1pt);
  \node[label,anchor=south west] at (2.06,1.58) {$E_{2}$};
  \draw[edge,draw=TwoAdic] (0,0) .. controls (0.68,1) and (3.32,1) .. (4,0);
  \draw[tip,TwoAdic,shift={(1.78,0.75)},rotate=0] (-2.5pt,2.1pt) -- (0,0) -- (-2.5pt,-2.1pt);
  \node[label,anchor=south west] at (2.06,0.83) {$E_{3}$};
  \draw[edge,draw=ThreeAdic] (0,0) .. controls (0.68,0) and (3.32,0) .. (4,0);
  \draw[tip,ThreeAdic,shift={(1.78,0)},rotate=0] (-2.5pt,2.1pt) -- (0,0) -- (-2.5pt,-2.1pt);
  \node[label,anchor=south west] at (2.06,0.08) {$E_{4}$};
  \draw[edge,draw=ThreeAdic] (0,0) .. controls (0.68,-1) and (3.32,-1) .. (4,0);
  \draw[tip,ThreeAdic,shift={(1.78,-0.75)},rotate=0] (-2.5pt,2.1pt) -- (0,0) -- (-2.5pt,-2.1pt);
  \node[label,anchor=north west] at (2.06,-0.83) {$E_{5}$};
  \draw[edge,draw=ThreeAdic] (0,0) .. controls (0.68,-2) and (3.32,-2) .. (4,0);
  \draw[tip,ThreeAdic,shift={(1.78,-1.5)},rotate=0] (-2.5pt,2.1pt) -- (0,0) -- (-2.5pt,-2.1pt);
  \node[label,anchor=north west] at (2.06,-1.58) {$E_{6}$};
  \draw[edge,draw=ThreeAdic] (0,0) .. controls (0.68,-3) and (3.32,-3) .. (4,0);
  \draw[tip,ThreeAdic,shift={(1.78,-2.25)},rotate=0] (-2.5pt,2.1pt) -- (0,0) -- (-2.5pt,-2.1pt);
  \node[label,anchor=north west] at (2.06,-2.33) {$E_{7}$};
  \node[vone] at (0,0) {};
  \node[vtwo] at (4,0) {};
  \fill[cell] (6.4,0.7) rectangle (7.9,2.2);
  \draw[mirror] (6.4,1.45) -- (7.9,1.45);
  \draw[side,draw=TwoAdic] (6.4,0.7) -- (6.4,2.2);
  \draw[tip,TwoAdic,shift={(6.4,1.23)},rotate=90] (-2.5pt,2.1pt) -- (0,0) -- (-2.5pt,-2.1pt);
  \node[label,anchor=east] at (6.31,1.45) {$E_{1}$};
  \draw[side,draw=ThreeAdic] (6.4,2.2) -- (7.9,2.2);
  \draw[tip,ThreeAdic,shift={(7.37,2.2)},rotate=180] (-2.5pt,2.1pt) -- (0,0) -- (-2.5pt,-2.1pt);
  \node[label,anchor=south] at (7.15,2.29) {$E_{4}$};
  \draw[side,draw=TwoAdic] (7.9,2.2) -- (7.9,0.7);
  \draw[tip,TwoAdic,shift={(7.9,1.67)},rotate=-90] (-2.5pt,2.1pt) -- (0,0) -- (-2.5pt,-2.1pt);
  \node[label,anchor=west] at (7.99,1.45) {$E_{3}$};
  \draw[side,draw=ThreeAdic] (7.9,0.7) -- (6.4,0.7);
  \draw[tip,ThreeAdic,shift={(6.93,0.7)},rotate=0] (-2.5pt,2.1pt) -- (0,0) -- (-2.5pt,-2.1pt);
  \node[label,anchor=north] at (7.15,0.61) {$E_{4}$};
  \node[cone] at (6.4,0.7) {};
  \node[ctwo] at (6.4,2.2) {};
  \node[cone] at (7.9,2.2) {};
  \node[ctwo] at (7.9,0.7) {};
  \node[inner sep=1.4pt,fill=SoftInk!5,rounded corners=1pt] at (7.15,1.45) {$F_{1}$};
  \fill[cell] (9.25,0.7) rectangle (10.75,2.2);
  \draw[mirror] (9.25,1.45) -- (10.75,1.45);
  \draw[side,draw=TwoAdic] (9.25,0.7) -- (9.25,2.2);
  \draw[tip,TwoAdic,shift={(9.25,1.23)},rotate=90] (-2.5pt,2.1pt) -- (0,0) -- (-2.5pt,-2.1pt);
  \node[label,anchor=east] at (9.16,1.45) {$E_{2}$};
  \draw[side,draw=ThreeAdic] (9.25,2.2) -- (10.75,2.2);
  \draw[tip,ThreeAdic,shift={(10.22,2.2)},rotate=180] (-2.5pt,2.1pt) -- (0,0) -- (-2.5pt,-2.1pt);
  \node[label,anchor=south] at (10,2.29) {$E_{5}$};
  \draw[side,draw=TwoAdic] (10.75,2.2) -- (10.75,0.7);
  \draw[tip,TwoAdic,shift={(10.75,1.67)},rotate=-90] (-2.5pt,2.1pt) -- (0,0) -- (-2.5pt,-2.1pt);
  \node[label,anchor=west] at (10.84,1.45) {$E_{3}$};
  \draw[side,draw=ThreeAdic] (10.75,0.7) -- (9.25,0.7);
  \draw[tip,ThreeAdic,shift={(9.78,0.7)},rotate=0] (-2.5pt,2.1pt) -- (0,0) -- (-2.5pt,-2.1pt);
  \node[label,anchor=north] at (10,0.61) {$E_{5}$};
  \node[cone] at (9.25,0.7) {};
  \node[ctwo] at (9.25,2.2) {};
  \node[cone] at (10.75,2.2) {};
  \node[ctwo] at (10.75,0.7) {};
  \node[inner sep=1.4pt,fill=SoftInk!5,rounded corners=1pt] at (10,1.45) {$F_{2}$};
  \fill[cell] (12.1,0.7) rectangle (13.6,2.2);
  \draw[mirror] (12.1,1.45) -- (13.6,1.45);
  \draw[side,draw=TwoAdic] (12.1,0.7) -- (12.1,2.2);
  \draw[tip,TwoAdic,shift={(12.1,1.23)},rotate=90] (-2.5pt,2.1pt) -- (0,0) -- (-2.5pt,-2.1pt);
  \node[label,anchor=east] at (12.01,1.45) {$E_{1}$};
  \draw[side,draw=ThreeAdic] (12.1,2.2) -- (13.6,2.2);
  \draw[tip,ThreeAdic,shift={(13.07,2.2)},rotate=180] (-2.5pt,2.1pt) -- (0,0) -- (-2.5pt,-2.1pt);
  \node[label,anchor=south] at (12.85,2.29) {$E_{6}$};
  \draw[side,draw=TwoAdic] (13.6,2.2) -- (13.6,0.7);
  \draw[tip,TwoAdic,shift={(13.6,1.67)},rotate=-90] (-2.5pt,2.1pt) -- (0,0) -- (-2.5pt,-2.1pt);
  \node[label,anchor=west] at (13.69,1.45) {$E_{2}$};
  \draw[side,draw=ThreeAdic] (13.6,0.7) -- (12.1,0.7);
  \draw[tip,ThreeAdic,shift={(12.63,0.7)},rotate=0] (-2.5pt,2.1pt) -- (0,0) -- (-2.5pt,-2.1pt);
  \node[label,anchor=north] at (12.85,0.61) {$E_{6}$};
  \node[cone] at (12.1,0.7) {};
  \node[ctwo] at (12.1,2.2) {};
  \node[cone] at (13.6,2.2) {};
  \node[ctwo] at (13.6,0.7) {};
  \node[inner sep=1.4pt,fill=SoftInk!5,rounded corners=1pt] at (12.85,1.45) {$F_{3}$};
  \fill[cell] (6.4,-2.2) rectangle (7.9,-0.7);
  \draw[mirror] (7.15,-2.2) -- (7.15,-0.7);
  \draw[side,draw=TwoAdic] (6.4,-2.2) -- (6.4,-0.7);
  \draw[tip,TwoAdic,shift={(6.4,-1.67)},rotate=90] (-2.5pt,2.1pt) -- (0,0) -- (-2.5pt,-2.1pt);
  \node[label,anchor=east] at (6.31,-1.45) {$E_{1}$};
  \draw[side,draw=ThreeAdic] (6.4,-0.7) -- (7.9,-0.7);
  \draw[tip,ThreeAdic,shift={(7.37,-0.7)},rotate=180] (-2.5pt,2.1pt) -- (0,0) -- (-2.5pt,-2.1pt);
  \node[label,anchor=south] at (7.15,-0.61) {$E_{5}$};
  \draw[side,draw=TwoAdic] (7.9,-0.7) -- (7.9,-2.2);
  \draw[tip,TwoAdic,shift={(7.9,-1.23)},rotate=-90] (-2.5pt,2.1pt) -- (0,0) -- (-2.5pt,-2.1pt);
  \node[label,anchor=west] at (7.99,-1.45) {$E_{1}$};
  \draw[side,draw=ThreeAdic] (7.9,-2.2) -- (6.4,-2.2);
  \draw[tip,ThreeAdic,shift={(6.93,-2.2)},rotate=0] (-2.5pt,2.1pt) -- (0,0) -- (-2.5pt,-2.1pt);
  \node[label,anchor=north] at (7.15,-2.29) {$E_{7}$};
  \node[cone] at (6.4,-2.2) {};
  \node[ctwo] at (6.4,-0.7) {};
  \node[cone] at (7.9,-0.7) {};
  \node[ctwo] at (7.9,-2.2) {};
  \node[inner sep=1.4pt,fill=SoftInk!5,rounded corners=1pt] at (7.15,-1.45) {$F_{4}$};
  \fill[cell] (9.25,-2.2) rectangle (10.75,-0.7);
  \draw[mirror] (10,-2.2) -- (10,-0.7);
  \draw[side,draw=TwoAdic] (9.25,-2.2) -- (9.25,-0.7);
  \draw[tip,TwoAdic,shift={(9.25,-1.67)},rotate=90] (-2.5pt,2.1pt) -- (0,0) -- (-2.5pt,-2.1pt);
  \node[label,anchor=east] at (9.16,-1.45) {$E_{2}$};
  \draw[side,draw=ThreeAdic] (9.25,-0.7) -- (10.75,-0.7);
  \draw[tip,ThreeAdic,shift={(10.22,-0.7)},rotate=180] (-2.5pt,2.1pt) -- (0,0) -- (-2.5pt,-2.1pt);
  \node[label,anchor=south] at (10,-0.61) {$E_{4}$};
  \draw[side,draw=TwoAdic] (10.75,-0.7) -- (10.75,-2.2);
  \draw[tip,TwoAdic,shift={(10.75,-1.23)},rotate=-90] (-2.5pt,2.1pt) -- (0,0) -- (-2.5pt,-2.1pt);
  \node[label,anchor=west] at (10.84,-1.45) {$E_{2}$};
  \draw[side,draw=ThreeAdic] (10.75,-2.2) -- (9.25,-2.2);
  \draw[tip,ThreeAdic,shift={(9.78,-2.2)},rotate=0] (-2.5pt,2.1pt) -- (0,0) -- (-2.5pt,-2.1pt);
  \node[label,anchor=north] at (10,-2.29) {$E_{7}$};
  \node[cone] at (9.25,-2.2) {};
  \node[ctwo] at (9.25,-0.7) {};
  \node[cone] at (10.75,-0.7) {};
  \node[ctwo] at (10.75,-2.2) {};
  \node[inner sep=1.4pt,fill=SoftInk!5,rounded corners=1pt] at (10,-1.45) {$F_{5}$};
  \fill[cell] (12.1,-2.2) rectangle (13.6,-0.7);
  \draw[mirror] (12.85,-2.2) -- (12.85,-0.7);
  \draw[side,draw=TwoAdic] (12.1,-2.2) -- (12.1,-0.7);
  \draw[tip,TwoAdic,shift={(12.1,-1.67)},rotate=90] (-2.5pt,2.1pt) -- (0,0) -- (-2.5pt,-2.1pt);
  \node[label,anchor=east] at (12.01,-1.45) {$E_{3}$};
  \draw[side,draw=ThreeAdic] (12.1,-0.7) -- (13.6,-0.7);
  \draw[tip,ThreeAdic,shift={(13.07,-0.7)},rotate=180] (-2.5pt,2.1pt) -- (0,0) -- (-2.5pt,-2.1pt);
  \node[label,anchor=south] at (12.85,-0.61) {$E_{6}$};
  \draw[side,draw=TwoAdic] (13.6,-0.7) -- (13.6,-2.2);
  \draw[tip,TwoAdic,shift={(13.6,-1.23)},rotate=-90] (-2.5pt,2.1pt) -- (0,0) -- (-2.5pt,-2.1pt);
  \node[label,anchor=west] at (13.69,-1.45) {$E_{3}$};
  \draw[side,draw=ThreeAdic] (13.6,-2.2) -- (12.1,-2.2);
  \draw[tip,ThreeAdic,shift={(12.63,-2.2)},rotate=0] (-2.5pt,2.1pt) -- (0,0) -- (-2.5pt,-2.1pt);
  \node[label,anchor=north] at (12.85,-2.29) {$E_{7}$};
  \node[cone] at (12.1,-2.2) {};
  \node[ctwo] at (12.1,-0.7) {};
  \node[cone] at (13.6,-0.7) {};
  \node[ctwo] at (13.6,-2.2) {};
  \node[inner sep=1.4pt,fill=SoftInk!5,rounded corners=1pt] at (12.85,-1.45) {$F_{6}$};
  \draw[edge,draw=TwoAdic] (3.9,-3.35) -- (4.45,-3.35);
  \node[label,anchor=west] at (4.55,-3.35) {$2$-adic edge};
  \draw[edge,draw=ThreeAdic] (7.2,-3.35) -- (7.75,-3.35);
  \node[label,anchor=west] at (7.85,-3.35) {$3$-adic edge};
\end{tikzpicture}
    \caption{The square complex $X_\Lambda$ of Section \ref{sec:q-disc5}. The covering $X_\Lambda \to X_\Gamma$ is the quotient by the involution given by reflection in the dashed lines. This involution inverts every edge of $X_\Lambda$.}
    \label{fig:q-disc5-covering}
\end{figure}

\begin{figure}[htbp]
    \centering
    \input{figures/rational_disc13_covering}
    \caption{The covering $X_\Lambda \to X_\Gamma$ of Section \ref{sec:q-disc13}.}
    \label{fig:q-disc13-covering}
\end{figure}

\begin{figure}[htbp]
    \centering
    \definecolor{SoftInk}{HTML}{2F3437}
\definecolor{TwoAdic}{HTML}{4F7D72}
\definecolor{ThreeAdic}{HTML}{A85D67}
\definecolor{VertexOne}{HTML}{4B5F9E}
\definecolor{VertexTwo}{HTML}{F2C94C}
\begin{tikzpicture}[x=1cm,y=1cm,
  every node/.style={font=\scriptsize,text=SoftInk},
  edge/.style={line width=1.05pt,line cap=round},
  side/.style={line width=0.95pt,line cap=round},
  tip/.style={line width=0.9pt,line cap=round,line join=round},
  cell/.style={fill=SoftInk!5},
  label/.style={inner sep=1.2pt},
  title/.style={font=\normalsize,text=SoftInk},
  vone/.style={circle,fill=VertexOne,draw=SoftInk!70,line width=0.3pt,inner sep=1.8pt},
  vtwo/.style={circle,fill=VertexTwo,draw=SoftInk!70,line width=0.3pt,inner sep=1.8pt},
  cone/.style={circle,fill=VertexOne,draw=SoftInk!70,line width=0.25pt,inner sep=1.2pt},
  ctwo/.style={circle,fill=VertexTwo,draw=SoftInk!70,line width=0.25pt,inner sep=1.2pt}]
  \node[title,anchor=east] at (-1.0,0) {$X_\Gamma$};
  \draw[edge,draw=TwoAdic] (0,0) .. controls (0.68,3) and (3.32,3) .. (4,0);
  \draw[tip,TwoAdic,shift={(2,2.25)},rotate=0] (-2.5pt,2.1pt) -- (0,0) -- (-2.5pt,-2.1pt);
  \node[label,anchor=south] at (2,2.35) {$E_{1}$};
  \draw[edge,draw=TwoAdic] (0,0) .. controls (0.68,2) and (3.32,2) .. (4,0);
  \draw[tip,TwoAdic,shift={(2,1.5)},rotate=0] (-2.5pt,2.1pt) -- (0,0) -- (-2.5pt,-2.1pt);
  \node[label,anchor=south] at (2,1.6) {$E_{2}$};
  \draw[edge,draw=TwoAdic] (0,0) .. controls (0.68,1) and (3.32,1) .. (4,0);
  \draw[tip,TwoAdic,shift={(2,0.75)},rotate=0] (-2.5pt,2.1pt) -- (0,0) -- (-2.5pt,-2.1pt);
  \node[label,anchor=south] at (2,0.85) {$E_{3}$};
  \draw[edge,draw=ThreeAdic] (0,0) .. controls (0.68,0) and (3.32,0) .. (4,0);
  \draw[tip,ThreeAdic,shift={(2,0)},rotate=0] (-2.5pt,2.1pt) -- (0,0) -- (-2.5pt,-2.1pt);
  \node[label,anchor=south] at (2,0.1) {$E_{4}$};
  \draw[edge,draw=ThreeAdic] (0,0) .. controls (0.68,-1) and (3.32,-1) .. (4,0);
  \draw[tip,ThreeAdic,shift={(2,-0.75)},rotate=0] (-2.5pt,2.1pt) -- (0,0) -- (-2.5pt,-2.1pt);
  \node[label,anchor=north] at (2,-0.85) {$E_{5}$};
  \draw[edge,draw=ThreeAdic] (0,0) .. controls (0.68,-2) and (3.32,-2) .. (4,0);
  \draw[tip,ThreeAdic,shift={(2,-1.5)},rotate=0] (-2.5pt,2.1pt) -- (0,0) -- (-2.5pt,-2.1pt);
  \node[label,anchor=north] at (2,-1.6) {$E_{6}$};
  \draw[edge,draw=ThreeAdic] (0,0) .. controls (0.68,-3) and (3.32,-3) .. (4,0);
  \draw[tip,ThreeAdic,shift={(2,-2.25)},rotate=0] (-2.5pt,2.1pt) -- (0,0) -- (-2.5pt,-2.1pt);
  \node[label,anchor=north] at (2,-2.35) {$E_{7}$};
  \node[vone] at (0,0) {};
  \node[vtwo] at (4,0) {};
  \fill[cell] (6.4,0.7) rectangle (7.9,2.2);
  \draw[side,draw=TwoAdic] (6.4,0.7) -- (6.4,2.2);
  \draw[tip,TwoAdic,shift={(6.4,1.45)},rotate=90] (-2.5pt,2.1pt) -- (0,0) -- (-2.5pt,-2.1pt);
  \node[label,anchor=east] at (6.31,1.45) {$E_{1}$};
  \draw[side,draw=ThreeAdic] (6.4,2.2) -- (7.9,2.2);
  \draw[tip,ThreeAdic,shift={(7.15,2.2)},rotate=180] (-2.5pt,2.1pt) -- (0,0) -- (-2.5pt,-2.1pt);
  \node[label,anchor=south] at (7.15,2.29) {$E_{4}$};
  \draw[side,draw=TwoAdic] (7.9,0.7) -- (7.9,2.2);
  \draw[tip,TwoAdic,shift={(7.9,1.45)},rotate=-90] (-2.5pt,2.1pt) -- (0,0) -- (-2.5pt,-2.1pt);
  \node[label,anchor=west] at (7.99,1.45) {$E_{2}$};
  \draw[side,draw=ThreeAdic] (6.4,0.7) -- (7.9,0.7);
  \draw[tip,ThreeAdic,shift={(7.15,0.7)},rotate=0] (-2.5pt,2.1pt) -- (0,0) -- (-2.5pt,-2.1pt);
  \node[label,anchor=north] at (7.15,0.61) {$E_{7}$};
  \node[cone] at (6.4,0.7) {};
  \node[ctwo] at (6.4,2.2) {};
  \node[cone] at (7.9,2.2) {};
  \node[ctwo] at (7.9,0.7) {};
  \node[inner sep=1.4pt,fill=SoftInk!5,rounded corners=1pt] at (7.15,1.45) {$F_{1}$};
  \fill[cell] (9.25,0.7) rectangle (10.75,2.2);
  \draw[side,draw=TwoAdic] (9.25,0.7) -- (9.25,2.2);
  \draw[tip,TwoAdic,shift={(9.25,1.45)},rotate=90] (-2.5pt,2.1pt) -- (0,0) -- (-2.5pt,-2.1pt);
  \node[label,anchor=east] at (9.16,1.45) {$E_{2}$};
  \draw[side,draw=ThreeAdic] (9.25,2.2) -- (10.75,2.2);
  \draw[tip,ThreeAdic,shift={(10,2.2)},rotate=180] (-2.5pt,2.1pt) -- (0,0) -- (-2.5pt,-2.1pt);
  \node[label,anchor=south] at (10,2.29) {$E_{5}$};
  \draw[side,draw=TwoAdic] (10.75,0.7) -- (10.75,2.2);
  \draw[tip,TwoAdic,shift={(10.75,1.45)},rotate=-90] (-2.5pt,2.1pt) -- (0,0) -- (-2.5pt,-2.1pt);
  \node[label,anchor=west] at (10.84,1.45) {$E_{3}$};
  \draw[side,draw=ThreeAdic] (9.25,0.7) -- (10.75,0.7);
  \draw[tip,ThreeAdic,shift={(10,0.7)},rotate=0] (-2.5pt,2.1pt) -- (0,0) -- (-2.5pt,-2.1pt);
  \node[label,anchor=north] at (10,0.61) {$E_{7}$};
  \node[cone] at (9.25,0.7) {};
  \node[ctwo] at (9.25,2.2) {};
  \node[cone] at (10.75,2.2) {};
  \node[ctwo] at (10.75,0.7) {};
  \node[inner sep=1.4pt,fill=SoftInk!5,rounded corners=1pt] at (10,1.45) {$F_{2}$};
  \fill[cell] (12.1,0.7) rectangle (13.6,2.2);
  \draw[side,draw=TwoAdic] (12.1,0.7) -- (12.1,2.2);
  \draw[tip,TwoAdic,shift={(12.1,1.45)},rotate=90] (-2.5pt,2.1pt) -- (0,0) -- (-2.5pt,-2.1pt);
  \node[label,anchor=east] at (12.01,1.45) {$E_{3}$};
  \draw[side,draw=ThreeAdic] (12.1,2.2) -- (13.6,2.2);
  \draw[tip,ThreeAdic,shift={(12.85,2.2)},rotate=180] (-2.5pt,2.1pt) -- (0,0) -- (-2.5pt,-2.1pt);
  \node[label,anchor=south] at (12.85,2.29) {$E_{6}$};
  \draw[side,draw=TwoAdic] (13.6,0.7) -- (13.6,2.2);
  \draw[tip,TwoAdic,shift={(13.6,1.45)},rotate=-90] (-2.5pt,2.1pt) -- (0,0) -- (-2.5pt,-2.1pt);
  \node[label,anchor=west] at (13.69,1.45) {$E_{2}$};
  \draw[side,draw=ThreeAdic] (12.1,0.7) -- (13.6,0.7);
  \draw[tip,ThreeAdic,shift={(12.85,0.7)},rotate=0] (-2.5pt,2.1pt) -- (0,0) -- (-2.5pt,-2.1pt);
  \node[label,anchor=north] at (12.85,0.61) {$E_{6}$};
  \node[cone] at (12.1,0.7) {};
  \node[ctwo] at (12.1,2.2) {};
  \node[cone] at (13.6,2.2) {};
  \node[ctwo] at (13.6,0.7) {};
  \node[inner sep=1.4pt,fill=SoftInk!5,rounded corners=1pt] at (12.85,1.45) {$F_{3}$};
  \fill[cell] (6.4,-2.2) rectangle (7.9,-0.7);
  \draw[side,draw=TwoAdic] (6.4,-2.2) -- (6.4,-0.7);
  \draw[tip,TwoAdic,shift={(6.4,-1.45)},rotate=90] (-2.5pt,2.1pt) -- (0,0) -- (-2.5pt,-2.1pt);
  \node[label,anchor=east] at (6.31,-1.45) {$E_{1}$};
  \draw[side,draw=ThreeAdic] (6.4,-0.7) -- (7.9,-0.7);
  \draw[tip,ThreeAdic,shift={(7.15,-0.7)},rotate=180] (-2.5pt,2.1pt) -- (0,0) -- (-2.5pt,-2.1pt);
  \node[label,anchor=south] at (7.15,-0.61) {$E_{5}$};
  \draw[side,draw=TwoAdic] (7.9,-2.2) -- (7.9,-0.7);
  \draw[tip,TwoAdic,shift={(7.9,-1.45)},rotate=-90] (-2.5pt,2.1pt) -- (0,0) -- (-2.5pt,-2.1pt);
  \node[label,anchor=west] at (7.99,-1.45) {$E_{1}$};
  \draw[side,draw=ThreeAdic] (6.4,-2.2) -- (7.9,-2.2);
  \draw[tip,ThreeAdic,shift={(7.15,-2.2)},rotate=0] (-2.5pt,2.1pt) -- (0,0) -- (-2.5pt,-2.1pt);
  \node[label,anchor=north] at (7.15,-2.29) {$E_{6}$};
  \node[cone] at (6.4,-2.2) {};
  \node[ctwo] at (6.4,-0.7) {};
  \node[cone] at (7.9,-0.7) {};
  \node[ctwo] at (7.9,-2.2) {};
  \node[inner sep=1.4pt,fill=SoftInk!5,rounded corners=1pt] at (7.15,-1.45) {$F_{4}$};
  \fill[cell] (9.25,-2.2) rectangle (10.75,-0.7);
  \draw[side,draw=TwoAdic] (9.25,-2.2) -- (9.25,-0.7);
  \draw[tip,TwoAdic,shift={(9.25,-1.45)},rotate=90] (-2.5pt,2.1pt) -- (0,0) -- (-2.5pt,-2.1pt);
  \node[label,anchor=east] at (9.16,-1.45) {$E_{3}$};
  \draw[side,draw=ThreeAdic] (9.25,-0.7) -- (10.75,-0.7);
  \draw[tip,ThreeAdic,shift={(10,-0.7)},rotate=180] (-2.5pt,2.1pt) -- (0,0) -- (-2.5pt,-2.1pt);
  \node[label,anchor=south] at (10,-0.61) {$E_{4}$};
  \draw[side,draw=TwoAdic] (10.75,-2.2) -- (10.75,-0.7);
  \draw[tip,TwoAdic,shift={(10.75,-1.45)},rotate=-90] (-2.5pt,2.1pt) -- (0,0) -- (-2.5pt,-2.1pt);
  \node[label,anchor=west] at (10.84,-1.45) {$E_{1}$};
  \draw[side,draw=ThreeAdic] (9.25,-2.2) -- (10.75,-2.2);
  \draw[tip,ThreeAdic,shift={(10,-2.2)},rotate=0] (-2.5pt,2.1pt) -- (0,0) -- (-2.5pt,-2.1pt);
  \node[label,anchor=north] at (10,-2.29) {$E_{7}$};
  \node[cone] at (9.25,-2.2) {};
  \node[ctwo] at (9.25,-0.7) {};
  \node[cone] at (10.75,-0.7) {};
  \node[ctwo] at (10.75,-2.2) {};
  \node[inner sep=1.4pt,fill=SoftInk!5,rounded corners=1pt] at (10,-1.45) {$F_{5}$};
  \fill[cell] (12.1,-2.2) rectangle (13.6,-0.7);
  \draw[side,draw=TwoAdic] (12.1,-2.2) -- (12.1,-0.7);
  \draw[tip,TwoAdic,shift={(12.1,-1.45)},rotate=90] (-2.5pt,2.1pt) -- (0,0) -- (-2.5pt,-2.1pt);
  \node[label,anchor=east] at (12.01,-1.45) {$E_{2}$};
  \draw[side,draw=ThreeAdic] (12.1,-0.7) -- (13.6,-0.7);
  \draw[tip,ThreeAdic,shift={(12.85,-0.7)},rotate=180] (-2.5pt,2.1pt) -- (0,0) -- (-2.5pt,-2.1pt);
  \node[label,anchor=south] at (12.85,-0.61) {$E_{4}$};
  \draw[side,draw=TwoAdic] (13.6,-2.2) -- (13.6,-0.7);
  \draw[tip,TwoAdic,shift={(13.6,-1.45)},rotate=-90] (-2.5pt,2.1pt) -- (0,0) -- (-2.5pt,-2.1pt);
  \node[label,anchor=west] at (13.69,-1.45) {$E_{3}$};
  \draw[side,draw=ThreeAdic] (12.1,-2.2) -- (13.6,-2.2);
  \draw[tip,ThreeAdic,shift={(12.85,-2.2)},rotate=0] (-2.5pt,2.1pt) -- (0,0) -- (-2.5pt,-2.1pt);
  \node[label,anchor=north] at (12.85,-2.29) {$E_{5}$};
  \node[cone] at (12.1,-2.2) {};
  \node[ctwo] at (12.1,-0.7) {};
  \node[cone] at (13.6,-0.7) {};
  \node[ctwo] at (13.6,-2.2) {};
  \node[inner sep=1.4pt,fill=SoftInk!5,rounded corners=1pt] at (12.85,-1.45) {$F_{6}$};
  \draw[edge,draw=TwoAdic] (3.9,-3.3) -- (4.45,-3.3);
  \node[label,anchor=west] at (4.55,-3.3) {$2$-adic edge};
  \draw[edge,draw=ThreeAdic] (7.2,-3.3) -- (7.75,-3.3);
  \node[label,anchor=west] at (7.85,-3.3) {$3$-adic edge};
\end{tikzpicture}
    \caption{The square complex $X_\Gamma$ of Section \ref{sec:qsqrt3}.}
    \label{fig:qsqrt3-complex}
\end{figure}

\begin{figure}[htbp]
    \centering
    \input{figures/qsqrt13_euler_one_covering}
    \caption{The covering $X_\Lambda \to X_\Gamma$ of Section \ref{sec:qsqrt13}.}
    \label{fig:qsqrt13-covering}
\end{figure}

\begin{figure}[htbp]
    \centering
    \definecolor{SoftInk}{HTML}{2F3437}
\definecolor{EdgeColour}{HTML}{4F7D72}
\definecolor{PrimeOne}{HTML}{4F7D72}
\definecolor{PrimeTwo}{HTML}{A85D67}
\begin{tikzpicture}[x=1cm,y=1cm,
  every node/.style={font=\scriptsize,text=SoftInk},
  edge/.style={line width=1.05pt,line cap=round},
  side/.style={line width=0.95pt,line cap=round},
  tip/.style={line width=0.9pt,line cap=round,line join=round},
  cell/.style={fill=SoftInk!5},
  mirror/.style={draw=SoftInk!75,line width=0.75pt,dash pattern=on 2.2pt off 1.6pt},
  label/.style={inner sep=1.2pt},
  vert/.style={circle,fill=SoftInk,inner sep=1.8pt},
  corner/.style={circle,fill=SoftInk,inner sep=1.2pt},
  orbi/.style={circle,fill=white,draw=SoftInk,line width=0.6pt,inner sep=1.35pt}]
  \draw[edge,draw=EdgeColour] (0,0) .. controls (-0.598,2.231) and (-2.231,0.598) .. (0,0);
  \draw[tip,draw=EdgeColour,shift={(-1.061,1.061)},rotate=-135] (-2.5pt,2.1pt) -- (0,0) -- (-2.5pt,-2.1pt);
  \node[label,anchor=315] at (-1.131,1.131) {$E_{1}$};
  \draw[edge,draw=EdgeColour] (0,0) -- (-1.358,-0.563);
  \node[orbi] at (-1.358,-0.563) {};
  \node[label,anchor=382.5] at (-1.469,-0.609) {$E_{2}$};
  \draw[edge,draw=EdgeColour] (0,0) -- (-0.563,-1.358);
  \node[orbi] at (-0.563,-1.358) {};
  \node[label,anchor=427.5] at (-0.609,-1.469) {$E_{3}$};
  \draw[edge,draw=EdgeColour] (0,0) -- (0.563,-1.358);
  \node[orbi] at (0.563,-1.358) {};
  \node[label,anchor=472.5] at (0.609,-1.469) {$E_{4}$};
  \draw[edge,draw=EdgeColour] (0,0) -- (1.358,-0.563);
  \node[orbi] at (1.358,-0.563) {};
  \node[label,anchor=517.5] at (1.469,-0.609) {$E_{5}$};
  \draw[edge,draw=EdgeColour] (0,0) -- (1.358,0.563);
  \node[orbi] at (1.358,0.563) {};
  \node[label,anchor=562.5] at (1.469,0.609) {$E_{6}$};
  \draw[edge,draw=EdgeColour] (0,0) -- (0.563,1.358);
  \node[orbi] at (0.563,1.358) {};
  \node[label,anchor=607.5] at (0.609,1.469) {$E_{7}$};
  \node[vert] at (0,0) {};
  \fill[cell] (2.688,0.55) rectangle (3.988,1.85);
  \draw[side,draw=EdgeColour] (2.688,0.55) -- (2.688,1.85);
  \node[orbi] at (2.688,1.2) {};
  \node[label,anchor=east] at (2.598,1.2) {$E_{6}$};
  \draw[side,draw=EdgeColour] (2.688,1.85) -- (3.988,1.85);
  \node[orbi] at (3.338,1.85) {};
  \node[label,anchor=south] at (3.338,1.94) {$E_{3}$};
  \draw[side,draw=EdgeColour] (3.988,1.85) -- (3.988,0.55);
  \node[orbi] at (3.988,1.2) {};
  \node[label,anchor=west] at (4.078,1.2) {$E_{2}$};
  \draw[side,draw=EdgeColour] (3.988,0.55) -- (2.688,0.55);
  \node[orbi] at (3.338,0.55) {};
  \node[label,anchor=north] at (3.338,0.46) {$E_{7}$};
  \node[corner] at (2.688,0.55) {};
  \node[corner] at (2.688,1.85) {};
  \node[corner] at (3.988,1.85) {};
  \node[corner] at (3.988,0.55) {};
  \node[inner sep=1.4pt,fill=SoftInk!5,rounded corners=1pt] at (3.338,1.2) {$F_{1}$};
  \fill[cell] (5.138,0.55) rectangle (6.438,1.85);
  \draw[mirror] (5.138,1.2) -- (6.438,1.2);
  \draw[mirror] (5.788,0.55) -- (5.788,1.85);
  \draw[side,draw=EdgeColour] (5.138,0.55) -- (5.138,1.85);
  \node[orbi] at (5.138,1.2) {};
  \node[label,anchor=east] at (5.048,1.2) {$E_{4}$};
  \draw[side,draw=EdgeColour] (5.138,1.85) -- (6.438,1.85);
  \node[orbi] at (5.788,1.85) {};
  \node[label,anchor=south] at (5.788,1.94) {$E_{3}$};
  \draw[side,draw=EdgeColour] (6.438,1.85) -- (6.438,0.55);
  \node[orbi] at (6.438,1.2) {};
  \node[label,anchor=west] at (6.528,1.2) {$E_{4}$};
  \draw[side,draw=EdgeColour] (6.438,0.55) -- (5.138,0.55);
  \node[orbi] at (5.788,0.55) {};
  \node[label,anchor=north] at (5.788,0.46) {$E_{3}$};
  \node[corner] at (5.138,0.55) {};
  \node[corner] at (5.138,1.85) {};
  \node[corner] at (6.438,1.85) {};
  \node[corner] at (6.438,0.55) {};
  \node[inner sep=1.4pt,fill=SoftInk!5,rounded corners=1pt] at (5.788,1.2) {$F_{2}$};
  \fill[cell] (7.588,0.55) rectangle (8.888,1.85);
  \draw[mirror] (8.238,0.55) -- (8.238,1.85);
  \draw[side,draw=EdgeColour] (7.588,0.55) -- (7.588,1.85);
  \draw[tip,draw=EdgeColour,shift={(7.588,1.2)},rotate=90] (-2.5pt,2.1pt) -- (0,0) -- (-2.5pt,-2.1pt);
  \node[label,anchor=east] at (7.498,1.2) {$E_{1}$};
  \draw[side,draw=EdgeColour] (7.588,1.85) -- (8.888,1.85);
  \node[orbi] at (8.238,1.85) {};
  \node[label,anchor=south] at (8.238,1.94) {$E_{3}$};
  \draw[side,draw=EdgeColour] (8.888,1.85) -- (8.888,0.55);
  \draw[tip,draw=EdgeColour,shift={(8.888,1.2)},rotate=90] (-2.5pt,2.1pt) -- (0,0) -- (-2.5pt,-2.1pt);
  \node[label,anchor=west] at (8.978,1.2) {$E_{1}$};
  \draw[side,draw=EdgeColour] (8.888,0.55) -- (7.588,0.55);
  \node[orbi] at (8.238,0.55) {};
  \node[label,anchor=north] at (8.238,0.46) {$E_{2}$};
  \node[corner] at (7.588,0.55) {};
  \node[corner] at (7.588,1.85) {};
  \node[corner] at (8.888,1.85) {};
  \node[corner] at (8.888,0.55) {};
  \node[inner sep=1.4pt,fill=SoftInk!5,rounded corners=1pt] at (8.238,1.2) {$F_{3}$};
  \fill[cell] (10.038,0.55) rectangle (11.338,1.85);
  \draw[mirror] (10.038,1.2) -- (11.338,1.2);
  \draw[mirror] (10.688,0.55) -- (10.688,1.85);
  \draw[side,draw=EdgeColour] (10.038,0.55) -- (10.038,1.85);
  \node[orbi] at (10.038,1.2) {};
  \node[label,anchor=east] at (9.948,1.2) {$E_{5}$};
  \draw[side,draw=EdgeColour] (10.038,1.85) -- (11.338,1.85);
  \node[orbi] at (10.688,1.85) {};
  \node[label,anchor=south] at (10.688,1.94) {$E_{2}$};
  \draw[side,draw=EdgeColour] (11.338,1.85) -- (11.338,0.55);
  \node[orbi] at (11.338,1.2) {};
  \node[label,anchor=west] at (11.428,1.2) {$E_{5}$};
  \draw[side,draw=EdgeColour] (11.338,0.55) -- (10.038,0.55);
  \node[orbi] at (10.688,0.55) {};
  \node[label,anchor=north] at (10.688,0.46) {$E_{2}$};
  \node[corner] at (10.038,0.55) {};
  \node[corner] at (10.038,1.85) {};
  \node[corner] at (11.338,1.85) {};
  \node[corner] at (11.338,0.55) {};
  \node[inner sep=1.4pt,fill=SoftInk!5,rounded corners=1pt] at (10.688,1.2) {$F_{4}$};
  \begin{scope}[shift={(1.225,0)}]
  \fill[cell] (2.688,-1.95) rectangle (3.988,-0.65);
  \draw[mirror] (2.688,-1.3) -- (3.988,-1.3);
  \draw[side,draw=EdgeColour] (2.688,-1.95) -- (2.688,-0.65);
  \node[orbi] at (2.688,-1.3) {};
  \node[label,anchor=east] at (2.598,-1.3) {$E_{6}$};
  \draw[side,draw=EdgeColour] (2.688,-0.65) -- (3.988,-0.65);
  \node[orbi] at (3.338,-0.65) {};
  \node[label,anchor=south] at (3.338,-0.56) {$E_{5}$};
  \draw[side,draw=EdgeColour] (3.988,-0.65) -- (3.988,-1.95);
  \node[orbi] at (3.988,-1.3) {};
  \node[label,anchor=west] at (4.078,-1.3) {$E_{4}$};
  \draw[side,draw=EdgeColour] (3.988,-1.95) -- (2.688,-1.95);
  \node[orbi] at (3.338,-1.95) {};
  \node[label,anchor=north] at (3.338,-2.04) {$E_{5}$};
  \node[corner] at (2.688,-1.95) {};
  \node[corner] at (2.688,-0.65) {};
  \node[corner] at (3.988,-0.65) {};
  \node[corner] at (3.988,-1.95) {};
  \node[inner sep=1.4pt,fill=SoftInk!5,rounded corners=1pt] at (3.338,-1.3) {$F_{5}$};
  \fill[cell] (5.138,-1.95) rectangle (6.438,-0.65);
  \draw[side,draw=EdgeColour] (5.138,-1.95) -- (5.138,-0.65);
  \node[orbi] at (5.138,-1.3) {};
  \node[label,anchor=east] at (5.048,-1.3) {$E_{7}$};
  \draw[side,draw=EdgeColour] (5.138,-0.65) -- (6.438,-0.65);
  \draw[tip,draw=EdgeColour,shift={(5.788,-0.65)},rotate=180] (-2.5pt,2.1pt) -- (0,0) -- (-2.5pt,-2.1pt);
  \node[label,anchor=south] at (5.788,-0.56) {$E_{1}$};
  \draw[side,draw=EdgeColour] (6.438,-0.65) -- (6.438,-1.95);
  \draw[tip,draw=EdgeColour,shift={(6.438,-1.3)},rotate=90] (-2.5pt,2.1pt) -- (0,0) -- (-2.5pt,-2.1pt);
  \node[label,anchor=west] at (6.528,-1.3) {$E_{1}$};
  \draw[side,draw=EdgeColour] (6.438,-1.95) -- (5.138,-1.95);
  \node[orbi] at (5.788,-1.95) {};
  \node[label,anchor=north] at (5.788,-2.04) {$E_{4}$};
  \node[corner] at (5.138,-1.95) {};
  \node[corner] at (5.138,-0.65) {};
  \node[corner] at (6.438,-0.65) {};
  \node[corner] at (6.438,-1.95) {};
  \node[inner sep=1.4pt,fill=SoftInk!5,rounded corners=1pt] at (5.788,-1.3) {$F_{6}$};
  \fill[cell] (7.588,-1.95) rectangle (8.888,-0.65);
  \draw[mirror] (7.588,-1.3) -- (8.888,-1.3);
  \draw[side,draw=EdgeColour] (7.588,-1.95) -- (7.588,-0.65);
  \node[orbi] at (7.588,-1.3) {};
  \node[label,anchor=east] at (7.498,-1.3) {$E_{6}$};
  \draw[side,draw=EdgeColour] (7.588,-0.65) -- (8.888,-0.65);
  \draw[tip,draw=EdgeColour,shift={(8.238,-0.65)},rotate=0] (-2.5pt,2.1pt) -- (0,0) -- (-2.5pt,-2.1pt);
  \node[label,anchor=south] at (8.238,-0.56) {$E_{1}$};
  \draw[side,draw=EdgeColour] (8.888,-0.65) -- (8.888,-1.95);
  \node[orbi] at (8.888,-1.3) {};
  \node[label,anchor=west] at (8.978,-1.3) {$E_{5}$};
  \draw[side,draw=EdgeColour] (8.888,-1.95) -- (7.588,-1.95);
  \draw[tip,draw=EdgeColour,shift={(8.238,-1.95)},rotate=0] (-2.5pt,2.1pt) -- (0,0) -- (-2.5pt,-2.1pt);
  \node[label,anchor=north] at (8.238,-2.04) {$E_{1}$};
  \node[corner] at (7.588,-1.95) {};
  \node[corner] at (7.588,-0.65) {};
  \node[corner] at (8.888,-0.65) {};
  \node[corner] at (8.888,-1.95) {};
  \node[inner sep=1.4pt,fill=SoftInk!5,rounded corners=1pt] at (8.238,-1.3) {$F_{7}$};
  \end{scope}
\end{tikzpicture}
    \caption{The complex $X_{\Gamma_{1}}$ of Section \ref{sec:qsqrt7}.}
    \label{fig:qsqrt7-gamma1}
\end{figure}

\begin{figure}[htbp]
    \centering
    \definecolor{SoftInk}{HTML}{2F3437}
\definecolor{EdgeColour}{HTML}{4F7D72}
\definecolor{PrimeOne}{HTML}{4F7D72}
\definecolor{PrimeTwo}{HTML}{A85D67}
\begin{tikzpicture}[x=1cm,y=1cm,
  every node/.style={font=\scriptsize,text=SoftInk},
  edge/.style={line width=1.05pt,line cap=round},
  side/.style={line width=0.95pt,line cap=round},
  tip/.style={line width=0.9pt,line cap=round,line join=round},
  cell/.style={fill=SoftInk!5},
  mirror/.style={draw=SoftInk!75,line width=0.75pt,dash pattern=on 2.2pt off 1.6pt},
  label/.style={inner sep=1.2pt},
  vert/.style={circle,fill=SoftInk,inner sep=1.8pt},
  corner/.style={circle,fill=SoftInk,inner sep=1.2pt},
  orbi/.style={circle,fill=white,draw=SoftInk,line width=0.6pt,inner sep=1.35pt}]
  \draw[edge,draw=PrimeOne] (0,0) -- (-0.563,1.358);
  \node[orbi] at (-0.563,1.358) {};
  \node[label,anchor=292.5] at (-0.609,1.469) {$E_{1}$};
  \draw[edge,draw=PrimeOne] (0,0) -- (-1.358,0.563);
  \node[orbi] at (-1.358,0.563) {};
  \node[label,anchor=337.5] at (-1.469,0.609) {$E_{2}$};
  \draw[edge,draw=PrimeOne] (0,0) -- (-1.358,-0.563);
  \node[orbi] at (-1.358,-0.563) {};
  \node[label,anchor=382.5] at (-1.469,-0.609) {$E_{3}$};
  \draw[edge,draw=PrimeOne] (0,0) -- (-0.563,-1.358);
  \node[orbi] at (-0.563,-1.358) {};
  \node[label,anchor=427.5] at (-0.609,-1.469) {$E_{4}$};
  \draw[edge,draw=PrimeTwo] (0,0) .. controls (0.598,-2.231) and (2.231,-0.598) .. (0,0);
  \draw[tip,draw=PrimeTwo,shift={(1.061,-1.061)},rotate=45] (-2.5pt,2.1pt) -- (0,0) -- (-2.5pt,-2.1pt);
  \node[label,anchor=495] at (1.131,-1.131) {$E_{5}$};
  \draw[edge,draw=PrimeTwo] (0,0) -- (1.358,0.563);
  \node[orbi] at (1.358,0.563) {};
  \node[label,anchor=562.5] at (1.469,0.609) {$E_{6}$};
  \draw[edge,draw=PrimeTwo] (0,0) -- (0.563,1.358);
  \node[orbi] at (0.563,1.358) {};
  \node[label,anchor=607.5] at (0.609,1.469) {$E_{7}$};
  \node[vert] at (0,0) {};
  \fill[cell] (2.688,0.55) rectangle (3.988,1.85);
  \draw[side,draw=PrimeOne] (2.688,0.55) -- (2.688,1.85);
  \node[orbi] at (2.688,1.2) {};
  \node[label,anchor=east] at (2.598,1.2) {$E_{4}$};
  \draw[side,draw=PrimeTwo] (2.688,1.85) -- (3.988,1.85);
  \node[orbi] at (3.338,1.85) {};
  \node[label,anchor=south] at (3.338,1.94) {$E_{6}$};
  \draw[side,draw=PrimeOne] (3.988,1.85) -- (3.988,0.55);
  \node[orbi] at (3.988,1.2) {};
  \node[label,anchor=west] at (4.078,1.2) {$E_{1}$};
  \draw[side,draw=PrimeTwo] (3.988,0.55) -- (2.688,0.55);
  \node[orbi] at (3.338,0.55) {};
  \node[label,anchor=north] at (3.338,0.46) {$E_{7}$};
  \node[corner] at (2.688,0.55) {};
  \node[corner] at (2.688,1.85) {};
  \node[corner] at (3.988,1.85) {};
  \node[corner] at (3.988,0.55) {};
  \node[inner sep=1.4pt,fill=SoftInk!5,rounded corners=1pt] at (3.338,1.2) {$F_{1}$};
  \fill[cell] (5.138,0.55) rectangle (6.438,1.85);
  \draw[mirror] (5.138,1.2) -- (6.438,1.2);
  \draw[side,draw=PrimeOne] (5.138,0.55) -- (5.138,1.85);
  \node[orbi] at (5.138,1.2) {};
  \node[label,anchor=east] at (5.048,1.2) {$E_{2}$};
  \draw[side,draw=PrimeTwo] (5.138,1.85) -- (6.438,1.85);
  \node[orbi] at (5.788,1.85) {};
  \node[label,anchor=south] at (5.788,1.94) {$E_{6}$};
  \draw[side,draw=PrimeOne] (6.438,1.85) -- (6.438,0.55);
  \node[orbi] at (6.438,1.2) {};
  \node[label,anchor=west] at (6.528,1.2) {$E_{3}$};
  \draw[side,draw=PrimeTwo] (6.438,0.55) -- (5.138,0.55);
  \node[orbi] at (5.788,0.55) {};
  \node[label,anchor=north] at (5.788,0.46) {$E_{6}$};
  \node[corner] at (5.138,0.55) {};
  \node[corner] at (5.138,1.85) {};
  \node[corner] at (6.438,1.85) {};
  \node[corner] at (6.438,0.55) {};
  \node[inner sep=1.4pt,fill=SoftInk!5,rounded corners=1pt] at (5.788,1.2) {$F_{2}$};
  \fill[cell] (7.588,0.55) rectangle (8.888,1.85);
  \draw[side,draw=PrimeOne] (7.588,0.55) -- (7.588,1.85);
  \node[orbi] at (7.588,1.2) {};
  \node[label,anchor=east] at (7.498,1.2) {$E_{1}$};
  \draw[side,draw=PrimeTwo] (7.588,1.85) -- (8.888,1.85);
  \draw[tip,draw=PrimeTwo,shift={(8.238,1.85)},rotate=0] (-2.5pt,2.1pt) -- (0,0) -- (-2.5pt,-2.1pt);
  \node[label,anchor=south] at (8.238,1.94) {$E_{5}$};
  \draw[side,draw=PrimeOne] (8.888,1.85) -- (8.888,0.55);
  \node[orbi] at (8.888,1.2) {};
  \node[label,anchor=west] at (8.978,1.2) {$E_{1}$};
  \draw[side,draw=PrimeTwo] (8.888,0.55) -- (7.588,0.55);
  \draw[tip,draw=PrimeTwo,shift={(8.238,0.55)},rotate=180] (-2.5pt,2.1pt) -- (0,0) -- (-2.5pt,-2.1pt);
  \node[label,anchor=north] at (8.238,0.46) {$E_{5}$};
  \node[corner] at (7.588,0.55) {};
  \node[corner] at (7.588,1.85) {};
  \node[corner] at (8.888,1.85) {};
  \node[corner] at (8.888,0.55) {};
  \node[orbi] at (8.238,1.2) {};
  \node[inner sep=1.4pt,fill=SoftInk!5,rounded corners=1pt] at (8.238,0.9) {$F_{3}$};
  \fill[cell] (2.688,-1.95) rectangle (3.988,-0.65);
  \draw[mirror] (2.688,-1.3) -- (3.988,-1.3);
  \draw[side,draw=PrimeOne] (2.688,-1.95) -- (2.688,-0.65);
  \node[orbi] at (2.688,-1.3) {};
  \node[label,anchor=east] at (2.598,-1.3) {$E_{4}$};
  \draw[side,draw=PrimeTwo] (2.688,-0.65) -- (3.988,-0.65);
  \draw[tip,draw=PrimeTwo,shift={(3.338,-0.65)},rotate=180] (-2.5pt,2.1pt) -- (0,0) -- (-2.5pt,-2.1pt);
  \node[label,anchor=south] at (3.338,-0.56) {$E_{5}$};
  \draw[side,draw=PrimeOne] (3.988,-0.65) -- (3.988,-1.95);
  \node[orbi] at (3.988,-1.3) {};
  \node[label,anchor=west] at (4.078,-1.3) {$E_{3}$};
  \draw[side,draw=PrimeTwo] (3.988,-1.95) -- (2.688,-1.95);
  \draw[tip,draw=PrimeTwo,shift={(3.338,-1.95)},rotate=180] (-2.5pt,2.1pt) -- (0,0) -- (-2.5pt,-2.1pt);
  \node[label,anchor=north] at (3.338,-2.04) {$E_{5}$};
  \node[corner] at (2.688,-1.95) {};
  \node[corner] at (2.688,-0.65) {};
  \node[corner] at (3.988,-0.65) {};
  \node[corner] at (3.988,-1.95) {};
  \node[inner sep=1.4pt,fill=SoftInk!5,rounded corners=1pt] at (3.338,-1.3) {$F_{4}$};
  \fill[cell] (5.138,-1.95) rectangle (6.438,-0.65);
  \draw[side,draw=PrimeOne] (5.138,-1.95) -- (5.138,-0.65);
  \node[orbi] at (5.138,-1.3) {};
  \node[label,anchor=east] at (5.048,-1.3) {$E_{3}$};
  \draw[side,draw=PrimeTwo] (5.138,-0.65) -- (6.438,-0.65);
  \draw[tip,draw=PrimeTwo,shift={(5.788,-0.65)},rotate=180] (-2.5pt,2.1pt) -- (0,0) -- (-2.5pt,-2.1pt);
  \node[label,anchor=south] at (5.788,-0.56) {$E_{5}$};
  \draw[side,draw=PrimeOne] (6.438,-0.65) -- (6.438,-1.95);
  \node[orbi] at (6.438,-1.3) {};
  \node[label,anchor=west] at (6.528,-1.3) {$E_{2}$};
  \draw[side,draw=PrimeTwo] (6.438,-1.95) -- (5.138,-1.95);
  \node[orbi] at (5.788,-1.95) {};
  \node[label,anchor=north] at (5.788,-2.04) {$E_{7}$};
  \node[corner] at (5.138,-1.95) {};
  \node[corner] at (5.138,-0.65) {};
  \node[corner] at (6.438,-0.65) {};
  \node[corner] at (6.438,-1.95) {};
  \node[inner sep=1.4pt,fill=SoftInk!5,rounded corners=1pt] at (5.788,-1.3) {$F_{5}$};
  \fill[cell] (7.588,-1.95) rectangle (8.888,-0.65);
  \draw[mirror] (7.588,-1.3) -- (8.888,-1.3);
  \draw[side,draw=PrimeOne] (7.588,-1.95) -- (7.588,-0.65);
  \node[orbi] at (7.588,-1.3) {};
  \node[label,anchor=east] at (7.498,-1.3) {$E_{4}$};
  \draw[side,draw=PrimeTwo] (7.588,-0.65) -- (8.888,-0.65);
  \draw[tip,draw=PrimeTwo,shift={(8.238,-0.65)},rotate=0] (-2.5pt,2.1pt) -- (0,0) -- (-2.5pt,-2.1pt);
  \node[label,anchor=south] at (8.238,-0.56) {$E_{5}$};
  \draw[side,draw=PrimeOne] (8.888,-0.65) -- (8.888,-1.95);
  \node[orbi] at (8.888,-1.3) {};
  \node[label,anchor=west] at (8.978,-1.3) {$E_{2}$};
  \draw[side,draw=PrimeTwo] (8.888,-1.95) -- (7.588,-1.95);
  \draw[tip,draw=PrimeTwo,shift={(8.238,-1.95)},rotate=0] (-2.5pt,2.1pt) -- (0,0) -- (-2.5pt,-2.1pt);
  \node[label,anchor=north] at (8.238,-2.04) {$E_{5}$};
  \node[corner] at (7.588,-1.95) {};
  \node[corner] at (7.588,-0.65) {};
  \node[corner] at (8.888,-0.65) {};
  \node[corner] at (8.888,-1.95) {};
  \node[inner sep=1.4pt,fill=SoftInk!5,rounded corners=1pt] at (8.238,-1.3) {$F_{6}$};
  \draw[edge,draw=PrimeOne] (2.3,-2.9) -- (2.85,-2.9);
  \node[label,anchor=west] at (2.95,-2.9) {$\fp_3$-adic edge};
  \draw[edge,draw=PrimeTwo] (5.3,-2.9) -- (5.85,-2.9);
  \node[label,anchor=west] at (5.95,-2.9) {$\fq_3$-adic edge};
\end{tikzpicture}
    \caption{The complex $X_{\Gamma_{2}}$ of Section \ref{sec:qsqrt7}.}
    \label{fig:qsqrt7-gamma2}
\end{figure}

\begin{figure}[htbp]
    \centering
    \definecolor{SoftInk}{HTML}{2F3437}
\definecolor{EdgeColour}{HTML}{4F7D72}
\definecolor{PrimeOne}{HTML}{4F7D72}
\definecolor{PrimeTwo}{HTML}{A85D67}
\begin{tikzpicture}[x=1cm,y=1cm,
  every node/.style={font=\scriptsize,text=SoftInk},
  edge/.style={line width=1.05pt,line cap=round},
  side/.style={line width=0.95pt,line cap=round},
  tip/.style={line width=0.9pt,line cap=round,line join=round},
  cell/.style={fill=SoftInk!5},
  mirror/.style={draw=SoftInk!75,line width=0.75pt,dash pattern=on 2.2pt off 1.6pt},
  label/.style={inner sep=1.2pt},
  vert/.style={circle,fill=SoftInk,inner sep=1.8pt},
  corner/.style={circle,fill=SoftInk,inner sep=1.2pt},
  orbi/.style={circle,fill=white,draw=SoftInk,line width=0.6pt,inner sep=1.35pt}]
  \draw[edge,draw=PrimeOne] (0,0) -- (-0.563,1.358);
  \node[orbi] at (-0.563,1.358) {};
  \node[label,anchor=292.5] at (-0.609,1.469) {$E_{1}$};
  \draw[edge,draw=PrimeOne] (0,0) -- (-1.358,0.563);
  \node[orbi] at (-1.358,0.563) {};
  \node[label,anchor=337.5] at (-1.469,0.609) {$E_{2}$};
  \draw[edge,draw=PrimeOne] (0,0) -- (-1.358,-0.563);
  \node[orbi] at (-1.358,-0.563) {};
  \node[label,anchor=382.5] at (-1.469,-0.609) {$E_{3}$};
  \draw[edge,draw=PrimeOne] (0,0) -- (-0.563,-1.358);
  \node[orbi] at (-0.563,-1.358) {};
  \node[label,anchor=427.5] at (-0.609,-1.469) {$E_{4}$};
  \draw[edge,draw=PrimeTwo] (0,0) .. controls (0.598,-2.231) and (2.231,-0.598) .. (0,0);
  \draw[tip,draw=PrimeTwo,shift={(1.061,-1.061)},rotate=45] (-2.5pt,2.1pt) -- (0,0) -- (-2.5pt,-2.1pt);
  \node[label,anchor=495] at (1.131,-1.131) {$E_{5}$};
  \draw[edge,draw=PrimeTwo] (0,0) -- (1.358,0.563);
  \node[orbi] at (1.358,0.563) {};
  \node[label,anchor=562.5] at (1.469,0.609) {$E_{6}$};
  \draw[edge,draw=PrimeTwo] (0,0) -- (0.563,1.358);
  \node[orbi] at (0.563,1.358) {};
  \node[label,anchor=607.5] at (0.609,1.469) {$E_{7}$};
  \node[vert] at (0,0) {};
  \fill[cell] (2.688,0.55) rectangle (3.988,1.85);
  \draw[side,draw=PrimeOne] (2.688,0.55) -- (2.688,1.85);
  \node[orbi] at (2.688,1.2) {};
  \node[label,anchor=east] at (2.598,1.2) {$E_{4}$};
  \draw[side,draw=PrimeTwo] (2.688,1.85) -- (3.988,1.85);
  \node[orbi] at (3.338,1.85) {};
  \node[label,anchor=south] at (3.338,1.94) {$E_{6}$};
  \draw[side,draw=PrimeOne] (3.988,1.85) -- (3.988,0.55);
  \node[orbi] at (3.988,1.2) {};
  \node[label,anchor=west] at (4.078,1.2) {$E_{1}$};
  \draw[side,draw=PrimeTwo] (3.988,0.55) -- (2.688,0.55);
  \node[orbi] at (3.338,0.55) {};
  \node[label,anchor=north] at (3.338,0.46) {$E_{7}$};
  \node[corner] at (2.688,0.55) {};
  \node[corner] at (2.688,1.85) {};
  \node[corner] at (3.988,1.85) {};
  \node[corner] at (3.988,0.55) {};
  \node[inner sep=1.4pt,fill=SoftInk!5,rounded corners=1pt] at (3.338,1.2) {$F_{1}$};
  \fill[cell] (5.138,0.55) rectangle (6.438,1.85);
  \draw[mirror] (5.138,1.2) -- (6.438,1.2);
  \draw[mirror] (5.788,0.55) -- (5.788,1.85);
  \draw[side,draw=PrimeOne] (5.138,0.55) -- (5.138,1.85);
  \node[orbi] at (5.138,1.2) {};
  \node[label,anchor=east] at (5.048,1.2) {$E_{2}$};
  \draw[side,draw=PrimeTwo] (5.138,1.85) -- (6.438,1.85);
  \node[orbi] at (5.788,1.85) {};
  \node[label,anchor=south] at (5.788,1.94) {$E_{6}$};
  \draw[side,draw=PrimeOne] (6.438,1.85) -- (6.438,0.55);
  \node[orbi] at (6.438,1.2) {};
  \node[label,anchor=west] at (6.528,1.2) {$E_{2}$};
  \draw[side,draw=PrimeTwo] (6.438,0.55) -- (5.138,0.55);
  \node[orbi] at (5.788,0.55) {};
  \node[label,anchor=north] at (5.788,0.46) {$E_{6}$};
  \node[corner] at (5.138,0.55) {};
  \node[corner] at (5.138,1.85) {};
  \node[corner] at (6.438,1.85) {};
  \node[corner] at (6.438,0.55) {};
  \node[inner sep=1.4pt,fill=SoftInk!5,rounded corners=1pt] at (5.788,1.2) {$F_{2}$};
  \fill[cell] (7.588,0.55) rectangle (8.888,1.85);
  \draw[mirror] (7.588,1.2) -- (8.888,1.2);
  \draw[side,draw=PrimeOne] (7.588,0.55) -- (7.588,1.85);
  \node[orbi] at (7.588,1.2) {};
  \node[label,anchor=east] at (7.498,1.2) {$E_{1}$};
  \draw[side,draw=PrimeTwo] (7.588,1.85) -- (8.888,1.85);
  \draw[tip,draw=PrimeTwo,shift={(8.238,1.85)},rotate=180] (-2.5pt,2.1pt) -- (0,0) -- (-2.5pt,-2.1pt);
  \node[label,anchor=south] at (8.238,1.94) {$E_{5}$};
  \draw[side,draw=PrimeOne] (8.888,1.85) -- (8.888,0.55);
  \node[orbi] at (8.888,1.2) {};
  \node[label,anchor=west] at (8.978,1.2) {$E_{1}$};
  \draw[side,draw=PrimeTwo] (8.888,0.55) -- (7.588,0.55);
  \draw[tip,draw=PrimeTwo,shift={(8.238,0.55)},rotate=180] (-2.5pt,2.1pt) -- (0,0) -- (-2.5pt,-2.1pt);
  \node[label,anchor=north] at (8.238,0.46) {$E_{5}$};
  \node[corner] at (7.588,0.55) {};
  \node[corner] at (7.588,1.85) {};
  \node[corner] at (8.888,1.85) {};
  \node[corner] at (8.888,0.55) {};
  \node[inner sep=1.4pt,fill=SoftInk!5,rounded corners=1pt] at (8.238,1.2) {$F_{3}$};
  \fill[cell] (10.038,0.55) rectangle (11.338,1.85);
  \draw[mirror] (10.038,1.2) -- (11.338,1.2);
  \draw[mirror] (10.688,0.55) -- (10.688,1.85);
  \draw[side,draw=PrimeOne] (10.038,0.55) -- (10.038,1.85);
  \node[orbi] at (10.038,1.2) {};
  \node[label,anchor=east] at (9.948,1.2) {$E_{3}$};
  \draw[side,draw=PrimeTwo] (10.038,1.85) -- (11.338,1.85);
  \node[orbi] at (10.688,1.85) {};
  \node[label,anchor=south] at (10.688,1.94) {$E_{6}$};
  \draw[side,draw=PrimeOne] (11.338,1.85) -- (11.338,0.55);
  \node[orbi] at (11.338,1.2) {};
  \node[label,anchor=west] at (11.428,1.2) {$E_{3}$};
  \draw[side,draw=PrimeTwo] (11.338,0.55) -- (10.038,0.55);
  \node[orbi] at (10.688,0.55) {};
  \node[label,anchor=north] at (10.688,0.46) {$E_{6}$};
  \node[corner] at (10.038,0.55) {};
  \node[corner] at (10.038,1.85) {};
  \node[corner] at (11.338,1.85) {};
  \node[corner] at (11.338,0.55) {};
  \node[inner sep=1.4pt,fill=SoftInk!5,rounded corners=1pt] at (10.688,1.2) {$F_{4}$};
  \begin{scope}[shift={(1.225,0)}]
  \fill[cell] (2.688,-1.95) rectangle (3.988,-0.65);
  \draw[mirror] (2.688,-1.3) -- (3.988,-1.3);
  \draw[side,draw=PrimeOne] (2.688,-1.95) -- (2.688,-0.65);
  \node[orbi] at (2.688,-1.3) {};
  \node[label,anchor=east] at (2.598,-1.3) {$E_{4}$};
  \draw[side,draw=PrimeTwo] (2.688,-0.65) -- (3.988,-0.65);
  \draw[tip,draw=PrimeTwo,shift={(3.338,-0.65)},rotate=180] (-2.5pt,2.1pt) -- (0,0) -- (-2.5pt,-2.1pt);
  \node[label,anchor=south] at (3.338,-0.56) {$E_{5}$};
  \draw[side,draw=PrimeOne] (3.988,-0.65) -- (3.988,-1.95);
  \node[orbi] at (3.988,-1.3) {};
  \node[label,anchor=west] at (4.078,-1.3) {$E_{3}$};
  \draw[side,draw=PrimeTwo] (3.988,-1.95) -- (2.688,-1.95);
  \draw[tip,draw=PrimeTwo,shift={(3.338,-1.95)},rotate=180] (-2.5pt,2.1pt) -- (0,0) -- (-2.5pt,-2.1pt);
  \node[label,anchor=north] at (3.338,-2.04) {$E_{5}$};
  \node[corner] at (2.688,-1.95) {};
  \node[corner] at (2.688,-0.65) {};
  \node[corner] at (3.988,-0.65) {};
  \node[corner] at (3.988,-1.95) {};
  \node[inner sep=1.4pt,fill=SoftInk!5,rounded corners=1pt] at (3.338,-1.3) {$F_{5}$};
  \fill[cell] (5.138,-1.95) rectangle (6.438,-0.65);
  \draw[side,draw=PrimeOne] (5.138,-1.95) -- (5.138,-0.65);
  \node[orbi] at (5.138,-1.3) {};
  \node[label,anchor=east] at (5.048,-1.3) {$E_{3}$};
  \draw[side,draw=PrimeTwo] (5.138,-0.65) -- (6.438,-0.65);
  \draw[tip,draw=PrimeTwo,shift={(5.788,-0.65)},rotate=180] (-2.5pt,2.1pt) -- (0,0) -- (-2.5pt,-2.1pt);
  \node[label,anchor=south] at (5.788,-0.56) {$E_{5}$};
  \draw[side,draw=PrimeOne] (6.438,-0.65) -- (6.438,-1.95);
  \node[orbi] at (6.438,-1.3) {};
  \node[label,anchor=west] at (6.528,-1.3) {$E_{2}$};
  \draw[side,draw=PrimeTwo] (6.438,-1.95) -- (5.138,-1.95);
  \node[orbi] at (5.788,-1.95) {};
  \node[label,anchor=north] at (5.788,-2.04) {$E_{7}$};
  \node[corner] at (5.138,-1.95) {};
  \node[corner] at (5.138,-0.65) {};
  \node[corner] at (6.438,-0.65) {};
  \node[corner] at (6.438,-1.95) {};
  \node[inner sep=1.4pt,fill=SoftInk!5,rounded corners=1pt] at (5.788,-1.3) {$F_{6}$};
  \fill[cell] (7.588,-1.95) rectangle (8.888,-0.65);
  \draw[mirror] (7.588,-1.3) -- (8.888,-1.3);
  \draw[side,draw=PrimeOne] (7.588,-1.95) -- (7.588,-0.65);
  \node[orbi] at (7.588,-1.3) {};
  \node[label,anchor=east] at (7.498,-1.3) {$E_{4}$};
  \draw[side,draw=PrimeTwo] (7.588,-0.65) -- (8.888,-0.65);
  \draw[tip,draw=PrimeTwo,shift={(8.238,-0.65)},rotate=0] (-2.5pt,2.1pt) -- (0,0) -- (-2.5pt,-2.1pt);
  \node[label,anchor=south] at (8.238,-0.56) {$E_{5}$};
  \draw[side,draw=PrimeOne] (8.888,-0.65) -- (8.888,-1.95);
  \node[orbi] at (8.888,-1.3) {};
  \node[label,anchor=west] at (8.978,-1.3) {$E_{2}$};
  \draw[side,draw=PrimeTwo] (8.888,-1.95) -- (7.588,-1.95);
  \draw[tip,draw=PrimeTwo,shift={(8.238,-1.95)},rotate=0] (-2.5pt,2.1pt) -- (0,0) -- (-2.5pt,-2.1pt);
  \node[label,anchor=north] at (8.238,-2.04) {$E_{5}$};
  \node[corner] at (7.588,-1.95) {};
  \node[corner] at (7.588,-0.65) {};
  \node[corner] at (8.888,-0.65) {};
  \node[corner] at (8.888,-1.95) {};
  \node[inner sep=1.4pt,fill=SoftInk!5,rounded corners=1pt] at (8.238,-1.3) {$F_{7}$};
  \end{scope}
  \draw[edge,draw=PrimeOne] (2.3,-2.9) -- (2.85,-2.9);
  \node[label,anchor=west] at (2.95,-2.9) {$\fp_3$-adic edge};
  \draw[edge,draw=PrimeTwo] (5.3,-2.9) -- (5.85,-2.9);
  \node[label,anchor=west] at (5.95,-2.9) {$\fq_3$-adic edge};
\end{tikzpicture}
    \caption{The complex $X_{\Gamma_{3}}$ of Section \ref{sec:qsqrt7}.}
    \label{fig:qsqrt7-gamma3}
\end{figure}

\begin{figure}[htbp]
    \centering
    \definecolor{SoftInk}{HTML}{2F3437}
\definecolor{EdgeColour}{HTML}{4F7D72}
\definecolor{PrimeOne}{HTML}{4F7D72}
\definecolor{PrimeTwo}{HTML}{A85D67}
\begin{tikzpicture}[x=1cm,y=1cm,
  every node/.style={font=\scriptsize,text=SoftInk},
  edge/.style={line width=1.05pt,line cap=round},
  side/.style={line width=0.95pt,line cap=round},
  tip/.style={line width=0.9pt,line cap=round,line join=round},
  cell/.style={fill=SoftInk!5},
  mirror/.style={draw=SoftInk!75,line width=0.75pt,dash pattern=on 2.2pt off 1.6pt},
  label/.style={inner sep=1.2pt},
  vert/.style={circle,fill=SoftInk,inner sep=1.8pt},
  corner/.style={circle,fill=SoftInk,inner sep=1.2pt},
  orbi/.style={circle,fill=white,draw=SoftInk,line width=0.6pt,inner sep=1.35pt}]
  \draw[edge,draw=EdgeColour] (0,0) .. controls (-0.598,2.231) and (-2.231,0.598) .. (0,0);
  \draw[tip,draw=EdgeColour,shift={(-1.061,1.061)},rotate=-135] (-2.5pt,2.1pt) -- (0,0) -- (-2.5pt,-2.1pt);
  \node[label,anchor=315] at (-1.131,1.131) {$E_{1}$};
  \draw[edge,draw=EdgeColour] (0,0) .. controls (-2.231,-0.598) and (-0.598,-2.231) .. (0,0);
  \draw[tip,draw=EdgeColour,shift={(-1.061,-1.061)},rotate=-45] (-2.5pt,2.1pt) -- (0,0) -- (-2.5pt,-2.1pt);
  \node[label,anchor=405] at (-1.131,-1.131) {$E_{2}$};
  \draw[edge,draw=EdgeColour] (0,0) .. controls (0.598,-2.231) and (2.231,-0.598) .. (0,0);
  \draw[tip,draw=EdgeColour,shift={(1.061,-1.061)},rotate=45] (-2.5pt,2.1pt) -- (0,0) -- (-2.5pt,-2.1pt);
  \node[label,anchor=495] at (1.131,-1.131) {$E_{3}$};
  \draw[edge,draw=EdgeColour] (0,0) -- (1.358,0.563);
  \node[orbi] at (1.358,0.563) {};
  \node[label,anchor=562.5] at (1.469,0.609) {$E_{4}$};
  \draw[edge,draw=EdgeColour] (0,0) -- (0.563,1.358);
  \node[orbi] at (0.563,1.358) {};
  \node[label,anchor=607.5] at (0.609,1.469) {$E_{5}$};
  \node[vert] at (0,0) {};
  \fill[cell] (2.688,0.55) rectangle (3.988,1.85);
  \draw[side,draw=EdgeColour] (2.688,0.55) -- (2.688,1.85);
  \draw[tip,draw=EdgeColour,shift={(2.688,1.2)},rotate=90] (-2.5pt,2.1pt) -- (0,0) -- (-2.5pt,-2.1pt);
  \node[label,anchor=east] at (2.598,1.2) {$E_{3}$};
  \draw[side,draw=EdgeColour] (2.688,1.85) -- (3.988,1.85);
  \node[orbi] at (3.338,1.85) {};
  \node[label,anchor=south] at (3.338,1.94) {$E_{4}$};
  \draw[side,draw=EdgeColour] (3.988,1.85) -- (3.988,0.55);
  \node[orbi] at (3.988,1.2) {};
  \node[label,anchor=west] at (4.078,1.2) {$E_{5}$};
  \draw[side,draw=EdgeColour] (3.988,0.55) -- (2.688,0.55);
  \draw[tip,draw=EdgeColour,shift={(3.338,0.55)},rotate=180] (-2.5pt,2.1pt) -- (0,0) -- (-2.5pt,-2.1pt);
  \node[label,anchor=north] at (3.338,0.46) {$E_{3}$};
  \node[corner] at (2.688,0.55) {};
  \node[corner] at (2.688,1.85) {};
  \node[corner] at (3.988,1.85) {};
  \node[corner] at (3.988,0.55) {};
  \node[inner sep=1.4pt,fill=SoftInk!5,rounded corners=1pt] at (3.338,1.2) {$F_{1}$};
  \fill[cell] (5.138,0.55) rectangle (6.438,1.85);
  \draw[mirror] (5.788,0.55) -- (5.788,1.85);
  \draw[side,draw=EdgeColour] (5.138,0.55) -- (5.138,1.85);
  \draw[tip,draw=EdgeColour,shift={(5.138,1.2)},rotate=90] (-2.5pt,2.1pt) -- (0,0) -- (-2.5pt,-2.1pt);
  \node[label,anchor=east] at (5.048,1.2) {$E_{1}$};
  \draw[side,draw=EdgeColour] (5.138,1.85) -- (6.438,1.85);
  \node[orbi] at (5.788,1.85) {};
  \node[label,anchor=south] at (5.788,1.94) {$E_{4}$};
  \draw[side,draw=EdgeColour] (6.438,1.85) -- (6.438,0.55);
  \draw[tip,draw=EdgeColour,shift={(6.438,1.2)},rotate=90] (-2.5pt,2.1pt) -- (0,0) -- (-2.5pt,-2.1pt);
  \node[label,anchor=west] at (6.528,1.2) {$E_{1}$};
  \draw[side,draw=EdgeColour] (6.438,0.55) -- (5.138,0.55);
  \node[orbi] at (5.788,0.55) {};
  \node[label,anchor=north] at (5.788,0.46) {$E_{5}$};
  \node[corner] at (5.138,0.55) {};
  \node[corner] at (5.138,1.85) {};
  \node[corner] at (6.438,1.85) {};
  \node[corner] at (6.438,0.55) {};
  \node[inner sep=1.4pt,fill=SoftInk!5,rounded corners=1pt] at (5.788,1.2) {$F_{2}$};
  \fill[cell] (7.588,0.55) rectangle (8.888,1.85);
  \draw[mirror] (8.238,0.55) -- (8.238,1.85);
  \draw[side,draw=EdgeColour] (7.588,0.55) -- (7.588,1.85);
  \draw[tip,draw=EdgeColour,shift={(7.588,1.2)},rotate=90] (-2.5pt,2.1pt) -- (0,0) -- (-2.5pt,-2.1pt);
  \node[label,anchor=east] at (7.498,1.2) {$E_{2}$};
  \draw[side,draw=EdgeColour] (7.588,1.85) -- (8.888,1.85);
  \node[orbi] at (8.238,1.85) {};
  \node[label,anchor=south] at (8.238,1.94) {$E_{5}$};
  \draw[side,draw=EdgeColour] (8.888,1.85) -- (8.888,0.55);
  \draw[tip,draw=EdgeColour,shift={(8.888,1.2)},rotate=90] (-2.5pt,2.1pt) -- (0,0) -- (-2.5pt,-2.1pt);
  \node[label,anchor=west] at (8.978,1.2) {$E_{2}$};
  \draw[side,draw=EdgeColour] (8.888,0.55) -- (7.588,0.55);
  \node[orbi] at (8.238,0.55) {};
  \node[label,anchor=north] at (8.238,0.46) {$E_{4}$};
  \node[corner] at (7.588,0.55) {};
  \node[corner] at (7.588,1.85) {};
  \node[corner] at (8.888,1.85) {};
  \node[corner] at (8.888,0.55) {};
  \node[inner sep=1.4pt,fill=SoftInk!5,rounded corners=1pt] at (8.238,1.2) {$F_{3}$};
  \begin{scope}[shift={(1.225,0)}]
  \fill[cell] (2.688,-1.95) rectangle (3.988,-0.65);
  \draw[side,draw=EdgeColour] (2.688,-1.95) -- (2.688,-0.65);
  \draw[tip,draw=EdgeColour,shift={(2.688,-1.3)},rotate=90] (-2.5pt,2.1pt) -- (0,0) -- (-2.5pt,-2.1pt);
  \node[label,anchor=east] at (2.598,-1.3) {$E_{3}$};
  \draw[side,draw=EdgeColour] (2.688,-0.65) -- (3.988,-0.65);
  \draw[tip,draw=EdgeColour,shift={(3.338,-0.65)},rotate=180] (-2.5pt,2.1pt) -- (0,0) -- (-2.5pt,-2.1pt);
  \node[label,anchor=south] at (3.338,-0.56) {$E_{2}$};
  \draw[side,draw=EdgeColour] (3.988,-0.65) -- (3.988,-1.95);
  \draw[tip,draw=EdgeColour,shift={(3.988,-1.3)},rotate=90] (-2.5pt,2.1pt) -- (0,0) -- (-2.5pt,-2.1pt);
  \node[label,anchor=west] at (4.078,-1.3) {$E_{2}$};
  \draw[side,draw=EdgeColour] (3.988,-1.95) -- (2.688,-1.95);
  \draw[tip,draw=EdgeColour,shift={(3.338,-1.95)},rotate=180] (-2.5pt,2.1pt) -- (0,0) -- (-2.5pt,-2.1pt);
  \node[label,anchor=north] at (3.338,-2.04) {$E_{1}$};
  \node[corner] at (2.688,-1.95) {};
  \node[corner] at (2.688,-0.65) {};
  \node[corner] at (3.988,-0.65) {};
  \node[corner] at (3.988,-1.95) {};
  \node[inner sep=1.4pt,fill=SoftInk!5,rounded corners=1pt] at (3.338,-1.3) {$F_{4}$};
  \fill[cell] (5.138,-1.95) rectangle (6.438,-0.65);
  \draw[side,draw=EdgeColour] (5.138,-1.95) -- (5.138,-0.65);
  \draw[tip,draw=EdgeColour,shift={(5.138,-1.3)},rotate=-90] (-2.5pt,2.1pt) -- (0,0) -- (-2.5pt,-2.1pt);
  \node[label,anchor=east] at (5.048,-1.3) {$E_{3}$};
  \draw[side,draw=EdgeColour] (5.138,-0.65) -- (6.438,-0.65);
  \draw[tip,draw=EdgeColour,shift={(5.788,-0.65)},rotate=0] (-2.5pt,2.1pt) -- (0,0) -- (-2.5pt,-2.1pt);
  \node[label,anchor=south] at (5.788,-0.56) {$E_{2}$};
  \draw[side,draw=EdgeColour] (6.438,-0.65) -- (6.438,-1.95);
  \draw[tip,draw=EdgeColour,shift={(6.438,-1.3)},rotate=90] (-2.5pt,2.1pt) -- (0,0) -- (-2.5pt,-2.1pt);
  \node[label,anchor=west] at (6.528,-1.3) {$E_{1}$};
  \draw[side,draw=EdgeColour] (6.438,-1.95) -- (5.138,-1.95);
  \draw[tip,draw=EdgeColour,shift={(5.788,-1.95)},rotate=0] (-2.5pt,2.1pt) -- (0,0) -- (-2.5pt,-2.1pt);
  \node[label,anchor=north] at (5.788,-2.04) {$E_{1}$};
  \node[corner] at (5.138,-1.95) {};
  \node[corner] at (5.138,-0.65) {};
  \node[corner] at (6.438,-0.65) {};
  \node[corner] at (6.438,-1.95) {};
  \node[inner sep=1.4pt,fill=SoftInk!5,rounded corners=1pt] at (5.788,-1.3) {$F_{5}$};
  \end{scope}
\end{tikzpicture}
    \caption{The complex $X_{\Gamma_{5}}$ of Section \ref{sec:qsqrt7}.}
    \label{fig:qsqrt7-gamma5}
\end{figure}

\begin{figure}[htbp]
    \centering
    \definecolor{SoftInk}{HTML}{2F3437}
\definecolor{EdgeColour}{HTML}{4F7D72}
\definecolor{PrimeOne}{HTML}{4F7D72}
\definecolor{PrimeTwo}{HTML}{A85D67}
\begin{tikzpicture}[x=1cm,y=1cm,
  every node/.style={font=\scriptsize,text=SoftInk},
  edge/.style={line width=1.05pt,line cap=round},
  side/.style={line width=0.95pt,line cap=round},
  tip/.style={line width=0.9pt,line cap=round,line join=round},
  cell/.style={fill=SoftInk!5},
  mirror/.style={draw=SoftInk!75,line width=0.75pt,dash pattern=on 2.2pt off 1.6pt},
  label/.style={inner sep=1.2pt},
  vert/.style={circle,fill=SoftInk,inner sep=1.8pt},
  corner/.style={circle,fill=SoftInk,inner sep=1.2pt},
  orbi/.style={circle,fill=white,draw=SoftInk,line width=0.6pt,inner sep=1.35pt}]
  \draw[edge,draw=EdgeColour] (0,0) .. controls (-0.598,2.231) and (-2.231,0.598) .. (0,0);
  \draw[tip,draw=EdgeColour,shift={(-1.061,1.061)},rotate=-135] (-2.5pt,2.1pt) -- (0,0) -- (-2.5pt,-2.1pt);
  \node[label,anchor=315] at (-1.131,1.131) {$E_{1}$};
  \draw[edge,draw=EdgeColour] (0,0) .. controls (-2.231,-0.598) and (-0.598,-2.231) .. (0,0);
  \draw[tip,draw=EdgeColour,shift={(-1.061,-1.061)},rotate=-45] (-2.5pt,2.1pt) -- (0,0) -- (-2.5pt,-2.1pt);
  \node[label,anchor=405] at (-1.131,-1.131) {$E_{2}$};
  \draw[edge,draw=EdgeColour] (0,0) .. controls (0.598,-2.231) and (2.231,-0.598) .. (0,0);
  \draw[tip,draw=EdgeColour,shift={(1.061,-1.061)},rotate=45] (-2.5pt,2.1pt) -- (0,0) -- (-2.5pt,-2.1pt);
  \node[label,anchor=495] at (1.131,-1.131) {$E_{3}$};
  \draw[edge,draw=EdgeColour] (0,0) -- (1.358,0.563);
  \node[orbi] at (1.358,0.563) {};
  \node[label,anchor=562.5] at (1.469,0.609) {$E_{4}$};
  \draw[edge,draw=EdgeColour] (0,0) -- (0.563,1.358);
  \node[orbi] at (0.563,1.358) {};
  \node[label,anchor=607.5] at (0.609,1.469) {$E_{5}$};
  \node[vert] at (0,0) {};
  \fill[cell] (2.688,0.55) rectangle (3.988,1.85);
  \draw[side,draw=EdgeColour] (2.688,0.55) -- (2.688,1.85);
  \draw[tip,draw=EdgeColour,shift={(2.688,1.2)},rotate=90] (-2.5pt,2.1pt) -- (0,0) -- (-2.5pt,-2.1pt);
  \node[label,anchor=east] at (2.598,1.2) {$E_{3}$};
  \draw[side,draw=EdgeColour] (2.688,1.85) -- (3.988,1.85);
  \node[orbi] at (3.338,1.85) {};
  \node[label,anchor=south] at (3.338,1.94) {$E_{4}$};
  \draw[side,draw=EdgeColour] (3.988,1.85) -- (3.988,0.55);
  \node[orbi] at (3.988,1.2) {};
  \node[label,anchor=west] at (4.078,1.2) {$E_{5}$};
  \draw[side,draw=EdgeColour] (3.988,0.55) -- (2.688,0.55);
  \draw[tip,draw=EdgeColour,shift={(3.338,0.55)},rotate=180] (-2.5pt,2.1pt) -- (0,0) -- (-2.5pt,-2.1pt);
  \node[label,anchor=north] at (3.338,0.46) {$E_{3}$};
  \node[corner] at (2.688,0.55) {};
  \node[corner] at (2.688,1.85) {};
  \node[corner] at (3.988,1.85) {};
  \node[corner] at (3.988,0.55) {};
  \node[inner sep=1.4pt,fill=SoftInk!5,rounded corners=1pt] at (3.338,1.2) {$F_{1}$};
  \fill[cell] (5.138,0.55) rectangle (6.438,1.85);
  \draw[side,draw=EdgeColour] (5.138,0.55) -- (5.138,1.85);
  \draw[tip,draw=EdgeColour,shift={(5.138,1.2)},rotate=90] (-2.5pt,2.1pt) -- (0,0) -- (-2.5pt,-2.1pt);
  \node[label,anchor=east] at (5.048,1.2) {$E_{1}$};
  \draw[side,draw=EdgeColour] (5.138,1.85) -- (6.438,1.85);
  \node[orbi] at (5.788,1.85) {};
  \node[label,anchor=south] at (5.788,1.94) {$E_{4}$};
  \draw[side,draw=EdgeColour] (6.438,1.85) -- (6.438,0.55);
  \draw[tip,draw=EdgeColour,shift={(6.438,1.2)},rotate=-90] (-2.5pt,2.1pt) -- (0,0) -- (-2.5pt,-2.1pt);
  \node[label,anchor=west] at (6.528,1.2) {$E_{2}$};
  \draw[side,draw=EdgeColour] (6.438,0.55) -- (5.138,0.55);
  \node[orbi] at (5.788,0.55) {};
  \node[label,anchor=north] at (5.788,0.46) {$E_{5}$};
  \node[corner] at (5.138,0.55) {};
  \node[corner] at (5.138,1.85) {};
  \node[corner] at (6.438,1.85) {};
  \node[corner] at (6.438,0.55) {};
  \node[inner sep=1.4pt,fill=SoftInk!5,rounded corners=1pt] at (5.788,1.2) {$F_{2}$};
  \fill[cell] (2.688,-1.95) rectangle (3.988,-0.65);
  \draw[side,draw=EdgeColour] (2.688,-1.95) -- (2.688,-0.65);
  \draw[tip,draw=EdgeColour,shift={(2.688,-1.3)},rotate=90] (-2.5pt,2.1pt) -- (0,0) -- (-2.5pt,-2.1pt);
  \node[label,anchor=east] at (2.598,-1.3) {$E_{3}$};
  \draw[side,draw=EdgeColour] (2.688,-0.65) -- (3.988,-0.65);
  \draw[tip,draw=EdgeColour,shift={(3.338,-0.65)},rotate=180] (-2.5pt,2.1pt) -- (0,0) -- (-2.5pt,-2.1pt);
  \node[label,anchor=south] at (3.338,-0.56) {$E_{2}$};
  \draw[side,draw=EdgeColour] (3.988,-0.65) -- (3.988,-1.95);
  \draw[tip,draw=EdgeColour,shift={(3.988,-1.3)},rotate=90] (-2.5pt,2.1pt) -- (0,0) -- (-2.5pt,-2.1pt);
  \node[label,anchor=west] at (4.078,-1.3) {$E_{2}$};
  \draw[side,draw=EdgeColour] (3.988,-1.95) -- (2.688,-1.95);
  \draw[tip,draw=EdgeColour,shift={(3.338,-1.95)},rotate=180] (-2.5pt,2.1pt) -- (0,0) -- (-2.5pt,-2.1pt);
  \node[label,anchor=north] at (3.338,-2.04) {$E_{1}$};
  \node[corner] at (2.688,-1.95) {};
  \node[corner] at (2.688,-0.65) {};
  \node[corner] at (3.988,-0.65) {};
  \node[corner] at (3.988,-1.95) {};
  \node[inner sep=1.4pt,fill=SoftInk!5,rounded corners=1pt] at (3.338,-1.3) {$F_{3}$};
  \fill[cell] (5.138,-1.95) rectangle (6.438,-0.65);
  \draw[side,draw=EdgeColour] (5.138,-1.95) -- (5.138,-0.65);
  \draw[tip,draw=EdgeColour,shift={(5.138,-1.3)},rotate=-90] (-2.5pt,2.1pt) -- (0,0) -- (-2.5pt,-2.1pt);
  \node[label,anchor=east] at (5.048,-1.3) {$E_{3}$};
  \draw[side,draw=EdgeColour] (5.138,-0.65) -- (6.438,-0.65);
  \draw[tip,draw=EdgeColour,shift={(5.788,-0.65)},rotate=0] (-2.5pt,2.1pt) -- (0,0) -- (-2.5pt,-2.1pt);
  \node[label,anchor=south] at (5.788,-0.56) {$E_{2}$};
  \draw[side,draw=EdgeColour] (6.438,-0.65) -- (6.438,-1.95);
  \draw[tip,draw=EdgeColour,shift={(6.438,-1.3)},rotate=90] (-2.5pt,2.1pt) -- (0,0) -- (-2.5pt,-2.1pt);
  \node[label,anchor=west] at (6.528,-1.3) {$E_{1}$};
  \draw[side,draw=EdgeColour] (6.438,-1.95) -- (5.138,-1.95);
  \draw[tip,draw=EdgeColour,shift={(5.788,-1.95)},rotate=0] (-2.5pt,2.1pt) -- (0,0) -- (-2.5pt,-2.1pt);
  \node[label,anchor=north] at (5.788,-2.04) {$E_{1}$};
  \node[corner] at (5.138,-1.95) {};
  \node[corner] at (5.138,-0.65) {};
  \node[corner] at (6.438,-0.65) {};
  \node[corner] at (6.438,-1.95) {};
  \node[inner sep=1.4pt,fill=SoftInk!5,rounded corners=1pt] at (5.788,-1.3) {$F_{4}$};
\end{tikzpicture}
    \caption{The complex $X_{\Gamma_{6}}$ of Section \ref{sec:qsqrt7}.}
    \label{fig:qsqrt7-gamma6}
\end{figure}

\begin{figure}[htbp]
    \centering
    \definecolor{SoftInk}{HTML}{2F3437}
\definecolor{EdgeColour}{HTML}{4F7D72}
\definecolor{PrimeOne}{HTML}{4F7D72}
\definecolor{PrimeTwo}{HTML}{A85D67}
\begin{tikzpicture}[x=1cm,y=1cm,
  every node/.style={font=\scriptsize,text=SoftInk},
  edge/.style={line width=1.05pt,line cap=round},
  side/.style={line width=0.95pt,line cap=round},
  tip/.style={line width=0.9pt,line cap=round,line join=round},
  cell/.style={fill=SoftInk!5},
  mirror/.style={draw=SoftInk!75,line width=0.75pt,dash pattern=on 2.2pt off 1.6pt},
  label/.style={inner sep=1.2pt},
  vert/.style={circle,fill=SoftInk,inner sep=1.8pt},
  corner/.style={circle,fill=SoftInk,inner sep=1.2pt},
  orbi/.style={circle,fill=white,draw=SoftInk,line width=0.6pt,inner sep=1.35pt}]
  \draw[edge,draw=EdgeColour] (0,0) .. controls (-0.598,2.231) and (-2.231,0.598) .. (0,0);
  \draw[tip,draw=EdgeColour,shift={(-1.061,1.061)},rotate=-135] (-2.5pt,2.1pt) -- (0,0) -- (-2.5pt,-2.1pt);
  \node[label,anchor=315] at (-1.131,1.131) {$E_{1}$};
  \draw[edge,draw=EdgeColour] (0,0) -- (-1.358,-0.563);
  \node[orbi] at (-1.358,-0.563) {};
  \node[label,anchor=382.5] at (-1.469,-0.609) {$E_{2}$};
  \draw[edge,draw=EdgeColour] (0,0) -- (-0.563,-1.358);
  \node[orbi] at (-0.563,-1.358) {};
  \node[label,anchor=427.5] at (-0.609,-1.469) {$E_{3}$};
  \draw[edge,draw=EdgeColour] (0,0) -- (0.563,-1.358);
  \node[orbi] at (0.563,-1.358) {};
  \node[label,anchor=472.5] at (0.609,-1.469) {$E_{4}$};
  \draw[edge,draw=EdgeColour] (0,0) -- (1.358,-0.563);
  \node[orbi] at (1.358,-0.563) {};
  \node[label,anchor=517.5] at (1.469,-0.609) {$E_{5}$};
  \draw[edge,draw=EdgeColour] (0,0) -- (1.358,0.563);
  \node[orbi] at (1.358,0.563) {};
  \node[label,anchor=562.5] at (1.469,0.609) {$E_{6}$};
  \draw[edge,draw=EdgeColour] (0,0) -- (0.563,1.358);
  \node[orbi] at (0.563,1.358) {};
  \node[label,anchor=607.5] at (0.609,1.469) {$E_{7}$};
  \node[vert] at (0,0) {};
  \fill[cell] (2.688,0.55) rectangle (3.988,1.85);
  \draw[side,draw=EdgeColour] (2.688,0.55) -- (2.688,1.85);
  \node[orbi] at (2.688,1.2) {};
  \node[label,anchor=east] at (2.598,1.2) {$E_{6}$};
  \draw[side,draw=EdgeColour] (2.688,1.85) -- (3.988,1.85);
  \draw[tip,draw=EdgeColour,shift={(3.338,1.85)},rotate=180] (-2.5pt,2.1pt) -- (0,0) -- (-2.5pt,-2.1pt);
  \node[label,anchor=south] at (3.338,1.94) {$E_{1}$};
  \draw[side,draw=EdgeColour] (3.988,1.85) -- (3.988,0.55);
  \draw[tip,draw=EdgeColour,shift={(3.988,1.2)},rotate=90] (-2.5pt,2.1pt) -- (0,0) -- (-2.5pt,-2.1pt);
  \node[label,anchor=west] at (4.078,1.2) {$E_{1}$};
  \draw[side,draw=EdgeColour] (3.988,0.55) -- (2.688,0.55);
  \node[orbi] at (3.338,0.55) {};
  \node[label,anchor=north] at (3.338,0.46) {$E_{7}$};
  \node[corner] at (2.688,0.55) {};
  \node[corner] at (2.688,1.85) {};
  \node[corner] at (3.988,1.85) {};
  \node[corner] at (3.988,0.55) {};
  \node[inner sep=1.4pt,fill=SoftInk!5,rounded corners=1pt] at (3.338,1.2) {$F_{1}$};
  \fill[cell] (5.138,0.55) rectangle (6.438,1.85);
  \draw[mirror] (5.138,1.2) -- (6.438,1.2);
  \draw[side,draw=EdgeColour] (5.138,0.55) -- (5.138,1.85);
  \node[orbi] at (5.138,1.2) {};
  \node[label,anchor=east] at (5.048,1.2) {$E_{2}$};
  \draw[side,draw=EdgeColour] (5.138,1.85) -- (6.438,1.85);
  \draw[tip,draw=EdgeColour,shift={(5.788,1.85)},rotate=180] (-2.5pt,2.1pt) -- (0,0) -- (-2.5pt,-2.1pt);
  \node[label,anchor=south] at (5.788,1.94) {$E_{1}$};
  \draw[side,draw=EdgeColour] (6.438,1.85) -- (6.438,0.55);
  \node[orbi] at (6.438,1.2) {};
  \node[label,anchor=west] at (6.528,1.2) {$E_{4}$};
  \draw[side,draw=EdgeColour] (6.438,0.55) -- (5.138,0.55);
  \draw[tip,draw=EdgeColour,shift={(5.788,0.55)},rotate=180] (-2.5pt,2.1pt) -- (0,0) -- (-2.5pt,-2.1pt);
  \node[label,anchor=north] at (5.788,0.46) {$E_{1}$};
  \node[corner] at (5.138,0.55) {};
  \node[corner] at (5.138,1.85) {};
  \node[corner] at (6.438,1.85) {};
  \node[corner] at (6.438,0.55) {};
  \node[inner sep=1.4pt,fill=SoftInk!5,rounded corners=1pt] at (5.788,1.2) {$F_{2}$};
  \fill[cell] (7.588,0.55) rectangle (8.888,1.85);
  \draw[mirror] (7.588,1.2) -- (8.888,1.2);
  \draw[side,draw=EdgeColour] (7.588,0.55) -- (7.588,1.85);
  \node[orbi] at (7.588,1.2) {};
  \node[label,anchor=east] at (7.498,1.2) {$E_{3}$};
  \draw[side,draw=EdgeColour] (7.588,1.85) -- (8.888,1.85);
  \draw[tip,draw=EdgeColour,shift={(8.238,1.85)},rotate=0] (-2.5pt,2.1pt) -- (0,0) -- (-2.5pt,-2.1pt);
  \node[label,anchor=south] at (8.238,1.94) {$E_{1}$};
  \draw[side,draw=EdgeColour] (8.888,1.85) -- (8.888,0.55);
  \node[orbi] at (8.888,1.2) {};
  \node[label,anchor=west] at (8.978,1.2) {$E_{5}$};
  \draw[side,draw=EdgeColour] (8.888,0.55) -- (7.588,0.55);
  \draw[tip,draw=EdgeColour,shift={(8.238,0.55)},rotate=0] (-2.5pt,2.1pt) -- (0,0) -- (-2.5pt,-2.1pt);
  \node[label,anchor=north] at (8.238,0.46) {$E_{1}$};
  \node[corner] at (7.588,0.55) {};
  \node[corner] at (7.588,1.85) {};
  \node[corner] at (8.888,1.85) {};
  \node[corner] at (8.888,0.55) {};
  \node[inner sep=1.4pt,fill=SoftInk!5,rounded corners=1pt] at (8.238,1.2) {$F_{3}$};
  \fill[cell] (10.038,0.55) rectangle (11.338,1.85);
  \draw[mirror] (10.038,1.2) -- (11.338,1.2);
  \draw[side,draw=EdgeColour] (10.038,0.55) -- (10.038,1.85);
  \node[orbi] at (10.038,1.2) {};
  \node[label,anchor=east] at (9.948,1.2) {$E_{6}$};
  \draw[side,draw=EdgeColour] (10.038,1.85) -- (11.338,1.85);
  \node[orbi] at (10.688,1.85) {};
  \node[label,anchor=south] at (10.688,1.94) {$E_{4}$};
  \draw[side,draw=EdgeColour] (11.338,1.85) -- (11.338,0.55);
  \node[orbi] at (11.338,1.2) {};
  \node[label,anchor=west] at (11.428,1.2) {$E_{2}$};
  \draw[side,draw=EdgeColour] (11.338,0.55) -- (10.038,0.55);
  \node[orbi] at (10.688,0.55) {};
  \node[label,anchor=north] at (10.688,0.46) {$E_{4}$};
  \node[corner] at (10.038,0.55) {};
  \node[corner] at (10.038,1.85) {};
  \node[corner] at (11.338,1.85) {};
  \node[corner] at (11.338,0.55) {};
  \node[inner sep=1.4pt,fill=SoftInk!5,rounded corners=1pt] at (10.688,1.2) {$F_{4}$};
  \begin{scope}[shift={(1.225,0)}]
  \fill[cell] (2.688,-1.95) rectangle (3.988,-0.65);
  \draw[mirror] (2.688,-1.3) -- (3.988,-1.3);
  \draw[side,draw=EdgeColour] (2.688,-1.95) -- (2.688,-0.65);
  \node[orbi] at (2.688,-1.3) {};
  \node[label,anchor=east] at (2.598,-1.3) {$E_{7}$};
  \draw[side,draw=EdgeColour] (2.688,-0.65) -- (3.988,-0.65);
  \node[orbi] at (3.338,-0.65) {};
  \node[label,anchor=south] at (3.338,-0.56) {$E_{2}$};
  \draw[side,draw=EdgeColour] (3.988,-0.65) -- (3.988,-1.95);
  \node[orbi] at (3.988,-1.3) {};
  \node[label,anchor=west] at (4.078,-1.3) {$E_{3}$};
  \draw[side,draw=EdgeColour] (3.988,-1.95) -- (2.688,-1.95);
  \node[orbi] at (3.338,-1.95) {};
  \node[label,anchor=north] at (3.338,-2.04) {$E_{2}$};
  \node[corner] at (2.688,-1.95) {};
  \node[corner] at (2.688,-0.65) {};
  \node[corner] at (3.988,-0.65) {};
  \node[corner] at (3.988,-1.95) {};
  \node[inner sep=1.4pt,fill=SoftInk!5,rounded corners=1pt] at (3.338,-1.3) {$F_{5}$};
  \fill[cell] (5.138,-1.95) rectangle (6.438,-0.65);
  \draw[mirror] (5.138,-1.3) -- (6.438,-1.3);
  \draw[side,draw=EdgeColour] (5.138,-1.95) -- (5.138,-0.65);
  \node[orbi] at (5.138,-1.3) {};
  \node[label,anchor=east] at (5.048,-1.3) {$E_{7}$};
  \draw[side,draw=EdgeColour] (5.138,-0.65) -- (6.438,-0.65);
  \node[orbi] at (5.788,-0.65) {};
  \node[label,anchor=south] at (5.788,-0.56) {$E_{5}$};
  \draw[side,draw=EdgeColour] (6.438,-0.65) -- (6.438,-1.95);
  \node[orbi] at (6.438,-1.3) {};
  \node[label,anchor=west] at (6.528,-1.3) {$E_{4}$};
  \draw[side,draw=EdgeColour] (6.438,-1.95) -- (5.138,-1.95);
  \node[orbi] at (5.788,-1.95) {};
  \node[label,anchor=north] at (5.788,-2.04) {$E_{5}$};
  \node[corner] at (5.138,-1.95) {};
  \node[corner] at (5.138,-0.65) {};
  \node[corner] at (6.438,-0.65) {};
  \node[corner] at (6.438,-1.95) {};
  \node[inner sep=1.4pt,fill=SoftInk!5,rounded corners=1pt] at (5.788,-1.3) {$F_{6}$};
  \fill[cell] (7.588,-1.95) rectangle (8.888,-0.65);
  \draw[mirror] (7.588,-1.3) -- (8.888,-1.3);
  \draw[side,draw=EdgeColour] (7.588,-1.95) -- (7.588,-0.65);
  \node[orbi] at (7.588,-1.3) {};
  \node[label,anchor=east] at (7.498,-1.3) {$E_{6}$};
  \draw[side,draw=EdgeColour] (7.588,-0.65) -- (8.888,-0.65);
  \node[orbi] at (8.238,-0.65) {};
  \node[label,anchor=south] at (8.238,-0.56) {$E_{3}$};
  \draw[side,draw=EdgeColour] (8.888,-0.65) -- (8.888,-1.95);
  \node[orbi] at (8.888,-1.3) {};
  \node[label,anchor=west] at (8.978,-1.3) {$E_{5}$};
  \draw[side,draw=EdgeColour] (8.888,-1.95) -- (7.588,-1.95);
  \node[orbi] at (8.238,-1.95) {};
  \node[label,anchor=north] at (8.238,-2.04) {$E_{3}$};
  \node[corner] at (7.588,-1.95) {};
  \node[corner] at (7.588,-0.65) {};
  \node[corner] at (8.888,-0.65) {};
  \node[corner] at (8.888,-1.95) {};
  \node[inner sep=1.4pt,fill=SoftInk!5,rounded corners=1pt] at (8.238,-1.3) {$F_{7}$};
  \end{scope}
\end{tikzpicture}
    \caption{The complex $X_{\Gamma_{9}}$ of Section \ref{sec:qsqrt7}.}
    \label{fig:qsqrt7-gamma9}
\end{figure}

\begin{figure}[htbp]
    \centering
    \definecolor{SoftInk}{HTML}{2F3437}
\definecolor{EdgeColour}{HTML}{4F7D72}
\definecolor{PrimeOne}{HTML}{4F7D72}
\definecolor{PrimeTwo}{HTML}{A85D67}
\begin{tikzpicture}[x=1cm,y=1cm,
  every node/.style={font=\scriptsize,text=SoftInk},
  edge/.style={line width=1.05pt,line cap=round},
  side/.style={line width=0.95pt,line cap=round},
  tip/.style={line width=0.9pt,line cap=round,line join=round},
  cell/.style={fill=SoftInk!5},
  mirror/.style={draw=SoftInk!75,line width=0.75pt,dash pattern=on 2.2pt off 1.6pt},
  label/.style={inner sep=1.2pt},
  vert/.style={circle,fill=SoftInk,inner sep=1.8pt},
  corner/.style={circle,fill=SoftInk,inner sep=1.2pt},
  orbi/.style={circle,fill=white,draw=SoftInk,line width=0.6pt,inner sep=1.35pt}]
  \draw[edge,draw=EdgeColour] (0,0) .. controls (-0.598,2.231) and (-2.231,0.598) .. (0,0);
  \draw[tip,draw=EdgeColour,shift={(-1.061,1.061)},rotate=-135] (-2.5pt,2.1pt) -- (0,0) -- (-2.5pt,-2.1pt);
  \node[label,anchor=315] at (-1.131,1.131) {$E_{1}$};
  \draw[edge,draw=EdgeColour] (0,0) .. controls (-2.231,-0.598) and (-0.598,-2.231) .. (0,0);
  \draw[tip,draw=EdgeColour,shift={(-1.061,-1.061)},rotate=-45] (-2.5pt,2.1pt) -- (0,0) -- (-2.5pt,-2.1pt);
  \node[label,anchor=405] at (-1.131,-1.131) {$E_{2}$};
  \draw[edge,draw=EdgeColour] (0,0) .. controls (0.598,-2.231) and (2.231,-0.598) .. (0,0);
  \draw[tip,draw=EdgeColour,shift={(1.061,-1.061)},rotate=45] (-2.5pt,2.1pt) -- (0,0) -- (-2.5pt,-2.1pt);
  \node[label,anchor=495] at (1.131,-1.131) {$E_{3}$};
  \draw[edge,draw=EdgeColour] (0,0) -- (1.358,0.563);
  \node[orbi] at (1.358,0.563) {};
  \node[label,anchor=562.5] at (1.469,0.609) {$E_{4}$};
  \draw[edge,draw=EdgeColour] (0,0) -- (0.563,1.358);
  \node[orbi] at (0.563,1.358) {};
  \node[label,anchor=607.5] at (0.609,1.469) {$E_{5}$};
  \node[vert] at (0,0) {};
  \fill[cell] (2.688,0.55) rectangle (3.988,1.85);
  \draw[side,draw=EdgeColour] (2.688,0.55) -- (2.688,1.85);
  \node[orbi] at (2.688,1.2) {};
  \node[label,anchor=east] at (2.598,1.2) {$E_{4}$};
  \draw[side,draw=EdgeColour] (2.688,1.85) -- (3.988,1.85);
  \draw[tip,draw=EdgeColour,shift={(3.338,1.85)},rotate=180] (-2.5pt,2.1pt) -- (0,0) -- (-2.5pt,-2.1pt);
  \node[label,anchor=south] at (3.338,1.94) {$E_{1}$};
  \draw[side,draw=EdgeColour] (3.988,1.85) -- (3.988,0.55);
  \draw[tip,draw=EdgeColour,shift={(3.988,1.2)},rotate=90] (-2.5pt,2.1pt) -- (0,0) -- (-2.5pt,-2.1pt);
  \node[label,anchor=west] at (4.078,1.2) {$E_{1}$};
  \draw[side,draw=EdgeColour] (3.988,0.55) -- (2.688,0.55);
  \node[orbi] at (3.338,0.55) {};
  \node[label,anchor=north] at (3.338,0.46) {$E_{5}$};
  \node[corner] at (2.688,0.55) {};
  \node[corner] at (2.688,1.85) {};
  \node[corner] at (3.988,1.85) {};
  \node[corner] at (3.988,0.55) {};
  \node[inner sep=1.4pt,fill=SoftInk!5,rounded corners=1pt] at (3.338,1.2) {$F_{1}$};
  \fill[cell] (5.138,0.55) rectangle (6.438,1.85);
  \draw[side,draw=EdgeColour] (5.138,0.55) -- (5.138,1.85);
  \draw[tip,draw=EdgeColour,shift={(5.138,1.2)},rotate=90] (-2.5pt,2.1pt) -- (0,0) -- (-2.5pt,-2.1pt);
  \node[label,anchor=east] at (5.048,1.2) {$E_{2}$};
  \draw[side,draw=EdgeColour] (5.138,1.85) -- (6.438,1.85);
  \draw[tip,draw=EdgeColour,shift={(5.788,1.85)},rotate=180] (-2.5pt,2.1pt) -- (0,0) -- (-2.5pt,-2.1pt);
  \node[label,anchor=south] at (5.788,1.94) {$E_{1}$};
  \draw[side,draw=EdgeColour] (6.438,1.85) -- (6.438,0.55);
  \draw[tip,draw=EdgeColour,shift={(6.438,1.2)},rotate=-90] (-2.5pt,2.1pt) -- (0,0) -- (-2.5pt,-2.1pt);
  \node[label,anchor=west] at (6.528,1.2) {$E_{3}$};
  \draw[side,draw=EdgeColour] (6.438,0.55) -- (5.138,0.55);
  \draw[tip,draw=EdgeColour,shift={(5.788,0.55)},rotate=180] (-2.5pt,2.1pt) -- (0,0) -- (-2.5pt,-2.1pt);
  \node[label,anchor=north] at (5.788,0.46) {$E_{1}$};
  \node[corner] at (5.138,0.55) {};
  \node[corner] at (5.138,1.85) {};
  \node[corner] at (6.438,1.85) {};
  \node[corner] at (6.438,0.55) {};
  \node[inner sep=1.4pt,fill=SoftInk!5,rounded corners=1pt] at (5.788,1.2) {$F_{2}$};
  \fill[cell] (2.688,-1.95) rectangle (3.988,-0.65);
  \draw[side,draw=EdgeColour] (2.688,-1.95) -- (2.688,-0.65);
  \node[orbi] at (2.688,-1.3) {};
  \node[label,anchor=east] at (2.598,-1.3) {$E_{4}$};
  \draw[side,draw=EdgeColour] (2.688,-0.65) -- (3.988,-0.65);
  \draw[tip,draw=EdgeColour,shift={(3.338,-0.65)},rotate=0] (-2.5pt,2.1pt) -- (0,0) -- (-2.5pt,-2.1pt);
  \node[label,anchor=south] at (3.338,-0.56) {$E_{3}$};
  \draw[side,draw=EdgeColour] (3.988,-0.65) -- (3.988,-1.95);
  \draw[tip,draw=EdgeColour,shift={(3.988,-1.3)},rotate=-90] (-2.5pt,2.1pt) -- (0,0) -- (-2.5pt,-2.1pt);
  \node[label,anchor=west] at (4.078,-1.3) {$E_{2}$};
  \draw[side,draw=EdgeColour] (3.988,-1.95) -- (2.688,-1.95);
  \draw[tip,draw=EdgeColour,shift={(3.338,-1.95)},rotate=180] (-2.5pt,2.1pt) -- (0,0) -- (-2.5pt,-2.1pt);
  \node[label,anchor=north] at (3.338,-2.04) {$E_{3}$};
  \node[corner] at (2.688,-1.95) {};
  \node[corner] at (2.688,-0.65) {};
  \node[corner] at (3.988,-0.65) {};
  \node[corner] at (3.988,-1.95) {};
  \node[inner sep=1.4pt,fill=SoftInk!5,rounded corners=1pt] at (3.338,-1.3) {$F_{3}$};
  \fill[cell] (5.138,-1.95) rectangle (6.438,-0.65);
  \draw[side,draw=EdgeColour] (5.138,-1.95) -- (5.138,-0.65);
  \node[orbi] at (5.138,-1.3) {};
  \node[label,anchor=east] at (5.048,-1.3) {$E_{5}$};
  \draw[side,draw=EdgeColour] (5.138,-0.65) -- (6.438,-0.65);
  \draw[tip,draw=EdgeColour,shift={(5.788,-0.65)},rotate=180] (-2.5pt,2.1pt) -- (0,0) -- (-2.5pt,-2.1pt);
  \node[label,anchor=south] at (5.788,-0.56) {$E_{2}$};
  \draw[side,draw=EdgeColour] (6.438,-0.65) -- (6.438,-1.95);
  \draw[tip,draw=EdgeColour,shift={(6.438,-1.3)},rotate=-90] (-2.5pt,2.1pt) -- (0,0) -- (-2.5pt,-2.1pt);
  \node[label,anchor=west] at (6.528,-1.3) {$E_{3}$};
  \draw[side,draw=EdgeColour] (6.438,-1.95) -- (5.138,-1.95);
  \draw[tip,draw=EdgeColour,shift={(5.788,-1.95)},rotate=0] (-2.5pt,2.1pt) -- (0,0) -- (-2.5pt,-2.1pt);
  \node[label,anchor=north] at (5.788,-2.04) {$E_{2}$};
  \node[corner] at (5.138,-1.95) {};
  \node[corner] at (5.138,-0.65) {};
  \node[corner] at (6.438,-0.65) {};
  \node[corner] at (6.438,-1.95) {};
  \node[inner sep=1.4pt,fill=SoftInk!5,rounded corners=1pt] at (5.788,-1.3) {$F_{4}$};
\end{tikzpicture}
    \caption{The complex $X_{\Gamma_{10}}$ of Section \ref{sec:qsqrt7}.}
    \label{fig:qsqrt7-gamma10}
\end{figure}

\begin{figure}[htbp]
    \centering
    \definecolor{SoftInk}{HTML}{2F3437}
\definecolor{EdgeColour}{HTML}{4F7D72}
\definecolor{PrimeOne}{HTML}{4F7D72}
\definecolor{PrimeTwo}{HTML}{A85D67}
\begin{tikzpicture}[x=1cm,y=1cm,
  every node/.style={font=\scriptsize,text=SoftInk},
  edge/.style={line width=1.05pt,line cap=round},
  side/.style={line width=0.95pt,line cap=round},
  tip/.style={line width=0.9pt,line cap=round,line join=round},
  cell/.style={fill=SoftInk!5},
  mirror/.style={draw=SoftInk!75,line width=0.75pt,dash pattern=on 2.2pt off 1.6pt},
  label/.style={inner sep=1.2pt},
  vert/.style={circle,fill=SoftInk,inner sep=1.8pt},
  corner/.style={circle,fill=SoftInk,inner sep=1.2pt},
  orbi/.style={circle,fill=white,draw=SoftInk,line width=0.6pt,inner sep=1.35pt}]
  \draw[edge,draw=PrimeOne] (0,0) -- (-0.563,1.358);
  \node[orbi] at (-0.563,1.358) {};
  \node[label,anchor=292.5] at (-0.609,1.469) {$E_{1}$};
  \draw[edge,draw=PrimeOne] (0,0) -- (-1.358,0.563);
  \node[orbi] at (-1.358,0.563) {};
  \node[label,anchor=337.5] at (-1.469,0.609) {$E_{2}$};
  \draw[edge,draw=PrimeOne] (0,0) -- (-1.358,-0.563);
  \node[orbi] at (-1.358,-0.563) {};
  \node[label,anchor=382.5] at (-1.469,-0.609) {$E_{3}$};
  \draw[edge,draw=PrimeOne] (0,0) -- (-0.563,-1.358);
  \node[orbi] at (-0.563,-1.358) {};
  \node[label,anchor=427.5] at (-0.609,-1.469) {$E_{4}$};
  \draw[edge,draw=PrimeTwo] (0,0) -- (0.563,-1.358);
  \node[orbi] at (0.563,-1.358) {};
  \node[label,anchor=472.5] at (0.609,-1.469) {$E_{5}$};
  \draw[edge,draw=PrimeTwo] (0,0) -- (1.358,-0.563);
  \node[orbi] at (1.358,-0.563) {};
  \node[label,anchor=517.5] at (1.469,-0.609) {$E_{6}$};
  \draw[edge,draw=PrimeTwo] (0,0) -- (1.358,0.563);
  \node[orbi] at (1.358,0.563) {};
  \node[label,anchor=562.5] at (1.469,0.609) {$E_{7}$};
  \draw[edge,draw=PrimeTwo] (0,0) -- (0.563,1.358);
  \node[orbi] at (0.563,1.358) {};
  \node[label,anchor=607.5] at (0.609,1.469) {$E_{8}$};
  \node[vert] at (0,0) {};
  \fill[cell] (2.688,0.55) rectangle (3.988,1.85);
  \draw[side,draw=PrimeOne] (2.688,0.55) -- (2.688,1.85);
  \node[orbi] at (2.688,1.2) {};
  \node[label,anchor=east] at (2.598,1.2) {$E_{4}$};
  \draw[side,draw=PrimeTwo] (2.688,1.85) -- (3.988,1.85);
  \node[orbi] at (3.338,1.85) {};
  \node[label,anchor=south] at (3.338,1.94) {$E_{5}$};
  \draw[side,draw=PrimeOne] (3.988,1.85) -- (3.988,0.55);
  \node[orbi] at (3.988,1.2) {};
  \node[label,anchor=west] at (4.078,1.2) {$E_{1}$};
  \draw[side,draw=PrimeTwo] (3.988,0.55) -- (2.688,0.55);
  \node[orbi] at (3.338,0.55) {};
  \node[label,anchor=north] at (3.338,0.46) {$E_{8}$};
  \node[corner] at (2.688,0.55) {};
  \node[corner] at (2.688,1.85) {};
  \node[corner] at (3.988,1.85) {};
  \node[corner] at (3.988,0.55) {};
  \node[inner sep=1.4pt,fill=SoftInk!5,rounded corners=1pt] at (3.338,1.2) {$F_{1}$};
  \fill[cell] (5.138,0.55) rectangle (6.438,1.85);
  \draw[mirror] (5.138,1.2) -- (6.438,1.2);
  \draw[mirror] (5.788,0.55) -- (5.788,1.85);
  \draw[side,draw=PrimeOne] (5.138,0.55) -- (5.138,1.85);
  \node[orbi] at (5.138,1.2) {};
  \node[label,anchor=east] at (5.048,1.2) {$E_{2}$};
  \draw[side,draw=PrimeTwo] (5.138,1.85) -- (6.438,1.85);
  \node[orbi] at (5.788,1.85) {};
  \node[label,anchor=south] at (5.788,1.94) {$E_{5}$};
  \draw[side,draw=PrimeOne] (6.438,1.85) -- (6.438,0.55);
  \node[orbi] at (6.438,1.2) {};
  \node[label,anchor=west] at (6.528,1.2) {$E_{2}$};
  \draw[side,draw=PrimeTwo] (6.438,0.55) -- (5.138,0.55);
  \node[orbi] at (5.788,0.55) {};
  \node[label,anchor=north] at (5.788,0.46) {$E_{5}$};
  \node[corner] at (5.138,0.55) {};
  \node[corner] at (5.138,1.85) {};
  \node[corner] at (6.438,1.85) {};
  \node[corner] at (6.438,0.55) {};
  \node[inner sep=1.4pt,fill=SoftInk!5,rounded corners=1pt] at (5.788,1.2) {$F_{2}$};
  \fill[cell] (7.588,0.55) rectangle (8.888,1.85);
  \draw[mirror] (8.238,0.55) -- (8.238,1.85);
  \draw[side,draw=PrimeOne] (7.588,0.55) -- (7.588,1.85);
  \node[orbi] at (7.588,1.2) {};
  \node[label,anchor=east] at (7.498,1.2) {$E_{1}$};
  \draw[side,draw=PrimeTwo] (7.588,1.85) -- (8.888,1.85);
  \node[orbi] at (8.238,1.85) {};
  \node[label,anchor=south] at (8.238,1.94) {$E_{7}$};
  \draw[side,draw=PrimeOne] (8.888,1.85) -- (8.888,0.55);
  \node[orbi] at (8.888,1.2) {};
  \node[label,anchor=west] at (8.978,1.2) {$E_{1}$};
  \draw[side,draw=PrimeTwo] (8.888,0.55) -- (7.588,0.55);
  \node[orbi] at (8.238,0.55) {};
  \node[label,anchor=north] at (8.238,0.46) {$E_{6}$};
  \node[corner] at (7.588,0.55) {};
  \node[corner] at (7.588,1.85) {};
  \node[corner] at (8.888,1.85) {};
  \node[corner] at (8.888,0.55) {};
  \node[inner sep=1.4pt,fill=SoftInk!5,rounded corners=1pt] at (8.238,1.2) {$F_{3}$};
  \fill[cell] (10.038,0.55) rectangle (11.338,1.85);
  \draw[mirror] (10.038,1.2) -- (11.338,1.2);
  \draw[mirror] (10.688,0.55) -- (10.688,1.85);
  \draw[side,draw=PrimeOne] (10.038,0.55) -- (10.038,1.85);
  \node[orbi] at (10.038,1.2) {};
  \node[label,anchor=east] at (9.948,1.2) {$E_{3}$};
  \draw[side,draw=PrimeTwo] (10.038,1.85) -- (11.338,1.85);
  \node[orbi] at (10.688,1.85) {};
  \node[label,anchor=south] at (10.688,1.94) {$E_{5}$};
  \draw[side,draw=PrimeOne] (11.338,1.85) -- (11.338,0.55);
  \node[orbi] at (11.338,1.2) {};
  \node[label,anchor=west] at (11.428,1.2) {$E_{3}$};
  \draw[side,draw=PrimeTwo] (11.338,0.55) -- (10.038,0.55);
  \node[orbi] at (10.688,0.55) {};
  \node[label,anchor=north] at (10.688,0.46) {$E_{5}$};
  \node[corner] at (10.038,0.55) {};
  \node[corner] at (10.038,1.85) {};
  \node[corner] at (11.338,1.85) {};
  \node[corner] at (11.338,0.55) {};
  \node[inner sep=1.4pt,fill=SoftInk!5,rounded corners=1pt] at (10.688,1.2) {$F_{4}$};
  \fill[cell] (2.688,-1.95) rectangle (3.988,-0.65);
  \draw[mirror] (2.688,-1.3) -- (3.988,-1.3);
  \draw[side,draw=PrimeOne] (2.688,-1.95) -- (2.688,-0.65);
  \node[orbi] at (2.688,-1.3) {};
  \node[label,anchor=east] at (2.598,-1.3) {$E_{4}$};
  \draw[side,draw=PrimeTwo] (2.688,-0.65) -- (3.988,-0.65);
  \node[orbi] at (3.338,-0.65) {};
  \node[label,anchor=south] at (3.338,-0.56) {$E_{7}$};
  \draw[side,draw=PrimeOne] (3.988,-0.65) -- (3.988,-1.95);
  \node[orbi] at (3.988,-1.3) {};
  \node[label,anchor=west] at (4.078,-1.3) {$E_{2}$};
  \draw[side,draw=PrimeTwo] (3.988,-1.95) -- (2.688,-1.95);
  \node[orbi] at (3.338,-1.95) {};
  \node[label,anchor=north] at (3.338,-2.04) {$E_{7}$};
  \node[corner] at (2.688,-1.95) {};
  \node[corner] at (2.688,-0.65) {};
  \node[corner] at (3.988,-0.65) {};
  \node[corner] at (3.988,-1.95) {};
  \node[inner sep=1.4pt,fill=SoftInk!5,rounded corners=1pt] at (3.338,-1.3) {$F_{5}$};
  \fill[cell] (5.138,-1.95) rectangle (6.438,-0.65);
  \draw[mirror] (5.788,-1.95) -- (5.788,-0.65);
  \draw[side,draw=PrimeOne] (5.138,-1.95) -- (5.138,-0.65);
  \node[orbi] at (5.138,-1.3) {};
  \node[label,anchor=east] at (5.048,-1.3) {$E_{2}$};
  \draw[side,draw=PrimeTwo] (5.138,-0.65) -- (6.438,-0.65);
  \node[orbi] at (5.788,-0.65) {};
  \node[label,anchor=south] at (5.788,-0.56) {$E_{6}$};
  \draw[side,draw=PrimeOne] (6.438,-0.65) -- (6.438,-1.95);
  \node[orbi] at (6.438,-1.3) {};
  \node[label,anchor=west] at (6.528,-1.3) {$E_{2}$};
  \draw[side,draw=PrimeTwo] (6.438,-1.95) -- (5.138,-1.95);
  \node[orbi] at (5.788,-1.95) {};
  \node[label,anchor=north] at (5.788,-2.04) {$E_{8}$};
  \node[corner] at (5.138,-1.95) {};
  \node[corner] at (5.138,-0.65) {};
  \node[corner] at (6.438,-0.65) {};
  \node[corner] at (6.438,-1.95) {};
  \node[inner sep=1.4pt,fill=SoftInk!5,rounded corners=1pt] at (5.788,-1.3) {$F_{6}$};
  \fill[cell] (7.588,-1.95) rectangle (8.888,-0.65);
  \draw[mirror] (8.238,-1.95) -- (8.238,-0.65);
  \draw[side,draw=PrimeOne] (7.588,-1.95) -- (7.588,-0.65);
  \node[orbi] at (7.588,-1.3) {};
  \node[label,anchor=east] at (7.498,-1.3) {$E_{3}$};
  \draw[side,draw=PrimeTwo] (7.588,-0.65) -- (8.888,-0.65);
  \node[orbi] at (8.238,-0.65) {};
  \node[label,anchor=south] at (8.238,-0.56) {$E_{7}$};
  \draw[side,draw=PrimeOne] (8.888,-0.65) -- (8.888,-1.95);
  \node[orbi] at (8.888,-1.3) {};
  \node[label,anchor=west] at (8.978,-1.3) {$E_{3}$};
  \draw[side,draw=PrimeTwo] (8.888,-1.95) -- (7.588,-1.95);
  \node[orbi] at (8.238,-1.95) {};
  \node[label,anchor=north] at (8.238,-2.04) {$E_{8}$};
  \node[corner] at (7.588,-1.95) {};
  \node[corner] at (7.588,-0.65) {};
  \node[corner] at (8.888,-0.65) {};
  \node[corner] at (8.888,-1.95) {};
  \node[inner sep=1.4pt,fill=SoftInk!5,rounded corners=1pt] at (8.238,-1.3) {$F_{7}$};
  \fill[cell] (10.038,-1.95) rectangle (11.338,-0.65);
  \draw[mirror] (10.038,-1.3) -- (11.338,-1.3);
  \draw[side,draw=PrimeOne] (10.038,-1.95) -- (10.038,-0.65);
  \node[orbi] at (10.038,-1.3) {};
  \node[label,anchor=east] at (9.948,-1.3) {$E_{4}$};
  \draw[side,draw=PrimeTwo] (10.038,-0.65) -- (11.338,-0.65);
  \node[orbi] at (10.688,-0.65) {};
  \node[label,anchor=south] at (10.688,-0.56) {$E_{6}$};
  \draw[side,draw=PrimeOne] (11.338,-0.65) -- (11.338,-1.95);
  \node[orbi] at (11.338,-1.3) {};
  \node[label,anchor=west] at (11.428,-1.3) {$E_{3}$};
  \draw[side,draw=PrimeTwo] (11.338,-1.95) -- (10.038,-1.95);
  \node[orbi] at (10.688,-1.95) {};
  \node[label,anchor=north] at (10.688,-2.04) {$E_{6}$};
  \node[corner] at (10.038,-1.95) {};
  \node[corner] at (10.038,-0.65) {};
  \node[corner] at (11.338,-0.65) {};
  \node[corner] at (11.338,-1.95) {};
  \node[inner sep=1.4pt,fill=SoftInk!5,rounded corners=1pt] at (10.688,-1.3) {$F_{8}$};
  \draw[edge,draw=PrimeOne] (2.3,-2.9) -- (2.85,-2.9);
  \node[label,anchor=west] at (2.95,-2.9) {$\fp_3$-adic edge};
  \draw[edge,draw=PrimeTwo] (5.3,-2.9) -- (5.85,-2.9);
  \node[label,anchor=west] at (5.95,-2.9) {$\fq_3$-adic edge};
\end{tikzpicture}
    \caption{The complex $X_{\Gamma_{13}}$ of Section \ref{sec:qsqrt7}.}
    \label{fig:qsqrt7-gamma13}
\end{figure}

\begin{figure}[htbp]
    \centering
    \definecolor{SoftInk}{HTML}{2F3437}
\definecolor{EdgeColour}{HTML}{4F7D72}
\definecolor{PrimeOne}{HTML}{4F7D72}
\definecolor{PrimeTwo}{HTML}{A85D67}
\begin{tikzpicture}[x=1cm,y=1cm,
  every node/.style={font=\scriptsize,text=SoftInk},
  edge/.style={line width=1.05pt,line cap=round},
  side/.style={line width=0.95pt,line cap=round},
  tip/.style={line width=0.9pt,line cap=round,line join=round},
  cell/.style={fill=SoftInk!5},
  mirror/.style={draw=SoftInk!75,line width=0.75pt,dash pattern=on 2.2pt off 1.6pt},
  label/.style={inner sep=1.2pt},
  vert/.style={circle,fill=SoftInk,inner sep=1.8pt},
  corner/.style={circle,fill=SoftInk,inner sep=1.2pt},
  orbi/.style={circle,fill=white,draw=SoftInk,line width=0.6pt,inner sep=1.35pt}]
  \draw[edge,draw=PrimeOne] (0,0) .. controls (-0.598,2.231) and (-2.231,0.598) .. (0,0);
  \draw[tip,draw=PrimeOne,shift={(-1.061,1.061)},rotate=-135] (-2.5pt,2.1pt) -- (0,0) -- (-2.5pt,-2.1pt);
  \node[label,anchor=315] at (-1.131,1.131) {$E_{1}$};
  \draw[edge,draw=PrimeOne] (0,0) -- (-1.358,-0.563);
  \node[orbi] at (-1.358,-0.563) {};
  \node[label,anchor=382.5] at (-1.469,-0.609) {$E_{2}$};
  \draw[edge,draw=PrimeOne] (0,0) -- (-0.563,-1.358);
  \node[orbi] at (-0.563,-1.358) {};
  \node[label,anchor=427.5] at (-0.609,-1.469) {$E_{3}$};
  \draw[edge,draw=PrimeTwo] (0,0) .. controls (0.598,-2.231) and (2.231,-0.598) .. (0,0);
  \draw[tip,draw=PrimeTwo,shift={(1.061,-1.061)},rotate=45] (-2.5pt,2.1pt) -- (0,0) -- (-2.5pt,-2.1pt);
  \node[label,anchor=495] at (1.131,-1.131) {$E_{4}$};
  \draw[edge,draw=PrimeTwo] (0,0) -- (1.358,0.563);
  \node[orbi] at (1.358,0.563) {};
  \node[label,anchor=562.5] at (1.469,0.609) {$E_{5}$};
  \draw[edge,draw=PrimeTwo] (0,0) -- (0.563,1.358);
  \node[orbi] at (0.563,1.358) {};
  \node[label,anchor=607.5] at (0.609,1.469) {$E_{6}$};
  \node[vert] at (0,0) {};
  \fill[cell] (2.688,0.55) rectangle (3.988,1.85);
  \draw[side,draw=PrimeOne] (2.688,0.55) -- (2.688,1.85);
  \node[orbi] at (2.688,1.2) {};
  \node[label,anchor=east] at (2.598,1.2) {$E_{3}$};
  \draw[side,draw=PrimeTwo] (2.688,1.85) -- (3.988,1.85);
  \node[orbi] at (3.338,1.85) {};
  \node[label,anchor=south] at (3.338,1.94) {$E_{5}$};
  \draw[side,draw=PrimeOne] (3.988,1.85) -- (3.988,0.55);
  \node[orbi] at (3.988,1.2) {};
  \node[label,anchor=west] at (4.078,1.2) {$E_{2}$};
  \draw[side,draw=PrimeTwo] (3.988,0.55) -- (2.688,0.55);
  \node[orbi] at (3.338,0.55) {};
  \node[label,anchor=north] at (3.338,0.46) {$E_{6}$};
  \node[corner] at (2.688,0.55) {};
  \node[corner] at (2.688,1.85) {};
  \node[corner] at (3.988,1.85) {};
  \node[corner] at (3.988,0.55) {};
  \node[inner sep=1.4pt,fill=SoftInk!5,rounded corners=1pt] at (3.338,1.2) {$F_{1}$};
  \fill[cell] (5.138,0.55) rectangle (6.438,1.85);
  \draw[mirror] (5.788,0.55) -- (5.788,1.85);
  \draw[side,draw=PrimeOne] (5.138,0.55) -- (5.138,1.85);
  \draw[tip,draw=PrimeOne,shift={(5.138,1.2)},rotate=90] (-2.5pt,2.1pt) -- (0,0) -- (-2.5pt,-2.1pt);
  \node[label,anchor=east] at (5.048,1.2) {$E_{1}$};
  \draw[side,draw=PrimeTwo] (5.138,1.85) -- (6.438,1.85);
  \node[orbi] at (5.788,1.85) {};
  \node[label,anchor=south] at (5.788,1.94) {$E_{5}$};
  \draw[side,draw=PrimeOne] (6.438,1.85) -- (6.438,0.55);
  \draw[tip,draw=PrimeOne,shift={(6.438,1.2)},rotate=90] (-2.5pt,2.1pt) -- (0,0) -- (-2.5pt,-2.1pt);
  \node[label,anchor=west] at (6.528,1.2) {$E_{1}$};
  \draw[side,draw=PrimeTwo] (6.438,0.55) -- (5.138,0.55);
  \node[orbi] at (5.788,0.55) {};
  \node[label,anchor=north] at (5.788,0.46) {$E_{5}$};
  \node[corner] at (5.138,0.55) {};
  \node[corner] at (5.138,1.85) {};
  \node[corner] at (6.438,1.85) {};
  \node[corner] at (6.438,0.55) {};
  \node[inner sep=1.4pt,fill=SoftInk!5,rounded corners=1pt] at (5.788,1.2) {$F_{2}$};
  \fill[cell] (7.588,0.55) rectangle (8.888,1.85);
  \draw[side,draw=PrimeOne] (7.588,0.55) -- (7.588,1.85);
  \node[orbi] at (7.588,1.2) {};
  \node[label,anchor=east] at (7.498,1.2) {$E_{2}$};
  \draw[side,draw=PrimeTwo] (7.588,1.85) -- (8.888,1.85);
  \draw[tip,draw=PrimeTwo,shift={(8.238,1.85)},rotate=0] (-2.5pt,2.1pt) -- (0,0) -- (-2.5pt,-2.1pt);
  \node[label,anchor=south] at (8.238,1.94) {$E_{4}$};
  \draw[side,draw=PrimeOne] (8.888,1.85) -- (8.888,0.55);
  \node[orbi] at (8.888,1.2) {};
  \node[label,anchor=west] at (8.978,1.2) {$E_{2}$};
  \draw[side,draw=PrimeTwo] (8.888,0.55) -- (7.588,0.55);
  \draw[tip,draw=PrimeTwo,shift={(8.238,0.55)},rotate=180] (-2.5pt,2.1pt) -- (0,0) -- (-2.5pt,-2.1pt);
  \node[label,anchor=north] at (8.238,0.46) {$E_{4}$};
  \node[corner] at (7.588,0.55) {};
  \node[corner] at (7.588,1.85) {};
  \node[corner] at (8.888,1.85) {};
  \node[corner] at (8.888,0.55) {};
  \node[orbi] at (8.238,1.2) {};
  \node[inner sep=1.4pt,fill=SoftInk!5,rounded corners=1pt] at (8.238,0.9) {$F_{3}$};
  \begin{scope}[shift={(1.225,0)}]
  \fill[cell] (2.688,-1.95) rectangle (3.988,-0.65);
  \draw[side,draw=PrimeOne] (2.688,-1.95) -- (2.688,-0.65);
  \node[orbi] at (2.688,-1.3) {};
  \node[label,anchor=east] at (2.598,-1.3) {$E_{3}$};
  \draw[side,draw=PrimeTwo] (2.688,-0.65) -- (3.988,-0.65);
  \draw[tip,draw=PrimeTwo,shift={(3.338,-0.65)},rotate=0] (-2.5pt,2.1pt) -- (0,0) -- (-2.5pt,-2.1pt);
  \node[label,anchor=south] at (3.338,-0.56) {$E_{4}$};
  \draw[side,draw=PrimeOne] (3.988,-0.65) -- (3.988,-1.95);
  \draw[tip,draw=PrimeOne,shift={(3.988,-1.3)},rotate=-90] (-2.5pt,2.1pt) -- (0,0) -- (-2.5pt,-2.1pt);
  \node[label,anchor=west] at (4.078,-1.3) {$E_{1}$};
  \draw[side,draw=PrimeTwo] (3.988,-1.95) -- (2.688,-1.95);
  \draw[tip,draw=PrimeTwo,shift={(3.338,-1.95)},rotate=180] (-2.5pt,2.1pt) -- (0,0) -- (-2.5pt,-2.1pt);
  \node[label,anchor=north] at (3.338,-2.04) {$E_{4}$};
  \node[corner] at (2.688,-1.95) {};
  \node[corner] at (2.688,-0.65) {};
  \node[corner] at (3.988,-0.65) {};
  \node[corner] at (3.988,-1.95) {};
  \node[inner sep=1.4pt,fill=SoftInk!5,rounded corners=1pt] at (3.338,-1.3) {$F_{4}$};
  \fill[cell] (5.138,-1.95) rectangle (6.438,-0.65);
  \draw[side,draw=PrimeOne] (5.138,-1.95) -- (5.138,-0.65);
  \draw[tip,draw=PrimeOne,shift={(5.138,-1.3)},rotate=-90] (-2.5pt,2.1pt) -- (0,0) -- (-2.5pt,-2.1pt);
  \node[label,anchor=east] at (5.048,-1.3) {$E_{1}$};
  \draw[side,draw=PrimeTwo] (5.138,-0.65) -- (6.438,-0.65);
  \draw[tip,draw=PrimeTwo,shift={(5.788,-0.65)},rotate=0] (-2.5pt,2.1pt) -- (0,0) -- (-2.5pt,-2.1pt);
  \node[label,anchor=south] at (5.788,-0.56) {$E_{4}$};
  \draw[side,draw=PrimeOne] (6.438,-0.65) -- (6.438,-1.95);
  \draw[tip,draw=PrimeOne,shift={(6.438,-1.3)},rotate=90] (-2.5pt,2.1pt) -- (0,0) -- (-2.5pt,-2.1pt);
  \node[label,anchor=west] at (6.528,-1.3) {$E_{1}$};
  \draw[side,draw=PrimeTwo] (6.438,-1.95) -- (5.138,-1.95);
  \node[orbi] at (5.788,-1.95) {};
  \node[label,anchor=north] at (5.788,-2.04) {$E_{6}$};
  \node[corner] at (5.138,-1.95) {};
  \node[corner] at (5.138,-0.65) {};
  \node[corner] at (6.438,-0.65) {};
  \node[corner] at (6.438,-1.95) {};
  \node[inner sep=1.4pt,fill=SoftInk!5,rounded corners=1pt] at (5.788,-1.3) {$F_{5}$};
  \end{scope}
  \draw[edge,draw=PrimeOne] (2.3,-2.9) -- (2.85,-2.9);
  \node[label,anchor=west] at (2.95,-2.9) {$\fp_3$-adic edge};
  \draw[edge,draw=PrimeTwo] (5.3,-2.9) -- (5.85,-2.9);
  \node[label,anchor=west] at (5.95,-2.9) {$\fq_3$-adic edge};
\end{tikzpicture}
    \caption{The complex $X_{\Gamma_{14}}$ of Section \ref{sec:qsqrt7}.}
    \label{fig:qsqrt7-gamma14}
\end{figure}

\begin{figure}[htbp]
    \centering
    \input{figures/qsqrt7_lambda_1}
    \caption{The square complex $X_{\Lambda_1}$ of Section \ref{sec:qsqrt7}.}
    \label{fig:qsqrt7-lambda1}
\end{figure}

\begin{figure}[htbp]
    \centering
    \input{figures/qsqrt7_lambda_2}
    \caption{The square complex $X_{\Lambda_2}$ of Section \ref{sec:qsqrt7}.}
    \label{fig:qsqrt7-lambda2}
\end{figure}

\begin{figure}[htbp]
    \centering
    \input{figures/qsqrt7_lambda_3}
    \caption{The square complex $X_{\Lambda_3}$ of Section \ref{sec:qsqrt7}.}
    \label{fig:qsqrt7-lambda3}
\end{figure}

\begin{figure}[htbp]
    \centering
    \input{figures/qsqrt7_lambda_4}
    \caption{The square complex $X_{\Lambda_4}$ of Section \ref{sec:qsqrt7}.}
    \label{fig:qsqrt7-lambda4}
\end{figure}

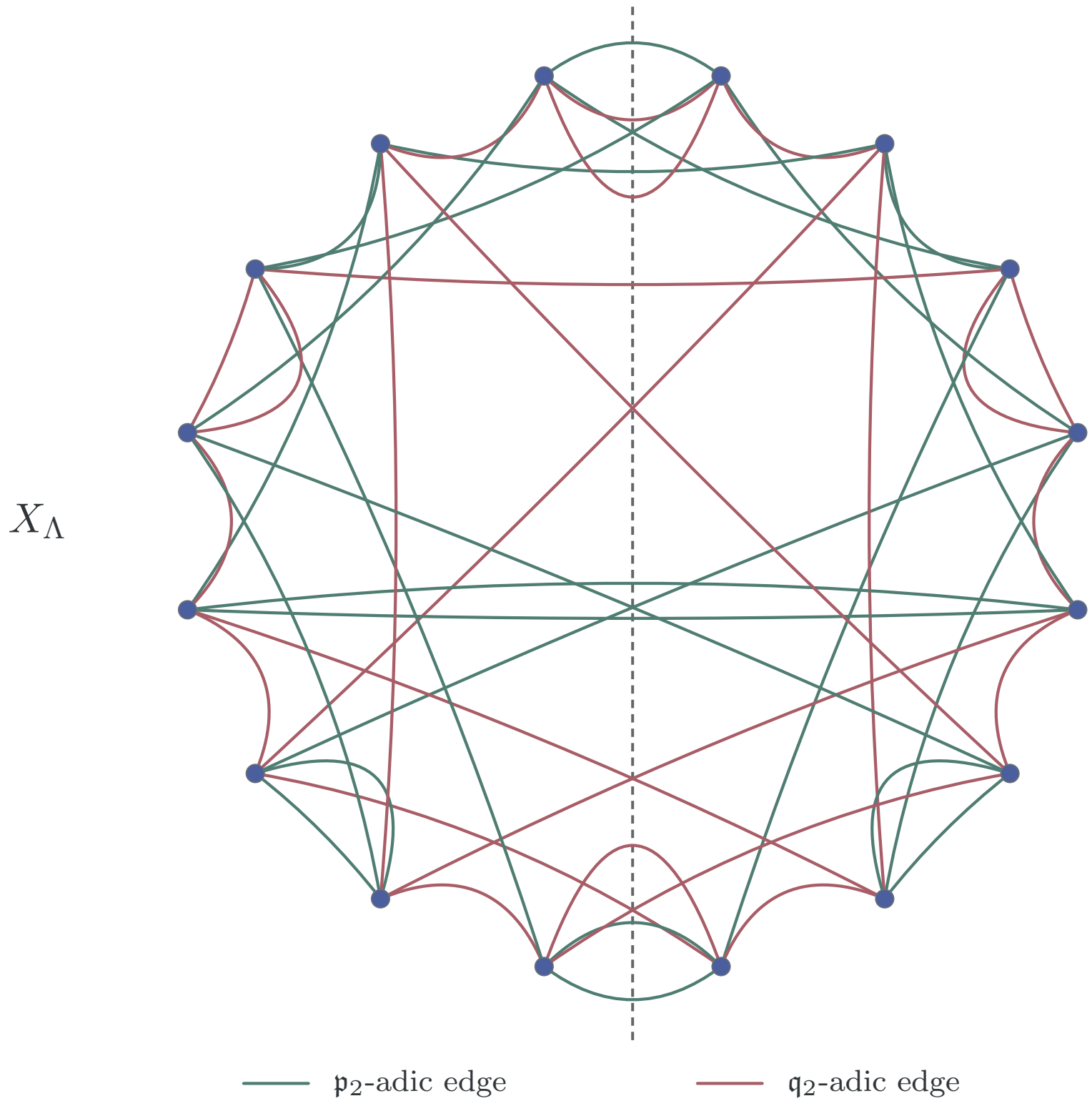
\begin{figure}[htbp]
    \centering
    \definecolor{SoftInk}{HTML}{2F3437}
\definecolor{PrimeOne}{HTML}{4F7D72}
\definecolor{PrimeTwo}{HTML}{A85D67}
\definecolor{VertexOne}{HTML}{4B5F9E}
\begin{tikzpicture}[x=1cm,y=1cm,
  every node/.style={font=\scriptsize,text=SoftInk},
  edge/.style={line width=0.8pt,line cap=round},
  vert/.style={circle,fill=VertexOne,draw=SoftInk!70,line width=0.3pt,inner sep=1.7pt},
  orbi/.style={circle,fill=white,draw=SoftInk,line width=0.6pt,inner sep=1.35pt},
  title/.style={font=\normalsize,text=SoftInk}]
  \node[title,anchor=east] at (-5.1,0) {$X_\Lambda$};
  \draw[SoftInk!75,line width=0.75pt,dash pattern=on 2.2pt off 1.6pt] (0,-4.80) -- (0,4.80);
  \draw[edge,draw=PrimeTwo] (4.119,0.819) .. controls (3.850,1.293) and (3.637,1.808) .. (3.492,2.333);
  \draw[edge,draw=PrimeTwo] (4.119,0.819) .. controls (2.972,0.929) and (2.758,1.444) .. (3.492,2.333);
  \draw[edge,draw=PrimeOne] (4.119,0.819) .. controls (1.565,-0.119) and (-1.023,-1.191) .. (-3.492,-2.333);
  \draw[edge,draw=PrimeOne] (4.119,0.819) .. controls (3.212,-0.472) and (2.605,-1.938) .. (2.333,-3.492);
  \draw[edge,draw=PrimeOne] (4.119,0.819) .. controls (2.787,1.665) and (1.665,2.787) .. (0.819,4.119);
  \draw[edge,draw=PrimeTwo] (4.119,0.819) .. controls (3.577,0.279) and (3.577,-0.279) .. (4.119,-0.819);
  \draw[edge,draw=PrimeOne] (3.492,2.333) .. controls (2.726,2.332) and (2.332,2.726) .. (2.333,3.492);
  \draw[edge,draw=PrimeOne] (3.492,2.333) .. controls (2.432,0.278) and (1.523,-1.916) .. (0.819,-4.119);
  \draw[edge,draw=PrimeOne] (3.492,2.333) .. controls (1.938,2.605) and (0.472,3.212) .. (-0.819,4.119);
  \draw[edge,draw=PrimeTwo] (3.492,2.333) .. controls (1.187,2.141) and (-1.187,2.141) .. (-3.492,2.333);
  \draw[edge,draw=PrimeTwo] (2.333,3.492) .. controls (0.489,1.492) and (-1.492,-0.489) .. (-3.492,-2.333);
  \draw[edge,draw=PrimeTwo] (2.333,3.492) .. controls (2.141,1.187) and (2.141,-1.187) .. (2.333,-3.492);
  \draw[edge,draw=PrimeTwo] (2.333,3.492) .. controls (1.626,3.198) and (1.111,3.411) .. (0.819,4.119);
  \draw[edge,draw=PrimeOne] (2.333,3.492) .. controls (2.605,1.938) and (3.212,0.472) .. (4.119,-0.819);
  \draw[edge,draw=PrimeOne] (2.333,3.492) .. controls (0.793,3.148) and (-0.793,3.148) .. (-2.333,3.492);
  \draw[edge,draw=PrimeTwo] (-3.492,-2.333) .. controls (-1.938,-2.605) and (-0.472,-3.212) .. (0.819,-4.119);
  \draw[edge,draw=PrimeTwo] (-3.492,-2.333) .. controls (-3.198,-1.626) and (-3.411,-1.111) .. (-4.119,-0.819);
  \draw[edge,draw=PrimeOne] (-3.492,-2.333) .. controls (-3.062,-2.668) and (-2.668,-3.062) .. (-2.333,-3.492);
  \draw[edge,draw=PrimeOne] (-3.492,-2.333) .. controls (-2.390,-1.996) and (-1.996,-2.390) .. (-2.333,-3.492);
  \draw[edge,draw=PrimeTwo] (2.333,-3.492) .. controls (1.626,-3.198) and (1.111,-3.411) .. (0.819,-4.119);
  \draw[edge,draw=PrimeTwo] (2.333,-3.492) .. controls (0.278,-2.432) and (-1.916,-1.523) .. (-4.119,-0.819);
  \draw[edge,draw=PrimeOne] (2.333,-3.492) .. controls (2.668,-3.062) and (3.062,-2.668) .. (3.492,-2.333);
  \draw[edge,draw=PrimeOne] (2.333,-3.492) .. controls (1.996,-2.390) and (2.390,-1.996) .. (3.492,-2.333);
  \draw[edge,draw=PrimeOne] (0.819,-4.119) .. controls (0.279,-4.528) and (-0.279,-4.528) .. (-0.819,-4.119);
  \draw[edge,draw=PrimeOne] (0.819,-4.119) .. controls (0.279,-3.577) and (-0.279,-3.577) .. (-0.819,-4.119);
  \draw[edge,draw=PrimeTwo] (0.819,-4.119) .. controls (0.279,-2.625) and (-0.279,-2.625) .. (-0.819,-4.119);
  \draw[edge,draw=PrimeOne] (0.819,4.119) .. controls (0.279,4.528) and (-0.279,4.528) .. (-0.819,4.119);
  \draw[edge,draw=PrimeTwo] (0.819,4.119) .. controls (0.279,3.577) and (-0.279,3.577) .. (-0.819,4.119);
  \draw[edge,draw=PrimeTwo] (0.819,4.119) .. controls (0.279,2.625) and (-0.279,2.625) .. (-0.819,4.119);
  \draw[edge,draw=PrimeOne] (0.819,4.119) .. controls (-0.472,3.212) and (-1.938,2.605) .. (-3.492,2.333);
  \draw[edge,draw=PrimeOne] (4.119,-0.819) .. controls (1.401,-0.925) and (-1.401,-0.925) .. (-4.119,-0.819);
  \draw[edge,draw=PrimeOne] (4.119,-0.819) .. controls (1.401,-0.492) and (-1.401,-0.492) .. (-4.119,-0.819);
  \draw[edge,draw=PrimeTwo] (4.119,-0.819) .. controls (1.916,-1.523) and (-0.278,-2.432) .. (-2.333,-3.492);
  \draw[edge,draw=PrimeTwo] (4.119,-0.819) .. controls (3.411,-1.111) and (3.198,-1.626) .. (3.492,-2.333);
  \draw[edge,draw=PrimeOne] (-4.119,-0.819) .. controls (-3.212,0.472) and (-2.605,1.938) .. (-2.333,3.492);
  \draw[edge,draw=PrimeTwo] (-4.119,-0.819) .. controls (-3.577,-0.279) and (-3.577,0.279) .. (-4.119,0.819);
  \draw[edge,draw=PrimeTwo] (-0.819,4.119) .. controls (-1.111,3.411) and (-1.626,3.198) .. (-2.333,3.492);
  \draw[edge,draw=PrimeOne] (-0.819,4.119) .. controls (-1.665,2.787) and (-2.787,1.665) .. (-4.119,0.819);
  \draw[edge,draw=PrimeTwo] (-0.819,-4.119) .. controls (-1.111,-3.411) and (-1.626,-3.198) .. (-2.333,-3.492);
  \draw[edge,draw=PrimeTwo] (-0.819,-4.119) .. controls (0.472,-3.212) and (1.938,-2.605) .. (3.492,-2.333);
  \draw[edge,draw=PrimeOne] (-0.819,-4.119) .. controls (-1.523,-1.916) and (-2.432,0.278) .. (-3.492,2.333);
  \draw[edge,draw=PrimeTwo] (-2.333,-3.492) .. controls (-2.141,-1.187) and (-2.141,1.187) .. (-2.333,3.492);
  \draw[edge,draw=PrimeOne] (-2.333,-3.492) .. controls (-2.605,-1.938) and (-3.212,-0.472) .. (-4.119,0.819);
  \draw[edge,draw=PrimeTwo] (3.492,-2.333) .. controls (1.492,-0.489) and (-0.489,1.492) .. (-2.333,3.492);
  \draw[edge,draw=PrimeOne] (3.492,-2.333) .. controls (1.023,-1.191) and (-1.565,-0.119) .. (-4.119,0.819);
  \draw[edge,draw=PrimeOne] (-2.333,3.492) .. controls (-2.332,2.726) and (-2.726,2.332) .. (-3.492,2.333);
  \draw[edge,draw=PrimeTwo] (-3.492,2.333) .. controls (-3.637,1.808) and (-3.850,1.293) .. (-4.119,0.819);
  \draw[edge,draw=PrimeTwo] (-3.492,2.333) .. controls (-2.758,1.444) and (-2.972,0.929) .. (-4.119,0.819);
  \node[vert] at (0.819,4.119) {};
  \node[vert] at (-0.819,4.119) {};
  \node[vert] at (2.333,3.492) {};
  \node[vert] at (-2.333,3.492) {};
  \node[vert] at (3.492,2.333) {};
  \node[vert] at (-3.492,2.333) {};
  \node[vert] at (4.119,0.819) {};
  \node[vert] at (-4.119,0.819) {};
  \node[vert] at (4.119,-0.819) {};
  \node[vert] at (-4.119,-0.819) {};
  \node[vert] at (3.492,-2.333) {};
  \node[vert] at (-3.492,-2.333) {};
  \node[vert] at (2.333,-3.492) {};
  \node[vert] at (-2.333,-3.492) {};
  \node[vert] at (0.819,-4.119) {};
  \node[vert] at (-0.819,-4.119) {};
  \draw[edge,draw=PrimeOne] (-3.6,-5.2) -- (-3.0,-5.2);
  \node[anchor=west] at (-2.9,-5.2) {$\fp_2$-adic edge};
  \draw[edge,draw=PrimeTwo] (0.6,-5.2) -- (1.2,-5.2);
  \node[anchor=west] at (1.3,-5.2) {$\fq_2$-adic edge};
\end{tikzpicture}
    \caption{The one-skeleton of $X_\Lambda$ of Section \ref{sec:qsqrt17}. The covering $X_\Lambda \to X_\Gamma$ is the quotient by the involution given by reflection in the dashed line. Some edges are sent to loops and some are inverted.}
    \label{fig:qsqrt17-covering}
\end{figure}

\begin{figure}[htbp]
    \centering
    \input{figures/qsqrt6_lambda_1}
    \caption{The square complex $X_{\Lambda_1}$ of Section \ref{sec:qsqrt6}.}
    \label{fig:qsqrt6-lambda1}
\end{figure}

\begin{figure}[htbp]
    \centering
    \input{figures/qsqrt6_lambda_2}
    \caption{The square complex $X_{\Lambda_2}$ of Section \ref{sec:qsqrt6}.}
    \label{fig:qsqrt6-lambda2}
\end{figure}

\clearpage

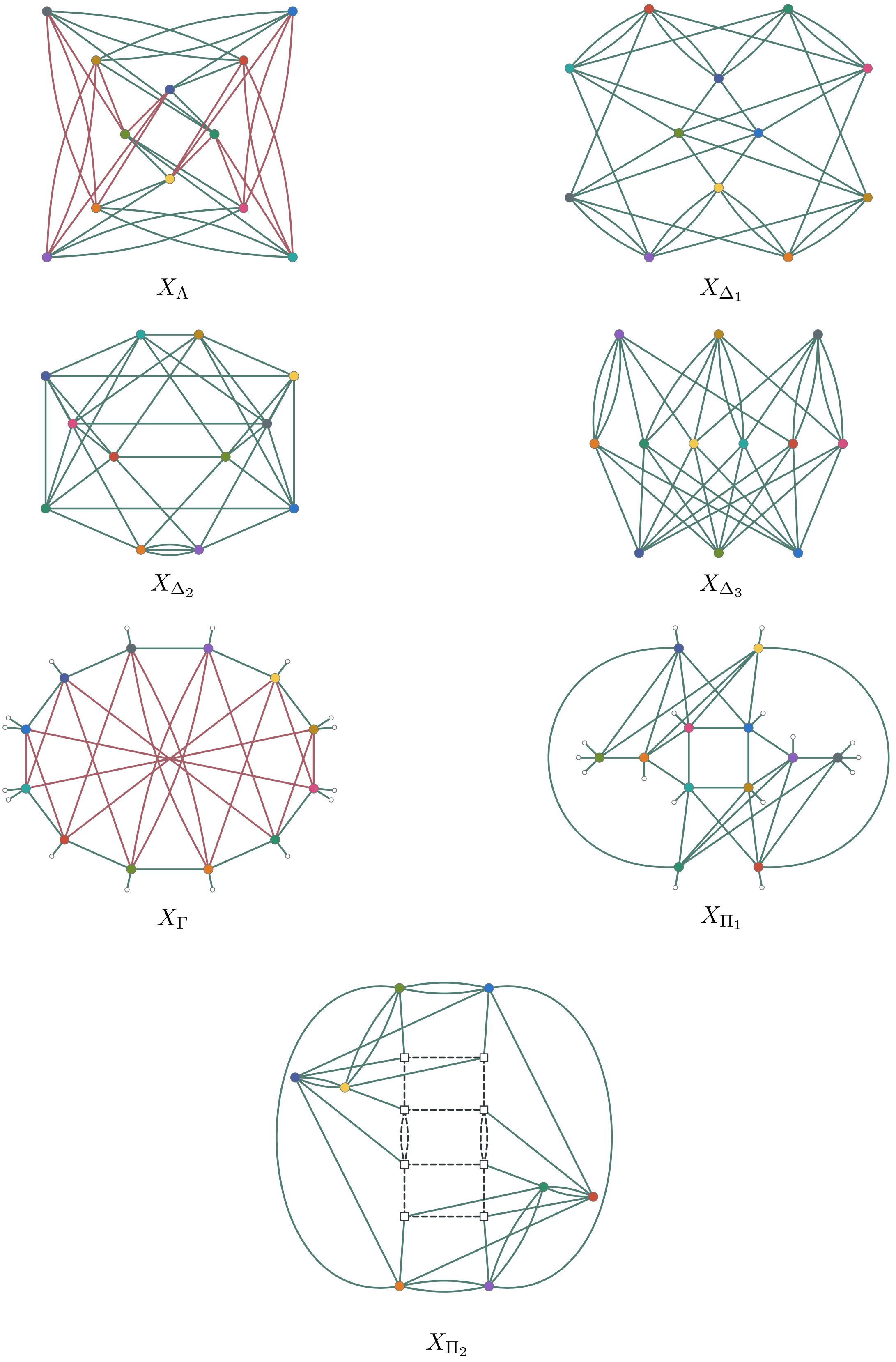
\begin{figure}[htbp]
    \centering
    \begin{minipage}{.48\textwidth}
        \centering
        \scalebox{0.8}{\definecolor{SoftInk}{HTML}{2F3437}%
\definecolor{TwoAdic}{HTML}{4F7D72}%
\definecolor{ThreeAdic}{HTML}{A85D67}%
\definecolor{ChiThreeVOne}{HTML}{4B5F9E}%
\definecolor{ChiThreeVTwo}{HTML}{F2C94C}%
\definecolor{ChiThreeVThree}{HTML}{2F8F6A}%
\definecolor{ChiThreeVFour}{HTML}{C64E3B}%
\definecolor{ChiThreeVFive}{HTML}{6D8F2F}%
\definecolor{ChiThreeVSix}{HTML}{2F76C8}%
\definecolor{ChiThreeVSeven}{HTML}{E07B28}%
\definecolor{ChiThreeVEight}{HTML}{8B5FBF}%
\definecolor{ChiThreeVNine}{HTML}{D84F82}%
\definecolor{ChiThreeVTen}{HTML}{2BA7A0}%
\definecolor{ChiThreeVEleven}{HTML}{B88920}%
\definecolor{ChiThreeVTwelve}{HTML}{5F6B73}%
\begin{tikzpicture}[x=1cm,y=1cm,
  vertex/.style={circle,draw=SoftInk!70,line width=0.3pt,inner sep=1.9pt},
  edge/.style={line width=1.05pt,line cap=round}]
  \coordinate (v1) at (0.0000,0.9000);
  \coordinate (v3) at (0.9000,0.0000);
  \coordinate (v2) at (0.0000,-0.9000);
  \coordinate (v5) at (-0.9000,-0.0000);
  \coordinate (v4) at (1.4849,1.4849);
  \coordinate (v9) at (1.4849,-1.4849);
  \coordinate (v7) at (-1.4849,-1.4849);
  \coordinate (v11) at (-1.4849,1.4849);
  \coordinate (v6) at (2.4749,2.4749);
  \coordinate (v10) at (2.4749,-2.4749);
  \coordinate (v8) at (-2.4749,-2.4749);
  \coordinate (v12) at (-2.4749,2.4749);
  \draw[edge,draw=TwoAdic] (v1) -- (v3);
  \draw[edge,draw=TwoAdic] (v1) -- (v4);
  \draw[edge,draw=TwoAdic] (v1) -- (v6);
  \draw[edge,draw=TwoAdic] (v2) -- (v5);
  \draw[edge,draw=TwoAdic] (v2) -- (v7);
  \draw[edge,draw=TwoAdic] (v2) -- (v8);
  \draw[edge,draw=TwoAdic] (v3) -- (v11);
  \draw[edge,draw=TwoAdic] (v3) -- (v12);
  \draw[edge,draw=TwoAdic] (v4) to[bend left=12] (v12);
  \draw[edge,draw=TwoAdic] (v4) to[bend right=12] (v12);
  \draw[edge,draw=TwoAdic] (v5) -- (v9);
  \draw[edge,draw=TwoAdic] (v5) -- (v10);
  \draw[edge,draw=TwoAdic] (v6) to[bend left=12] (v11);
  \draw[edge,draw=TwoAdic] (v6) to[bend right=12] (v11);
  \draw[edge,draw=TwoAdic] (v7) to[bend left=12] (v10);
  \draw[edge,draw=TwoAdic] (v7) to[bend right=12] (v10);
  \draw[edge,draw=TwoAdic] (v8) to[bend left=12] (v9);
  \draw[edge,draw=TwoAdic] (v8) to[bend right=12] (v9);
  \draw[edge,draw=ThreeAdic] (v1) -- (v5);
  \draw[edge,draw=ThreeAdic] (v1) -- (v7);
  \draw[edge,draw=ThreeAdic] (v1) -- (v8);
  \draw[edge,draw=ThreeAdic] (v2) -- (v3);
  \draw[edge,draw=ThreeAdic] (v2) -- (v4);
  \draw[edge,draw=ThreeAdic] (v2) -- (v6);
  \draw[edge,draw=ThreeAdic] (v3) -- (v9);
  \draw[edge,draw=ThreeAdic] (v3) -- (v10);
  \draw[edge,draw=ThreeAdic] (v4) to[bend left=12] (v10);
  \draw[edge,draw=ThreeAdic] (v4) to[bend right=12] (v10);
  \draw[edge,draw=ThreeAdic] (v5) -- (v11);
  \draw[edge,draw=ThreeAdic] (v5) -- (v12);
  \draw[edge,draw=ThreeAdic] (v6) to[bend left=12] (v9);
  \draw[edge,draw=ThreeAdic] (v6) to[bend right=12] (v9);
  \draw[edge,draw=ThreeAdic] (v7) to[bend left=12] (v12);
  \draw[edge,draw=ThreeAdic] (v7) to[bend right=12] (v12);
  \draw[edge,draw=ThreeAdic] (v8) to[bend left=12] (v11);
  \draw[edge,draw=ThreeAdic] (v8) to[bend right=12] (v11);
  \node[vertex,fill=ChiThreeVOne] at (v1) {};
  \node[vertex,fill=ChiThreeVTwo] at (v2) {};
  \node[vertex,fill=ChiThreeVThree] at (v3) {};
  \node[vertex,fill=ChiThreeVFour] at (v4) {};
  \node[vertex,fill=ChiThreeVFive] at (v5) {};
  \node[vertex,fill=ChiThreeVSix] at (v6) {};
  \node[vertex,fill=ChiThreeVSeven] at (v7) {};
  \node[vertex,fill=ChiThreeVEight] at (v8) {};
  \node[vertex,fill=ChiThreeVNine] at (v9) {};
  \node[vertex,fill=ChiThreeVTen] at (v10) {};
  \node[vertex,fill=ChiThreeVEleven] at (v11) {};
  \node[vertex,fill=ChiThreeVTwelve] at (v12) {};
\end{tikzpicture}}\\[2pt]
        $X_{\Lambda}$
    \end{minipage}\hfill
    \begin{minipage}{.48\textwidth}
        \centering
        \scalebox{0.8}{\definecolor{SoftInk}{HTML}{2F3437}%
\definecolor{TwoAdic}{HTML}{4F7D72}%
\definecolor{ThreeAdic}{HTML}{A85D67}%
\definecolor{ChiThreeVOne}{HTML}{4B5F9E}%
\definecolor{ChiThreeVTwo}{HTML}{F2C94C}%
\definecolor{ChiThreeVThree}{HTML}{2F8F6A}%
\definecolor{ChiThreeVFour}{HTML}{C64E3B}%
\definecolor{ChiThreeVFive}{HTML}{6D8F2F}%
\definecolor{ChiThreeVSix}{HTML}{2F76C8}%
\definecolor{ChiThreeVSeven}{HTML}{E07B28}%
\definecolor{ChiThreeVEight}{HTML}{8B5FBF}%
\definecolor{ChiThreeVNine}{HTML}{D84F82}%
\definecolor{ChiThreeVTen}{HTML}{2BA7A0}%
\definecolor{ChiThreeVEleven}{HTML}{B88920}%
\definecolor{ChiThreeVTwelve}{HTML}{5F6B73}%
\begin{tikzpicture}[x=1cm,y=1cm,
  vertex/.style={circle,draw=SoftInk!70,line width=0.3pt,inner sep=1.9pt},
  edge/.style={line width=1.05pt,line cap=round}]
  \coordinate (v1) at (0,1.1);
  \coordinate (v2) at (0,-1.1);
  \coordinate (v3) at (1.4,2.5);
  \coordinate (v4) at (-1.4,2.5);
  \coordinate (v7) at (1.4,-2.5);
  \coordinate (v8) at (-1.4,-2.5);
  \coordinate (v9) at (3,1.3);
  \coordinate (v10) at (-3,1.3);
  \coordinate (v11) at (3,-1.3);
  \coordinate (v12) at (-3,-1.3);
  \coordinate (v5) at (-0.8,0);
  \coordinate (v6) at (0.8,0);
  \draw[edge,draw=TwoAdic] (v1) to[bend left=12] (v3);
  \draw[edge,draw=TwoAdic] (v1) to[bend right=12] (v3);
  \draw[edge,draw=TwoAdic] (v1) to[bend left=12] (v4);
  \draw[edge,draw=TwoAdic] (v1) to[bend right=12] (v4);
  \draw[edge,draw=TwoAdic] (v1) -- (v5);
  \draw[edge,draw=TwoAdic] (v1) -- (v6);
  \draw[edge,draw=TwoAdic] (v2) -- (v5);
  \draw[edge,draw=TwoAdic] (v2) -- (v6);
  \draw[edge,draw=TwoAdic] (v2) to[bend left=12] (v7);
  \draw[edge,draw=TwoAdic] (v2) to[bend right=12] (v7);
  \draw[edge,draw=TwoAdic] (v2) to[bend left=12] (v8);
  \draw[edge,draw=TwoAdic] (v2) to[bend right=12] (v8);
  \draw[edge,draw=TwoAdic] (v3) to[bend left=12] (v9);
  \draw[edge,draw=TwoAdic] (v3) to[bend right=12] (v9);
  \draw[edge,draw=TwoAdic] (v3) -- (v10);
  \draw[edge,draw=TwoAdic] (v3) -- (v11);
  \draw[edge,draw=TwoAdic] (v4) -- (v9);
  \draw[edge,draw=TwoAdic] (v4) to[bend left=12] (v10);
  \draw[edge,draw=TwoAdic] (v4) to[bend right=12] (v10);
  \draw[edge,draw=TwoAdic] (v4) -- (v12);
  \draw[edge,draw=TwoAdic] (v5) -- (v9);
  \draw[edge,draw=TwoAdic] (v5) -- (v10);
  \draw[edge,draw=TwoAdic] (v5) -- (v11);
  \draw[edge,draw=TwoAdic] (v5) -- (v12);
  \draw[edge,draw=TwoAdic] (v6) -- (v9);
  \draw[edge,draw=TwoAdic] (v6) -- (v10);
  \draw[edge,draw=TwoAdic] (v6) -- (v11);
  \draw[edge,draw=TwoAdic] (v6) -- (v12);
  \draw[edge,draw=TwoAdic] (v7) -- (v9);
  \draw[edge,draw=TwoAdic] (v7) to[bend left=12] (v11);
  \draw[edge,draw=TwoAdic] (v7) to[bend right=12] (v11);
  \draw[edge,draw=TwoAdic] (v7) -- (v12);
  \draw[edge,draw=TwoAdic] (v8) -- (v10);
  \draw[edge,draw=TwoAdic] (v8) -- (v11);
  \draw[edge,draw=TwoAdic] (v8) to[bend left=12] (v12);
  \draw[edge,draw=TwoAdic] (v8) to[bend right=12] (v12);
  \node[vertex,fill=ChiThreeVOne] at (v1) {};
  \node[vertex,fill=ChiThreeVTwo] at (v2) {};
  \node[vertex,fill=ChiThreeVThree] at (v3) {};
  \node[vertex,fill=ChiThreeVFour] at (v4) {};
  \node[vertex,fill=ChiThreeVFive] at (v5) {};
  \node[vertex,fill=ChiThreeVSix] at (v6) {};
  \node[vertex,fill=ChiThreeVSeven] at (v7) {};
  \node[vertex,fill=ChiThreeVEight] at (v8) {};
  \node[vertex,fill=ChiThreeVNine] at (v9) {};
  \node[vertex,fill=ChiThreeVTen] at (v10) {};
  \node[vertex,fill=ChiThreeVEleven] at (v11) {};
  \node[vertex,fill=ChiThreeVTwelve] at (v12) {};
\end{tikzpicture}}\\[2pt]
        $X_{\Delta_1}$
    \end{minipage}\\[12pt]
    \begin{minipage}{.48\textwidth}
        \centering
        \scalebox{0.8}{\definecolor{SoftInk}{HTML}{2F3437}%
\definecolor{TwoAdic}{HTML}{4F7D72}%
\definecolor{ChiThreeVOne}{HTML}{4B5F9E}%
\definecolor{ChiThreeVTwo}{HTML}{F2C94C}%
\definecolor{ChiThreeVThree}{HTML}{2F8F6A}%
\definecolor{ChiThreeVFour}{HTML}{C64E3B}%
\definecolor{ChiThreeVFive}{HTML}{6D8F2F}%
\definecolor{ChiThreeVSix}{HTML}{2F76C8}%
\definecolor{ChiThreeVSeven}{HTML}{E07B28}%
\definecolor{ChiThreeVEight}{HTML}{8B5FBF}%
\definecolor{ChiThreeVNine}{HTML}{D84F82}%
\definecolor{ChiThreeVTen}{HTML}{2BA7A0}%
\definecolor{ChiThreeVEleven}{HTML}{B88920}%
\definecolor{ChiThreeVTwelve}{HTML}{5F6B73}%
\begin{tikzpicture}[x=1cm,y=1cm,
  vertex/.style={circle,draw=SoftInk!70,line width=0.3pt,inner sep=1.9pt},
  edge/.style={line width=1.05pt,line cap=round}]
  \coordinate (v1) at (-2.5000,1.3333);
  \coordinate (v2) at (2.5000,1.3333);
  \coordinate (v3) at (-2.5000,-1.3333);
  \coordinate (v4) at (-1.1250,-0.2917);
  \coordinate (v5) at (1.1250,-0.2917);
  \coordinate (v6) at (2.5000,-1.3333);
  \coordinate (v7) at (-0.5833,-2.1667);
  \coordinate (v8) at (0.5833,-2.1667);
  \coordinate (v9) at (-1.9583,0.3750);
  \coordinate (v10) at (-0.5833,2.1667);
  \coordinate (v11) at (0.5833,2.1667);
  \coordinate (v12) at (1.9583,0.3750);
  \draw[edge,draw=TwoAdic] (v1) -- (v2);
  \draw[edge,draw=TwoAdic] (v1) -- (v3);
  \draw[edge,draw=TwoAdic] (v1) -- (v4);
  \draw[edge,draw=TwoAdic] (v1) -- (v7);
  \draw[edge,draw=TwoAdic] (v1) -- (v9);
  \draw[edge,draw=TwoAdic] (v1) -- (v10);
  \draw[edge,draw=TwoAdic] (v2) -- (v5);
  \draw[edge,draw=TwoAdic] (v2) -- (v6);
  \draw[edge,draw=TwoAdic] (v2) -- (v8);
  \draw[edge,draw=TwoAdic] (v2) -- (v11);
  \draw[edge,draw=TwoAdic] (v2) -- (v12);
  \draw[edge,draw=TwoAdic] (v3) -- (v4);
  \draw[edge,draw=TwoAdic] (v3) -- (v6);
  \draw[edge,draw=TwoAdic] (v3) -- (v7);
  \draw[edge,draw=TwoAdic] (v3) -- (v9);
  \draw[edge,draw=TwoAdic] (v3) -- (v10);
  \draw[edge,draw=TwoAdic] (v4) -- (v5);
  \draw[edge,draw=TwoAdic] (v4) -- (v8);
  \draw[edge,draw=TwoAdic] (v4) -- (v9);
  \draw[edge,draw=TwoAdic] (v4) -- (v11);
  \draw[edge,draw=TwoAdic] (v5) -- (v6);
  \draw[edge,draw=TwoAdic] (v5) -- (v7);
  \draw[edge,draw=TwoAdic] (v5) -- (v10);
  \draw[edge,draw=TwoAdic] (v5) -- (v12);
  \draw[edge,draw=TwoAdic] (v6) -- (v8);
  \draw[edge,draw=TwoAdic] (v6) -- (v11);
  \draw[edge,draw=TwoAdic] (v6) -- (v12);
  \draw[edge,draw=TwoAdic] (v7) to[bend left=18] (v8);
  \draw[edge,draw=TwoAdic] (v7) -- (v8);
  \draw[edge,draw=TwoAdic] (v7) to[bend right=18] (v8);
  \draw[edge,draw=TwoAdic] (v9) -- (v10);
  \draw[edge,draw=TwoAdic] (v9) -- (v11);
  \draw[edge,draw=TwoAdic] (v9) -- (v12);
  \draw[edge,draw=TwoAdic] (v10) -- (v11);
  \draw[edge,draw=TwoAdic] (v10) -- (v12);
  \draw[edge,draw=TwoAdic] (v11) -- (v12);
  \node[vertex,fill=ChiThreeVOne] at (v1) {};
  \node[vertex,fill=ChiThreeVTwo] at (v2) {};
  \node[vertex,fill=ChiThreeVThree] at (v3) {};
  \node[vertex,fill=ChiThreeVFour] at (v4) {};
  \node[vertex,fill=ChiThreeVFive] at (v5) {};
  \node[vertex,fill=ChiThreeVSix] at (v6) {};
  \node[vertex,fill=ChiThreeVSeven] at (v7) {};
  \node[vertex,fill=ChiThreeVEight] at (v8) {};
  \node[vertex,fill=ChiThreeVNine] at (v9) {};
  \node[vertex,fill=ChiThreeVTen] at (v10) {};
  \node[vertex,fill=ChiThreeVEleven] at (v11) {};
  \node[vertex,fill=ChiThreeVTwelve] at (v12) {};
\end{tikzpicture}}\\[2pt]
        $X_{\Delta_2}$
    \end{minipage}\hfill
    \begin{minipage}{.48\textwidth}
        \centering
        \scalebox{0.8}{\definecolor{SoftInk}{HTML}{2F3437}%
\definecolor{TwoAdic}{HTML}{4F7D72}%
\definecolor{ChiThreeVOne}{HTML}{4B5F9E}%
\definecolor{ChiThreeVTwo}{HTML}{F2C94C}%
\definecolor{ChiThreeVThree}{HTML}{2F8F6A}%
\definecolor{ChiThreeVFour}{HTML}{C64E3B}%
\definecolor{ChiThreeVFive}{HTML}{6D8F2F}%
\definecolor{ChiThreeVSix}{HTML}{2F76C8}%
\definecolor{ChiThreeVSeven}{HTML}{E07B28}%
\definecolor{ChiThreeVEight}{HTML}{8B5FBF}%
\definecolor{ChiThreeVNine}{HTML}{D84F82}%
\definecolor{ChiThreeVTen}{HTML}{2BA7A0}%
\definecolor{ChiThreeVEleven}{HTML}{B88920}%
\definecolor{ChiThreeVTwelve}{HTML}{5F6B73}%
\begin{tikzpicture}[x=1cm,y=1cm,
  vertex/.style={circle,draw=SoftInk!70,line width=0.3pt,inner sep=1.9pt},
  edge/.style={line width=1.05pt,line cap=round}]
  \coordinate (v1) at (-1.6000,-2.2000);
  \coordinate (v2) at (-0.5000,0.0000);
  \coordinate (v3) at (-1.5000,0.0000);
  \coordinate (v4) at (1.5000,0.0000);
  \coordinate (v5) at (0.0000,-2.2000);
  \coordinate (v6) at (1.6000,-2.2000);
  \coordinate (v7) at (-2.5000,0.0000);
  \coordinate (v8) at (-2.0000,2.2000);
  \coordinate (v9) at (2.5000,0.0000);
  \coordinate (v10) at (0.5000,0.0000);
  \coordinate (v11) at (0.0000,2.2000);
  \coordinate (v12) at (2.0000,2.2000);
  \draw[edge,draw=TwoAdic] (v1) -- (v2);
  \draw[edge,draw=TwoAdic] (v1) -- (v3);
  \draw[edge,draw=TwoAdic] (v1) -- (v4);
  \draw[edge,draw=TwoAdic] (v1) -- (v7);
  \draw[edge,draw=TwoAdic] (v1) -- (v9);
  \draw[edge,draw=TwoAdic] (v1) -- (v10);
  \draw[edge,draw=TwoAdic] (v2) -- (v5);
  \draw[edge,draw=TwoAdic] (v2) -- (v6);
  \draw[edge,draw=TwoAdic] (v2) -- (v8);
  \draw[edge,draw=TwoAdic] (v2) -- (v11);
  \draw[edge,draw=TwoAdic] (v2) -- (v12);
  \draw[edge,draw=TwoAdic] (v3) -- (v5);
  \draw[edge,draw=TwoAdic] (v3) -- (v6);
  \draw[edge,draw=TwoAdic] (v3) -- (v8);
  \draw[edge,draw=TwoAdic] (v3) to[bend left=12] (v11);
  \draw[edge,draw=TwoAdic] (v3) to[bend right=12] (v11);
  \draw[edge,draw=TwoAdic] (v4) -- (v5);
  \draw[edge,draw=TwoAdic] (v4) -- (v6);
  \draw[edge,draw=TwoAdic] (v4) -- (v8);
  \draw[edge,draw=TwoAdic] (v4) to[bend left=12] (v12);
  \draw[edge,draw=TwoAdic] (v4) to[bend right=12] (v12);
  \draw[edge,draw=TwoAdic] (v5) -- (v7);
  \draw[edge,draw=TwoAdic] (v5) -- (v9);
  \draw[edge,draw=TwoAdic] (v5) -- (v10);
  \draw[edge,draw=TwoAdic] (v6) -- (v7);
  \draw[edge,draw=TwoAdic] (v6) -- (v9);
  \draw[edge,draw=TwoAdic] (v6) -- (v10);
  \draw[edge,draw=TwoAdic] (v7) to[bend left=18] (v8);
  \draw[edge,draw=TwoAdic] (v7) -- (v8);
  \draw[edge,draw=TwoAdic] (v7) to[bend right=18] (v8);
  \draw[edge,draw=TwoAdic] (v9) -- (v11);
  \draw[edge,draw=TwoAdic] (v9) to[bend left=12] (v12);
  \draw[edge,draw=TwoAdic] (v9) to[bend right=12] (v12);
  \draw[edge,draw=TwoAdic] (v10) to[bend left=12] (v11);
  \draw[edge,draw=TwoAdic] (v10) to[bend right=12] (v11);
  \draw[edge,draw=TwoAdic] (v10) -- (v12);
  \node[vertex,fill=ChiThreeVOne] at (v1) {};
  \node[vertex,fill=ChiThreeVTwo] at (v2) {};
  \node[vertex,fill=ChiThreeVThree] at (v3) {};
  \node[vertex,fill=ChiThreeVFour] at (v4) {};
  \node[vertex,fill=ChiThreeVFive] at (v5) {};
  \node[vertex,fill=ChiThreeVSix] at (v6) {};
  \node[vertex,fill=ChiThreeVSeven] at (v7) {};
  \node[vertex,fill=ChiThreeVEight] at (v8) {};
  \node[vertex,fill=ChiThreeVNine] at (v9) {};
  \node[vertex,fill=ChiThreeVTen] at (v10) {};
  \node[vertex,fill=ChiThreeVEleven] at (v11) {};
  \node[vertex,fill=ChiThreeVTwelve] at (v12) {};
\end{tikzpicture}}\\[2pt]
        $X_{\Delta_3}$
    \end{minipage}\\[12pt]
    \begin{minipage}{.48\textwidth}
        \centering
        \scalebox{0.8}{\definecolor{SoftInk}{HTML}{2F3437}%
\definecolor{TwoAdic}{HTML}{4F7D72}%
\definecolor{ThreeAdic}{HTML}{A85D67}%
\definecolor{ChiThreeVOne}{HTML}{4B5F9E}%
\definecolor{ChiThreeVTwo}{HTML}{F2C94C}%
\definecolor{ChiThreeVThree}{HTML}{2F8F6A}%
\definecolor{ChiThreeVFour}{HTML}{C64E3B}%
\definecolor{ChiThreeVFive}{HTML}{6D8F2F}%
\definecolor{ChiThreeVSix}{HTML}{2F76C8}%
\definecolor{ChiThreeVSeven}{HTML}{E07B28}%
\definecolor{ChiThreeVEight}{HTML}{8B5FBF}%
\definecolor{ChiThreeVNine}{HTML}{D84F82}%
\definecolor{ChiThreeVTen}{HTML}{2BA7A0}%
\definecolor{ChiThreeVEleven}{HTML}{B88920}%
\definecolor{ChiThreeVTwelve}{HTML}{5F6B73}%
\begin{tikzpicture}[x=1cm,y=1cm,
  vertex/.style={circle,draw=SoftInk!70,line width=0.3pt,inner sep=1.9pt},
  fixed/.style={rectangle,draw=SoftInk,fill=white,line width=0.6pt,inner sep=2.2pt},
  edge/.style={line width=1.05pt,line cap=round},
  stab/.style={line width=1.05pt,line cap=round,densely dashed},
  half/.style={circle,draw=SoftInk!70,fill=white,line width=0.4pt,inner sep=0.9pt}]
  \coordinate (v1) at (-2.121,1.626);
  \coordinate (v2) at (2.121,1.626);
  \coordinate (v3) at (2.121,-1.626);
  \coordinate (v4) at (-2.121,-1.626);
  \coordinate (v5) at (-0.776,-2.222);
  \coordinate (v6) at (-2.898,0.595);
  \coordinate (v7) at (0.776,-2.222);
  \coordinate (v8) at (0.776,2.222);
  \coordinate (v9) at (2.898,-0.595);
  \coordinate (v10) at (-2.898,-0.595);
  \coordinate (v11) at (2.898,0.595);
  \coordinate (v12) at (-0.776,2.222);
  \draw[edge,draw=TwoAdic] (v1) -- (v6);
  \draw[edge,draw=TwoAdic] (v1) -- (v12);
  \draw[edge,draw=TwoAdic] (v2) -- (v8);
  \draw[edge,draw=TwoAdic] (v2) -- (v11);
  \draw[edge,draw=TwoAdic] (v3) -- (v7);
  \draw[edge,draw=TwoAdic] (v3) -- (v9);
  \draw[edge,draw=TwoAdic] (v4) -- (v5);
  \draw[edge,draw=TwoAdic] (v4) -- (v10);
  \draw[edge,draw=TwoAdic] (v5) -- (v7);
  \draw[edge,draw=TwoAdic] (v8) -- (v12);
  \draw[edge,draw=ThreeAdic] (v1) -- (v3);
  \draw[edge,draw=ThreeAdic] (v1) -- (v5);
  \draw[edge,draw=ThreeAdic] (v1) -- (v10);
  \draw[edge,draw=ThreeAdic] (v2) -- (v4);
  \draw[edge,draw=ThreeAdic] (v2) -- (v7);
  \draw[edge,draw=ThreeAdic] (v2) -- (v9);
  \draw[edge,draw=ThreeAdic] (v3) -- (v8);
  \draw[edge,draw=ThreeAdic] (v3) -- (v11);
  \draw[edge,draw=ThreeAdic] (v4) -- (v6);
  \draw[edge,draw=ThreeAdic] (v4) -- (v12);
  \draw[edge,draw=ThreeAdic] (v5) to[bend left=10] (v8);
  \draw[edge,draw=ThreeAdic] (v5) to[bend right=10] (v8);
  \draw[edge,draw=ThreeAdic] (v6) -- (v9);
  \draw[edge,draw=ThreeAdic] (v6) -- (v10);
  \draw[edge,draw=ThreeAdic] (v7) to[bend left=10] (v12);
  \draw[edge,draw=ThreeAdic] (v7) to[bend right=10] (v12);
  \draw[edge,draw=ThreeAdic] (v9) -- (v11);
  \draw[edge,draw=ThreeAdic] (v10) -- (v11);
  \draw[edge,draw=TwoAdic] (v1) -- ++(127.476:0.42) node[half] {};
  \draw[edge,draw=TwoAdic] (v2) -- ++(52.524:0.42) node[half] {};
  \draw[edge,draw=TwoAdic] (v3) -- ++(-52.524:0.42) node[half] {};
  \draw[edge,draw=TwoAdic] (v4) -- ++(-127.476:0.42) node[half] {};
  \draw[edge,draw=TwoAdic] (v5) -- ++(-101.609:0.42) node[half] {};
  \draw[edge,draw=TwoAdic] (v6) -- ++(146.736:0.42) node[half] {};
  \draw[edge,draw=TwoAdic] (v6) -- ++(174.736:0.42) node[half] {};
  \draw[edge,draw=TwoAdic] (v7) -- ++(-78.391:0.42) node[half] {};
  \draw[edge,draw=TwoAdic] (v8) -- ++(78.391:0.42) node[half] {};
  \draw[edge,draw=TwoAdic] (v9) -- ++(-33.264:0.42) node[half] {};
  \draw[edge,draw=TwoAdic] (v9) -- ++(-5.264:0.42) node[half] {};
  \draw[edge,draw=TwoAdic] (v10) -- ++(-174.736:0.42) node[half] {};
  \draw[edge,draw=TwoAdic] (v10) -- ++(-146.736:0.42) node[half] {};
  \draw[edge,draw=TwoAdic] (v11) -- ++(5.264:0.42) node[half] {};
  \draw[edge,draw=TwoAdic] (v11) -- ++(33.264:0.42) node[half] {};
  \draw[edge,draw=TwoAdic] (v12) -- ++(101.609:0.42) node[half] {};
  \node[vertex,fill=ChiThreeVOne] at (v1) {};
  \node[vertex,fill=ChiThreeVTwo] at (v2) {};
  \node[vertex,fill=ChiThreeVThree] at (v3) {};
  \node[vertex,fill=ChiThreeVFour] at (v4) {};
  \node[vertex,fill=ChiThreeVFive] at (v5) {};
  \node[vertex,fill=ChiThreeVSix] at (v6) {};
  \node[vertex,fill=ChiThreeVSeven] at (v7) {};
  \node[vertex,fill=ChiThreeVEight] at (v8) {};
  \node[vertex,fill=ChiThreeVNine] at (v9) {};
  \node[vertex,fill=ChiThreeVTen] at (v10) {};
  \node[vertex,fill=ChiThreeVEleven] at (v11) {};
  \node[vertex,fill=ChiThreeVTwelve] at (v12) {};
\end{tikzpicture}}\\[2pt]
        $X_{\Gamma}$
    \end{minipage}\hfill
    \begin{minipage}{.48\textwidth}
        \centering
        \scalebox{0.8}{\definecolor{SoftInk}{HTML}{2F3437}%
\definecolor{TwoAdic}{HTML}{4F7D72}%
\definecolor{ThreeAdic}{HTML}{A85D67}%
\definecolor{ChiThreeVOne}{HTML}{4B5F9E}%
\definecolor{ChiThreeVTwo}{HTML}{F2C94C}%
\definecolor{ChiThreeVThree}{HTML}{2F8F6A}%
\definecolor{ChiThreeVFour}{HTML}{C64E3B}%
\definecolor{ChiThreeVFive}{HTML}{6D8F2F}%
\definecolor{ChiThreeVSix}{HTML}{2F76C8}%
\definecolor{ChiThreeVSeven}{HTML}{E07B28}%
\definecolor{ChiThreeVEight}{HTML}{8B5FBF}%
\definecolor{ChiThreeVNine}{HTML}{D84F82}%
\definecolor{ChiThreeVTen}{HTML}{2BA7A0}%
\definecolor{ChiThreeVEleven}{HTML}{B88920}%
\definecolor{ChiThreeVTwelve}{HTML}{5F6B73}%
\begin{tikzpicture}[x=1cm,y=1cm,
  vertex/.style={circle,draw=SoftInk!70,line width=0.3pt,inner sep=1.9pt},
  edge/.style={line width=1.05pt,line cap=round},
  half/.style={circle,draw=SoftInk!70,fill=white,line width=0.4pt,inner sep=0.9pt}]
  \coordinate (v1) at (-0.8,2.2);
  \coordinate (v2) at (0.8,2.2);
  \coordinate (v3) at (-0.8,-2.2);
  \coordinate (v4) at (0.8,-2.2);
  \coordinate (v5) at (-2.4,0);
  \coordinate (v6) at (0.6,0.6);
  \coordinate (v7) at (-1.5,0);
  \coordinate (v8) at (1.5,0);
  \coordinate (v9) at (-0.6,0.6);
  \coordinate (v10) at (-0.6,-0.6);
  \coordinate (v11) at (0.6,-0.6);
  \coordinate (v12) at (2.4,0);
  \draw[edge,draw=TwoAdic] (v1) .. controls (-4.3,2.4) and (-4.3,-2.4) .. (v3);
  \draw[edge,draw=TwoAdic] (v1) -- (v5);
  \draw[edge,draw=TwoAdic] (v1) -- (v6);
  \draw[edge,draw=TwoAdic] (v1) -- (v7);
  \draw[edge,draw=TwoAdic] (v1) -- (v9);
  \draw[edge,draw=TwoAdic] (v2) .. controls (4.3,2.4) and (4.3,-2.4) .. (v4);
  \draw[edge,draw=TwoAdic] (v2) -- (v5);
  \draw[edge,draw=TwoAdic] (v2) -- (v6);
  \draw[edge,draw=TwoAdic] (v2) -- (v7);
  \draw[edge,draw=TwoAdic] (v2) -- (v9);
  \draw[edge,draw=TwoAdic] (v3) -- (v8);
  \draw[edge,draw=TwoAdic] (v3) -- (v10);
  \draw[edge,draw=TwoAdic] (v3) -- (v11);
  \draw[edge,draw=TwoAdic] (v3) -- (v12);
  \draw[edge,draw=TwoAdic] (v4) -- (v8);
  \draw[edge,draw=TwoAdic] (v4) -- (v10);
  \draw[edge,draw=TwoAdic] (v4) -- (v11);
  \draw[edge,draw=TwoAdic] (v4) -- (v12);
  \draw[edge,draw=TwoAdic] (v5) -- (v7);
  \draw[edge,draw=TwoAdic] (v6) -- (v8);
  \draw[edge,draw=TwoAdic] (v6) -- (v9);
  \draw[edge,draw=TwoAdic] (v6) -- (v11);
  \draw[edge,draw=TwoAdic] (v7) -- (v9);
  \draw[edge,draw=TwoAdic] (v7) -- (v10);
  \draw[edge,draw=TwoAdic] (v8) -- (v11);
  \draw[edge,draw=TwoAdic] (v8) -- (v12);
  \draw[edge,draw=TwoAdic] (v9) -- (v10);
  \draw[edge,draw=TwoAdic] (v10) -- (v11);
  \draw[edge,draw=TwoAdic] (v5) -- ++(135:0.42) node[half] {};
  \draw[edge,draw=TwoAdic] (v5) -- ++(180:0.42) node[half] {};
  \draw[edge,draw=TwoAdic] (v5) -- ++(225:0.42) node[half] {};
  \draw[edge,draw=TwoAdic] (v12) -- ++(45:0.42) node[half] {};
  \draw[edge,draw=TwoAdic] (v12) -- ++(0:0.42) node[half] {};
  \draw[edge,draw=TwoAdic] (v12) -- ++(315:0.42) node[half] {};
  \draw[edge,draw=TwoAdic] (v1) -- ++(100:0.42) node[half] {};
  \draw[edge,draw=TwoAdic] (v2) -- ++(80:0.42) node[half] {};
  \draw[edge,draw=TwoAdic] (v3) -- ++(260:0.42) node[half] {};
  \draw[edge,draw=TwoAdic] (v4) -- ++(280:0.42) node[half] {};
  \draw[edge,draw=TwoAdic] (v7) -- ++(270:0.42) node[half] {};
  \draw[edge,draw=TwoAdic] (v8) -- ++(90:0.42) node[half] {};
  \draw[edge,draw=TwoAdic] (v9) -- ++(135:0.42) node[half] {};
  \draw[edge,draw=TwoAdic] (v11) -- ++(315:0.42) node[half] {};
  \draw[edge,draw=TwoAdic] (v6) -- ++(45:0.42) node[half] {};
  \draw[edge,draw=TwoAdic] (v10) -- ++(225:0.42) node[half] {};
  \node[vertex,fill=ChiThreeVOne] at (v1) {};
  \node[vertex,fill=ChiThreeVTwo] at (v2) {};
  \node[vertex,fill=ChiThreeVThree] at (v3) {};
  \node[vertex,fill=ChiThreeVFour] at (v4) {};
  \node[vertex,fill=ChiThreeVFive] at (v5) {};
  \node[vertex,fill=ChiThreeVSix] at (v6) {};
  \node[vertex,fill=ChiThreeVSeven] at (v7) {};
  \node[vertex,fill=ChiThreeVEight] at (v8) {};
  \node[vertex,fill=ChiThreeVNine] at (v9) {};
  \node[vertex,fill=ChiThreeVTen] at (v10) {};
  \node[vertex,fill=ChiThreeVEleven] at (v11) {};
  \node[vertex,fill=ChiThreeVTwelve] at (v12) {};
\end{tikzpicture}}\\[2pt]
        $X_{\Pi_1}$
    \end{minipage}\\[12pt]
    \begin{minipage}{.48\textwidth}
        \centering
        \scalebox{0.8}{\definecolor{SoftInk}{HTML}{2F3437}%
\definecolor{TwoAdic}{HTML}{4F7D72}%
\definecolor{ThreeAdic}{HTML}{A85D67}%
\definecolor{ChiThreeVOne}{HTML}{4B5F9E}%
\definecolor{ChiThreeVTwo}{HTML}{F2C94C}%
\definecolor{ChiThreeVThree}{HTML}{2F8F6A}%
\definecolor{ChiThreeVFour}{HTML}{C64E3B}%
\definecolor{ChiThreeVFive}{HTML}{6D8F2F}%
\definecolor{ChiThreeVSix}{HTML}{2F76C8}%
\definecolor{ChiThreeVSeven}{HTML}{E07B28}%
\definecolor{ChiThreeVEight}{HTML}{8B5FBF}%
\definecolor{ChiThreeVNine}{HTML}{D84F82}%
\definecolor{ChiThreeVTen}{HTML}{2BA7A0}%
\definecolor{ChiThreeVEleven}{HTML}{B88920}%
\definecolor{ChiThreeVTwelve}{HTML}{5F6B73}%
\begin{tikzpicture}[x=1cm,y=1cm,
  vertex/.style={circle,draw=SoftInk!70,line width=0.3pt,inner sep=1.9pt},
  fixed/.style={rectangle,draw=SoftInk,fill=white,line width=0.6pt,inner sep=2.2pt},
  edge/.style={line width=1.05pt,line cap=round},
  stab/.style={line width=1.05pt,line cap=round,densely dashed},
  half/.style={circle,draw=SoftInk!70,fill=white,line width=0.4pt,inner sep=0.9pt}]
  \coordinate (v1) at (-0.8,1.6);
  \coordinate (v2) at (0.8,1.6);
  \coordinate (v3) at (-0.8,-1.6);
  \coordinate (v4) at (0.8,-1.6);
  \coordinate (v5) at (-3,1.2);
  \coordinate (v6) at (-2,1);
  \coordinate (v7) at (2,-1);
  \coordinate (v8) at (3,-1.2);
  \coordinate (v9) at (-0.8,0.55);
  \coordinate (v10) at (0.8,0.55);
  \coordinate (v11) at (-0.8,-0.55);
  \coordinate (v12) at (0.8,-0.55);
  \coordinate (v13) at (-0.9,3);
  \coordinate (v14) at (0.9,3);
  \coordinate (v15) at (-0.9,-3);
  \coordinate (v16) at (0.9,-3);
  \draw[stab,draw=SoftInk] (v1) -- (v2);
  \draw[edge,draw=TwoAdic] (v1) -- (v5);
  \draw[stab,draw=SoftInk] (v1) -- (v9);
  \draw[edge,draw=TwoAdic] (v1) -- (v13);
  \draw[edge,draw=TwoAdic] (v2) -- (v6);
  \draw[stab,draw=SoftInk] (v2) -- (v10);
  \draw[edge,draw=TwoAdic] (v2) -- (v14);
  \draw[stab,draw=SoftInk] (v3) -- (v4);
  \draw[edge,draw=TwoAdic] (v3) -- (v7);
  \draw[stab,draw=SoftInk] (v3) -- (v11);
  \draw[edge,draw=TwoAdic] (v3) -- (v15);
  \draw[edge,draw=TwoAdic] (v4) -- (v8);
  \draw[stab,draw=SoftInk] (v4) -- (v12);
  \draw[edge,draw=TwoAdic] (v4) -- (v16);
  \draw[edge,draw=TwoAdic] (v5) to[bend left=12] (v6);
  \draw[edge,draw=TwoAdic] (v5) to[bend right=12] (v6);
  \draw[edge,draw=TwoAdic] (v5) -- (v11);
  \draw[edge,draw=TwoAdic] (v5) -- (v14);
  \draw[edge,draw=TwoAdic] (v5) -- (v15);
  \draw[edge,draw=TwoAdic] (v6) -- (v9);
  \draw[edge,draw=TwoAdic] (v6) to[bend left=12] (v13);
  \draw[edge,draw=TwoAdic] (v6) to[bend right=12] (v13);
  \draw[edge,draw=TwoAdic] (v7) to[bend left=12] (v8);
  \draw[edge,draw=TwoAdic] (v7) to[bend right=12] (v8);
  \draw[edge,draw=TwoAdic] (v7) -- (v12);
  \draw[edge,draw=TwoAdic] (v7) to[bend left=12] (v16);
  \draw[edge,draw=TwoAdic] (v7) to[bend right=12] (v16);
  \draw[edge,draw=TwoAdic] (v8) -- (v10);
  \draw[edge,draw=TwoAdic] (v8) -- (v14);
  \draw[edge,draw=TwoAdic] (v8) -- (v15);
  \draw[stab,draw=SoftInk] (v9) -- (v10);
  \draw[stab,draw=SoftInk] (v9) to[bend left=12] (v11);
  \draw[stab,draw=SoftInk] (v9) to[bend right=12] (v11);
  \draw[stab,draw=SoftInk] (v10) to[bend left=12] (v12);
  \draw[stab,draw=SoftInk] (v10) to[bend right=12] (v12);
  \draw[stab,draw=SoftInk] (v11) -- (v12);
  \draw[edge,draw=TwoAdic] (v13) to[bend left=12] (v14);
  \draw[edge,draw=TwoAdic] (v13) to[bend right=12] (v14);
  \draw[edge,draw=TwoAdic] (v13) .. controls (-4.2,3.6) and (-4.2,-3.6) .. (v15);
  \draw[edge,draw=TwoAdic] (v14) .. controls (4.2,3.6) and (4.2,-3.6) .. (v16);
  \draw[edge,draw=TwoAdic] (v15) to[bend left=12] (v16);
  \draw[edge,draw=TwoAdic] (v15) to[bend right=12] (v16);
  \node[vertex,fill=ChiThreeVOne] at (v5) {};
  \node[vertex,fill=ChiThreeVTwo] at (v6) {};
  \node[vertex,fill=ChiThreeVThree] at (v7) {};
  \node[vertex,fill=ChiThreeVFour] at (v8) {};
  \node[vertex,fill=ChiThreeVFive] at (v13) {};
  \node[vertex,fill=ChiThreeVSix] at (v14) {};
  \node[vertex,fill=ChiThreeVSeven] at (v15) {};
  \node[vertex,fill=ChiThreeVEight] at (v16) {};
  \node[fixed] at (v1) {};
  \node[fixed] at (v2) {};
  \node[fixed] at (v3) {};
  \node[fixed] at (v4) {};
  \node[fixed] at (v9) {};
  \node[fixed] at (v10) {};
  \node[fixed] at (v11) {};
  \node[fixed] at (v12) {};
\end{tikzpicture}}\\[2pt]
        $X_{\Pi_2}$
    \end{minipage}
    \caption{The one-skeleta of the quotients of $T_\fs$ by the lattices of Section \ref{sec:qsqrt33}. The two edge colors of $X_\Gamma$ and $X_\Lambda$ are the two tree directions; the other lattices interchange the two trees. In $X_\Gamma$ and $X_{\Pi_1}$, an edge ending in a small open circle is the image of an edge inverted by an involution. In $X_{\Pi_2}$, the squares and dashed edges are the vertices and edges with stabilizer of order $2$.}
    \label{fig:qsqrt33-one-skeleta}
\end{figure}

\bibliography{references}

\end{document}